\documentclass[11pt,letterpaper]{article}
\RequirePackage[backend=biber,bibstyle=authoryear,citestyle=authoryear,
natbib=true,sorting=nyt,sortcites=true,maxnames=10,minnames=10,
hyperref=true,abbreviate=false,date=long,isbn=true,eprint=true,
uniquename=init,giveninits=true]{biblatex}
\usepackage{hyphenat,graphicx}
\usepackage[colorlinks,linktocpage,breaklinks]{hyperref}
\makeatletter
\def\@filecolor{blue}
\def\@linkcolor{blue}
\def\@citecolor{blue}
\def\@urlcolor{blue}
\let\@old@citep\citep
\let\@old@citet\citet
\let\@old@citeauthor\citeauthor
\def\citep{\@old@citep*}
\def\citet{\@old@citet*}
\def\citeauthor{\@old@citeauthor*}
\let\cite\citep
\makeatother
\usepackage[top=1.5in,bottom=1.2in,left=1.15in,right=1.15in,letterpaper]%
  {geometry}
\makeatletter
\def\ps@headings{%
    \let\@mkboth\@gobbletwo
    \def\@oddhead{\hss\scshape\shorttitle\hss\reset@font\rmfamily\thepage}
    \def\@evenhead{\reset@font\rmfamily\thepage\hss\scshape\shortauthors\hss}
    \let\@oddfoot\@empty\let\@evenfoot\@empty}
\@twosidetrue
\makeatother
\usepackage[neveradjust]{paralist}
\makeatletter
\@definecounter{@ADL@savenum}
\newcommand\savenum{\setcounter{@ADL@savenum}%
  {\the\@nameuse{c@\@listctr}}}
\newcommand\resumenum{\setcounter{\@listctr}{\arabic{@ADL@savenum}}}
\makeatother
\usepackage[reqno,centertags]{amsmath}
\usepackage{amssymb,xspace,twoopt,xy,mathrsfs}
\xyoption{arrow}
\xyoption{matrix}
\xyoption{curve}
\usepackage[mathcal]{eucal}
\makeatletter
\DeclareFontFamily{OMX}{MnSymbolE}{}
\DeclareSymbolFont{MnLargeSymbols}{OMX}{MnSymbolE}{m}{n}
\SetSymbolFont{MnLargeSymbols}{bold}{OMX}{MnSymbolE}{b}{n}
\DeclareFontShape{OMX}{MnSymbolE}{m}{n}{
    <-6>  MnSymbolE5
   <6-7>  MnSymbolE6
   <7-8>  MnSymbolE7
   <8-9>  MnSymbolE8
   <9-10> MnSymbolE9
  <10-12> MnSymbolE10
  <12->   MnSymbolE12
}{}
\DeclareFontShape{OMX}{MnSymbolE}{b}{n}{
    <-6>  MnSymbolE-Bold5
   <6-7>  MnSymbolE-Bold6
   <7-8>  MnSymbolE-Bold7
   <8-9>  MnSymbolE-Bold8
   <9-10> MnSymbolE-Bold9
  <10-12> MnSymbolE-Bold10
  <12->   MnSymbolE-Bold12
}{}
\let\llangle\@undefined
\let\rrangle\@undefined
\DeclareMathDelimiter{\llangle}{\mathopen}%
                     {MnLargeSymbols}{'164}{MnLargeSymbols}{'164}
\DeclareMathDelimiter{\rrangle}{\mathclose}%
                     {MnLargeSymbols}{'171}{MnLargeSymbols}{'171}
\DeclareRobustCommand\widecheck[1]{{\mathpalette\@widecheck{#1}}}
\def\@widecheck#1#2{%
    \setbox\z@\hbox{\m@th$#1#2$}%
    \setbox\tw@\hbox{\m@th$#1%
       \widehat{%
          \vrule\@width\z@\@height\ht\z@
          \vrule\@height\z@\@width\wd\z@}$}%
    \dp\tw@-\ht\z@
    \@tempdima\ht\z@ \advance\@tempdima2\ht\tw@ \divide\@tempdima\thr@@
    \setbox\tw@\hbox{%
       \raise\@tempdima\hbox{\scalebox{1}[-1]{\lower\@tempdima\box
\tw@}}}%
    {\ooalign{\box\tw@ \cr \box\z@}}}
\makeatother
\let\sinfty\infty
\input cyracc.def
\let\subsetneq\subset
\let\subset\subseteq

\newcommand\real{\mathbb{R}}
\newcommand\realp{\real_{>0}}
\newcommand\realnn{\real_{\ge0}}
\newcommand\realn{\real_{<0}}
\newcommand\realnp{\real_{\le0}}
\newcommand\complex{\mathbb{C}}
\newcommand\integer{\mathbb{Z}}
\newcommand\integerp{\integer_{>0}}
\newcommand\integernn{\integer_{\ge0}}
\renewcommand\Re{\operatorname{Re}}
\newcommand\eul{\textup{e}}
\newcommand\imag{\textup{i}}
\newcommand\scirc{\raise1pt\hbox{$\,\scriptstyle\circ\,$}}
\newcommand\sscirc{\hbox{$\,\scriptscriptstyle\circ\,$}}
\newcommand\eqdef{\triangleq}
\newcommand\subscr[2]{#1_{\textup{#2}}}

\newcommand\ol[1]{\overline{#1}}

\newcommand\disjointunion{\mathop{\overset{\sscirc}{\cup}}}

\makeatletter
\newcommand\slnorm{\lvert}
\newcommand\srnorm{\rvert}
\newcommand\snorm[1]{\slnorm #1\srnorm}

\newcommand\dlnorm{\lVert}
\newcommand\drnorm{\rVert}
\newcommand\dnorm[1]{\dlnorm #1\drnorm}
\newcommand\adnorm[1]{\left\dlnorm #1\right\drnorm}
\newcommand\setdef[2]{\{#1\;|\enspace#2\}}
\newcommand\asetdef[2]{\left\{#1\immediate\vphantom{#2}\;\right|
  \left.\immediate\vphantom{#1}\enspace#2\right\}}
\newcommand\ifam[1]{(#1)}
\newcommand\natpair[2]{\langle#1;#2\rangle}
\newcommand\inprod[2]{\langle#1,#2\rangle}
\newcommand\vecspan{\@ifnextchar[{\@ADL@normalvecspan}
  {\@ADL@@normalvecspan}}
\newcommand\avecspan{\@ifnextchar[{\@ADL@vecspan}{\@ADL@@vecspan}}
\def\@ADL@vecspan[#1]#2{\operatorname{span}_{#1}\left(#2\right)}
\newcommand\@ADL@@vecspan[1]{\operatorname{span}\left(#1\right)}
\def\@ADL@normalvecspan[#1]#2{\operatorname{span}_{#1}(#2)}
\newcommand\@ADL@@normalvecspan[1]{\operatorname{span}(#1)}
\makeatother
\newcommand\ie{i.e.,}
\newcommand\eg{e.g.,}
\newcommand\cf{cf.}

\newcommand\resp{resp.}
\newcommand\der{\boldsymbol{d}}
\newcommand\mathupper[1]{\textup{#1}}
\renewcommand\d[1]{{\normalfont\textrm{d}}#1}
\newcommand\codim{\operatorname{codim}}
\newcommand\image{\operatorname{image}}
\newcommand\rank{\operatorname{rank}}
\newcommand\id{\operatorname{id}}
\newcommand\interior{\operatorname{cl}}
\newcommand\closure{\operatorname{cl}}

\newcommand\affhull{\operatorname{aff}}
\newcommand\cohull{\operatorname{conv}}
\newcommand\deriv[2]{\frac{{\normalfont\mathupper{d}}#1}
  {{\normalfont\mathupper{d}}#2}}
\newcommand\derivatzero[2]{\frac{{\normalfont\mathupper{d}}#1}
  {{\normalfont\mathupper{d}}#2}\Big|_{#2=0}}

\newcommand\tderivatzero[2]{{\textstyle\frac{{\normalfont\mathupper{d}}#1}
  {{\normalfont\mathupper{d}}#2}}\big|_{#2=0}}
\newcommand\pderiv[2]{\frac{\partial#1}{\partial#2}}
\makeatletter
\newcommand\linder{\@ifnextchar[{\@ADL@rlinder}{\@ADL@linder}}
\def\@ADL@rlinder[#1]#2{\boldsymbol{D}^{#1}#2}
\newcommand\@ADL@linder[1]{\boldsymbol{D}#1}
\newcommand\plinder{\@ifnextchar[{\@ADL@rplinder}{\@ADL@plinder}}
\def\@ADL@rplinder[#1]#2#3{\boldsymbol{D}^{#1}_{#2}#3}
\newcommand\@ADL@plinder[2]{\boldsymbol{D}_{#1}#2}
\makeatother
\newcommand\map[3]{#1\colon#2\rightarrow#3}

\newcommand\mapdef[5]{\begin{aligned}
  #1\colon&\begin{aligned}[t]#2\end{aligned}\rightarrow
  \begin{aligned}[t]#3\end{aligned}\\&\begin{aligned}[t]#4\end{aligned}
  \mapsto\begin{aligned}[t]#5\end{aligned}\end{aligned}}
\newcommand\mapschar{\mathupper{C}}
\newcommand\C{\mapschar}

\makeatletter
\newcommand\lin{\@ifnextchar[{\@ADL@klinmap}{\@ADL@linmap}}
\def\@ADL@klinmap[#1]#2#3{\def\@tempa{#3}%
  \ifx\@tempa\@empty\mathupper{L}^{#1}(#2;#2)\else%
  \mathupper{L}^{#1}(#2;#3)\fi}
\def\@ADL@linmap#1#2{\def\@tempa{#2}\ifx\@tempa\@empty\mathupper{L}(#1;#1)
  \else\mathupper{L}(#1;#2)\fi}
\makeatother
\newcommandtwoopt\sections[3][\sinfty][\null]{\Gamma^{#1}_{#2}(#3)}

\newcommand\man[1]{\mathsf{#1}}
\newcommand\zsec[1]{Z(#1)}
\makeatletter
\newcommand\tb{\@ifnextchar[{\@ADL@tbarg}{\@ADL@tb}}
\def\@ADL@tbarg[#1]#2{\man{T}_{#1}#2}
\newcommand\@ADL@tb[1]{\man{T}#1}

\newcommand\ctb{\@ifnextchar[{\@ADL@ctbarg}{\@ADL@ctb}}
\def\@ADL@ctbarg[#1]#2{\man{T}^*_{#1}#2}
\newcommand\@ADL@ctb[1]{\man{T}^*#1}

\newcommand\func{\@ifnextchar[{\@ADL@crfuncs}{\@ADL@cinftyfuncs}}
\def\@ADL@crfuncs[#1]#2{\mapschar^{#1}(#2)}
\newcommand\@ADL@cinftyfuncs[1]{\mapschar^\infty(#1)}
\newcommand\sfunc{\@ifnextchar[{\@ADL@crsfuncs}{\@ADL@cinftysfuncs}}
\def\@ADL@crsfuncs[#1]#2{\mathscr{C}^{#1}_{#2}}
\newcommand\@ADL@cinftysfuncs[1]{\mathscr{C}^\sinfty_{#1}}
\newcommand\mappings{\@ifnextchar[{\@ADL@crmappings}{\@ADL@cinftymappings}}
\def\@ADL@crmappings[#1]#2#3{\mapschar^{#1}(#2;#3)}
\newcommand\@ADL@cinftymappings[2]{\mapschar^\infty(#1;#2)}
\newcommand\gfunc{\@ifnextchar[{\@ADL@crgfuncs}{\@ADL@cinftygfuncs}}
\def\@ADL@crgfuncs[#1]#2#3{\mathscr{C}^{#1}_{#2,#3}}
\newcommand\@ADL@cinftygfuncs[2]{\mathscr{C}^\sinfty_{#1,#2}}
\newcommand\gsections{\@ifnextchar[{\@ADL@rgsections}{\@ADL@gsections}}
\def\@ADL@rgsections[#1]#2#3{\mathscr{G}^{#1}_{#2,#3}}
\newcommand\@ADL@gsections[2]{\mathscr{G}^\infty_{#1,#2}}
\newcommand\ssections{\@ifnextchar[{\@ADL@rssections}{\@ADL@ssections}}
\def\@ADL@rssections[#1]#2{\mathscr{G}^{#1}_{#2}}
\newcommand\@ADL@ssections[1]{\mathscr{G}^\infty_{#1}}
\newcommand\tf{\@ifnextchar[{\@ADL@tfarg}{\@ADL@tf}}
\def\@ADL@tfarg[#1]#2{\@ifnextchar[{\@ADL@@tfargr{#1}{#2}}
  {\@ADL@@tfarg{#1}{#2}}}
\def\@ADL@@tfargr#1#2[#3]{T^{#3}_{#1}#2}
\newcommand\@ADL@@tfarg[2]{T_{#1}#2}
\newcommand\@ADL@tf[1]{\@ifnextchar[{\@ADL@@tfr{#1}}{\@ADL@@tf{#1}}}
\def\@ADL@@tfr#1[#2]{T^{#2}#1}
\newcommand\@ADL@@tf[1]{T#1}
\newcommand\ctf{\@ifnextchar[{\@ADL@ctfarg}{\@ADL@ctf}}
\def\@ADL@ctfarg[#1]#2{T^*_{#1}#2}
\newcommand\@ADL@ctf[1]{T^*#1}
\newcommand\lieder[2]{\def\@tempa{#2}\ifx\@tempa\@empty%
  \boldsymbol{\mathscr{L}}_{#1}\else\boldsymbol{\mathscr{L}}_{#1}#2\fi}
\makeatletter
\def\interval{\@ifnextchar({\@ADL@openleftint}{\@ADL@closedleftint}}
\def\@ADL@openleftint(#1,#2{(#1,#2%
  \@ifnextchar){\@ADL@openrightint}{\@ADL@closedrightint}}
\def\@ADL@closedleftint[#1,#2{[#1,#2%
  \@ifnextchar){\@ADL@openrightint}{\@ADL@closedrightint}}
\def\@ADL@openrightint){)}
\def\@ADL@closedrightint]{]}
\makeatother
\usepackage{theorem}
\newtheorem{theorem}{Theorem}[section]
\newtheorem{proposition}[theorem]{Proposition}
\newtheorem{lemma}[theorem]{Lemma}
\newtheorem{corollary}[theorem]{Corollary}
\newtheorem{prooflemma}{Lemma}[theorem]

{\theorembodyfont{\normalfont\rmfamily}
\newtheorem{definition}[theorem]{Definition}
\newtheorem{remark}[theorem]{Remark}
\newtheorem{example}[theorem]{Example}
\newtheorem{examples}[theorem]{Examples}}
\makeatletter
\newcommand\pushright{\protect\@ADL@pushright}
\newcommand\@ADL@pushright[1]{{\ifvmode\null\hfill{#1}\par\else\ifmmode%
  \@ADLmaths@pushright{\hbox{#1}}\else\ifinner\@ADLhbox@pushright{#1}%
  \else\@ADLparag@pushright{#1}\fi\fi\fi}}
\newcommand\@ADLmaths@pushright[1]{{\ifinner\@ADLhbox@pushright{#1}\else%
  \tag*{$#1$}\fi}}
\newcommand\@ADLparag@pushright[1]{{\parfillskip=0pt\widowpenalty=10000%
  \displaywidowpenalty=10000\finalhyphendemerits=0\@ADLhbox@pushright#1\par}}
\newcommand\@ADLhbox@pushright{\unskip\nobreak\hfil\penalty50\hskip.2em%
  \null\hfill\hfill}
\newenvironment{proof}{\trivlist\item[\hskip\labelsep\textit{Proof:}\/]%
  \@ADLsave@set@qed\xspace\normalfont\rmfamily}
  {\qed\@ADLrestore@qed\endtrivlist}
\newif\if@ADL@qed\@ADL@qedfalse
\newcommand\qed{\protect\@ADL@qed{$\blacksquare$}}
\newcommand\@ADL@qed[1]{\if@ADL@qed\global\@ADL@qedfalse%
  \pushright{#1}\else\ifhmode\ifinner\else\par\fi\fi\fi}
\newcommand\@ADLrestore@qed{\global\let\if@ADL@qed\@ADLsaved@ifqed}
\newcommand\@ADLsave@set@qed{\let\@ADLsaved@ifqed
  \if@ADL@qed\global\@ADL@qedtrue}

\newenvironment{subproof}{\trivlist\item[\hskip\labelsep\textit{Proof:}\/]%
  \@ADLsave@set@subqed\normalfont\rmfamily}
  {\subqed\@ADLrestore@subqed\endtrivlist}
\newif\if@ADL@subqed\@ADL@subqedfalse
\newcommand\subqed{\protect\@ADL@subqed{$\blacktriangledown$}}
\newcommand\@ADL@subqed[1]{\if@ADL@subqed\global\@ADL@subqedfalse%
  \pushright{#1}\else\ifhmode\ifinner\else\par\fi\fi\fi}
\newcommand\@ADLrestore@subqed{\global\let\if@ADL@subqed\@ADLsaved@ifsubqed}
\newcommand\@ADLsave@set@subqed{\let\@ADLsaved@ifsubqed
  \if@ADL@subqed\global\@ADL@subqedtrue}
\newcommand\eqsubqed{\tag*{\subqed}}
\newif\if@ADL@oprocend\@ADL@oprocendfalse
\newcommand\oprocend{\@ADLsave@set@oprocend
  \protect\@ADL@oprocend{$\bullet$}\@ADLrestore@oprocend}
\newcommand\@ADL@oprocend[1]{\if@ADL@oprocend\global\@ADL@oprocendfalse%
  \pushright{#1}\else\ifhmode\ifinner\else\par\fi\fi\fi}
\newcommand\@ADLrestore@oprocend{\global
  \let\if@ADL@oprocend\@ADLsaved@ifoprocend}
\newcommand\@ADLsave@set@oprocend{\let\@ADLsaved@ifoprocend\if@ADL@oprocend%
  \global\@ADL@oprocendtrue}
\newcommand\eqoprocend{\tag*{\oprocend}}
\makeatother
\newenvironment{keywords}{\quote\small\textbf{Keywords.}}{\endquote}
\newenvironment{AMS}{\quote\small\textbf{AMS Subject Classifications (2010).}}
   {\endquote}
\newcommand\defn[1]{{\normalfont\bfseries\emph{\mathversion{bold}#1}}}
\numberwithin{equation}{section}
\newcommand\pldblref[2]{\mbox{\ref{#1}(\ref{#2})}}
\newcommand\enumdblref[2]{\mbox{\ref{#1}--\ref{#2}}}
\usepackage{xy}
\xyoption{arrow}
\xyoption{matrix}
\newcommand\mat[1]{\boldsymbol{#1}}
\newcommand\vect[1]{\boldsymbol{#1}}
\newcommand\alg[1]{\mathsf{#1}}
\newcommand\dist[1]{\mathsf{#1}}
\renewcommand\man[1]{\mathsf{#1}}
\newcommand\ts[1]{\mathcal{#1}}
\newcommand\ms[1]{\mathcal{#1}}
\newcommand\nbhd[1]{\mathcal{#1}}
\newcommand\sD{\mathscr{D}}
\newcommand\sF{\mathscr{F}}
\newcommand\sG{\mathscr{G}}
\newcommand\sI{\mathscr{I}}
\newcommand\sL{\mathscr{L}}
\newcommand\sM{\mathscr{M}}
\newcommand\sN{\mathscr{N}}
\newcommand\sO{\mathscr{O}}
\newcommand\sR{\mathscr{R}}

\newcommand\sU{\mathscr{U}}
\newcommand\sV{\mathscr{U}}
\newcommand\sX{\mathscr{X}}
\newcommand\sY{\mathscr{Y}}
\newcommand\Cbdd{\subscr{\C}{bdd}}
\newcommand\sCbdd{\subscr{\mathscr{C}}{bdd}}
\newcommand\Et{\mathupper{Et}}

\newcommand\Diff{\operatorname{Diff}}
\newcommand\Orb{\operatorname{Orb}}
\newcommand\grank{\operatorname{grank}}

\newcommand\metric{\mathbb{G}}
\newcommand\slieder[3]{\mathscr{L}^{#1}_{#2}#3}
\newcommand\modgen[1]{\langle#1\rangle}
\newcommand\cochain[2]{\mathupper{C}^{#2}_{#1}}
\newcommand\cohomker[2]{\mathupper{Z}^{#2}_{#1}}
\newcommand\cohomim[2]{\mathupper{B}^{#2}_{#1}}
\newcommand\cohom[2]{\mathupper{H}^{#2}_{#1}}
\makeatletter
\def\flow#1{\@ifnextchar[{\@tflow{#1}}{\@flow{#1}}}
\def\@tflow#1[#2]#3{\Phi^{#1}_{#2,#3}}
\def\@flow#1#2{\Phi^{#1}_{#2}}
\newcommand\oball{\@ifnextchar[{\@ADL@oballarg}{\@ADL@oballnoarg}}
\def\@ADL@oballarg[#1]#2#3{\mathsf{B}_{#1}(#2,#3)}
\newcommand\@ADL@oballnoarg[2]{\mathsf{B}(#1,#2)}
\newcommand\osections{\@ifnextchar[{\@ADL@rosections}{\@ADL@osections}}
\def\@ADL@rosections[#1]#2{\overline{\Gamma}^{#1}(#2)}
\newcommand\@ADL@osections[1]{\overline{\Gamma}^\sinfty(#1)}
\newcommand\ogsections{\@ifnextchar[{\@ADL@rogsections}{\@ADL@ogsections}}
\def\@ADL@rogsections[#1]#2#3{\ol{\mathscr{G}}^{#1}_{#2,#3}}
\newcommand\@ADL@ogsections[2]{\ol{\mathscr{G}}^\sinfty_{#1,#2}}
\newcommand\ossections{\@ifnextchar[{\@ADL@rossections}{\@ADL@ossections}}
\def\@ADL@rossections[#1]#2{\ol{\mathscr{G}}^{#1}_{#2}}
\newcommand\@ADL@ossections[1]{\ol{\mathscr{G}}^\sinfty_{#1}}
\def\jet{\@ifnextchar[{\@ADL@jetbase}{\@ADL@jetnobase}}
\def\@ADL@jetbase[#1]#2#3{\man{J}_{#1}^{#2}#3}
\def\@ADL@jetnobase#1#2{\man{J}^{#1}#2}
\makeatother
\title{Generalised subbundles and distributions:\\A comprehensive
review\thanks{Research supported in part by a grant from the Natural Sciences
and Engineering Research Council of Canada}}
\author{Andrew D.\ Lewis\thanks{Professor, Department of Mathematics and
Statistics, Queen's University, Kingston, ON K7L 3N6, Canada,
email:~\texttt{andrew.lewis@queensu.ca}}}
\newcommand\shorttitle{Generalised subbundles and distributions}
\newcommand\shortauthors{A.\ D.\ Lewis}
\date{2014/07/21}
\begin{document}
\maketitle

\begin{abstract}
Distributions,~\ie~subsets of tangent bundles formed by piecing together
subspaces of tangent spaces, are commonly encountered in the theory and
application of differential geometry.  Indeed, the theory of distributions is
a fundamental part of mechanics and control theory.

The theory of distributions is presented in a systematic way, and
self-contained proofs are given of some of the major results.  Parts of the
theory are presented in the context of generalised subbundles of vector
bundles.  Special emphasis is placed on understanding the r\^ole of sheaves
and understanding the distinctions between the smooth or finitely
differentiable cases and the real analytic case.  The Orbit Theorem and
applications, including Frobenius's Theorem and theorems on the equivalence
of families of vector fields, are considered in detail.  Examples illustrate
the phenomenon that can occur with generalised subbundles and distributions.
\end{abstract}

\tableofcontents

\section{Introduction}

Distributions arise naturally in differential geometry (for example, in the
characterisation of a Poisson manifold as being a disjoint union of its
symplectic leaves) and in many applications of differential geometry, such as
mechanics (for example, as arising from nonholonomic constraints in
mechanics), and control theory (for example, in characterisations of orbits
of families of vector fields).  As such, distributions have been widely
studied and much is known about them.  However, as is often the case with
objects that are widely used, there has arisen certain conventions for
handling distributions that are viable ``a lot of the time,'' but which, in
fact, do not have a basis in the general theory.  For just one example, when
distributions are used, there is often an unstated assumption of the
distribution and all distributions arising from it having locally constant
rank.  However, this assumption is not always valid, and interesting
phenomenon arise when it does not hold,~\eg~the general class of abnormal
sub-Riemannian minimisers described by \citet{WSL/HJS:94} rely on the
distribution generated by certain brackets being singular.  For this reason,
it seems that there may be some benefit to assembling the basic theory of
distributions in one place, with complete proofs of important results, and
this is what we do in this paper.  The objective is to present, in one place,
precise and general definitions and results relating to these definitions.
Many of these definitions and results are most naturally presented in the
setting of generalised subbundles of general vector bundles, and not just
tangent bundles.  Thus a substantial part of the paper deals with this.

A principal function of the paper is that of a review paper.  However, there
are also some other contributions which we now outline.

One of the contributions of the paper is to explicate clearly the r\^ole of
sheaf theory in handling distributions.  This contribution arises in four
ways:~(1)~in understanding the difficulties of dealing with real analytic
global sections of vector bundles, and the manner in which such global
sections arise;~(2)~in understanding clearly the r\^ole of analyticity in
results which require certain modules to be finitely generated;~(3)~in
understanding some of the algebraic properties of generalised
subbundles;~(4)~in properly characterising some local constructions involving
distributions,~\eg~invariance under vector fields and flows.  To carry out
these objectives in full detail requires a great deal to be drawn from sheaf
theory, particularly as it relates to complex differential geometry.  In this
paper we make these connections explicit, we believe for the first time.

Apart from this explication of sheaves in the theory of generalised subbundles and distributions, another contribution of the paper is to clarify certain ``folklorish'' parts of the theory,~\ie~things which are known to be true, but for which it is difficult (and for the author, impossible) to put together complete proofs in the existing literature.  Here are a few such bits of folklore:~(1)~the fact that the rank of a real analytic generalised subbundle attains its maximal rank on the complement of an analytic set (proved here as Proposition~\ref{prop:analytic-rank});~(2)~the Serre\textendash{}Swan Theorem for general vector bundles (proved here as Theorem~\ref{the:swan});~(3)~the Noetherian properties of real analytic germs (proved here as Proposition~\ref{prop:analytic-germs});~(4)~the r\^ole of these Noetherian properties in the local finite generation of modules (proved here as Theorem~\ref{the:analytic-fingen});~(5)~all variants of the Orbit Theorem (\ie~all combinations of fixed-time and arbitrary-time, and finitely generated and non-finitely generated versions);~(6)~a counterexample that shows that smooth involutive distributions need not be integrable.  The proofs of many of these facts are a matter of putting together known results; nonetheless, this has not been done to the best knowledge of the author.

There are also some new results in the paper.  In
Section~\ref{subsec:patchy-subsheaves} we introduce the class of ``patchy''
subsheaves of the sheaf of sections of a vector bundle.  In
Proposition~\ref{prop:smooth-patchy} and Theorem~\ref{the:patchy-analytic} we
show that the subsheaf of sections of a generalised subbundle is patchy.  We
subsequently show in Corollary~\ref{cor:analytic-coherent} that patchy real
analytic subsheaves are coherent, and this provides for these subsheaves
access to the machinery of the cohomology of coherent real analytic sheaves
presented in Section~\ref{subsec:coherent-sheaves}\@.  Other new results are
presented in Section~\ref{subsec:invariance}\@, where we consider the notions
of invariance of subsheaves of vector fields under diffeomorphisms and vector
fields, and discuss the relationship between invariance under a vector field
and invariance under the flow of the same vector field.  Here the topology of
stalks of sheaves discussed in Section~\ref{subsec:stalk-topology} features
prominently, and in the smooth case we reveal the r\^ole of the Whitney
Spectral Theorem.

What follows is an outline of the paper.  Many of the properties of
distributions, particularly their smoothness, are prescribed locally.  This
raises the immediate question as to whether these local constructions give
rise to meaningful global objects.  The best way to systematically handle
such an approach is via the use of sheaves.  Thus in
Section~\ref{sec:sheaves} we introduce the elements of sheaf theory that we
will require.  Many of the elementary constructions and results concerning
distributions are just as easily done with vector bundles, rather than
specifically with tangent bundles.  For this reason, we devote a significant
part of the paper to definitions and results for vector bundles.  In
Section~\ref{sec:vector-bundles} we give the basic definitions and properties
of generalised subbundles,~\ie~assignments of a subspace of each fibre of a
vector bundle.  The set of sections of a generalised subbundle is a submodule
of the module of sections.  This simple observation, appropriately parsed in
terms of sheaf theory, gives rise to some interesting algebraic structure for
generalised subbundles.  In Section~\ref{sec:algebra} we present some of this
algebraic theory.  In Section~\ref{sec:global-sections} we study the
important question of global generators of generalised subbundles.  This
provides an instance of the importance of sheaf theory in the study of
generalised subbundles.  Next in the paper we turn particularly to
distributions,~\ie~generalised subbundles of tangent bundles.  Here the
additional structure of vector fields having a flow, and all the consequences
of this, play a r\^ole.  In particular, in Section~\ref{sec:differential} we
look at invariant distributions and constructions related to the Lie bracket
of vector fields.  An important contribution of control theory to
differential geometry is the Orbit Theorem.  We study this theorem in some
detail in Section~\ref{sec:orbit-theorem}\@, in particular giving equal
emphasis to the so-called Fixed-time Orbit Theorem, something not normally
done (but see~\cite[\S2.4]{VJ:97}).  Related to the Orbit Theorem, but not
equivalent to it, is Frobenius's Theorem, which we present in
Section~\ref{sec:frobenius}\@.

\subsection*{Notation}

Let us establish the notation we use in the paper.

We write $A\subset B$ if $A$ is a subset of $B$\@, allowing that $A=B$\@.  If
$A$ is a strict subset of $B$ we denote this by $A\subsetneq B$\@.  The image
of a map $\map{f}{A}{B}$ is denoted by $\image(f)$\@.  The symbol
``$\eqdef$'' means ``is defined to be equal to.''

By $\integer$\@, $\integernn$\@, and $\integerp$ we denote the sets of
integers, nonnegative integers, and positive integers, respectively.  By
$\real$ and $\complex$ we denote the sets of real and complex numbers,
respectively.  By $\realp$\@, $\realnn$\@, $\realn$\@, and $\realnp$ we
denote the sets of positive real numbers, nonnegative real numbers, negative
real numbers, and nonpositive real numbers, respectively.  By $\real^n$ we
denote $n$-dimensional real Euclidean space (with $\complex^n$ being the
complex analogue) and by $\real^{m\times n}$ we denote the set of real
$m\times n$ matrices.

If $\alg{V}$ is a $\real$-vector space and if $S\subset\alg{V}$\@, by
$\vecspan[\real]{S}$\@, $\cohull(S)$\@, and $\affhull(S)$ we denote the
linear hull,~\ie~the linear span, the convex hull, and the affine hull of
$S$\@.  The kernel of a linear map $\map{L}{\alg{U}}{\alg{V}}$ is denoted by
$\ker(L)$\@.  Similarly, if $\alg{M}$ is a module over a commutative ring
$\alg{R}$ and if $S\subset\alg{M}$\@, then $\vecspan[\alg{R}]{S}$ is the
module generated by $S$\@.

By $\dnorm{\cdot}$ we denote the standard Euclidean norm on $\real^n$\@.  By
$\oball{r}{\vect{x}}$ we denote the open ball of radius $r$ centred at
$\vect{x}\in\real^n$\@.  If $\map{\vect{f}}{\nbhd{U}}{\real^m}$ is a
differentiable map from an open subset $\nbhd{U}\subset\real^n$\@, the $r$th
derivative of $\vect{f}$ at $\vect{x}\in\nbhd{U}$ is denoted by
$\linder[r]{\vect{f}}(\vect{x})$\@.  If $\nbhd{U}_j\subset\real^{n_j}$\@,
$j\in\{1,\dots,k\}$\@, and if
\begin{equation*}
\map{\vect{f}}{\nbhd{U}_1\times\dots\times\nbhd{U}_k}{\real^m}
\end{equation*}
is differentiable, we denote by
$\plinder{j}{\vect{f}}(\vect{x}_1,\dots,\vect{x}_k)$ the $j$th partial
derivative of $\vect{x}$ at $(\vect{x}_1,\dots,\vect{x}_k)\in
\nbhd{U}_1\times\dots\times\nbhd{U}_k$\@,~\ie~the derivative at $\vect{x}_j$
of the map
\begin{equation*}
\vect{x}'_j\mapsto\vect{f}(\vect{x}_1,\dots,\vect{x}'_j,\dots,\vect{x}_k).
\end{equation*}

For the most part, we follow the differential geometric notations and
conventions of~\cite{RA/JEM/TSR:88}\@.  The tangent bundle of a manifold
$\man{M}$ is denoted by $\tb{\man{M}}$ and $\tb[x]{\man{M}}$ denotes the
tangent space at $x$\@.  The cotangent bundle is denoted by $\ctb{\man{M}}$
and $\ctb[x]{\man{M}}$ denotes the cotangent space at $x$\@.  If
$\map{f}{\man{M}}{\man{N}}$ is a differentiable map between manifolds, we
denote its derivative by $\map{\tf{f}}{\tb{\man{M}}}{\tb{\man{N}}}$\@, with
$\tf[x]{f}$ denoting the restriction of this derivative to
$\tb[x]{\man{M}}$\@.  The differential of a differentiable function
$\map{f}{\man{M}}{\real}$ is denoted by
$\map{\der{f}}{\man{M}}{\ctb{\man{M}}}$\@.

By $\man{E}_x$ we denote the fibre at $x$ of a vector bundle
$\map{\pi}{\man{E}}{\man{M}}$\@.  If $\map{\pi}{\man{E}}{\man{M}}$ is a
vector bundle and if $\nbhd{U}\subset\man{M}$ is open, then
$\man{E}|\nbhd{U}$ denotes the restriction of this vector bundle to
$\nbhd{U}$\@.

We shall speak of geometric objects of class $\C^r$ for
$r\in\integernn\cup\{\infty,\omega\}$\@, with objects of class $\C^\infty$
being infinitely differentiable and objects of class $\C^\omega$ being real
analytic.  We shall often use language like, ``let $\man{M}$ be a manifold of
class $\C^\infty$ or class $\C^\omega$\@, as required.''  By this we mean
that the reader should ascribe the attributes of smoothness or real
analyticity as needed to make sense of the ensuing statements.  For
$r\in\integernn\cup\{\infty,\omega\}$\@, we denote the set of functions of
class $\C^r$ on a manifold $\man{M}$ by $\func[r]{\man{M}}$ and the set of
sections of class $\C^r$ of a smooth or real analytic (as is required) vector
bundle $\map{\pi}{\man{E}}{\man{M}}$ by $\sections[r]{\man{E}}$\@.

By $\lieder{X}{f}$ we denote the Lie derivative of the function $f$ by the
vector field $X$\@.  By $\flow{X}{t}$ we denote the flow of a vector field
$X$ on $\man{M}$\@.  Thus $t\mapsto\flow{X}{t}(x)$ is the integral curve of
$X$ through $x$\@.

\section{A wee bit of sheaf theory}\label{sec:sheaves}

Our Definition~\ref{def:subbundle} for what we mean by a generalised
subbundle with prescribed smoothness will be local.  One of the obvious
questions arising from this sort of construction is whether there are in fact
any globally defined vector fields taking values in the generalised
subbundle.  Sheaf theory, such as we study in this section, is designed to
answer such questions.  Sheaf theory is a large and complex subject, and here
we only develop in a limited way those facets of the theory that we will use.
There are various relevant references here, including
\cite{MK/PS:90,BRT:76}\@.  The discussion in \cite[Chapter~5]{FWW:83} is also
useful and concise.  The presentation of sheaf cohomology in~\cite{SR:05} is
at a level appropriate for someone with a good background in differential
geometry.  Coherent analytic sheaves are the subject of the book of
\citet{HG/RR:84}\@.

\subsection{Presheaves and sheaves of sets}

Although we will be interested almost exclusively in this paper with sheaves
of rings and modules, it is convenient to first define sheaves of sets.  The
starting point for the definition is that of presheaves.
\begin{definition}
Let $\man{M}$ be a smooth manifold.  A \defn{presheaf of sets} over $\man{M}$
is an assignment to each open set $\nbhd{U}\subset\man{M}$ a set
$F(\nbhd{U})$ and, to each pair of open sets
$\nbhd{V},\nbhd{U}\subset\man{M}$ with $\nbhd{V}\subset\nbhd{U}$\@, a map
$\map{r_{\nbhd{U},\nbhd{V}}}{F(\nbhd{U})}{F(\nbhd{V})}$ called the
\defn{restriction map}\@, with these assignments having the following
properties:
\begin{compactenum}[(i)]
\item $r_{\nbhd{U},\nbhd{U}}$ is the identity map;
\item if\/ $\nbhd{W},\nbhd{V},\nbhd{U}\subset\man{M}$ are open with\/
$\nbhd{W}\subset\nbhd{V}\subset\nbhd{U}$\@, then\/
$r_{\nbhd{U},\nbhd{W}}=r_{\nbhd{V},\nbhd{W}}\scirc r_{\nbhd{U},\nbhd{V}}$\@.
\end{compactenum}
We shall frequently use a single symbol, like $\sF$\@, to refer to a
presheaf, with the understanding that
$\sF=\ifam{F(\nbhd{U})}_{\nbhd{U}\,\textrm{open}}$\@, and that the
restriction maps are understood.  An element $s\in F(\nbhd{U})$ is called a
\defn{section over $\nbhd{U}$} and an element of $F(\man{M})$ is called a
\defn{global section}\@.\oprocend
\end{definition}

Presheaves can be restricted to open sets.
\begin{definition}
Let $\sF=\ifam{F(\nbhd{U})}_{\nbhd{U}\,\textrm{open}}$ be a presheaf of sets
over a smooth or real analytic manifold $\man{M}$\@.  If
$\nbhd{U}\subset\man{M}$ is open, then we denote by $\sF|\nbhd{U}$ the
\defn{restriction} of $\sF$ to $\nbhd{U}$\@, which is the presheaf over
$\nbhd{U}$ whose sections over $\nbhd{V}\subset\nbhd{U}$ are simply
$F(\nbhd{V})$\@.\oprocend
\end{definition}

Let us give the examples of presheaves that will be of interest to us here.
\begin{examples}\label{eg:presheaves}
Let $\man{M}$ be a smooth or real analytic manifold, as is required, let
$\map{\pi}{\man{E}}{\man{M}}$ be a smooth or real analytic vector bundle, as
is required, and let $r\in\integernn\cup\{\infty,\omega\}$\@.
\begin{compactenum}
\item Let us denote by $\sfunc[r]{\man{M}}$ the presheaf over $\man{M}$ for
which the sections over an open subset $\nbhd{U}\subset\man{M}$ is the set
$\func[r]{\nbhd{U}}$ of functions of class $\C^r$ on $\nbhd{U}$\@.  The
restriction maps are the natural restrictions of functions.  This presheaf we
call the \defn{presheaf of functions of class $\C^r$} on $\man{M}$\@.
\item Let us denote by $\ssections[r]{\man{E}}$ the presheaf over $\man{M}$
whose sections over an open subset $\nbhd{U}\subset\man{M}$ is the set
$\sections[r]{\nbhd{U}}$ of sections of $\man{E}|\nbhd{U}$ of class $\C^r$\@.
The restriction maps, again, are the natural restrictions.  This presheaf we
call the \defn{presheaf of sections of $\man{E}$ of class $\C^r$}\@.\oprocend
\end{compactenum}
\end{examples}

The notion of a presheaf has built into it a global character; for example,
the specification of global sections is part of the definition.  The power of
sheaf theory, however, is that it gives a framework for extending local
constructions to global ones.  (A good model to have in mind is using
analytic continuation to patch together locally defined holomorphic functions
to arrive at a globally defined holomorphic function.)  In order to do this
in a self-consistent way, one must place some conditions on the presheaves
one uses.  In this way we arrive at the notion of a sheaf, as in the
following definition.
\begin{definition}
Let $\man{M}$ be a smooth manifold and suppose that we have a presheaf
$\sF=\ifam{F(\nbhd{U})}_{\nbhd{U}\,\textrm{open}}$ of sets with restriction
maps $r_{\nbhd{U},\nbhd{V}}$ for $\nbhd{U},\nbhd{V}\subset\man{M}$ open and
satisfying $\nbhd{V}\subset\nbhd{U}$\@.
\begin{compactenum}[(i)]
\item \label{pl:sheafdef1} The presheaf $\sF$ is \defn{separated} when, if
$\nbhd{U}\subset\man{M}$ is open, if $\ifam{\nbhd{U}_a}_{a\in A}$ is an open
covering of $\nbhd{U}$\@, and if $s,t\in F(\nbhd{U})$ satisfy
$r_{\nbhd{U},\nbhd{U}_a}(s)=r_{\nbhd{U},\nbhd{U}_a}(t)$ for every $a\in A$\@,
then $s=t$\@;
\item \label{pl:sheafdef2} The presheaf $\sF$ has the \defn{gluing property}
when, if $\nbhd{U}\subset\man{M}$ is open, if $\ifam{\nbhd{U}_a}_{a\in A}$ is
an open covering of $\nbhd{U}$\@, and if, for each $a\in A$\@, there exists
$s_a\in F(\nbhd{U}_a)$ with the family $\ifam{s_a}_{a\in A}$ satisfying
\begin{equation*}
r_{\nbhd{U}_{a_1},\nbhd{U}_{a_1}\cap\nbhd{U}_{a_2}}(s_{a_1})=
r_{\nbhd{U}_{a_2},\nbhd{U}_{a_1}\cap\nbhd{U}_{a_2}}(s_{a_2})
\end{equation*}
for each\/ $a_1,a_2\in A$\@, then there exists\/ $s\in F(\nbhd{U})$ such
that\/ $s_a=r_{\nbhd{U},\nbhd{U}_a}(s)$ for each\/ $a\in A$\@.
\item The presheaf $\sF$ is a \defn{sheaf of sets} if it is separated and has
the gluing property.\oprocend
\end{compactenum}
\end{definition}

One fairly easily verifies that the presheaves $\sfunc[r]{\man{M}}$ and
$\ssections[r]{\man{E}}$ from Example~\ref{eg:presheaves} are, in fact,
sheaves.  Examples of presheaves failing to be sheaves arise, of course, by
failing either of the conditions~\eqref{pl:sheafdef1} or~\eqref{pl:sheafdef2}
of the definition.  Presheaves failing condition~\eqref{pl:sheafdef1} do not
often arise in settings such as that in this paper, and we refer the reader
to the references for a discussion of this phenomenon.  However, presheaves
failing to satisfy the gluing property~\eqref{pl:sheafdef2} can arise, and the
failure of a presheaf to satisfy this property is one that one is often
forced to deal with.  As an elementary example of a presheaf failing to
satisfy the gluing property, let $\man{M}$ be a noncompact smooth manifold
and denote by $\sCbdd^\infty(\man{M})$ the presheaf whose sections over an
open subset $\nbhd{U}$ is the set $\Cbdd^\infty(\nbhd{U})$ of bounded smooth
functions on $\nbhd{U}$\@.  Because it is possible to patch together locally
defined bounded smooth functions to arrive at a globally defined unbounded
function (we leave the straightforward construction of a counterexample to
the reader), $\sCbdd^\infty(\man{M})$ is not a sheaf.

\subsection{Germs and \'etale spaces}\label{subsec:etale}

Associated to every presheaf is a topological space that captures the local
behaviour of the presheaf.  To construct this space, the notion of a germ is
essential.  We work with a presheaf $\sF$ of sets on a manifold $\man{M}$\@.
We let $x\in\man{M}$ and let $\sN_x$ be the collection of open subsets of
$\man{M}$ containing $x$\@.  We define an equivalence relation in
$\ifam{F(\nbhd{U})}_{\nbhd{U}\in\sN_x}$ by saying that $s_1\in F(\nbhd{U}_1)$
and $s_2\in F(\nbhd{U}_2)$ are \defn{equivalent} if there exists
$\nbhd{V}\in\sN_x$ such that $\nbhd{V}\subset\nbhd{U}_1$\@,
$\nbhd{V}\subset\nbhd{U}_2$\@, and
$r_{\nbhd{U}_1,\nbhd{V}}(s_1)=r_{\nbhd{U}_2,\nbhd{V}}(s_2)$\@.  The
equivalence class of a section $s\in F(\nbhd{U})$ we denote by
$r_{\nbhd{U},x}(s)$\@, by $[(s,\nbhd{U})]_x$\@, or simply by $[s]_x$ if we
are able to forget about the neighbourhood on which $s$ is defined.  With
this construction, we make the following definition.
\begin{definition}
Let $\man{M}$ be a smooth manifold and let
$\sF=\ifam{F(\nbhd{U})}_{\nbhd{U}\,\textrm{open}}$ be a presheaf of sets over
$\man{M}$\@.  For $x\in\man{M}$\@, the \defn{stalk} of $\sF$ at $x$ is the
set of equivalence classes under the equivalence relation defined above, and
is denoted by $\sF_x$\@.  The equivalence class $r_{\nbhd{U},x}(s)$ of a
section $s\in F(\nbhd{U})$ is called the \defn{germ} of $s$ at
$x$\@.\oprocend
\end{definition}

The germs at $x$ of the presheaves $\sfunc[r]{\man{M}}$ and
$\ssections[r]{\man{E}}$ from Example~\ref{eg:presheaves} are denoted by
$\gfunc[r]{x}{\man{M}}$ and $\gsections[r]{x}{\man{M}}$\@, respectively.

With stalks at hand, we can make another useful construction associated with
a presheaf.
\begin{definition}
Let $\man{M}$ be a smooth manifold and let
$\sF=\ifam{F(\nbhd{U})}_{\nbhd{U}\,\textrm{open}}$ be a presheaf of sets over
$\man{M}$\@.  The \defn{\'etale space} of $\sF$ is the disjoint union of
the stalks of $\sF$\@:
\begin{equation*}
\Et(\sF)=\disjointunion_{x\in\man{M}}\sF_x.
\end{equation*}
The \defn{\'etale topology} on $\Et(\sF)$ is that topology whose basis
consists of subsets of the form
\begin{equation*}
\nbhd{B}(\nbhd{U},s)=\setdef{r_{\nbhd{U},x}(s)}{x\in\nbhd{U}},\qquad
\nbhd{U}\subset\man{M}\,\textrm{open},\ s\in F(\nbhd{U}).
\end{equation*}
By $\map{\pi_{\sF}}{\Et(\sF)}{\man{M}}$ we denote the canonical projection
$\pi_{\sF}(r_{\nbhd{U},x}(s))=x$ which we call the \defn{\'etale
projection}\@.\oprocend
\end{definition}

One verifies the following properties of \'etale spaces, the proofs for
which we refer to the referenced texts.
\begin{proposition}
Let\/ $\man{M}$ be a smooth manifold with\/
$\sF=\ifam{F(\nbhd{U})}_{\nbhd{U}\,\textrm{open}}$ a presheaf of sets over\/
$\man{M}$\@.  The \'etale topology on\/ $\Et(\sF)$ has the following
properties:
\begin{compactenum}[(i)]
\item \label{pl:etale1} the sets\/ $\nbhd{B}(\nbhd{U},s)$\@,\/
$\nbhd{U}\subset\man{M}$ open,\/ $s\in F(\nbhd{U})$\@, form a basis for a
topology;
\item \label{pl:etale2} the projection\/ $\pi_{\sF}$ is a local
homeomorphism,~\ie~about every\/ $[s]_x\in\Et(\sF)$ there exists a
neighbourhood\/ $\nbhd{O}\subset\Et(\sF)$ such that\/ $\pi_{\sF}|\nbhd{O}$ is
a homeomorphism onto its image.
\end{compactenum}
\end{proposition}

The way in which one should think of the \'etale topology is depicted in
Figure~\ref{fig:etale}\@.
\begin{figure}[htbp]
\centering
\includegraphics{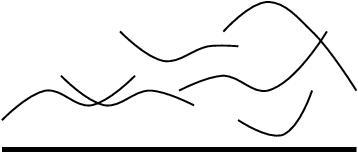}
\caption{How to think of open sets in the \'etale
topology}\label{fig:etale}
\end{figure}%
The point is that open sets in the \'etale topology can be thought of as the
``graphs'' of local sections.  It is a fun exercise to show that the \'etale
topologies for the \'etale spaces of the sheaves $\sfunc[r]{\man{M}}$ and
$\ssections[r]{\man{E}}$ from Example~\ref{eg:presheaves} are Hausdorff if
and only if $r=\omega$\@.

There is a natural notion of a local section of the \'etale space of a
presheaf.
\begin{definition}
Let $\man{M}$ be a smooth manifold and let
$\sF=\ifam{F(\nbhd{U})}_{\nbhd{U}\,\textrm{open}}$ be a presheaf of sets over
$\man{M}$\@.  For $\nbhd{U}\subset\man{M}$ open, a \defn{section} of
$\Et(\sF)$ over $\nbhd{U}$ is a continuous mapping
$\map{\sigma}{\nbhd{U}}{\Et(\sF)}$ with the property that
$\pi_{\sF}\scirc\sigma=\id_{\nbhd{U}}$\@.  The set of sections of $\Et(\sF)$
over $\nbhd{U}$ we denote by $\sections[]{\nbhd{U};\Et(\sF)}$\@.\oprocend
\end{definition}

Note that $\ifam{\sections[]{\nbhd{U};\Et(\sF)}}_{\nbhd{U}\,\textrm{open}}$
is a presheaf if we use the natural restriction maps,~\ie~the set theoretic
restrictions.  This presheaf can be verified to always be a sheaf.  Moreover,
if $\sF$ is itself a sheaf, then there exists a natural isomorphism from
$\sF$ to the presheaf of local sections of $\Et(\sF)$\@.  Explicitly, $s\in
F(\nbhd{U})$ is mapped by this natural isomorphism to the local section
$x\mapsto[s]_x$\@.

The upshot of this section is the following.  A sheaf $\sF$ is in natural
correspondence with the local sections of its \'etale space $\Et(\sF)$\@.
In particular, the attributes of a sheaf $\sF$ are determined by the germs
used in constructing its \'etale space.  Said otherwise, a presheaf that is
a sheaf is determined by its germs.  For this reason, we shall adopt the
usual convention and abandon the distinction between a sheaf and its
\'etale space, and write $\sF$ for both the presheaf and its \'etale
space.

\subsection{Sheaves of rings and modules}

The sheaves in which we are most interested are the sheaf
$\sfunc[r]{\man{M}}$ of functions and the sheaf $\ssections[r]{\man{E}}$ of
sections of a vector bundle.  Just like the set $\sections[r]{\man{E}}$ of
sections is a module over the ring $\func[r]{\man{M}}$\@, the corresponding
sheaves inherit some algebraic structure.  Let us first give a general
definition.
\begin{definition}
\begin{compactenum}[(i)]
\item A \defn{presheaf of rings} over a smooth manifold $\man{M}$ is a
presheaf $\sR=\ifam{R(\nbhd{U})}_{\nbhd{U}\,\textrm{open}}$ whose local
sections are rings and for which the restriction maps
$\map{r_{\nbhd{U},\nbhd{V}}}{R(\nbhd{U})}{R(\nbhd{V})}$\@,
$\nbhd{U},\nbhd{V}\subset\man{M}$ open, $\nbhd{V}\subset\nbhd{U}$\@, are
homomorphisms of rings.
\item If $\sR=\ifam{R(\nbhd{U})}_{\nbhd{U}\,\textrm{open}}$ is a presheaf of
rings over a smooth manifold $\man{M}$\@, a \defn{presheaf of $\sR$-modules}
is a presheaf $\sF=\ifam{F(\nbhd{U})}_{\nbhd{U}\,\textrm{open}}$ of sets such
that $F(\nbhd{U})$ is a module over $R(\nbhd{U})$\@ and such that the
restriction maps $r^{\sR}_{\nbhd{U},\nbhd{V}}$ and
$r^{\sF}_{\nbhd{U},\nbhd{V}}$ satisfy
\begin{align*}
r^{\sF}_{\nbhd{U},\nbhd{V}}(s+t)=r^{\sF}_{\nbhd{U},\nbhd{V}}(s)+
r^{\sF}_{\nbhd{U},\nbhd{V}}(t),\qquad&s,t\in F(\nbhd{U}),\\\eqoprocend
r^{\sF}_{\nbhd{U},\nbhd{V}}(f\,s)=r^{\sR}_{\nbhd{U},\nbhd{V}}(f)
r^{\sF}_{\nbhd{U},\nbhd{V}}(s),\qquad& f\in R(\nbhd{U}),\ s\in F(\nbhd{U}).
\end{align*}
\end{compactenum}
\end{definition}

Of course, if $\map{\pi}{\man{E}}{\man{M}}$ is a smooth or real analytic
vector bundle, as required, and if $r\in\integernn\cup\{\infty,\omega\}$\@,
then $\ssections[r]{\man{E}}$ is a sheaf of modules over the sheaf of rings
$\sfunc[r]{\man{M}}$\@.

\subsection{Morphisms and subsheaves}

Next we introduce maps between presheaves.
\begin{definition}\label{def:morphism}
Let $\man{M}$ be a smooth manifold and let
$\sF=\ifam{F(\nbhd{U})}_{\nbhd{U}\,\textrm{open}}$ and
$\sG=\ifam{G(\nbhd{U})}_{\nbhd{U}\,\textrm{open}}$ be presheaves of sets over
$\man{M}$\@.  A \defn{morphism} of the presheaves $\sF$ and $\sG$ is an
assignment to each open set $\nbhd{U}\subset\man{M}$ a map
$\map{\Phi_{\nbhd{U}}}{F(\nbhd{U})}{G(\nbhd{U})}$ such that the diagram
\begin{equation}\label{eq:presheaf-morphism}
\xymatrix{{F(\nbhd{U})}\ar[r]^{\Phi_{\nbhd{U}}}
\ar[d]_{r_{\nbhd{U},\nbhd{V}}}&
{G(\nbhd{U})}\ar[d]^{r_{\nbhd{U},\nbhd{V}}}\\
{F(\nbhd{V})}\ar[r]_{\Phi_{\nbhd{V}}}&{G(\nbhd{V})}}
\end{equation}
commutes for every open $\nbhd{U},\nbhd{V}\subset\man{M}$ with
$\nbhd{V}\subset\nbhd{U}$\@.  We shall often use the abbreviation
$\Phi=\ifam{\Phi_{\nbhd{U}}}_{\nbhd{U}\,\textrm{open}}$\@.  If $\sF$ and
$\sG$ are sheaves, $\Phi$ is called a \defn{morphism} of sheaves.\oprocend
\end{definition}

Our interest is mainly in, not morphisms between presheaves of sets, but in
morphisms between presheaves of $\sfunc[r]{\man{M}}$-modules.
\begin{definition}
Let $\man{M}$ be a smooth or real analytic manifold as required, let
$r\in\integernn\cup\{\infty,\omega\}$\@, and let
$\sF=\ifam{F(\nbhd{U})}_{\nbhd{U}\,\textrm{open}}$ and
$\sG=\ifam{G(\nbhd{U})}_{\nbhd{U}\,\textrm{open}}$ be sheaves of
$\sfunc[r]{\man{M}}$-modules.  A morphism
$\ifam{\Phi_{\nbhd{U}}}_{\nbhd{U}\,\textrm{open}}$ of the presheaves $\sF$
and $\sG$ is a \defn{morphism of $\sfunc[r]{\man{M}}$-modules} if
$\map{\Phi_{\nbhd{U}}}{F(\nbhd{U})}{G(\nbhd{U})}$ is a homomorphism of
$\func[r]{\nbhd{U}}$ modules for each open set
$\nbhd{U}\subset\man{M}$\@.\oprocend
\end{definition}

Another form of morphism, one that maps from one manifold to another, will
also be useful for us.  In order to state the definition, we need some
notation.  We let $\man{M}$ be a manifold and let
$\sF=\ifam{F(\nbhd{U})}_{\nbhd{U}\,\textrm{open}}$ be a presheaf over
$\man{M}$\@.  Let $A\subset\man{M}$\@.  Let $\nbhd{U},\nbhd{V}\subset\man{M}$
be neighbourhoods of $A$\@.  Sections $s\in F(\nbhd{U})$ and $t\in
F(\nbhd{V})$ are \defn{equivalent} if there exists a neighbourhood
$\nbhd{W}\subset\nbhd{U}\cap\nbhd{V}$ of $A$ such that
$r_{\nbhd{U},\nbhd{W}}(s)=r_{\nbhd{V},\nbhd{W}}(t)$\@.  Let $\sF_A$ denote
the set of equivalence classes under this equivalence relation.  Let us
denote an equivalence class by $[(s,\nbhd{U})]_A$ or by $[s]_A$ if the subset
$\nbhd{U}$ is of no consequence.  Restriction maps can be defined between
such sets of equivalence classes as well.  Thus we let $A,B\subset\man{M}$ be
subsets for which $A\subset B$\@.  If $[(s,\nbhd{U})]_B\in\sF_B$ then, since
$\nbhd{U}$ is also a neighbourhood of $A$\@, $[(s,\nbhd{U})]_B\in\sF_A$\@,
and we denote by $r_{B,A}([(s,\nbhd{U})]_B)$ the equivalence class in
$\sF_A$\@.  One can readily verify that these restriction maps are
well-defined.
\begin{definition}\label{def:di-ii}
Let $r\in\integernn\cup\{\infty,\omega\}$\@, let $\man{M}$ and $\man{N}$ be
smooth or real analytic manifolds, as required, let
$\Phi\in\mappings[r]{\man{M}}{\man{N}}$ be a $\C^r$-map, and let
$\sF=\ifam{F(\nbhd{U})}_{\nbhd{U}\,\textrm{open}}$ be a presheaf of
$\sfunc[r]{\man{M}}$-modules and
$\sG=\ifam{G(\nbhd{V})}_{\nbhd{V}\,\textrm{open}}$ be a presheaf of
$\sfunc[r]{\man{N}}$-modules.
\begin{compactenum}[(i)]
\item The \defn{direct image presheaf} of $\sF$ by $\Phi$ is the presheaf
$\Phi_*\sF$ on $\man{N}$ given by
$\Phi_*\sF(\nbhd{V})=F(\Phi^{-1}(\nbhd{V}))$ for $\nbhd{V}\subset\man{N}$
open.  If $r_{\nbhd{U},\nbhd{V}}$ denote the restriction maps for $\sF$\@,
the restriction maps $\Phi_*r_{\nbhd{U},\nbhd{V}}$ for $\Phi_*\sF$ satisfy,
for $\nbhd{U},\nbhd{V}\subset\man{N}$ open with $\nbhd{V}\subset\nbhd{U}$\@,
\begin{equation*}
\Phi_*r_{\nbhd{U},\nbhd{V}}(s)=
r_{\Phi^{-1}(\nbhd{U}),\Phi^{-1}(\nbhd{V})}(s)
\end{equation*}
for $s\in\Phi_*\sF(\nbhd{U})=F(\Phi^{-1}(\nbhd{U}))$\@.
\item The \defn{inverse image presheaf} of $\sF$ by $\Phi$ is the presheaf
$\Phi^{-1}\sF$ over $\man{M}$ defined by
$\Phi^{-1}\sF(\nbhd{U})=\sF_{\Phi(\nbhd{U})}$\@.  The restriction maps for
$\Phi^{-1}\sF$ are defined by $\Phi^{-1}r_{\nbhd{U},\nbhd{V}}([s])=
r_{\Phi(\nbhd{U}),\Phi(\nbhd{V})}([s])$\@.\oprocend
\end{compactenum}
\end{definition}

If $\sF$ is a sheaf, one readily verifies that $\Phi_*\sF$ is also a sheaf.
From this one readily deduces that, if $\Phi$ is a diffeomorphism, then
$\Phi^{-1}\sG$ is a sheaf if $\sG$ is a sheaf.  The asymmetry in the notation
for the direct and inverse image (\ie~the fact that $\Phi^{-1}$ seems like it
should be $\Phi^*$) is explained by the fact that $\Phi^*$ is the notation
used in a related, but not exactly identical, situation.  We refer
to~\cite[\S7.3]{JLT:02} for details.

We can also talk about subsheaves in a more or less obvious way.
\begin{definition}
Let $\man{M}$ be a smooth manifold, and let
$\sF=\ifam{F(\nbhd{U})}_{\nbhd{U}\,\textrm{open}}$ and
$\sG=\ifam{G(\nbhd{U})}_{\nbhd{U}\,\textrm{open}}$ be presheaves of sets over
$\man{M}$\@.  The presheaf $\sF$ is a \defn{subpresheaf} of $\sG$ if, for
each open set $\nbhd{U}\subset\man{M}$\@, $F(\nbhd{U})$ is a subset of
$G(\nbhd{U})$ and if the inclusion maps
$\map{i_{\sF,\nbhd{U}}}{F(\nbhd{U})}{G(\nbhd{U})}$\@,
$\nbhd{U}\subset\man{M}$ open\@, define a morphism
$i_{\sF}=\ifam{i_{\sF,\nbhd{U}}}_{\nbhd{U}\,\textrm{open}}$ of presheaves of
sets.  If $\sF$ and $\sG$ are sheaves, we say that $\sF$ is a
\defn{presheaf} of $\sG$\@.\oprocend
\end{definition}

Of course, if one replaces ``presheaf of sets'' with ``presheaf of
$\sfunc[r]{\man{M}}$-modules'' in the above definition one arrives at the
notions of a subpresheaf and subsheaf of $\sfunc[r]{\man{M}}$-modules.

\subsection{Coherent real analytic sheaves}\label{subsec:coherent-sheaves}

As we have mentioned several times, one of the uses of sheaf theory is that
it allows one to systematically address global existence questions.  One such
question is the following.  Suppose that we are given a vector bundle
$\map{\pi}{\man{E}}{\man{M}}$\@.  Clearly, about any $x\in\man{M}$ there are
many local sections.  One can legitimately ask whether global sections are
plentiful.  If the vector bundle is smooth, then one can use constructions
with cutoff functions and partitions of unity to construct a section having
any ``reasonable'' property.  (In sheaf language, this is because the sheaf
of smooth sections has a property called ``softness.'')  However, if the
vector bundle is real analytic, the question is not so easy to answer.  To
motivate this a little further, let us recall that a holomorphic vector
bundle over a compact base has few global sections; precisely, the space of
global sections is finite-dimensional over $\complex$\@.  For example, the
dimension of the $\complex$-vector space of sections of the so-called
tautological bundle over complex projective space is
zero~\cite[page~133]{KES/LK/PK/WT:00}\@.  This immediately makes one think an
analogous situation likely holds for real analytic vector bundles.  However,
this is not so, but the reasons for this are not trivial.  In this section we
outline some of the historical developments leading to a few main results
that we shall make use of.

First of all, let us deal with the fact that real analytic manifolds are not
exactly analogous to holomorphic manifolds.  In fact, \cite{HG:58} shows that
real analytic manifolds are analogous to the class of holomorphic manifolds
known as Stein manifolds.  Stein manifolds, unlike general holomorphic
manifolds, possess many holomorphic functions.  For example, about any point
in a Stein manifold, one can find globally defined functions that, in a
neighbourhood of the point, form the components of a holomorphic coordinate
chart.  Thus we have some hope that the question about the plenitude of
global sections of a real analytic vector bundle has an answer unlike that
for holomorphic vector bundles over a compact base.

However, there is still much work to be done.  In the holomorphic case, the
big result here, proved by~\citet{HC:51} and known as ``Cartan's Theorem~A''
(there is also a ``Theorem~B'' which we will get to in time), has as a
consequence that the module of germs at $x$ of sections of a holomorphic
vector bundle over a Stein base is generated by germs of global sections.
In~\cite{HC:57} these holomorphic results are extended to the real analytic
case.  Thus, as a consequence of Cartan's Theorem~A in the real analytic
case, the module $\gsections[\omega]{x}{\man{E}}$ is generated by germs of
global sections.  However, Cartan's results extend far beyond sheaves of
sections of vector bundles to the setting of coherent analytic sheaves.  We
shall actually access these more general results, so in this section we give
the definitions and sketch the results to which we shall subsequently make
reference.

We begin with the notion of a locally finitely generated sheaf of modules.
\begin{definition}\label{def:locfingen}
Let $\man{M}$ be a smooth or real analytic manifold, as required, let
$r\in\integernn\cup\{\infty,\omega\}$\@, and let
$\sF=\ifam{F(\nbhd{U})}_{\nbhd{U}\,\textrm{open}}$ be a sheaf of
$\sfunc[r]{\man{M}}$-modules.  The sheaf $\sF$ is \defn{locally finitely
generated} if, for each $x_0\in\man{M}$\@, there exists a neighbourhood
$\nbhd{U}$ of $x_0$ and sections $s_1,\dots,s_k\in F(\nbhd{U})$ such that
$[s_1]_x,\dots,[s_k]_x$ generate the $\gfunc[r]{x}{\man{M}}$-module $\sF_x$
for every $x\in\nbhd{U}$\@.\oprocend
\end{definition}

Next we turn to the other property required of a coherent sheaf.  We let
$\man{M}$ be a smooth or real analytic manifold, as required, let
$r\in\integernn\cup\{\infty,\omega\}$\@, and let
$\sF=\ifam{F(\nbhd{U})}_{\nbhd{U}\,\textrm{open}}$ be a sheaf of
$\sfunc[r]{\man{M}}$-modules.  Let $\nbhd{U}\subset\man{M}$ be open and let
$s_1,\dots,s_k\in F(\nbhd{U})$\@.  We define a morphism
$\varrho(s_1,\dots,s_k)$ of sheaves from $(\sfunc[r]{\nbhd{U}})^k$ to
$\sF|\nbhd{U}$ by defining it stalkwise:
\begin{equation*}
\varrho(s_1,\dots,s_k)_x([f^1]_x,\dots,[f^k]_x)=
\sum_{j=1}^k[f^j]_x[s_j]_x,\qquad x\in\nbhd{U}.
\end{equation*}
The kernel $\ker(\varrho(s_1,\dots,s_k))$ of this morphism we call the
\defn{sheaf of relations} of the sections $s_1,\dots,s_k$ over $\nbhd{U}$\@.

With the preceding construction, we now make the following definition.
\begin{definition}
Let $\man{M}$ be a smooth or real analytic manifold, as required, let
$r\in\integernn\cup\{\infty,\omega\}$\@, and let
$\sF=\ifam{F(\nbhd{U})}_{\nbhd{U}\,\textrm{open}}$ be a sheaf of
$\sfunc[r]{\man{M}}$-modules.  The sheaf $\sF$ is \defn{coherent}
\begin{compactenum}[(i)]
\item if it is locally finitely generated and
\item if, for every open $\nbhd{U}\subset\man{M}$ and $s_1,\dots,s_k\in
F(\nbhd{U})$\@, $\ker(\varrho(s_1,\dots,s_k))$ is locally finitely
generated.\oprocend
\end{compactenum}
\end{definition}

Now we can define the objects of interest to us.
\begin{definition}
Let $\man{M}$ be a real analytic manifold.  A \defn{coherent real analytic
sheaf} is a coherent sheaf $\sF$ of
$\sfunc[\omega]{\man{M}}$-modules.\oprocend
\end{definition}

We can give an important example of a coherent real analytic sheaf.
\begin{theorem}\label{the:oka}
If\/ $\map{\pi}{\man{E}}{\man{M}}$ is a real analytic vector bundle then\/
$\ssections[\omega]{\man{E}}$ is coherent.
\begin{proof}
Proofs of this result in the holomorphic case can be found in many texts on
several complex variables,~see~\cite[\eg][Theorem~6.4.1]{LH:66}\@.  The
proofs are all lengthy inductive arguments based on the Weierstrass
Preparation Theorem.  The Weierstrass Preparation Theorem in the real
analytic case is given in~\cite[Theorem~6.1.3]{SGK/HRP:02}\@.  With this
version of the theorem, the standard holomorphic proofs of Oka's Theorem
apply to the real analytic case.
\end{proof}
\end{theorem}

Now we are in a position to state an important result concerning global
sections of coherent real analytic sheaves.
\begin{theorem}\label{the:theoremA}
Let\/ $\man{M}$ be a paracompact Hausdorff real analytic manifold and let\/
$\sF$ be a coherent real analytic sheaf.  Then, for\/ $x\in\man{M}$\@, the\/
$\gfunc[\omega]{x}{\man{M}}$-module\/ $\sF_x$ is generated by germs of global
sections of\/ $\sF$\@.
\begin{proof}
The holomorphic case,~\ie~for coherent complex analytic sheaves over Stein
manifolds, was first proved by \citet{HC:51}\@.  Proofs of this result are
often found in texts on several complex
variables~\cite[\eg][Theorem~7.2.8]{LH:66}\@.  The real analytic case we
state here was proved by \citet{HC:57} using the fact that real analytic
manifolds can be, in an appropriate sense, be approximated by Stein
manifolds.  In~\cite{HC:57} the theorems are stated for real analytic
submanifolds of $\real^n$\@.  However, by the real analytic embedding theorem
of \citet{HG:58}\@, this assumption holds for any paracompact Hausdorff real
analytic manifold.
\end{proof}
\end{theorem}

We note that the sheaf $\ssections[r]{\man{E}}$ of $\C^r$-sections of a 
smooth vector bundle $\map{\pi}{\man{E}}{\man{M}}$ is \emph{not} coherent. 
Thus coherence is really an analytic tool, and indeed is the device one uses 
to compensate for the fact that one does not have analytic partitions of 
unity.

\subsection{Real analytic spaces}

The study of analytic spaces is fundamental in the theory of complex analytic
geometry.  As such, references for this theory abound, with~\cite{JLT:02} as
an example.  In the real analytic case, there are fewer references.  The
study of these spaces was really initiated with the work of \citet{HC:57} and
of \citet{HW/FB:59}\@.  Monographs on results in this area are those of
\citet{RN:66} and \citet{FG/PM/AT:86}\@.  Since this theory has some
important repercussions for us, we shall review it here.

We begin with the definition.
\begin{definition}
If $\man{M}$ is a real analytic manifold, a \defn{real analytic space} in
$\man{M}$ is a subset $\man{S}\subset\man{M}$ such that, for each
$x_0\in\man{S}$\@, there is a neighbourhood $\nbhd{U}$ of $x_0$ and
$f_1,\dots,f_k\in\func[\omega]{\nbhd{U}}$ such that
\begin{equation*}
\man{S}\cap\nbhd{U}=\setdef{x\in\nbhd{U}}{f_1(x)=\dots=f_k(x)}.
\end{equation*}
The \defn{ideal sheaf} of a real analytic space $\man{S}$ is the subsheaf
$\sI_{\man{S}}=\ifam{I_{\man{S}}(\nbhd{U})}_{\nbhd{U}\,\textrm{open}}$ of $\sfunc[\omega]{\man{M}}$ defined by
\begin{equation*}\eqoprocend
I_{\man{S}}(\nbhd{U})=\setdef{f\in\func[\omega]{\nbhd{U}}}
{f(x)=0\ \textrm{for all}\ x\in\man{S}\cap\nbhd{U}}.
\end{equation*}
\end{definition}

Characterisations of the ideal sheaf are important in the theory of analytic
spaces.  As an example of a question of interest, one wonders whether, given
an analytic space, its ideal sheaf is generated by globally defined
functions.  One is not surprised to learn that coherence of the ideal sheaf
is important for answering such questions, and indeed \citet{HC:50} proves
that the ideal sheaf of a complex analytic space is coherent.  However, in
the real analytic case, it is no longer true that ideal sheaves are coherent,
as the following example of \citeauthor{HC:57} shows.
\begin{example}
We take $\man{M}=\real^3$ with $(x_1,x_2,x_3)$ the standard coordinates.
Consider the analytic function
\begin{equation*}
\mapdef{f}{\real^3}{\real}{(x_1,x_2,x_3)}{x_3(x_1^2+x_2^2)-x_2^3.}
\end{equation*}
In Figure~\ref{fig:umbrella}
\begin{figure}[htbp]
\centering
\includegraphics[trim={0.45cm 1.4cm 0.25cm 0.6cm},
clip,width=0.5\hsize]{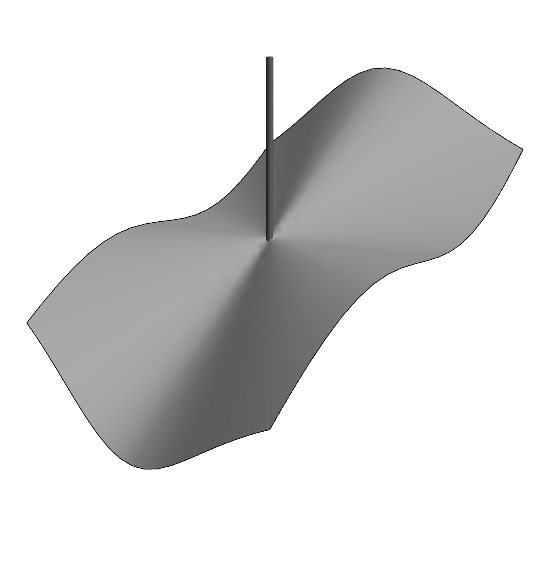}
\caption{Cartan's umbrella}\label{fig:umbrella}
\end{figure}%
we show the $0$-level set of $f$\@, which is thus an analytic space that we
denote by $\man{C}$\@.  We claim that $\sI_{\man{C}}$ is not coherent.  To
show this, we first look at the stalk of $\sI_{\man{C}}$ at
$\vect{0}=(0,0,0)$\@.  We claim that $[f]_{\vect{0}}$ generates this stalk.

To see now that $\sI_{\man{C}}$ is not locally finitely generated, let
$\nbhd{U}$ be a neighbourhood of $\vect{0}$ and let
$\vect{x}_0=(0,0,a)\in\nbhd{U}$ with $a\not=0$\@.  For a sufficiently small
neighbourhood $\nbhd{V}$ of $\vect{x}_0$ (specifically, require that
$\nbhd{V}$ be a ball not containing $\vect{0}$) define
$g\in\func[\omega]{\nbhd{V}}$ by $g(x_1,x_2,x_3)=x_1$\@, and note that
$[(g,\nbhd{V})]_{\vect{x}_0}\in\sI_{\man{C},\vect{x}_0}$\@.  We claim that
$[(g,\nbhd{V})]_{\vect{x}_0}$ is not in the ideal generated by
$[f]_{\vect{x}_0}$\@.  Indeed, note that if $g=hf$ in some neighbourhood of
$\vect{x}_0$\@, then we must have
\begin{equation*}
h(x_1,x_2,x_3)=\frac{x_3(x_1^2+x_2^2)-x_2^3}{x_1},
\end{equation*}
which can easily seen to not be real analytic in any neighbourhood of
$\vect{x}_0$\@.  Finally, by Lemma~\ref{lem:local-generators} we can now
conclude that $\sI_{\man{C}}$ is not locally finitely generated.

A consequence of this is that real analytic functions on $\man{C}$ cannot
generally be extended to real analytic functions away from $\man{C}$\@.
Consider, for example, the function
\begin{equation*}
g(x_1,x_2,x_3)=\frac{x_1}{x_1^2+x_2^2+(x_3-1)^2}.
\end{equation*}
One readily checks that $g$ be extended to a real analytic function in a
neighbourhood of any point, but cannot be extended to a real analytic
function on a neighbourhood of $\man{C}$\@.\oprocend
\end{example}

\subsection{The beginnings of sheaf cohomology}

Sheaf cohomology is a powerful tool for dealing systematically with the
problems concerning the ``local to global'' passage.  The theory has a
reputation for being difficult to learn.  This is in some sense true, but is
also exacerbated by many treatments of the subject which provide a purely
category theory based treatment which is difficult for a beginner to
penetrate.  Here we sketch the beginnings of sheaf cohomology in a fairly
concrete manner, and state a weak form of Cartan's Theorem~B that we shall
subsequently use.  A readable, but not comprehensive, introduction to sheaf
cohomology may be found in the book of \citet{SR:05}\@.  A less readable (by
non-experts) but comprehensive account may be found in the book of
\citet{MK/PS:90}\@.

We let $r\in\integernn\cup\{\infty,\omega\}$\@, let $\man{M}$ be a smooth or
real analytic manifold, as required, and let
$\sF=\ifam{F(\nbhd{U})}_{\nbhd{U}\,\textrm{open}}$ be a sheaf of
$\sfunc[r]{\man{M}}$-modules.  We suppose that we are given an open cover
$\sU=\ifam{\nbhd{U}_a}_{a\in A}$ for $\man{M}$\@, and we let
$\cochain{}{0}(\sU;\sF)$ be the set of all sections over all open sets in
$\sU$\@.  Thus an element of $\cochain{}{0}(\sU;\sF)$ is a family
$\ifam{s_a}_{a\in A}$ where $s_a\in F(\nbhd{U}_a)$\@.  Now let
$\cohomker{}{0}(\sU;\sF)$ be the elements of $\cochain{}{0}(\sU;\sF)$ that
agree on their intersection.  Thus $\ifam{s_a}_{a\in
A}\in\cohomker{}{0}(\sU;\sF)$ if the restrictions of $s_a$ and $s_b$ to
$\nbhd{U}_a\cap\nbhd{U}_b$ agree whenever
$\nbhd{U}_a\cap\nbhd{U}_b\not=\emptyset$\@.  Since $\sF$ is a sheaf, if
$\ifam{s_a}_{a\in A}\in\cohomker{}{0}(\sU;\sF)$ then there exists a unique
$s\in\sections[]{\man{M};\sF}$ such that the restriction of $s$ to
$\nbhd{U}_a$ agrees with $s_a$ for each $a\in A$\@.  Thus we naturally
identify $\cohomker{}{0}(\sU;\sF)$ with $\sections[]{\man{M};\sF}$\@.  Let us
also define $\cohomim{}{0}(\sU;\sF)=0$ by convention.  We take
$\cohom{}{0}(\sU;\sF)=\cohomker{}{0}(\sU;\sF)/\cohomim{}{0}(\sU;\sF)$ so that
$\cohom{}{0}(\sU;\sF)$ is naturally identified with
$\sections[]{\man{M};\sF}$\@.  This is the \defn{zeroth cohomology group} of
$\sF$ for the cover $\sU$\@.

The preceding constructions are related to restricting global sections to
sets from the open cover.  Now we restrict further.  Let
$\cochain{}{1}(\sU;\sF)$ be the set of sections over
$\nbhd{U}_a\cap\nbhd{U}_b$\@, $a,b\in A$\@.  Thus an element of
$\cochain{}{1}(\sU;\sF)$ is a family $\ifam{s_{ab}}_{a,b\in A}$ such that
$s_{ab}\in F(\nbhd{U}_a\cap\nbhd{U}_b)$\@.  Given $\ifam{s_a}_{a\in
A}\in\cochain{}{0}(\sU;\sF)$ we have an induced element
$\ifam{s_{ab}}_{a,b\in A}\in\cochain{}{1}(\sU;\sF)$ defined by
$s_{ab}=s_b-s_a$ (to keep things simple, we omit the restrictions which are
really required here).  Let us denote by $\cohomim{}{1}(\sU;\sF)$ the
elements of $\cochain{}{1}(\sU;\sF)$ obtained in this way.  We define
$\cohomker{}{1}(\sU;\sF)\subset\cochain{}{1}(\sU;\sF)$ by an algebraic
condition that is vacuously satisfied by elements of
$\cohomim{}{1}(\sU;\sF)$\@.  If we take $\cohomker{}{1}(\sU;\sF)$ to be those
elements $\ifam{s_{ab}}_{a,b\in A}\in\cochain{}{1}(\sU;\sF)$ for which
\begin{equation*}
s_{bc}-s_{ac}+s_{ab}=0,\qquad a,b,c\in A,
\end{equation*}
(again, we omit restrictions for brevity), we can see that this condition is
satisfied by elements of $\cohomim{}{1}(\sU;\sF)$\@.  We thus define
$\cohom{}{1}(\sU;\sF)=\cohomker{}{1}(\sU;\sF)/\cohomim{}{1}(\sU;\sF)$\@,
which is the \defn{first cohomology group} of $\sF$ for the cover $\sU$\@.
The vanishing of the first cohomology group is intimately connected with the
capacity of the sheaf to support the patching together of local constructions
to form a global construction.

The path from the preceding rather elementary constructions to higher
cohomology groups now typically proceeds in one of two equivalent directions.
One can continue with the constructions with open covers and prove, that for
suitable open covers, one arrives at a cover-independent theory.  This gives
what is known as \v{C}ech cohomology.  This approach can often be used to
compute the cohomology of a concrete sheaf.  Another approach, more abstract
and so more difficult to understand, realises cohomology groups as ``the
right derived functors for the global section functor.''  About this we shall
say nothing more, but refer to the references.

For us, the following results are useful.  The first result is useful in the
smooth case.
\begin{theorem}\label{the:smooth-cohom}
If\/ $r\in\integernn\cup\{\infty\}$\@, if\/ $\man{M}$ is a smooth paracompact
Hausdorff manifold, if\/ $\sF$ is a sheaf of\/ $\sfunc[r]{\man{M}}$-modules,
and if\/ $\sU=\ifam{\nbhd{U}_a}_{a\in A}$ is an open cover of\/ $\man{M}$\@,
then\/ $\cohom{}{1}(\sU;\sF)=0$\@.
\begin{proof}
The sheaf of rings $\sfunc[r]{\man{M}}$ is easily shown to have the property
of ``fineness''~\cite[Definition~3.3]{ROW:08}\@; this amounts to the fact
that a smooth manifold possesses a $\C^r$-partition of
unity~\cite[Theorem~5.5.7]{RA/JEM/TSR:88}\@.  One then can
show~\cite[Proposition~3.5]{ROW:08} that fine sheaves have the property of
``softness''~\cite[Definition~3.1]{ROW:08}\@.  A sheaf of modules over a soft
sheaf of rings can be shown to be soft~\cite[Lemma~3.16]{ROW:08}\@.  Finally,
the cohomology of soft sheaves may be shown to vanish at orders larger than
zero~\cite[Theorem~3.11]{ROW:08}\@.
\end{proof}
\end{theorem}

In the real analytic case, the following result is the one we shall find
useful.
\begin{theorem}
If\/ $\man{M}$ is a paracompact Hausdorff real analytic manifold, if\/ $\sF$
is a coherent sheaf of\/ $\sfunc[\omega]{\man{M}}$-modules, and if\/
$\sU=\ifam{\nbhd{U}_a}_{a\in A}$ is an open cover of\/ $\man{M}$\@, then\/
$\cohom{}{1}(\sU;\sF)=0$\@.
\begin{proof}
The history of the proof is rather like that for Cartan's Theorem~A given above.
\end{proof}
\end{theorem}

\subsection{Topologies on stalks of sheaves of
sections}\label{subsec:stalk-topology}

We will require topologies on the stalks $\gsections[r]{x}{\man{E}}$ of the
sheaf of sections of a vector bundle $\map{\pi}{\man{E}}{\man{M}}$ of class
$\C^r$\@.  This is done differently in the cases $r=\infty$ and $r=\omega$\@.

\subsubsection{The smooth case}

Let $\map{\pi}{\man{E}}{\man{M}}$ be a smooth vector bundle.  Without loss of
generality (since we are only topologising stalks, so all constructions need
only be local) we suppose that $\man{M}$ is Hausdorff and paracompact.  If
$\nbhd{U}\subset\man{M}$ is open, we recall the weak topology on
$\sections[\infty]{\man{E}|\nbhd{U}}$~\cite{PWM:80}\@.  This is most easily
described by assigning a smooth vector bundle metric $\metric$ to
$\man{E}$\@.  Thus $\metric_x$ is an inner product on
$\man{E}_x$\@.\footnote{The construction of $\metric$ in the smooth case
follows from standard arguments using partitions of unity,~\cf~the proof of
the existence of a Riemannian metric on a smooth, paracompact, Hausdorff
manifold~\cite[Corollary~5.5.13]{RA/JEM/TSR:88}\@.}  We let $\dnorm{\cdot}_x$
denote the induced norm on $\man{E}_x$\@.  Note that the $r$th jet bundle
$\map{\pi^r}{\jet{r}{\man{E}}}{\man{M}}$ of $\map{\pi}{\man{E}}{\man{M}}$ is
a vector bundle~\cite[\S12.17]{IK/PWM/JS:93}\@.  Thus we may define a vector
bundle metric on this bundle that we denote by $\metric^r$\@.  The
corresponding norm on the fibre over $x$ we denote by $\dnorm{\cdot}^r_x$\@.
For $K\subset\nbhd{U}$ compact and for $r\in\integernn$ we can then define a
seminorm $\dnorm{\cdot}_{r,K}$ on $\sections[\infty]{\man{E}|\nbhd{U}}$ by
\begin{equation*}
\dnorm{\xi}_{r,K}=\sup\setdef{\dnorm{j_r\xi(x)}^r_x}{x\in K}.
\end{equation*}
If $\ifam{K_j}_{j\in\integerp}$ is a sequence of compact sets such that
$\nbhd{U}=\cup_{j\in\integerp}K_j$ (by~\cite[Lemma~2.76]{CDA/KCB:06}) then
the locally convex topology defined by the family of seminorms
$\dnorm{\cdot}_{r,K}$\@, $r\in\integernn$\@, $K\subset\nbhd{U}$ compact, is
the same as the locally convex topology defined by the countable family of
seminorms $\dnorm{\cdot}_{r,K_j}$\@, $r\in\integernn$\@, $j\in\integerp$\@.
Moreover, this topology can be easily verified to be Hausdorff and complete.
Thus $\sections[\infty]{\man{E}|\nbhd{U}}$ is a Fr\'echet space with this
topology.  The topology can also be shown to be independent of the choices of
the vector bundle metrics $\metric^r$\@, $r\in\integernn$\@.

Now let $x\in\man{M}$ and let $\sN_x$ be the set of neighbourhoods of $x$\@,
noting that $\sN_x$ is a directed set under inclusion.  Note that
$\gsections[\infty]{x}{\man{E}}$ is the direct limit (in the category of
$\real$-vector spaces) of
$\ifam{\sections[\infty]{\man{E}|\nbhd{U}}}_{\nbhd{U}\in\sN_x}$ with respect
to the mappings $r_{\nbhd{U},x}$\@.  If $\ifam{\nbhd{U}_j}_{j\in\integerp}$
is a sequence of neighbourhoods of $x$ such that
$\nbhd{U}_{j+1}\subset\nbhd{U}_j$ and such that
$\cap_{j\in\integerp}\nbhd{U}_j=\{x\}$\@, then this family is cofinal in
$\sN_x$ and so the resulting direct limit topology on
$\gsections[\infty]{x}{\man{E}}$ induced by the mappings
$r_{\nbhd{U}_j,x}$\@, $j\in\integerp$\@, gives
$\gsections[\infty]{x}{\man{E}}$ the structure of an (LF)-space;
see~\cite[\S19.5]{GK:69}\@.

We will be interested in closed submodules of
$\gsections[\infty]{x}{\man{E}}$\@.  These can be described with the aid of
Whitney's Spectral Theorem~\cite{HW:48}\@.  By $j_\infty\xi(x)$ we denote the
infinite jet of a section $\xi$ at $x$\@.
\begin{theorem}
Let\/ $\map{\pi}{\man{E}}{\man{M}}$ be a smooth vector bundle with bounded
fibre dimension and with\/ $\man{M}$ smooth, second countable, and
Hausdorff.  If\/ $\sM\subset\sections[\infty]{\man{E}}$ is a submodule, then
the closure of\/ $\sM$ in the weak topology on\/ $\sections[\infty]{\man{E}}$
is
\begin{equation*}
\closure(\sM)=\setdef{\xi\in\sections[\infty]{\man{E}}}
{j_\infty\xi(x)\in\setdef{j_\infty\eta(x)}{\eta\in\sM}\ \textrm{for each}\
x\in\man{M}}.
\end{equation*}
\end{theorem}

With this, we can characterise closed submodules of the stalks
$\gsections[\infty]{x}{\man{E}}$\@.
\begin{proposition}\label{prop:closed-smooth}
Let\/ $\map{\pi}{\man{E}}{\man{M}}$ be a smooth vector bundle and let\/
$x\in\man{M}$\@.  If\/
$\sF_x\subset\gsections[\infty]{x}{\man{E}}$ is a submodule, then the closure
of\/ $\sF_x$ in the (LF)-topology on\/ $\gsections[\infty]{x}{\man{E}}$ is
\begin{multline*}
\closure(\sF_x)=\{[(\xi,\nbhd{U})]_x\in\gsections[\infty]{x}{\man{E}}|\enspace
\textrm{there exists a neighbourhood}\ \nbhd{V}\subset\nbhd{U}\
\textrm{of}\ x\ \textrm{such that}\\
j_\infty\xi(y)\in\setdef{j_\infty\eta(y)}{\eta\in r^{-1}_{\nbhd{V},x}(\sF_x)}\
\textrm{for every}\ y\in\nbhd{V}\}.
\end{multline*}
\begin{proof}
First let $[(\xi,\nbhd{U})]_x\in\gsections[\infty]{x}{\man{E}}$ be such that
there exists a neighbourhood $\nbhd{V}\subset\nbhd{U}$ of $x$ for which
\begin{equation*}
j_\infty\xi(y)\in\setdef{j_\infty\eta(y)}{\eta\in r_{\nbhd{V},x}^{-1}(\sF_x)}
\end{equation*}
for every $y\in\nbhd{V}$\@.  Let $\ifam{\nbhd{U}_j}_{j\in\integerp}$ be a
sequence of neighbourhoods of $x$ such that $\nbhd{U}_{j+1}\subset\nbhd{U}_j$
and such that $\cap_{j\in\integerp}\nbhd{U}_j=\{x\}$\@, and note that
$r_{\nbhd{U}_j,x}^{-1}(\sF_x)$ is a submodule of
$\sections[\infty]{\man{E}|\nbhd{U}_j}$ for each $j\in\integerp$\@.  Let
$N\in\integerp$ be sufficiently large that $\nbhd{U}_j\subset\nbhd{V}$ for
all $j\ge N$\@.  Then
\begin{equation*}
j_\infty\xi(y)\in\setdef{j_\infty\eta(y)}{\eta\in r_{\nbhd{U}_j,x}^{-1}(\sF_x)}
\end{equation*}
for all $y\in\nbhd{U}_j$ and $j\ge N$\@.  Thus
$r_{\nbhd{U},\nbhd{U}_j}(\xi)\in\closure(r_{\nbhd{U}_j,x}^{-1}(\sF_j))$ for
$j\ge N$ by the Whitney Spectral Theorem.  It then follows that, for every
$j\ge N$\@,
\begin{equation*}
r_{\nbhd{U}_j,x}(r_{\nbhd{U},\nbhd{U}_j}(\xi))=r_{\nbhd{U},x}(\xi)=[(\xi,\nbhd{U})]_x
\in\closure(\sF_x),
\end{equation*}
\cf~\cite[\S19.5]{GK:69}\@.

Next let $[(\xi,\nbhd{U})]_x$ be in the closure of $\sF_x$\@.  It follows
that $r_{\nbhd{U},x}^{-1}([(\xi,\nbhd{U})]_x)$ is in the closure of
$r_{\nbhd{U},x}^{-1}(\sF_x)$\@.  Therefore, by the Whitney Spectral Theorem,
\begin{equation*}
j_\infty\xi(y)\in\setdef{j_\infty\eta(y)}{\eta\in r^{-1}_{\nbhd{U},x}(\sF_x)}
\end{equation*}
for all $y\in\nbhd{U}$\@, giving the desired characterisation of
$\closure(\sF_x)$\@.
\end{proof}
\end{proposition}

The following example shows that there can be finitely generated submodules
of germs of sections that are not closed.  We refer to~\cite{BGR:70} for
further discussion along these lines.
\begin{example}
We take $\man{M}=\real$ with coordinate $x$\@.  We take
$f\in\func[\infty]{\real}$ to be defined by
\begin{equation*}
f(x)=\begin{cases}\eul^{-1/x^4},&x\not=0,\\0,&x=0.\end{cases}
\end{equation*}
We claim that the submodule $\modgen{f}_{\nbhd{U}}$ of
$\func[\infty]{\nbhd{U}}$ generated by $f|\nbhd{U}$ is not closed for any
neighbourhood $\nbhd{U}$ of $0$\@.  Indeed, let $\nbhd{U}$ be a neighbourhood
of $0$\@.  Note that the function $g\in\func[\infty]{\nbhd{U}}$ defined by
\begin{equation*}
g(x)=\begin{cases}\eul^{-1/x^2},&x\not=0,\\0,&x=0\end{cases}
\end{equation*}
has the property that its Taylor series at $0$ agrees with that of $f$\@:
both Taylor series are identically zero.  However, since
\begin{equation*}
\lim_{x\to0}\frac{g(x)}{f(x)}=\infty,
\end{equation*}
there is no function $h\in\func[\infty]{\nbhd{U}}$ such that $g=hf$\@, and so
$g\not\in\modgen{f}_{\nbhd{U}}$\@.  Moreover, by
Proposition~\ref{prop:closed-smooth}\@, the argument shows that the module
$\modgen{[f]_0}$ generated by $[f]_0$ is not closed.\oprocend
\end{example}

\subsubsection{The real analytic case}

The topology on the stalks of the sheaf of sections of a real analytic vector
bundle is more difficult to describe than the smooth case.  In the real
analytic case, we must first extend real analytic objects to holomorphic
objects on a complexification of the vector bundle.

We let $\map{\pi}{\man{E}}{\man{M}}$ be a real analytic vector bundle,
supposing that $\man{M}$ is paracompact and Hausdorff.  As in the smooth
case, our definition of the appropriate topologies is facilitated by the
introduction of a vector bundle metric on $\man{E}$\@.  Let us be sure we
understand how to construct such a metric in the real analytic case.  First
of all, by the real analytic embedding theorem of \citet{HG:58}\@, we
analytically embed $\man{E}$ into $\real^N$ for sufficiently large $N$\@.
Then $\man{M}$ and the fibres $\man{E}_x$ of $\man{E}$ are real analytic
submanifolds of $\real^N$\@.  Moreover, $\man{E}$ is naturally isomorphic to
the normal bundle of $\man{M}$ in $\man{E}$\@.  Using the Euclidean inner
product $\metric_x$ on the fibres of the normal bundle at $x\in\man{M}$\@, we
define a real analytic vector bundle metric $\metric$ on $\man{E}$\@.
Following \citet{HW/FB:59}\@, we can regard $\man{M}$ as the real part of a
corresponding holomorphic manifold $\ol{\man{M}}$\@.  Moreover, using our
observation above that $\man{E}$ is isomorphic to the normal bundle of
$\man{M}\subset\man{E}\subset\real^N$\@, we can extend $\man{E}$ to be a
holomorphic vector bundle $\map{\ol{\pi}}{\ol{\man{E}}}{\ol{\man{M}}}$\@.
The restriction of $\ol{\man{E}}$ to $\man{M}\subset\ol{\man{M}}$ agrees with
the complexification $\man{E}\otimes_{\real}\complex$ of $\man{E}$\@, and we
denote this by $\man{E}^{\complex}$\@.  That is,
\begin{equation*}
\man{E}^{\complex}\eqdef\ol{\man{E}}|\man{M}\simeq\man{E}\otimes_{\real}\complex.
\end{equation*}
By shrinking the holomorphic extension $\ol{\man{M}}$ if necessary, we can
suppose that $\metric$ extends to a Hermitian metric $\ol{\metric}$ on the
fibres of $\ol{\man{E}}$\@.  Denote by $\dnorm{\cdot}_{\ol{\metric}}$ the
norm induced on the fibres of $\ol{\man{E}}$\@.

With these holomorphic extensions at hand, we can now begin to describe the
topology on the stalks of $\ssections[\omega]{\man{E}}$\@.  As in the smooth
case, we first describe the topology of $\sections[\omega]{\man{E}|\nbhd{U}}$
for open sets $\nbhd{U}\subset\man{M}$\@.  The construction of this topology
is done in a few steps.
\begin{compactenum}
\item \emph{Topologise the holomorphic sections:} We first consider the
complexification.  Thus let $\ol{\nbhd{U}}\subset\ol{\man{M}}$ be open and
let $\sections[\textup{hol}]{\ol{\man{E}}|\ol{\nbhd{U}}}$ denote the
holomorphic sections of $\ol{\man{E}}|\ol{\nbhd{U}}$\@.  Define a topology on
$\sections[0]{\ol{\man{E}}|\ol{\nbhd{U}}}$ as that defined by the family of
seminorms $\dnorm{\cdot}_K$\@, $K\subset\ol{\nbhd{U}}$ compact, given by
\begin{equation*}
\dnorm{\ol{\xi}}_K=\sup\setdef{\dnorm{\ol{\xi}(z)}_{\ol{\metric}}}{z\in K}.
\end{equation*}
This defines the compact-open, or weak $\C^0$-, topology on
$\sections[0]{\ol{\man{E}}|\ol{\nbhd{U}}}$\@.  It is well-known that
$\sections[\textup{hol}]{\ol{\man{E}}|\ol{\nbhd{U}}}$ is a \emph{closed}
subspace of $\sections[0]{\ol{\man{E}}|\ol{\nbhd{U}}}$ with this
topology~\cite[Theorem~V.B.5]{RCG/HR:65}\@.  Thus
$\sections[\textup{hol}]{\ol{\man{E}}|\ol{\nbhd{U}}}$ has a natural Fr\'echet
topology.
\item \emph{Restrict from holomorphic to real analytic sections of
$\man{E}^{\complex}$:} Now let $\nbhd{U}\subset\man{M}$ be open and let
$\ol{\nbhd{U}}\subset\ol{\man{M}}$ be a neighbourhood of $\nbhd{U}$\@.  If
$\ol{\xi}\in\sections[\textup{hol}]{\ol{\man{E}}|\ol{\nbhd{U}}}$ then
$\ol{\xi}|\nbhd{U}\in\sections[\omega]{\man{E}^{\complex}|\nbhd{U}}$\@,~\cf\
\cite[Corollary~2.3.7]{SGK:92}\@.  Conversely, if $\nbhd{U}\subset\man{M}$ is
open and if $\xi\in\sections[\omega]{\man{E}^{\complex}|\nbhd{U}}$\@, then
there exists a neighbourhood $\ol{\nbhd{U}}$ of $\nbhd{U}$ and
$\ol{\xi}\in\sections[\textup{hol}]{\ol{\man{E}}|\ol{\nbhd{U}}}$ such that
$\ol{\xi}|\nbhd{U}=\xi$\@.  (To see this, it is easiest to think in terms of
Taylor series.  The Taylor series for $\xi$ will be a series with coordinates
$(x^1,\dots,x^n)$ for $\man{M}$ as indeterminates and with complex
coefficients.  The Taylor series for the corresponding section $\ol{\xi}$
will be given by replacing the real indeterminates $(x^1,\dots,x^n)$ with
complex indeterminates $(z^1,\dots,z^n)$ representing coordinates for
$\ol{\man{M}}$\@.)  Moreover, $\ol{\xi}$ is unique in that any two such
extensions will agree on any connected neighbourhood of $\nbhd{U}$\@.  For
$\nbhd{U}\subset\man{M}$ open and $\ol{\nbhd{U}}$ a neighbourhood of
$\nbhd{U}$ in $\ol{\man{M}}$\@, let us denote by
\begin{equation}\label{eq:analytic-rho}
\map{\ol{\rho}_{\ol{\nbhd{U}},\nbhd{U}}}
{\sections[\textup{hol}]{\ol{\man{E}}|\ol{\nbhd{U}}}}
{\sections[\omega]{\man{E}^{\complex}|\nbhd{U}}}
\end{equation}
the restriction map, noting that this map is injective if $\ol{\nbhd{U}}$ is
connected.  Moreover, if $\nbhd{U}\subset\man{M}$ is open and if
$\sN_{\nbhd{U}}$ denotes the set of neighbourhoods of $\nbhd{U}$ in
$\ol{\man{M}}$\@, then
\begin{equation*}
\sections[\omega]{\man{E}^{\complex}|\nbhd{U}}=
\cup_{\ol{\nbhd{U}}\in\sN_{\nbhd{U}}}\image(\ol{\rho}_{\ol{\nbhd{U}},\nbhd{U}}).
\end{equation*}
\item \label{enum:real-part} \emph{Projection to real analytic sections of
$\man{E}$\@:} We let
$\map{\Re_{\nbhd{U}}}{\man{E}^{\complex}|\nbhd{U}}{\man{E}|\nbhd{U}}$ be the
projection onto the real part of the fibres.  Let us also abuse notation
slightly and let
\begin{equation*}
\map{\Re_{\nbhd{U}}}{\sections[\omega]{\man{E}^{\complex}|\nbhd{U}}}
{\sections[\omega]{\man{E}|\nbhd{U}}}
\end{equation*}
be the induced map on sections.  We then define
\begin{equation*}
\map{\rho_{\ol{\nbhd{U}},\nbhd{U}}\eqdef\Re_{\nbhd{U}}\scirc
\ol{\rho}_{\ol{\nbhd{U}},\nbhd{U}}}{\sections[\textup{hol}]
{\ol{\man{E}}|\ol{\nbhd{U}}}}
{\sections[\omega]{\man{E}|\nbhd{U}}}.
\end{equation*}
Let us see how this homomorphism reacts with the module structures on the
domain and codomain.  We let $\nbhd{U}\subset\man{M}$ be open and let
$\ol{\nbhd{U}}\in\sN_{\nbhd{U}}$\@.  Making a slight abuse of notation, we
denote by
\begin{equation*}
\map{\rho_{\ol{\nbhd{U}},\nbhd{U}}}{\func[\textup{hol}]{\ol{\nbhd{U}}}}
{\func[\omega]{\nbhd{U}}}
\end{equation*}
the restriction map, where $\func[\textup{hol}]{\ol{\man{E}}|\ol{\nbhd{U}}}$
denotes the holomorphic functions on $\ol{\nbhd{U}}$\@.  We claim that if
$\sM\subset\sections[\omega]{\man{E}|\nbhd{U}}$ is a submodule over
$\func[\omega]{\nbhd{U}}$ then $\rho_{\ol{\nbhd{U}},\nbhd{U}}^{-1}(\sM)$ is a
submodule over
$\rho_{\ol{\nbhd{U}},\nbhd{U}}^{-1}(\func[\omega]{\nbhd{U}})$\@.  First of
all, let us show that
$\rho_{\ol{\nbhd{U}},\nbhd{U}}^{-1}(\func[\omega]{\nbhd{U}})$ is a ring.  If
$\ol{f}$ and $\ol{g}$ are such that $\rho_{\ol{\nbhd{U}},\nbhd{U}}(\ol{f})=f$
and $\rho_{\ol{\nbhd{U}},\nbhd{U}}(\ol{g})=g$ for
$f,g\in\func[\omega]{\nbhd{U}}$\@, then clearly
\begin{equation*}
\rho_{\ol{\nbhd{U}},\nbhd{U}}(\ol{f}+\ol{g})=f+g,\quad
\rho_{\ol{\nbhd{U}},\nbhd{U}}(\ol{f}\ol{g})=fg,
\end{equation*}
and so $\ol{f}+\ol{g},\ol{f}\ol{g}\in
\rho^{-1}_{\ol{\nbhd{U}},\nbhd{U}}(\func[\omega]{\nbhd{U}})$\@.  Now let
$\ol{\xi}$ and $\ol{\eta}$ be such that
$\rho_{\ol{\nbhd{U}},\nbhd{U}}(\ol{\xi})=\xi$ and
$\rho_{\ol{\nbhd{U}},\nbhd{U}}(\ol{\eta})=\eta$ for $\xi,\eta\in\sM$\@.  It
is clear that
\begin{equation*}
\rho_{\ol{\nbhd{U}},\nbhd{U}}(\ol{\xi}+\ol{\eta})=\xi+\eta
\end{equation*}
and so $\ol{\xi}+\ol{\eta}\in\rho_{\ol{\nbhd{U}},\nbhd{U}}^{-1}(\sM)$\@.  Also,
since the restriction of the product $\ol{f}\,\ol{\xi}$ agrees with the
products of the restrictions of $\ol{f}$ and $\ol{\xi}$\@, we have
$\rho_{\ol{\nbhd{U}},\nbhd{U}}(\ol{f}\ol{\xi})=f\xi$ and so
$\ol{f}\,\ol{\xi}\in\rho_{\ol{\nbhd{U}},\nbhd{U}}^{-1}(\sM)$\@, giving our claim.
\item \emph{Use the inductive limit topology:} Now let
$\nbhd{U}\subset\man{M}$ be open and let $\sN_{\nbhd{U}}$ be the set of
neighbourhoods of $\nbhd{U}$ in $\ol{\man{M}}$\@.  Note that $\sN_{\nbhd{U}}$
is a directed set under inclusion.  The topology on
$\sections[\omega]{\man{E}^{\complex}|\nbhd{U}}$ is the inductive limit
topology with respect to the family of restriction
mappings~\eqref{eq:analytic-rho}\@,~\ie~the finest topology for which all of
these maps are continuous.  Thus
$\sections[\omega]{\man{E}^{\complex}|\nbhd{U}}$ is an inductive limit of
Fr\'echet spaces, but it need not be an (LF)-space since this limit need not
be countable.  One verifies that $\sections[\omega]{\man{E}|\nbhd{U}}$ is a
closed subspace of $\sections[\omega]{\man{E}^{\complex}|\nbhd{U}}$\@.  Thus
the induced topology gives us the desired topology on
$\sections[\omega]{\man{E}|\nbhd{U}}$\@.  This topology has some not so
friendly properties; we refer to the work of \citet{AM:66} and the
discussions of this in~\cite[\S2.6]{SGK/HRP:02} and~\cite{PD/DV:00a} for
details.
\end{compactenum}

Let us consider some related constructions with stalks, as this notation will
be of use to us in the proof of Theorem~\ref{the:analytic-closed-module}
below.  With our notation above, we have the sheaf
$\ssections[\textup{hol}]{\ol{\man{E}}}$ of
$\sfunc[\textup{hol}]{\ol{\man{M}}}$-modules.  We also have the sheaves
$\ssections[\omega]{\man{E}}$ and $\ssections[\omega]{\man{E}^{\complex}}$ of
$\sfunc[\omega]{\man{M}}$-modules.  Our constructions above ensure that
restriction of germs to $\man{M}\subset\ol{\man{M}}$ gives a bijection from
$\gsections[\textup{hol}]{\ol{x}}{\ol{\man{E}}}$ to
$\gsections[\omega]{x}{\man{E}^{\complex}}$ for each $x\in\man{M}$\@, where
$\ol{x}$ denotes the image of $x\in\man{M}$ in $\ol{\man{M}}$\@.  We thus
have a homomorphism
\begin{equation*}
\map{\ol{\rho}_{\ol{x},x}}{\ssections[\textup{hol}]{\ol{\man{E}}}}
{\ssections[\omega]{\man{E}^{\complex}}}
\end{equation*}
of $\real$-vector spaces (topologies will be considered shortly).  By taking
real parts, we further get a homomorphism
\begin{equation*}
\map{\rho_{\ol{x},x}}{\ssections[\textup{hol}]{\ol{\man{E}}}}
{\ssections[\omega]{\man{E}}}.
\end{equation*}
We may argue as in step~\ref{enum:real-part} above that, if
$\sF_x\subset\gsections[\omega]{x}{\man{E}}$ is a
$\gfunc[\omega]{x}{\man{M}}$-submodule, then $\rho_{\ol{x},x}^{-1}(\sF_x)$ is
a $\gfunc[\textup{hol}]{\ol{x}}{\ol{\man{M}}}$-submodule of
$\gsections[\textup{hol}]{\ol{x}}{\ol{\man{E}}}$\@.

Finally, we topologise the stalks of $\ssections[\textup{hol}]{\ol{\man{E}}}$
and $\ssections[\omega]{\man{E}}$\@.  Let $x\in\man{M}$ and let
$\ol{x}\in\ol{\man{M}}$\@.  Let $\sN_x$ and $\sN_{\ol{x}}$ denote the
families of neighbourhoods of $x$ and $\ol{x}$ in $\man{M}$ and
$\ol{\man{M}}$\@, respectively.  Note that these are both directed sets under
inclusion.  The topologies on $\gsections[\omega]{x}{\man{E}}$ and
$\gsections[\textup{hol}]{\ol{x}}{\ol{\man{E}}}$ are then the inductive limit
topologies with respect to the families of mappings $r_{\nbhd{U},x}$\@,
$\nbhd{U}\in\sN_x$\@, and $r_{\ol{\nbhd{U}},\ol{x}}$\@,
$\ol{\nbhd{U}}\in\sN_{\ol{x}}$\@, respectively.

Next we need to describe the closed submodules of
$\gsections[\omega]{x}{\man{E}}$\@.
\begin{theorem}\label{the:analytic-closed-module}
Let\/ $\map{\pi}{\man{E}}{\man{M}}$ be a real analytic vector bundle and
let\/ $x\in\man{M}$\@.  If\/ $\sF_x\subset\gsections[\omega]{x}{\man{E}}$ is
a submodule, then it is closed in the inductive limit topology on\/
$\gsections[\omega]{x}{\man{E}}$\@.
\begin{proof}
We first claim that the map\/ $\rho_{\ol{x},x}$ is continuous.  By
Proposition~2 from Section~4.1 of~\cite{AG:73} it suffices to show that, for
any neighbourhood $\ol{\nbhd{U}}$ of $\ol{x}$ in $\ol{\man{M}}$\@, the
composition
\begin{equation*}
\map{\rho_{\ol{x},x}\scirc r_{\ol{\nbhd{U}},\ol{x}}}
{\sections[\textup{hol}]{\ol{\man{E}}|\ol{\nbhd{U}}}}
{\gsections[\omega]{x}{\man{E}}}
\end{equation*}
is continuous.  Note that the diagram
\begin{equation}\label{eq:analytic-stalk}
\xymatrix{{\sections[\textup{hol}]{\ol{\man{E}}|\ol{\nbhd{U}}}}
\ar[r]^{\rho_{\ol{\nbhd{U}},\nbhd{U}}}\ar[d]_{r_{\ol{\nbhd{U}},\ol{x}}}&
{\sections[\omega]{\man{E}|\nbhd{U}}}\ar[d]^{r_{\nbhd{U},x}}\\
{\gsections[\textup{hol}]{\ol{x}}{\ol{\man{E}}}}\ar[r]_{\rho_{\ol{x},x}}&
{\gsections[\omega]{x}{\man{E}}}}
\end{equation}
commutes for every $\ol{\nbhd{U}}\in\sN_{\ol{x}}$\@, where
$\nbhd{U}=\ol{\nbhd{U}}\cap\man{M}$\@.  Therefore, it suffices to show that 
\begin{equation*}
\map{r_{\nbhd{U},x}\scirc\rho_{\ol{\nbhd{U}},\nbhd{U}}}
{\sections[\textup{hol}]{\ol{\man{E}}|\ol{\nbhd{U}}}}
{\gsections[\omega]{x}{\man{E}}}
\end{equation*}
is continuous for every $\nbhd{U}\in\sN_x$ and
$\ol{\nbhd{U}}\in\sN_{\nbhd{U}}$\@.  However, again by Proposition~2 from
Section~4.1 of~\cite{AG:73}\@, the homomorphisms
$\rho_{\ol{\nbhd{U}},\nbhd{U}}$ and $r_{\nbhd{U},x}$ are continuous, and so
the claim follows.

Now let $[(\xi,\nbhd{U})]_x$ be in the closure of $\sF_x$ in
$\gsections[\omega]{x}{\man{E}}$\@.  Let $\ol{\nbhd{U}}\in\sN_{\nbhd{U}}$ be
such that $\xi$ extends to a section $\ol{\xi}$ of
$\ol{\man{E}}|\ol{\nbhd{U}}$\@.  Then
$r^{-1}_{\nbhd{U},x}([(\xi,\nbhd{U})]_x)$ is in the closure of
$r^{-1}_{\nbhd{U},x}(\sF_x)$ in $\sections[\omega]{\man{E}|\nbhd{U}}$ and
$\rho^{-1}_{\ol{\nbhd{U}},\nbhd{U}}(r^{-1}_{\nbhd{U},x}([(\xi,\nbhd{U})]_x))$
is in the closure of
$\rho^{-1}_{\ol{\nbhd{U}},\nbhd{U}}(r^{-1}_{\nbhd{U},x}(\sF_x))$ in
$\sections[\textup{hol}]{\ol{\man{E}}|\ol{\nbhd{U}}}$\@.  By the
commutativity of the diagram~\eqref{eq:analytic-stalk}\@, this implies that
$r^{-1}_{\ol{\nbhd{U}},\ol{x}}([(\ol{\xi},\ol{\nbhd{U}})]_{\ol{x}})$ is in
the closure of $r_{\ol{\nbhd{U}},\nbhd{U}}^{-1}(\ol{\sF}_{\ol{x}})$ in
$\osections[\textup{hol}]{\ol{\man{E}}|\ol{\nbhd{U}}}$\@.  Consequently,
$r^{-1}_{\ol{\nbhd{U}},\ol{x}}([(\ol{\xi},\ol{\nbhd{U}})]_{\ol{x}})$ is in
the closure of $r_{\ol{\nbhd{U}},\nbhd{U}}^{-1}(\rho_{\ol{x},x}^{-1}(\sF_x))$
in $\sections[\textup{hol}]{\ol{\man{E}}|\ol{\nbhd{U}}}$\@.  As we argued
before the statement of the theorem, $\rho_{\ol{x},x}^{-1}(\sF_x)$ is a
$\gfunc[\textup{hol}]{\ol{x}}{\ol{\man{M}}}$-submodule of
$\gsections[\textup{hol}]{\ol{x}}{\ol{\man{E}}}$\@.  Therefore,
by~\cite[Proposition~11.2.2]{JLT:02}\@,
$r_{\ol{\nbhd{U}},\nbhd{U}}^{-1}(\rho_{\ol{x},x}^{-1}(\sF_x))$ is closed in
$\sections[\textup{hol}]{\ol{\man{E}}|\ol{\nbhd{U}}}$ and so contains
$r^{-1}_{\ol{\nbhd{U}},\ol{x}}([(\ol{\xi},\ol{\nbhd{U}})]_{\ol{x}})$\@.
Consequently,
$\rho^{-1}_{\ol{\nbhd{U}},\nbhd{U}}(r^{-1}_{\nbhd{U},x}(\sF_x))$ contains
$\rho^{-1}_{\ol{\nbhd{U}},\nbhd{U}}(r^{-1}_{\nbhd{U},x}([(\xi,\nbhd{U})]_x))$\@,
and so $[(\xi,\nbhd{U})]_x\in\sF_x$\@, as desired.
\end{proof}
\end{theorem}

Note that, unsurprisingly given the topology on
$\gsections[\omega]{x}{\man{E}}$\@, the proof relies in an essential way on
the holomorphic analogue of the theorem.  This holomorphic analogue is an
essential ingredient in the proofs of Cartan's Theorems~A and~B.

\section{Generalised subbundles of vector bundles}\label{sec:vector-bundles}

In this section we introduce the major player in this paper in the general
setting of vector bundles; later in the paper we shall specialise to tangent
bundles.  We also introduce the connections between these constructions and
subsheaves of the sheaf of sections of a vector bundle.

\subsection{Generalised subbundles}

We begin by giving the definitions we shall use throughout the paper.
\begin{definition}\label{def:subbundle}
Let $\map{\pi}{\man{E}}{\man{M}}$ be a vector bundle of class $\C^\infty$ or
$\C^\omega$\@, as is required.  A \defn{generalised subbundle} of $\man{E}$
is a subset $\dist{F}\subset\man{E}$ such that, for each $x\in\man{M}$\@, the
subset $\dist{F}_x=\dist{F}\cap\man{E}_x$ is a subspace (and so, in
particular, is nonempty).  The subspace $\dist{F}_x$ is the \defn{fibre} of
$\dist{F}$ at $x$\@.  Associated with the notion of a generalised subbundle
we have the following.
\begin{compactenum}[(i)]
\item A generalised subbundle $\dist{F}$ is of \defn{class $\C^r$}\@,
$r\in\integernn\cup\{\infty,\omega\}$\@, if, for each $x_0\in\man{M}$\@,
there exists a neighbourhood $\nbhd{N}$ of $x_0$ and a family
$\ifam{\xi_j}_{j\in J}$ of $\C^r$-sections, called \defn{local generators}\@,
of $\man{E}|\nbhd{N}$ such that
\begin{equation*}
\dist{F}_x=\vecspan[\real]{\xi_j(x)\mid\enspace j\in J}
\end{equation*}
for each $x\in\nbhd{N}$\@.
\item A generalised subbundle $\dist{F}$ of class $\C^r$\@,
$r\in\integernn\cup\{\infty,\omega\}$\@, is \defn{locally finitely generated}
if, for each $x_0\in\man{M}$\@, there exists a neighbourhood $\nbhd{N}$ of
$x_0$ and a family $\ifam{\xi_1,\dots,\xi_k}$ of $\C^r$-sections, called
\defn{local generators}\@, of $\man{E}|\nbhd{N}$ such that
\begin{equation*}
\dist{F}_x=\vecspan[\real]{\xi_1(x),\dots,\xi_k(x)}
\end{equation*}
for each $x\in\nbhd{N}$\@.
\item A generalised subbundle $\dist{F}$ of class $\C^r$\@,
$r\in\integernn\cup\{\infty,\omega\}$\@, is \defn{finitely generated} if
there exists a family $\ifam{\xi_1,\dots,\xi_k}$ of $\C^r$-sections, called
\defn{generators}\@, of $\man{E}$ such that
\begin{equation*}
\dist{F}_x=\vecspan[\real]{\xi_1(x),\dots,\xi_k(x)}
\end{equation*}
for each $x\in\man{M}$\@.
\end{compactenum}
The nonnegative integer $\dim(\dist{F}_x)$ is called the \defn{rank} of
$\dist{F}$ at $x$ and is sometimes denoted by $\rank(\dist{F}_x)$\@.\oprocend
\end{definition}

Let us also give a few related standard definitions that we shall use.
\begin{definition}
Let $\map{\pi}{\man{E}}{\man{M}}$ be a vector bundle of class $\C^\infty$ or
$\C^\omega$\@, as is required, let $r\in\integernn\cup\{\infty,\omega\}$\@,
and let $\dist{F}\subset\man{E}$ be a $\C^r$-generalised subbundle.  If
$\nbhd{U}\subset\man{M}$ is open, the \defn{restriction} of $\dist{F}$ to
$\nbhd{U}$ is
\begin{equation*}\eqoprocend
\dist{F}|\nbhd{U}=\cup_{x\in\nbhd{U}}\dist{F}_x.
\end{equation*}
\end{definition}

\begin{definition}
Let $\map{\pi}{\man{E}}{\man{M}}$ be a vector bundle of class $\C^\infty$ or
$\C^\omega$\@, as is required, let $r\in\integernn\cup\{\infty,\omega\}$\@,
and let $\dist{F}\subset\man{E}$ be a $\C^r$-generalised subbundle.  If
$\nbhd{U}\subset\man{M}$ is open, a \defn{local section} of $\dist{F}$ over
$\nbhd{U}$ is a section $\map{\xi}{\nbhd{U}}{\man{E}}$ such that
$\xi(x)\in\dist{F}_x$ for every $x\in\nbhd{U}$\@.  A local section $\xi$ of
$\man{F}$ is of \defn{class $\C^k$}\@, $k\le r$\@, if it is of class $\C^k$
as a local section of $\man{E}$\@.  The set of local sections of $\dist{F}$
over $\nbhd{U}$ of class $\C^k$ is denoted by
$\sections[k]{\dist{F}|\nbhd{U}}$\@, or simply by $\sections[k]{\dist{F}}$
when $\nbhd{U}=\man{M}$\@.\oprocend
\end{definition}

Of particular interest to us in the paper are generalised subbundles of
tangent bundles.  These we give a special name.
\begin{definition}
Let $\man{M}$ be a smooth or real analytic manifold, as required.  A
\defn{distribution} on $\man{M}$ is a generalised subbundle $\dist{D}$ of
$\tb{\man{M}}$\@, and a distribution $\dist{D}$ is of class $\C^r$\@,
$r\in\integernn\cup\{\infty,\omega\}$\@, if it is of class $\C^r$ as a
generalised subbundle.\oprocend
\end{definition}

For the first few sections of the paper we shall focus on generalised
subbundles of vector bundles, turning especially to distributions in
Section~\ref{sec:differential}\@.

\subsection{Regular and singular points of generalised subbundles}

One of the complications ensuing from the notion of a generalised subbundle
arises if the dimensions of the subspaces $\dist{F}_x$\@, $x\in\man{M}$\@,
are not locally constant.  The following definition associates some language
with this.
\begin{definition}
Let $\dist{F}$ be a generalised subbundle of a vector bundle
$\map{\pi}{\man{E}}{\man{M}}$\@.  A point $x_0\in\man{M}$
\begin{compactenum}[(i)]
\item is a \defn{regular point} for $\dist{F}$ if there exists a
neighbourhood $\nbhd{N}$ of $x_0$ such that
$\rank(\dist{F}_x)=\rank(\dist{F}_{x_0})$ for every $x\in\nbhd{N}$ and
\item is a \defn{singular point} for $\dist{F}$ if it is not a regular point
for $\dist{F}$\@.
\end{compactenum}
A generalised subbundle $\dist{F}$ is \defn{regular} if every point in
$\man{M}$ is a regular point for $\dist{F}$\@, and is \defn{singular}
otherwise.\oprocend
\end{definition}

Although our definition of regular and singular points is made for arbitrary
generalised subbundles, these definitions only have real value in the case
when the generalised subbundle has some smoothness.  A regular generalised
subbundle of class $\C^r$\@, $r\in\integernn\cup\{\infty,\omega\}$\@, is
often called \defn{subbundle} of $\man{E}$ of class $\C^r$\@.

For continuous generalised subbundles one can make some statements about the
character of rank and the character of the set of regular and singular
points.  In the following result, if $\dist{F}$ is a generalised subbundle of
a vector bundle $\map{\pi}{\man{E}}{\man{M}}$\@, then we denote by
$\map{\rank_{\dist{F}}}{\man{M}}{\integernn}$ the function defined by
$\rank_{\dist{F}}(x)=\rank(\dist{F}_x)$\@.
\begin{proposition}\label{prop:rank-semicont}
If\/ $\dist{F}$ is a generalised subbundle of class\/ $\C^0$ of a vector
bundle\/ $\map{\pi}{\man{E}}{\man{M}}$\@, then the function\/
$\rank_{\dist{F}}$ is lower semicontinuous and the set of regular points of\/
$\dist{F}$ is open and dense.
\begin{proof}
Let $a\in\real$ and let $x_0\in\rank_{\dist{F}}^{-1}(\interval(a,\infty))$\@.
Thus $k\eqdef\rank_{\dist{F}}(x_0)>a$\@.  This means that there are $k$
sections $\xi_1,\dots,\xi_k$ of class $\C^0$ defined in a neighbourhood
$\nbhd{N}$ of $x_0$ such that
$\dist{F}_{x_0}=\vecspan[\real]{\xi_1(x_0),\dots,\xi_k(x_0)}$\@.  Now choose
a vector bundle chart $(\nbhd{V},\psi)$ for $\man{E}$ about $x_0$ so that the
local sections $\xi_1,\dots,\xi_k$ have local representatives
\begin{equation*}
\vect{x}\mapsto(\vect{x},\vect{\xi}_j(\vect{x})),\qquad
j\in\{1,\dots,k\}.
\end{equation*}
Let $(\nbhd{U},\phi)$ be the induced chart for $\man{M}$ and let
$\vect{x}_0=\phi(x_0)$\@.  The vectors
$\ifam{\vect{\xi}_1(\vect{x}_0),\dots,\vect{\xi}_k(\vect{x}_0)}$ are then
linearly independent.  Therefore, there exist $j_1,\dots,j_k\in\{1,\dots,n\}$
(supposing that $n$ is the dimension of $\man{M}$) such that the matrix
\begin{equation*}
\begin{bmatrix}\vect{\xi}^{j_1}_1(\vect{x}_0)&\cdots&
\vect{\xi}^{j_1}_k(\vect{x}_0)\\
\vdots&\ddots&\vdots\\\vect{\xi}^{j_k}_1(\vect{x}_0)&\cdots&
\vect{\xi}^{j_k}_k(\vect{x}_0)\end{bmatrix}
\end{equation*}
has nonzero determinant, where $\vect{\xi}^{j_l}_i(\vect{x}_0)$ is the
$j_l$th component of $\vect{\xi}_i$\@, $i,l\in\{1,\dots,k\}$\@.  By
continuity of the determinant there exists a neighbourhood $\nbhd{U}'$ of
$\vect{x}_0$ such that the matrix
\begin{equation*}
\begin{bmatrix}\vect{\xi}^{j_1}_1(\vect{x})&\cdots&
\vect{\xi}^{j_1}_k(\vect{x})\\
\vdots&\ddots&\vdots\\\vect{\xi}^{j_k}_1(\vect{x})&\cdots&
\vect{\xi}^{j_k}_k(\vect{x})\end{bmatrix}
\end{equation*}
has nonzero determinant for every $\vect{x}\in\nbhd{U}'$\@.  Thus the vectors
$\ifam{\vect{\xi}_1(\vect{x}),\dots,\vect{\xi}_k(\vect{x})}$ are linearly
independent for every $\vect{x}\in\nbhd{U}'$\@.  Therefore, the local
sections $\xi_1,\dots,\xi_k$ are linearly independent on
$\phi^{-1}(\nbhd{U}')$\@, and so
$\phi^{-1}(\nbhd{U}')\subset\rank_{\dist{F}}^{-1}(\interval(a,\infty))$ which
gives lower semicontinuity of $\rank_{\dist{F}}$\@.

Let us denote by $R_{\dist{F}}$ the set of regular points of $\dist{F}$ and
let $x_0\in R_{\dist{F}}$\@.  Then, by definition of $R_{\dist{F}}$\@, there
exists a neighbourhood $\nbhd{U}$ of $x_0$ such that $\nbhd{U}\subset
R_{\dist{F}}$\@.  Thus $R_{\dist{F}}$ is open.  Now let $x_0\in\man{M}$ and
let $\nbhd{U}$ be a connected neighbourhood of $x_0$\@.  Since the function
$\rank_{\dist{F}}$ is locally bounded, there exists a least integer $N$ such
that $\rank_{\dist{F}}(x)\le N$ for each $x\in\nbhd{U}$\@.  Moreover, since
$\rank_{\dist{F}}$ is integer-valued, there exists $x'\in\nbhd{U}$ such that
$\rank_{\dist{F}}(x')=N$\@.  Now, by lower semicontinuity of
$\rank_{\dist{F}}$\@, there exists a neighbourhood $\nbhd{U}'$ of $x'$ such
that $\rank_{\dist{F}}(x)\ge N$ for all $x\in\nbhd{U}'$\@.  By definition of
$N$ we also have $\rank_{\dist{F}}(x)\le N$ for each $x\in\nbhd{U}'$\@.  Thus
$x'\in R_{\dist{F}}$\@, and so $x_0\in\closure(R_{\dist{F}})$\@.  Therefore,
$R_{\dist{F}}$ is dense.
\end{proof}
\end{proposition}

For real analytic generalised subbundles, one can say much more about the set
of regular points.
\begin{proposition}\label{prop:analytic-rank}
If\/ $\dist{F}$ is a generalised subbundle of class\/ $\C^\omega$ of a real
analytic vector bundle\/ $\map{\pi}{\man{E}}{\man{M}}$ over a paracompact
analytic base\/ $\man{M}$\@, then the following statements hold:
\begin{compactenum}[(i)]
\item \label{pl:analytic-rank1} the set of singular points for\/ $\dist{F}$
is a locally analytic set,~\ie~for each\/ $x\in\man{M}$ there exists a
neighbourhood\/ $\nbhd{U}$ of\/ $x$ and real analytic functions\/
$f_1,\dots,f_k\in\func[\omega]{\nbhd{U}}$ such that the set of singular
points of\/ $\dist{F}$ in\/ $\nbhd{U}$ is given by\/
$\cap_{j=1}^kf_j^{-1}(0)$\@;
\item \label{pl:analytic-rank2} if\/ $x_1,x_2\in\man{M}$ are regular points
for\/ $\dist{F}$ in the same connected component of\/ $\man{M}$\@, then\/
$\rank_{\dist{F}}(x_1)=\rank_{\dist{F}}(x_2)$\@.
\end{compactenum}
\begin{proof}
\eqref{pl:analytic-rank1} Let $x_0\in\man{M}$\@.  By
Corollary~\ref{cor:subbunfingen} we assume that $(\nbhd{U},\phi)$ is a
coordinate chart about $x_0$ with $\nbhd{U}$ connected and that we have real
analytic sections $\xi_1,\dots,\xi_k$ of $\dist{E}|\nbhd{U}$ such that
\begin{equation*}
\dist{F}_x=\vecspan[\real]{\xi_1(x),\dots,\xi_k(x)},\qquad x\in\nbhd{U}.
\end{equation*}
Denote by $\vect{x}\mapsto(\vect{x},\vect{\xi}_j(\vect{x}))$\@,
$j\in\{1,\dots,k\}$\@, the local representatives of the sections
$\xi_1,\dots,\xi_k$\@.  Define $\map{\mat{\xi}}{\nbhd{U}}{\real^{n\times k}}$
by
\begin{equation*}
\mat{\xi}(x)=\left[\begin{array}{c|c|c}\vect{\xi}_1(x)&\cdots&
\vect{\xi}_k(x)\end{array}\right].
\end{equation*}
Let us denote
\begin{equation*}
\grank(\mat{\xi})=\max\setdef{\rank(\mat{\xi}(x))}{x\in\nbhd{U}}.
\end{equation*}
For $m\in\{1,\dots,\grank(\mat{\xi})\}$ define
\begin{equation*}
\nbhd{U}_m=\setdef{x\in\nbhd{U}}{\rank(\mat{\xi}(x))<m}.
\end{equation*}
Fix $m\in\{1,\dots,\grank(\mat{\xi})\}$\@.  We claim that if
$\nbhd{U}_m\not=\emptyset$ then it is closed with empty interior.  Indeed,
$x\in\nbhd{U}_m$ if and only if the determinants of all $m\times m$
submatrices of $\mat{\xi}(x)$ vanish.  Thus $\nbhd{U}_m$ is analytic.  Note
that the set of points where the determinant of a fixed $m\times m$ submatrix
vanishes is closed, being the preimage of $0\in\real$ under the real analytic
(and so continuous) determinant function.  Thus $\nbhd{U}_m$ is the
intersection of a finite collection of closed sets, and so is closed.
Suppose that $\nbhd{U}_m$ has a nonempty interior.  By the Identity Theorem
(this is proved in the holomorphic case in~\cite[Theorem~A.3]{RCG:90a}\@, and
the same proof applies to the real analytic case) and connectedness of
$\nbhd{U}$ it follows that all the determinants of all $m\times m$
submatrices of $\mat{\xi}$ vanish on $\nbhd{U}$\@.  This contradicts the fact
that $m<\grank(\mat{\xi})$ and the definition of $\grank(\mat{\xi})$\@.

To complete this part of the theorem, we claim that $x\in\nbhd{U}$ is a
singular point for $\dist{F}|\nbhd{U}$ if and only if
$x\in\nbhd{U}_{\grank(\mat{\xi})}$\@.  The assertion is trivial if
$\nbhd{U}_{\grank(\mat{\xi})}=\emptyset$\@, so we suppose otherwise.  If
$x\in\nbhd{U}_{\grank(\mat{\xi})}$ then the fact that
$\nbhd{U}_{\grank(\mat{\xi})}$ has empty interior ensures that $x$ is not a
regular point for $\dist{F}|\nbhd{U}$\@.  Conversely, if
$x\not\in\nbhd{U}_{\grank(\mat{\xi})}$ then
$\rank(\dist{F}_x)=\grank(\mat{\xi})$ and by
Proposition~\ref{prop:rank-semicont} we conclude that $x$ is a regular point
for $\dist{F}|\nbhd{U}$\@.

\eqref{pl:analytic-rank2} We suppose that $\man{M}$ is connected so,
consequently, $\man{M}$ has a well-defined dimension.  By our above
constructions, for each $x\in\man{M}$ let $(\nbhd{U}_x,\phi_x)$ be a chart
for $\man{M}$ about $x$ and define
\begin{equation*}
\grank_{\dist{F}}(x)=\max\setdef{\rank(\dist{F}_y)}{y\in\nbhd{U}_x}.
\end{equation*}
Let $x_1,x_2\in\man{M}$ and since $\man{M}$ is path connected
by~\cite[Proposition~1.6.7]{RA/JEM/TSR:88}\@, let
$\map{\gamma}{\interval[0,1]}{\man{M}}$ be a continuous curve for which
$\gamma(0)=x_1$ and $\gamma(1)=x_2$\@.  Define
$\map{\phi}{\interval[0,1]}{\integernn}$ by
\begin{equation*}
\phi(s)=\grank_{\dist{F}}(\gamma(s)).
\end{equation*}
We claim that $\phi$ is locally constant.  Indeed, let
$s_0\in\interval[0,1]$ and let $\delta\in\realp$ be such that
$\gamma(s)\in\nbhd{U}_{\gamma(s_0)}$ for
$s\in\interval[0,1]\cap\interval({s_0-\delta},{s_0+\delta})$\@.  From the
first part of the proof, the sets
\begin{equation*}
\setdef{x\in\nbhd{U}_{\gamma(s_0)}}{\rank(\dist{F}_x)<\phi(s_0)}
\end{equation*}
and
\begin{equation*}
\setdef{x\in\nbhd{U}_{\gamma(s)}}{\rank(\dist{F}_x)<\phi(s)}
\end{equation*}
are closed subsets with empty interior in $\nbhd{U}_{\gamma(s_0)}$ and
$\nbhd{U}_{\gamma(s)}$\@, respectively.  It follows that the sets
\begin{equation*}
\setdef{x\in\nbhd{U}_{\gamma(s)}\cap\nbhd{U}_{\gamma(s_0)}}
{\rank(\dist{F}_x)<\phi(s_0)}
\end{equation*}
and
\begin{equation*}
\setdef{x\in\nbhd{U}_{\gamma(s)}\cap\nbhd{U}_{\gamma(s_0)}}
{\rank(\dist{F}_x)<\phi(s)}
\end{equation*}
are closed with empty interior in
$\nbhd{U}_{\gamma(s)}\cap\nbhd{U}_{\gamma(s_0)}$\@.  Therefore, for
\begin{multline*}
x\in\nbhd{U}_{\gamma(s)}\cap\nbhd{U}_{\gamma(s_0)}\setminus
\left(\setdef{x\in\nbhd{U}_{\gamma(s)}\cap\nbhd{U}_{\gamma(s_0)}}
{\rank(\dist{F}_x)<\phi(s_0)}\right.\\\left.\cup
\setdef{x\in\nbhd{U}_{\gamma(s)}\cap\nbhd{U}_{\gamma(s_0)}}
{\rank(\dist{F}_x)<\phi(s)}\right),
\end{multline*}
$\phi(s)=\rank(\dist{F}_x)=\phi(s_0)$\@.  Thus $\phi$ is indeed locally
constant, and so constant since $\interval[0,1]$ is connected.
\end{proof}
\end{proposition}

The notions of regularity and singularity bear on the character of local
generators for a generalised subbundle.
\begin{proposition}\label{prop:local-generators}
Let\/ $\man{M}$ be a manifold of class\/ $\C^\infty$ or\/ $\C^\omega$\@, as
is required, let\/ $r\in\integernn\cup\{\infty,\omega\}$\@, and let\/
$\dist{F}$ be a\/ $\C^r$-generalised subbundle on\/ $\man{M}$\@.  Then, for
each\/ $x_0\in\man{M}$\@, there exists a neighbourhood\/ $\nbhd{N}$ of\/
$x_0$ and local generators\/ $\ifam{\xi_1,\dots,\xi_k}$ for\/ $\dist{F}$ on\/
$\nbhd{N}$ with the following properties:
\begin{compactenum}[(i)]
\item $\ifam{\xi_1(x_0),\dots,\xi_m(x_0)}$ is a basis for\/
$\dist{F}_{x_0}$\@;
\item $\xi_{m+1}(x_0)=\dots=\xi_k(x_0)=0_{x_0}$\@.
\end{compactenum}
In particular, if\/ $x_0$ is a regular point for\/ $\dist{F}$\@, the
sections\/ $\ifam{\xi_1,\dots,\xi_m}$ are local generators for\/ $\dist{F}$
in some neighbourhood (possibly smaller than\/ $\nbhd{N}$) of\/ $x_0$\@.
\begin{proof}
Let $\ifam{\eta_1,\dots,\eta_k}$ be local generators for $\dist{F}$ defined
on a neighbourhood $\nbhd{N}$ of $x_0$\@.  We can assume there are finitely
many of these by Theorem~\ref{the:global-span} in the case when
$r\in\integernn\cup\{\infty\}$ and by Corollary~\ref{cor:subbunfingen} in
case $r=\omega$\@.  We may rearrange the sections
$\ifam{\eta_1,\dots,\eta_k}$ so that $\ifam{\eta_1(x_0),\dots,\eta_m(x_0)}$
forms a basis for $\dist{F}_{x_0}$\@.  We then let
$\ifam{\vect{v}_{m+1},\dots,\vect{v}_k}\subset\real^k$ be a basis for the
kernel of the linear map $\map{L_{x_0}}{\real^k}{\dist{F}_{x_0}}$ defined by
\begin{equation*}
L_{x_0}(\vect{v})=\sum_{j=1}^kv_j\eta_j(x_0).
\end{equation*}
Define an invertible $k\times k$ matrix $\mat{R}$ by
\begin{equation*}
\mat{R}=\left[\begin{array}{c|c|c|c|c|c}
\vect{e}_1&\cdots&\vect{e}_m&\vect{v}_{m+1}&\cdots&\vect{v}_k
\end{array}\right],
\end{equation*}
where $\vect{e}_j\in\real^k$\@, $j\in\{1,\dots,k\}$\@, is the $j$th standard
basis vector, and define $\xi_1,\dots,\xi_k$ by
\begin{equation*}
\xi_j=\sum_{l=1}^kR_{lj}\eta_l,\qquad j\in\{1,\dots,k\}.
\end{equation*}
Then $\xi_j=\eta_j$ for $j\in\{1,\dots,k\}$\@, and
$\xi_j(x_0)=L_{x_0}(\vect{v}_j)=0_{x_0}$ for $j\in\{m+1,\dots,k\}$\@.  This
gives the first conclusion of the proposition.

For the second, if $x_0$ is a regular point for $\dist{F}$\@, then let
$\dist{F}'$ be the generalised subbundle on $\nbhd{N}$ generated by
$\ifam{\xi_1,\dots,\xi_m}$\@.  Then, for $x$ in a neighbourhood $\nbhd{N}'$
of $x_0$\@, $\rank_{\dist{F}'}(x)=m$ by lower semicontinuity of
$\rank_{\dist{F}'}$\@.  Since $\dist{F}'\subset\dist{F}$ it follows that
$\dist{F}'|\nbhd{N}=\dist{F}|\nbhd{N}$\@.
\end{proof}
\end{proposition}

\subsection{Generalised subbundles and subsheaves of sections}

One of the features of this paper is that the r\^ole of sheaves in
distribution theory is made explicit and we illustrate how the theory leads
to nontrivial fundamental results, particularly for real analytic
distributions.  In this section we make connections with sheaf theory to our
definitions for generalised subbundles in the preceding section.  We shall
see that there is a relationship between generalised subbundles of
$\map{\pi}{\man{E}}{\man{M}}$ and subsheaves of $\ssections[r]{\man{E}}$\@,
but the relationship is not always a perfect correspondence.

Let us first consider the natural subsheaf arising from a generalised
subbundle.
\begin{definition}
Let $\map{\pi}{\man{E}}{\man{M}}$ be a smooth or real analytic vector bundle,
as is required, let $r\in\integernn\cup\{\infty,\omega\}$\@, and let
$\dist{F}\subset\man{E}$ be a generalised subbundle of class $\C^r$\@.  The
\defn{sheaf of sections} of $\dist{F}$ is the sheaf $\ssections[r]{\dist{F}}$
whose local sections over the open set $\nbhd{U}\subset\man{M}$ is the set of
sections of $\dist{F}|\nbhd{U}$\@.\oprocend
\end{definition}

Clearly $\ssections[r]{\dist{F}}$ is a subsheaf of
$\sfunc[r]{\man{M}}$-modules of $\ssections[r]{\man{E}}$\@.  Note that a
$\C^r$-generalised subbundle $\dist{F}$, by definition, is constructed using
local generators.  It is clear that each local generator is a local section
of the sheaf $\ssections[r]{\dist{F}}$\@.  What is not generally true is that
the local generators for the distribution generate the stalks of the sheaf
$\ssections[r]{\dist{F}}$\@, as the following example shows.
\begin{example}\label{eg:!stalk-generate}
We take $\man{M}=\real$ and $\man{E}=\real\times\real$ with the vector bundle
projection $\map{\pi}{\man{E}}{\man{M}}$ being projection onto the first
factor: $\pi(x,v)=x$\@.  We let $r\in\integerp\cup\{\infty,\omega\}$\@.  For
$x\in\real$ we consider the neighbourhood $\nbhd{N}_x=\real$ of $x$ and on
$\nbhd{N}_x$ we take the local (in fact, global in this case) section
$\xi_x(y)=(y,y^2)$\@.  Thus, in this case, the neighbourhood of $x$ and the
local generators defined on this neighbourhood are the same for each $x$\@.
This cannot be expected to be the usual situation, but in this simple example
it turns out to be possible.  In any case, this data then defines a
$\C^r$-generalised subbundle $\dist{F}$ by
\begin{equation*}
\dist{F}_x=\begin{cases}\{x\}\times\real,&x\not=0,\\
\{x\}\times\{0\},&x=0.\end{cases}
\end{equation*}
Associated to this generalised subbundle we have the subsheaf
$\ssections[r]{\dist{F}}$ of sections of $\dist{F}$\@.  We can describe this
subsheaf explicitly.  Indeed, we claim that if $\nbhd{U}\subset\real$ is
open, then any section of $\dist{F}|\nbhd{U}$ is a
$\func[r]{\nbhd{U}}$-multiple of the section $\xi\colon x\mapsto(x,x)$\@.
This is clear if $0\not\in\nbhd{U}$ since $\xi$ is then a nonvanishing
section of a one-dimensional vector bundle.  So we need only consider the
case when $0\in\nbhd{U}$\@.  In this case, $\nbhd{U}$\@, being an open subset
of $\real$\@, is a union of open intervals, one of which contains $0$\@.  On
the other open intervals any section is clearly a multiple of $\xi$\@, again
since $\xi$ is a nonvanishing section of a one-dimensional vector bundle.
Thus we need only show that, if $a<0<b$\@, then any $\C^r$-section $\eta$ of
$\man{E}|\interval(a,b)$ is a multiple of $\xi$\@.  In the following
calculation, we identify $\eta$ with a $\real$-valued function on
$\interval(a,b)$\@.  For $x\in\interval(a,b)$ compute
\begin{equation*}
\eta(x)=\int_0^x\eta'(t)\,\d{t}=x\int_0^1\eta'(xs)\,\d{s}.
\end{equation*}
The function $\bar{\eta}\colon x\mapsto\int_0^1\eta'(xs)\,\d{s}$ is of class
$\C^r$ and so $\eta=\id_{\interval(a,b)}\cdot\bar{\eta}$\@.  Thus every
$\C^r$-section $\eta$ on $\interval(a,b)$ vanishing at $0$ is a product of
$\id_{\interval(a,b)}$ and a function of class
$\C^r$\@,~\ie~$\eta(x)=x\bar{\eta}(x)$\@.  This characterises local sections
of $\ssections[r]{\dist{F}}$ over $\nbhd{U}$ as multiples of $\xi$\@, as
claimed.

Note, however, that the generator $\xi_0$ for $\dist{F}$ about $x=0$ does not
generate the stalk $\gsections[r]{0}{\dist{F}}$ since, for example, the germ
of the local section $y\mapsto(y,y)$ is not a
$\gfunc[r]{0}{\man{M}}$-multiple of $[\xi_0]_0$\@.\oprocend
\end{example}

By our above considerations, we can proceed from generalised subbundles to
subsheaves, keeping in mind the caution that the local generators for the
generalised subbundles are not necessarily generators for stalks of the
sheaf.  We can also proceed from subsheaves to generalised subbundles of
$\man{E}$\@.
\begin{definition}\label{def:sheaf->subbundle}
Let $\map{\pi}{\man{E}}{\man{M}}$ be a smooth or real analytic vector bundle,
as is required, let $r\in\integernn\cup\{\infty,\omega\}$\@, and let $\sF$ be
a subsheaf of $\sfunc[r]{\man{M}}$-modules of $\ssections[r]{\man{E}}$\@.
The \defn{generalised subbundle generated by $\sF$} is defined by
\begin{equation*}
\dist{F}(\sF)_x=\vecspan[\real]{\xi(x)|\enspace[\xi]_x\in\sF_x}.
\end{equation*}
If\/ $\sF\subset\sections[r]{\tb{\man{M}}}$ is a subsheaf of vector fields,
we will use the notation $\dist{D}(\sF)$ to denote the \defn{distribution
generated by $\sF$}\@,~\ie
\begin{equation*}\eqoprocend
\dist{D}(\sF)_x=\vecspan[\real]{X(x)|\enspace[X]_x\in\sF_x}.
\end{equation*}
\end{definition}

There are two obvious questions related to this inclusion that come to mind
at this point.
\begin{compactenum}
\item Is every generalised subbundle generated by some subsheaf of
$\ssections[r]{\man{E}}$ of class $\C^r$\@?
\item It is clear that $\sF\subset\ssections[r]{\dist{F}(\sF)}$\@.  If some
generalised subbundle $\dist{F}$ is generated by a subsheaf $\sF$\@, is it
true that $\sF=\ssections[r]{\dist{F}}$\@?
\end{compactenum}
Let us address these questions in order.

The answer to the first question is, ``Yes,'' in the smooth case.  In fact,
in the smooth case one can say much more.
\begin{proposition}\label{prop:subsheaf->subbundle}
Let\/ $r\in\integernn\cup\{\infty\}$ and let\/ $\map{\pi}{\man{E}}{\man{M}}$
be a smooth vector bundle with\/ $\man{M}$ paracompact and Hausdorff.  If\/
$\sF=\ifam{F(\nbhd{U})}_{\nbhd{U}\,\textrm{open}}$ is a subsheaf of\/
$\ssections[r]{\man{E}}$ then the map
\begin{equation*}
F(\man{M})\ni\xi\mapsto[\xi]_x\in\sF_x
\end{equation*}
is surjective for each\/ $x\in\man{M}$\@.  In particular, the generalised
subbundle\/ $\dist{F}(\sF)$ generated by\/ $\sF$ is of class\/ $\C^r$\@.
\begin{proof}
Let $x_0\in\man{M}$ and let $[(\xi,\nbhd{U})]_{x_0}\in\sF_{x_0}$\@.  Let
$\nbhd{V}\subset\nbhd{W}$ be relatively compact neighbourhoods of $x_0$ such
that
\begin{equation*}
\closure(\nbhd{V})\subset\nbhd{W}\subset\closure(\nbhd{W})\subset\nbhd{U},
\end{equation*}
and, by the Tietze Extension
Theorem~\cite[Proposition~5.5.8]{RA/JEM/TSR:88}\@, let
$f\in\func[\infty]{\man{M}}$ be such that $f$ takes the value $1$ on
$\nbhd{V}$ and vanishes on the complement of $\nbhd{W}$\@.  We then define
\begin{equation*}
\ol{\xi}(x)=\begin{cases}f(x)\xi(x),&x\in\nbhd{U},\\
0,&\textrm{otherwise},\end{cases}
\end{equation*}
We claim that $\ol{\xi}\in F(\man{M})$\@.  Let $x\in\man{M}$\@.  If
$x\in\closure(\nbhd{W})$ then let $\nbhd{U}_x\subset\nbhd{U}$ be a
neighbourhood of $x$ and let $\xi_x=\ol{\xi}|\nbhd{U}_x$\@, noting that
$\xi_x\in F(\nbhd{U}_x)$\@.  If $x\in\man{M}\setminus\closure(\nbhd{W})$ then
let $\nbhd{U}_x$ be a neighbourhood of $x$ such that
$\nbhd{U}_x\subset\man{M}\setminus\closure(\nbhd{W})$ and let $\xi_x$ be the
zero section of $\nbhd{U}_x$\@.  Note that $\xi_x\in F(\nbhd{U}_x)$\@.  Since
$\sF$ is a sheaf, it follows that there exists $\eta\in F(\man{M})$ such that
$\eta|\nbhd{U}_x=\xi_x$ for every $x\in\man{M}$\@.  Again since $\sF$ is a
sheaf, $\eta=\ol{\xi}$\@, showing that $\ol{\xi}\in F(\man{M})$\@, as
desired.  Finally, since $[(\xi,\nbhd{U})]_{x_0}=[(\ol{\xi},\man{M})]_x$\@,
the first assertion of the proposition follows.

The second assertion follows since the first assertion implies that we can
choose global smooth generators for $\dist{F}(\sF)$\@.
\end{proof}
\end{proposition}

Now let us show that the analogue of the preceding result does not hold in
the real analytic case.
\begin{example}\label{eg:subsheaf!subbundle}
We take $\man{M}=\real$ and consider the trivial vector bundle
$\man{E}=\real\times\real$ with projection $\pi(x,v)=x$\@.  We take
$S=\{0\}\cup\setdef{\frac{1}{j}}{j\in\integerp}$ and define a subsheaf
$\sF_S$ of $\ssections[\omega]{\man{E}}$ by
\begin{equation*}
\sF_S(\nbhd{U})=\setdef{\xi\in\sections[\omega]{\man{E}|\nbhd{U}}}
{\xi(x)=0\ \textrm{for all}\ x\in\nbhd{U}\cap S}.
\end{equation*}
It is clear that
\begin{equation*}
\dist{F}(\sF_S)_x=\begin{cases}\{x\}\times\real,&x\not\in S,\\
\{(x,0)\},&x\in S.\end{cases}
\end{equation*}
By the Identity Theorem for real analytic
functions~\cite[Corollary~1.2.6]{SGK/HRP:02}\@, the stalk $\sF_{S,0}$ of
$\sF_S$ at $0$ consists of germs of functions that vanish in some
neighbourhood of $0$\@,~\ie~$\sF_{S,0}=\{0\}$\@.  More precisely, if
$\nbhd{N}$ is a connected neighbourhood of $0$\@, then
$\sF_S(\nbhd{N})=\{0\}$\@.  This precludes the existence of real analytic
local generators in any neighbourhood of $0$\@, and so $\dist{F}(\sF_S)$ is
not real analytic.\oprocend
\end{example}

This addresses the first of our questions.  As to the second question
concerning the relationship between $\sF$ and
$\ssections[r]{\dist{F}(\sF)}$\@, we consider the following example.
\begin{example}\label{eg:bad-generator}
We take $\man{M}=\real$ and $\man{E}=\real\times\real$ with the vector bundle
projection $\map{\pi}{\man{E}}{\man{M}}$ being projection onto the first
factor.  We consider the real analytic subsheaf $\sF$ generated by the global
section $x\mapsto(x,x^2)$\@.  That is, the local sections of $\sF$ over
$\nbhd{U}\subset\real$ are the $\sfunc[\omega]{\nbhd{U}}$-multiples of this
section restricted to $\nbhd{U}$\@.  As in
Example~\ref{eg:!stalk-generate}\@, the generalised subbundle generated by
this subsheaf is
\begin{equation*}
\dist{F}(\sF)_x=\begin{cases}\{x\}\times\real,&x\not=0,\\
\{x\}\times\{0\},&x=0.\end{cases}
\end{equation*}
The section $x\mapsto(x,x)$ is a section of
$\sections[\omega]{\dist{F}(\sF)}$ that is not a global section of $\sF$\@,
showing that $\sF\subsetneq\ssections[\omega]{\dist{F}(\sF)}$\@.\oprocend
\end{example}

\subsection{Generators for submodules and subsheaves of sections}

One often wishes to define a submodule of sections, or a subsheaf of
sections, of a vector bundle by using generators.  In this section we give
the notation for doing this, and establish a few elementary consequences of
these constructions.

Let us consider first the nonsheaf situation where we have families of
globally defined sections.
\begin{definition}\label{def:generated}
Let $r\in\integernn\cup\{\infty,\omega\}$ and let\/
$\map{\pi}{\man{E}}{\man{M}}$ be a smooth or real analytic manifold, as
required.  If\/ $\sX=\ifam{\xi_j}_{j\in J}$ is a family of sections of\/
$\man{E}$\@,
\begin{compactenum}[(i)]
\item the \defn{generalised subbundle generated by $\sX$} is
\begin{equation*}
\dist{F}(\sX)_x=\vecspan[\real]{\xi_j(x)|\enspace j\in J}
\end{equation*}
and
\item the \defn{module of sections generated by $\sX$} is
\begin{equation*}
\modgen{\sX}=\vecspan[{\func[r]{\man{M}}}]{\xi_j|\enspace j\in J}.
\end{equation*}
\end{compactenum}
In the special case $\man{E}=\tb{\man{M}}$ and where $\sX=\ifam{X_j}_{j\in
J}$ is a family of vector fields, the \defn{distribution generated by $\sX$}
is denoted by $\dist{D}(\sX)$\@.\oprocend
\end{definition}

The following elementary result is often useful to make certain assumptions
without losing generality.
\begin{proposition}\label{prop:subbundlespan1}
Let\/ $r\in\integernn\cup\{\infty,\omega\}$\@, let\/
$\map{\pi}{\man{E}}{\man{M}}$ be a\/ smooth or real analytic vector bundle,
as required, and let\/ $\sX\subset\sections[r]{\man{E}}$\@.  Then the
generalised subbundles
\begin{compactenum}[(i)]
\item $\dist{F}(\sX)$\@,
\item $\dist{F}(\vecspan[\real]{\sX})$\@, and
\item $\dist{F}(\modgen{\sX})$
\end{compactenum}
agree.
\begin{proof}
Since we clearly have
\begin{equation*}
\dist{F}(\sX)\subset\dist{F}(\vecspan[\real]{\sX})\subset
\dist{F}(\modgen{\sX}),
\end{equation*}
it suffices to show that $\dist{F}(\modgen{\sX})=\dist{F}(\sX)$\@.  We have
\begin{align*}
\dist{F}(\modgen{\sX}&)_x\\
=&\;\vecspan[\real]{(f^1\xi_1+\dots+f^k\xi_k)(x)|\enspace
k\in\integerp,\ f^1,\dots,f^k\in\func[r]{\man{M}},\
\xi_1,\dots,\xi_k\in\sX}\\
=&\;\vecspan[\real]{a^1\xi_1(x)+\dots+a^k\xi_k(x)|\enspace k\in\integerp,\
a^1,\dots,a^k\in\real,\ \xi_1,\dots,\xi_k\in\sX}\\
=&\;\vecspan[\real]{\xi(x)|\enspace \xi\in\sX}=\dist{F}(\sX)_x,
\end{align*}
giving the result.
\end{proof}
\end{proposition}

Let us now turn to sheaves generated by families of locally defined sections.
\begin{definition}\label{def:setsheaves}
Let $r\in\integernn\cup\{\infty,\omega\}$ and let
$\map{\pi}{\man{E}}{\man{M}}$ be a smooth or real analytic manifold, as
required.  Let $\sX=\ifam{X(\nbhd{U})}_{\nbhd{U}\,\textrm{open}}$ be a
subsheaf of sets of the sheaf $\ssections[r]{\man{E}}$\@\textemdash{}\ie~an
assignment to each open set $\nbhd{U}\subset\man{M}$ a subset
$X(\nbhd{U})\subset\sections[r]{\man{E}|\nbhd{U}}$ with the assignment
satisfying $X(\nbhd{V})=r_{\nbhd{U},\nbhd{V}}(X(\nbhd{U}))$ for every pair of
open sets $\nbhd{U},\nbhd{V}$ for which $\nbhd{V}\subset\nbhd{V}$\@.
\begin{compactenum}[(i)]
\item The \defn{sheaf of $\sfunc[r]{\man{M}}$-modules generated by $\sX$} is
the subsheaf
$\modgen{\sX}=\ifam{\modgen{X}(\nbhd{U})}_{\nbhd{U}\,\textrm{open}}$ of
$\sfunc[r]{\man{M}}$-modules defined by
$\modgen{X}(\nbhd{U})=\modgen{X(\nbhd{U})}$\@.
\item The \defn{generalised subbundle generated by $\sX$} is defined by
\begin{equation*}
\dist{F}(\sX)_x=\vecspan[\real]{\xi(x)|\enspace[\xi]_x\in\sX_x}.
\end{equation*}
\item If $\man{E}=\tb{\man{M}}$ and so $\sX$ is a subsheaf of sets of
$\ssections[r]{\tb{\man{M}}}$\@, then the \defn{distribution generated by
$\sX$} is defined by $\dist{D}(\sX)=\dist{F}(\sX)$\@.\oprocend
\end{compactenum}
\end{definition}

We can also adapt Proposition~\ref{prop:subbundlespan1} to the sheaf setting.
\begin{proposition}\label{prop:subbundlespan2}
Let\/ $r\in\integernn\cup\{\infty,\omega\}$\@, let\/
$\map{\pi}{\man{E}}{\man{M}}$ be a\/ smooth or real analytic vector bundle,
as required, and let\/ $\sX=\ifam{X(\nbhd{U})}_{\nbhd{U}\,\textrm{open}}$ be
a subsheaf of sets of $\ssections[r]{\man{E}}$\@.  Then, the generalised
subbundles
\begin{compactenum}[(i)]
\item $\dist{F}(\sX)$ and
\item $\dist{F}(\modgen{\sX})$
\end{compactenum}
agree.
\begin{proof}
This follows immediately from Proposition~\ref{prop:subbundlespan1}\@.
\end{proof}
\end{proposition}

Of course, it is also true that both generalised subbundles in the statement
of the preceding result agree with the generalised subbundle whose fibre at
$x$ is given by
\begin{equation*}
\vecspan[\real]{\xi(x)|\enspace[\xi]_x\in\sX_x},
\end{equation*}
the setting for this being to think of $\ssections[r]{\man{E}}$ as a sheaf of
$\real$-vector spaces.  We leave it to the reader to develop the attendant
definitions as we shall not make use of any of these.

Note that the specification of a family $\sX=\ifam{\xi_j}_{j\in J}$ of global
sections of $\map{\pi}{\man{E}}{\man{M}}$ also prescribes a subsheaf of
$\ssections[r]{\man{E}}$\@.
\begin{definition}\label{def:module->sheaf}
Let $r\in\integernn\cup\{\infty,\omega\}$\@, let
$\map{\pi}{\man{E}}{\man{M}}$ be a smooth or real analytic vector bundle, as
required, and let $\sX\subset\sections[r]{\man{E}}$\@.
\begin{compactenum}[(i)]
\item For $\nbhd{U}\subset\man{M}$ open, the \defn{restriction} of $\sX$ to
$\nbhd{U}$ is
\begin{equation*}
\sX|\nbhd{U}=\setdef{\xi|\nbhd{U}}{\xi\in\sX}.
\end{equation*}
\item The \defn{subsheaf generated by $\sX$} is the subsheaf
$\sF_{\sX}=\ifam{F_{\sX}(\nbhd{U})}_{\nbhd{U}\,\textrm{open}}$ of
$\ssections[r]{\man{E}}$ defined by
\begin{equation*}\eqoprocend
F_{\sX}(\nbhd{U})=
\vecspan[{\func[r]{\nbhd{U}}}]{\sX|\nbhd{U}}.
\end{equation*}
\end{compactenum}
\end{definition}

Thus prescribing the data of a subsheaf is more general than prescribing a
globally defined submodule.
\begin{definition}\label{def:finite-module}
Let $r\in\integernn\cup\{\infty,\omega\}$\@, let
$\map{\pi}{\man{E}}{\man{M}}$ be a smooth or real analytic vector bundle, as
required, and let $\sM\subset\sections[r]{\man{E}}$ be a submodule of
sections.
\begin{compactenum}[(i)]
\item The submodule $\sM$ is \defn{finitely generated} if it is finitely
generated in the usual sense.
\item The submodule $\sM$ is \defn{locally finitely generated} if, for each
$x\in\man{M}$\@, there exists a neighbourhood $\nbhd{U}$ of $x$ such that the
submodule $F_{\sM}(\nbhd{U})$ of the $\func[r]{\nbhd{U}}$-module
$\sections[r]{\man{E}|\nbhd{U}}$ is finitely generated.\oprocend
\end{compactenum}
\end{definition}

We have on hand now a few different situations where we can apply the notion

of being locally finitely generated.  We refer to Section~\ref{subsec:fingen}
for a discussion of the property of local finite generation for various
objects.

\subsection{Patchy subsheaves}\label{subsec:patchy-subsheaves}

Example~\ref{eg:!stalk-generate} indicates that generators for a generalised
subbundle do not generally serve as generators for stalks of the associated
subsheaf of sections.  In this section we describe a construction where local
generators \emph{can} be used to define a subsheaf.  This requires that we
place some compatibility conditions on local generators.
\begin{proposition}\label{prop:patchy-exist}
Let\/ $r\in\integernn\cup\{\infty,\omega\}$ and let\/
$\map{\pi}{\man{E}}{\man{M}}$ be a smooth or real analytic vector bundle, as
is required.  Consider the following data:
\begin{compactenum}[(i)]
\item an open cover\/ $\sU=\ifam{\nbhd{U}_a}_{a\in A}$ for\/ $\man{M}$\@;
\item \label{pl:patchy-exist2} for each\/ $a\in A$\@, a family\/
$\sX_a=\ifam{\xi_b}_{b\in B_a}$ of\/ $\C^r$-sections of\/
$\man{E}|\nbhd{U}_a$ such that, if\/
$\nbhd{U}_{a_1}\cap\nbhd{U}_{a_2}\not=\emptyset$\@, then
\begin{equation*}
\modgen{\sX_{a_1}|\nbhd{U}_{a_1}\cap\nbhd{U}_{a_2}}=
\modgen{\sX_{a_2}|\nbhd{U}_{a_1}\cap\nbhd{U}_{a_2}};
\end{equation*}
\item for each\/ $a\in A$\@, the sheaf\/
$\sF_a=\ifam{F_a(\nbhd{U})}_{\nbhd{U}\subset\nbhd{U}_a\,\textrm{open}}$ of\/
$\sfunc[r]{\nbhd{U}_a}$-modules given by
\begin{equation*}
F_a(\nbhd{U})=\modgen{\sX_a|\nbhd{U}}.
\end{equation*}
\end{compactenum}
Then there exists a unique subsheaf\/
$\sF_{\sU}=\ifam{F_{\sU}(\nbhd{U})}_{\nbhd{U}\,\textrm{open}}$ of\/
$\ssections[r]{\man{E}}$ with the property that\/
$\sF_{\sU}|\nbhd{U}_a=\sF_a$ for each\/ $a\in A$\@.
\begin{proof}
For $\nbhd{U}\subset\man{M}$ open we define
\begin{multline*}
F(\nbhd{U})=\left\{\ifam{\xi_a}_{a\in A}\right|\enspace
\xi_a\in F_a(\nbhd{U}\cap\nbhd{U}_a),\ a\in A,\\\left.
\xi_{a_1}|\nbhd{U}\cap\nbhd{U}_{a_1}\cap\nbhd{U}_{a_2}=
\xi_{a_2}|\nbhd{U}\cap\nbhd{U}_{a_1}\cap\nbhd{U}_{a_2},\ a_1,a_2\in A\right\}.
\end{multline*}
For $\nbhd{U},\nbhd{V}\subset\man{M}$ open and satisfying
$\nbhd{V}\subset\nbhd{U}$\@, define
$\map{r_{\nbhd{U},\nbhd{V}}}{F(\nbhd{U})}{F(\nbhd{V})}$ by
\begin{equation*}
r_{\nbhd{U},\nbhd{V}}(\ifam{\xi_a}_{a\in A})=
\ifam{\xi_a|\nbhd{V}\cap\nbhd{U}_a}_{a\in A}.
\end{equation*}
We will verify that $\sF=\ifam{F(\nbhd{U})}_{\nbhd{U}\,\textrm{open}}$ is a
sheaf.

Let $\nbhd{W}\subset\man{M}$ be open and let $\ifam{\nbhd{W}_i}_{i\in I}$ be
an open cover for $\nbhd{W}$\@.  Let $\xi,\eta\in F(\nbhd{W})$ satisfy
$r_{\nbhd{W},\nbhd{W}_i}(\xi)=r_{\nbhd{W},\nbhd{W}_i}(\eta)$ for each $i\in
I$\@.  We write $\xi=\ifam{\xi_a}_{a\in A}$ and $\eta=\ifam{\eta_a}_{a\in A}$
and note that we have
\begin{equation*}
\xi_a|\nbhd{W}_i\cap\nbhd{U}_a=\eta_a|\nbhd{W}_i\cap\nbhd{U}_a,
\qquad a\in A,\ i\in I.
\end{equation*}
Since $\sF_a$ is separated, $\xi_a=\eta_a$ for each $a\in A$ and so
$\xi=\eta$\@.

Let $\nbhd{W}\in\sO$ and let $\ifam{\nbhd{W}_i}_{i\in I}$ be an open cover
for $\nbhd{W}$\@.  For each $i\in I$ let $\xi_i\in F(\nbhd{W}_i)$ and suppose
that $r_{\nbhd{W}_i,\nbhd{W}_i\cap\nbhd{W}_j}(\xi_i)=
r_{\nbhd{W}_j,\nbhd{W}_i\cap\nbhd{W}_j}(\xi_j)$ for each $i,j\in I$\@.  We
write $\xi_i=\ifam{\xi_{i,a}}_{a\in A}$\@, $i\in I$\@, and note that
\begin{equation*}
\xi_{i,a}|\nbhd{W}_i\cap\nbhd{W}_j\cap\nbhd{U}_a=
\xi_{j,a}|\nbhd{W}_i\cap\nbhd{W}_j\cap\nbhd{U}_a,\qquad
i,j\in I,\ a\in A.
\end{equation*}
Since $\sF_a$ satisfies the gluing property, there exists $\xi_a\in
F_a(\nbhd{W}\cap\nbhd{U}_a)$ such that
\begin{equation*}
\xi_a|\nbhd{W}_i\cap\nbhd{U}_a=\xi_{i,a},\qquad i\in I,\ a\in A.
\end{equation*}
Let us define $\xi=\ifam{\xi_a}_{a\in A}$\@.  By \eqref{pl:patchy-exist2} we
have
\begin{equation*}
\xi_{i,a_1}|\nbhd{W}_i\cap\nbhd{U}_{a_1}\cap\nbhd{U}_{a_2}=
\xi_{i,a_2}|\nbhd{W}_i\cap\nbhd{U}_{a_1}\cap\nbhd{U}_{a_2},\qquad i\in A,\ a_1,a_2\in A.
\end{equation*}
Thus
\begin{equation*}
s_{a_1}|\nbhd{W}\cap\nbhd{U}_{a_1}\cap\nbhd{U}_{a_2}=
s_{a_2}|\nbhd{W}\cap\nbhd{U}_{a_1}\cap\nbhd{U}_{a_2},\qquad a_1,a_2\in A,
\end{equation*}
and so $s$ as constructed is an element of $F(\nbhd{W})$\@.

The preceding shows that $\sF$ is a sheaf of sets.  To verify that it is a
sheaf $\sfunc[r]{\man{M}}$-modules we define the algebraic operations in the
obvious way by defining them in $F(\nbhd{U})$ by
\begin{equation*}
\ifam{\xi_a}_{a\in A}+\ifam{\eta_a}_{a\in A}=\ifam{\xi_a+\eta_a}_{a\in A},\qquad
f\left(\ifam{\xi_a}_{a\in A}\right)=
\ifam{(f|\nbhd{U}\cap\nbhd{U}_a)\xi_a\cdot t_a}_{a\in A},
\end{equation*}
respectively.  One easily verifies that these operations are well-defined,
and that the restriction morphisms for $\sF$ are $\sfunc[r]{\man{M}}$-module
homomorphisms.

The uniqueness assertion of the proposition follows from the fact that our
constructions obviously prescribe the stalks of $\sF$\@, and so uniquely
prescribe $\sF$ since it is a sheaf.
\end{proof}
\end{proposition}

Let us give a name to the sheaf defined as in the preceding result.
\begin{definition}\label{def:patchy-subsheaf}
Let $r\in\integernn\cup\{\infty,\omega\}$ and let
$\map{\pi}{\man{E}}{\man{M}}$ be a smooth or real analytic vector bundle, as
is required.  Consider the data $\sU=\ifam{\nbhd{U}_a}_{a\in A}$ and, for
each $a\in A$\@, $\sX_a$ as in Proposition~\ref{prop:patchy-exist}\@.  We
call the subsheaf $\sF_{\sU}$ of $\sfunc[r]{\man{M}}$-modules defined as in
Proposition~\ref{prop:patchy-exist} the \defn{patchy subsheaf} defined by the
above data.\oprocend
\end{definition}

Note that the specification of a submodule $\sM\subset\sections[r]{\man{M}}$ 
is a particular case of a patchy subsheaf corresponding to the cover of
$\man{M}$ by the single open set $\man{M}$\@.  Thus the specification of a 
submodule of $\sections[r]{\man{E}}$ is a special case of specifying a 
subsheaf of $\ssections[r]{\man{E}}$\@.

Given our definition of a generalised subbundle, the motivation for studying
patchy sheaves is clear.  Specifically, the notions of smoothness we provide
for generalised subbundles are defined in terms of an open cover of the
manifold, on each subset of which there are sections of the prescribed
smoothness that generate the generalised subbundle.  The patchy condition is
simply one of ensuring that there is some compatibility on overlapping open
sets.  What is not clear is the extent to which a given subsheaf is patchy,
and, if a subsheaf is patchy, what are the implications of this.

As to the patchiness of subsheaf, let us first consider the smooth case.
$\sF_{\sU}\subset\ssections[r]{\dist{F}_{\sU}}$\@, the opposite inclusion
will not generally hold.
\begin{proposition}\label{prop:smooth-patchy}
If\/ $r\in\integernn\cup\{\infty\}$\@, if\/ $\map{\pi}{\man{E}}{\man{M}}$ is
a smooth vector bundle, and if\/ $\sF$ is a subsheaf of\/
$\ssections[r]{\man{E}}$\@, then\/ $\sF$ is a patchy subsheaf.
\begin{proof}
We will prove that we can take as our patchy subsheaf data the following:
\begin{enumerate}
\item the trivial open cover $\sU=\ifam{\man{M}}$\@;
\item for the single element $\man{M}$ of the open cover, the family
$F(\man{M})$ of global sections.
\end{enumerate}
Thus we are in the situation of Definition~\ref{def:module->sheaf} where we
are considering a sheaf defined by a submodule of global sections.  To show
that $\sF=\sF_{\sU}$\@, it is sufficient by~\cite[Proposition~II.1.1]{RH:77a}
to prove that the stalk $\sF_{x_0}$ is generated by global sections for every
$x_0\in\man{M}$\@.  This, however, has been shown in
Proposition~\ref{prop:subsheaf->subbundle}\@.
\end{proof}
\end{proposition}

In the real analytic case, however, there are subsheaves that are not patchy.
\begin{example}\label{eg:!patchy}
We consider $\man{M}=\real$\@, $\man{E}=\real\times\real$\@, and
$\pi(x,v)=x$\@.  Let
\begin{equation*}
S=\setdef{\tfrac{1}{j}}{j\in\integerp}\cup\{0\}.
\end{equation*}
Consider the presheaf $\sI_S=\ifam{I_S(\nbhd{U})}_{\nbhd{U}\,\textrm{open}}$
given by
\begin{equation*}
I_S(\nbhd{U})=\setdef{f\in\sections[\omega]{\nbhd{U}}}{f(x)=0\ \textrm{for}\
x\in\nbhd{U}\cap S}.
\end{equation*}
One can easily verify that $\sI_S$ is a subsheaf of
$\ssections[\omega]{\man{E}}$\@.  We claim $\sI_S$ is not patchy.  The
sneakiest way to see this is as follows.  Note that $\sI_{S,0}=\{0\}$ since
every section defined on a connected neighbourhood of $0$ and vanishing on
$S$ must be identically zero, as we saw in
Example~\ref{eg:subsheaf!subbundle}\@.  However, note that if $x\not=0$ then
$\sI_{S,x}\not=\{0\}$ and so, by Lemma~\ref{lem:local-generators}\@, it
follows that $\sI_S$ cannot be locally finitely generated.  By
Theorem~\ref{the:analytic-fingen} it then follows that $\sI_S$ is not
patchy.\oprocend
\end{example}

However, in the real analytic case, one does have the following result,
establishing a correspondence of sorts between patchy real analytic
subsheaves and real analytic distributions.
\begin{theorem}\label{the:patchy-analytic}
If\/ $\map{\pi}{\man{E}}{\man{M}}$ is a real analytic vector bundle, then the
following statements hold:
\begin{compactenum}[(i)]
\item \label{pl:patchy-analytic1} if\/ $\dist{F}$ is a real analytic
generalised subbundle of\/ $\man{E}$\@, then\/ $\ssections[\omega]{\dist{F}}$
is patchy;
\item \label{pl:patchy-analytic2} if\/ $\sF$ is a patchy subsheaf of\/
$\ssections[\omega]{\man{E}}$\@, then\/ $\dist{F}(\sF)$ is a real analytic
generalised subbundle.
\end{compactenum}
\begin{proof}
\eqref{pl:patchy-analytic1}

\newcommand\sH{\mathscr{H}}
\newcommand\sQ{\mathscr{Q}}

\begin{compactenum}
\item Fix some point $p\in\man{M}$\@.
\item Let $r_1<r_2<\dots<r_s$ be the ranks of the $\dist{F}_x$ (the images of
the sections of $\dist{F}$ in the fibres of $\dist{E}$)
\item Let $\man{M}_k$ be the (closed) set of points in the manifold $\man{M}$
where the rank of $\dist{F}_x$ is $\le r_k$ (so $\man{M}_{k-1}$ is contained
in $\man{M}_k$).  We only really care about those $\man{M}_k$ which contain
$p$\@, so let's ignore any that don't, and just assume that $p$ is in all of
them, i.e., that $\man{M}_1,\man{M}_2,\dots,\man{M}_s$ all contain $p$\@.
\item We now work inductively (descending on $k$), constructing subsheaves of
$\ssections[\omega]{\man{E}}$ that are equal to $\sF$ away from
$\man{M}_{k-1}$.  At each step we'll check that the resulting sheaf is
finitely generated near $p$\@.  At the last stage we'll be out of subsets
$\man{M}$\@, and the subsheaf constructed will be $\sG$\@.
\item Let $\sG_s$ be the subsheaf of sections of $\man{E}$ such that their
restriction to $\man{M}_s\setminus\man{M}_{s-1}$ lies in the rank $r_s$
bundle corresponding to the fibres of $\dist{F}$\@.  Thinking in coordinates
shows that this condition can be expressed as linear conditions on the
sections of $\man{E}$\@, and that the linear conditions extend across
$\man{M}_{s-1}$\@. (The linear conditions are the vanishing of certain
$(r_s+1)\times(r_s+1)$ minors, where the matrix in question has rows local
generators of $\dist{F}$ and one more row for a section of $\man{E}$\@,
thought of as a variable.  The columns are local coordinates of $\man{E}$\@.
The condition that the minors are zero is linear in the coordinates of the
section of $\man{E}$).

That means that we can express $\sG_s$ as the kernel of a map between locally
free sheaves.  Since the structure sheaf is coherent, the kernel is finitely
generated.  Let $\dist{H}_s$ be the trivial bundle of the rank equal to the
number of local generators of $\sG_s$\@, and $\man{H}_s\to\man{E}$ the map
whose image is $\sG_s$\@.
\item Now let $\sG_{s-1}$ be the subsheaf of sections of $\sG_s$ whose
restrictions to $\man{M}_{s-1}\setminus\man{M}_{s-2}$ lie in the rank
$r_{s-1}$ vector bundle corresponding to the fibres of $\dist{F}$\@.  It
suffices to show that the subsheaf of $\ssections[\omega]{\dist{H}_s}$ which
maps to $\sG_{s-1}$ is locally finitely generated, since then its image
$\sG_{s-1}$ will be too.

Let $\sQ_s$ be the restriction of $\ssections[\omega]{\dist{H}_s}$ to
$\man{M}_s$\@, and note that $\sQ_s$ is a quotient of
$\ssections[\omega]{\dist{H}_s}$\@. By the same local coordinates idea used
in part 4, we'll be able to express the condition that the sections of
$\sQ_s$ lie in the rank $r_{s-1}$ bundle corresponding to the sections of
$\dist{F}$ as linear conditions on the sections of $\sQ_s$ (restricted to
$\man{M}_s$), and so conclude that we have local finite generation (of a
particular subsheaf of $\sQ_s$).  The subsheaf of
$\ssections[\omega]{\dist{H}_s}$ we want is those sections which, upon
restriction to $\man{M}_s$, lie in the (locally finitely generated) subsheaf
of $\sQ_s$ constructed above.  This subsheaf of
$\ssections[\omega]{\dist{H}_s}$ should also be finitely generated (since
$\sQ_s$ is a quotient of $\ssections[\omega]{\dist{H}_s}$\@, this must follow
from a general principle, but I should think about the argument for a bit).

Let $\man{H}_{s-1}$ be the trivial bundle whose rank is equal to the number
of generators of the subsheaf of $\ssections[\omega]{\dist{H}_s}$ we just
constructed, and $\man{H}_{s-1}\to\man{H}_s\to\man{E}$ the induced map, with
image $\dist{G}_{s-1}$\@.

\item Now repeat step 5 until we get down to the bottom of the list, having constructed $\sG_1$ which is $\sG$, and a trivial bundle $\man{H}_1$ with a map $\man{H}_1\to\man{E}$ whose image is $\man{G}_1=\man{G}$.
\end{compactenum}

\eqref{pl:patchy-analytic2} Suppose that $\sF$ is the patchy subsheaf
corresponding to the data of an open cover $\sU=\ifam{\nbhd{U}_a}_{a\in A}$
with corresponding subsets $\ifam{\sX_a}_{a\in A}$ of real analytic sections
over the subsets of the open cover.  Let $x_0\in\man{M}$ and let $a\in A$ be
such that $x_0\in\nbhd{U}_a$\@.  For $x\in\nbhd{U}_a$ we have
\begin{equation*}
\dist{F}(\sF)_x=\setdef{\xi(x)}{[\xi]_x\in\sF_x}=
\setdef{\xi(x)}{\xi\in\modgen{\sX_a}},
\end{equation*}
showing that $\dist{F}(\sF)$ is generated, as a generalised subbundle, by
locally defined real analytic sections, and so is real analytic.
\end{proof}
\end{theorem}

The preceding results and example should be interpreted as follows.  The
result tells us that the attribute of being patchy is not interesting in the
smooth or finitely differentiable case.  The example tells us that the
attribute of being patchy is at least not vacuous in the real analytic case.
It may not be apparent at this point that patchiness is something useful to
study.  However, as we shall see in Corollary~\ref{cor:analytic-coherent}\@,
patchy real analytic sheaves are coherent and so Cartan's Theorems~A and~B
are available to be used for these sheaves.

\section{Algebraic constructions associated to generalised
subbundles}\label{sec:algebra}

In this section we consider some algebraic constructions associated with
generalised subbundles.  It is often the case that such considerations are
phrased in terms of the ring structure of functions on a manifold and the
corresponding module structure of the sections of a vector bundle.  This is
more naturally carried out using sheaves, and so much of what we say in this
section is based on sheaves.

\subsection{From stalks of a sheaf to fibres}

Let $r\in\{\infty,\omega\}$ and let $\map{\pi}{\man{E}}{\man{M}}$ be a vector
bundle of class $\C^r$\@.  The stalk of the sheaf $\ssections[r]{\man{E}}$ at
$x\in\man{M}$ is the set $\gsections[r]{x}{\man{E}}$ of germs of sections
which is a module over the ring $\gfunc[r]{x}{\man{M}}$ of germs of
functions.  The stalk is \emph{not} the same as the fibre $\man{E}_x$\@,
however, the fibre can be obtained from the stalk, and in this section we see
how this is done.  We shall couch this in a brief general algebraic
construction, just to add colour.

Recall that if $\alg{R}$ is a commutative unit ring, if
$\alg{I}\subset\alg{R}$ is an ideal, and if $\alg{A}$ is a unital
$\alg{R}$-module, $\alg{I}\alg{A}$ is the submodule of $\alg{A}$ generated by
elements of the form $rv$ where $r\in\alg{I}$ and $v\in\alg{A}$\@.
\begin{proposition}\label{prop:local-vector-space}
Let\/ $\alg{R}$ be a commutative unit ring that is local,~\ie~possess a
unique maximal ideal\/ $\mathfrak{m}$\@, and let\/ $\alg{A}$ be a unital\/
$\alg{R}$-module.  Then\/ $\alg{A}/\mathfrak{m}\alg{A}$ is a vector space
over\/ $\alg{R}/\mathfrak{m}$\@.  Moreover, this vector space is naturally
isomorphic to\/ $(\alg{R}/\mathfrak{m})\otimes_{\alg{R}}\alg{A}$\@.
\begin{proof}
We first prove that $\alg{R}/\mathfrak{m}$ is a field.  Denote by
$\map{\pi_{\mathfrak{m}}}{\alg{R}}{\alg{R}/\mathfrak{m}}$ the canonical
projection.  Let $\alg{I}\subset\alg{R}/\mathfrak{m}$ be an ideal.  We claim
that
\begin{equation*}
\tilde{\alg{I}}=\setdef{r\in\alg{R}}{\pi_{\mathfrak{m}}(r)\in\alg{I}}
\end{equation*}
is an ideal in $\alg{R}$\@.  Indeed, let $r_1,r_2\in\tilde{\alg{I}}$ and note
that $\pi_{\mathfrak{m}}(r_1-r_2)=
\pi_{\mathfrak{m}}(r_1)-\pi_{\mathfrak{m}}(r_2)\in\alg{I}$ since
$\pi_{\mathfrak{m}}$ is a ring homomorphism and since $\alg{I}$ is an ideal.
Thus $r_1-r_2\in\tilde{\alg{I}}$\@.  Now let $r\in\tilde{\alg{I}}$ and
$s\in\alg{R}$ and note that
$\pi_{\mathfrak{m}}(sr)=\pi_{\mathfrak{m}}(s)\pi_{\mathfrak{m}}(r)
\in\alg{I}$\@, again since $\pi_{\mathfrak{m}}$ is a ring homomorphism and
since $\alg{I}$ is an ideal.  Thus $\tilde{\alg{I}}$ is an ideal.  Clearly
$\mathfrak{m}\subset\tilde{\alg{I}}$ so that either
$\tilde{\alg{I}}=\mathfrak{m}$ or $\tilde{\alg{I}}=\alg{R}$\@.  In the first
case $\alg{I}=\{0_{\alg{R}}+\mathfrak{m}\}$ and in the second case
$\alg{I}=\alg{R}/\mathfrak{m}$\@.  Thus the only ideals of
$\alg{R}/\mathfrak{m}$ are $\{0_{\alg{R}}+\mathfrak{m}\}$ and
$\alg{R}/\mathfrak{m}$\@.  To see that this implies that
$\alg{R}/\mathfrak{m}$ is a field, let
$r+\mathfrak{m}\in\alg{R}/\mathfrak{m}$ be nonzero and consider the ideal
$(r+\mathfrak{m})$\@.  Since $(r+\mathfrak{m})$ is nontrivial we must have
$(r+\mathfrak{m})=\alg{R}/\mathfrak{m}$\@.  In particular,
$1=(r+\mathfrak{m})(s+\mathfrak{m})$ for some
$s+\mathfrak{m}\in\alg{R}/\mathfrak{m}$\@, and so $r+\mathfrak{m}$ is a unit.

Now we show that $\alg{A}/\mathfrak{m}\alg{A}$ is a vector space over
$\alg{R}/\mathfrak{m}$\@.  This amounts to showing that the natural vector
space operations
\begin{equation*}
(u+\mathfrak{m}\alg{A})+(v+\mathfrak{m}\alg{A})=
u+v+\mathfrak{m}\alg{A},\quad(r+\mathfrak{m})(u+\mathfrak{m}\alg{A})=
ru+\mathfrak{m}\alg{A}
\end{equation*}
make sense.  The only possible issue is with scalar multiplication, so
suppose that
\begin{equation*}
r+\mathfrak{m}=s+\mathfrak{m},\quad
u+\mathfrak{m}\alg{A}=v+\mathfrak{m}\alg{A}
\end{equation*}
so that $s=r+a$ for $a\in\mathfrak{m}$ and $v=u+w$ for
$w\in\mathfrak{m}\alg{A}$\@.  Then
\begin{equation*}
sv=(r+a)(u+w)=ru+au+rw+aw,
\end{equation*}
and we observe that $au,rw,aw\in\mathfrak{m}\alg{A}$\@, and so the
sensibility of scalar multiplication is proved.

For the last assertion, note that we have the exact sequence
\begin{equation*}
\xymatrix{{0}\ar[r]&{\mathfrak{m}}\ar[r]&{\alg{R}}\ar[r]&
{\alg{R}/\mathfrak{m}}\ar[r]&{0}}
\end{equation*}
By right exactness of the tensor product~\cite[Proposition~IV.5.4]{TWH:74}
this gives the exact sequence
\begin{equation*}
\xymatrix{{\mathfrak{m}\otimes_{\alg{R}}\alg{A}}\ar[r]&
{\alg{A}}\ar[r]&{(\alg{R}/\mathfrak{m})\otimes_{\alg{R}}\alg{A}}\ar[r]&{0}}
\end{equation*}
noting that $\alg{R}\otimes_{\alg{R}}\alg{A}\simeq\alg{A}$\@.  By this
isomorphism, the image of $\mathfrak{m}\otimes_{\alg{R}}\alg{A}$ in $\alg{A}$
is simply generated by elements of the form $rv$ for $r\in\mathfrak{m}$ and
$v\in\alg{A}$\@.  That is to say, the image of
$\mathfrak{m}\otimes_{\alg{R}}\alg{A}$ in $\alg{A}$ is simply
$\mathfrak{m}\alg{A}$\@.  Thus we have the induced commutative diagram
\begin{equation*}
\xymatrix{&{\mathfrak{m}\otimes_{\alg{R}}\alg{A}}\ar[r]\ar[d]&
{\alg{A}}\ar[r]\ar[r]\ar[d]&{(\alg{R}/\mathfrak{m})\otimes_{\alg{R}}\alg{A}}
\ar[r]\ar@{-->}[d]&{0}\\
{0}\ar[r]&{\mathfrak{m}\alg{A}}\ar[r]&\alg{A}\ar[r]&
{\alg{A}/\mathfrak{m}\alg{A}}\ar[r]&{0}}
\end{equation*}
with exact rows.  We claim that there is an induced mapping as indicated by
the dashed arrow, and that this mapping is an isomorphism.  To define the
mapping, let $\alpha\in(\alg{R}/\mathfrak{m})\otimes_{\alg{R}}\alg{A}$ and
let $v\in\alg{A}$ project to $\alpha$\@.  The image of $\beta$ is then taken
to be $v+\mathfrak{m}\alg{A}$\@.  It is a straightforward exercise to show
that this mapping is well-defined and is an isomorphism, using exactness of
the diagram.
\end{proof}
\end{proposition}

With this simple algebraic construction as background, we can then indicate
how to recover the fibres of a vector bundle from the stalks of its sheaf of
sections.  To do this, the notation
\begin{equation*}
\mathfrak{m}_x=\setdef{[f]_x\in\gfunc[r]{x}{\man{M}}}{f(x)=0}
\end{equation*}
will be useful for $r\in\integernn\cup\{\infty,\omega\}$\@.  Algebraically,
$\mathfrak{m}_x$ is the unique maximal ideal of the local ring
$\gfunc[r]{x}{\man{M}}$~\cite{JANV/JBSS:03,JN:03}\@.
\begin{proposition}\label{prop:stalk-fibre}
Let\/ $r\in\{\infty,\omega\}$ and let\/ $\map{\pi}{\man{E}}{\man{M}}$ be a
vector bundle of class\/ $\C^r$\@.  For\/ $x\in\man{M}$ let\/
$\mathfrak{m}_x$ denote the unique maximal ideal in\/
$\gfunc[r]{x}{\man{M}}$\@.  Then the following statements hold:
\begin{compactenum}[(i)]
\item \label{pl:stalk-fibre1} the field
$\gfunc[r]{x}{\man{M}}/\mathfrak{m}_x$ is isomorphic to\/ $\real$ via the
isomorphism
\begin{equation*}
[f]_x+\mathfrak{m}_x\mapsto f(x);
\end{equation*}
\item \label{pl:stalk-fibre2} the\/
$\gfunc[r]{x}{\man{M}}/\mathfrak{m}_x$-vector space\/
$\gsections[r]{x}{\man{E}}/\mathfrak{m}_x\gsections[r]{x}{\man{E}}$ is
isomorphic to\/ $\man{E}_x$ via the isomorphism
\begin{equation*}
[\xi]_x+\mathfrak{m}_k\gsections[r]{x}{\man{E}}\mapsto\xi(x);
\end{equation*}
\item \label{pl:stalk-fibre3} the map from\/
$(\gfunc[r]{x}{\man{M}}/\mathfrak{m}_x)\otimes_{\gfunc[r]{x}{\man{M}}}
\gsections[r]{x}{\man{E}}$ to\/ $\man{E}_x$ defined by
\begin{equation*}
([f]_x+\mathfrak{m}_x)\otimes[\xi]_x\mapsto f(x)\xi(x)
\end{equation*}
is an isomorphism of\/ $\real$-vector spaces.
\end{compactenum}
\begin{proof}
\eqref{pl:stalk-fibre1} The map is clearly a homomorphism of fields.  To show
that it is surjective, if $a\in\real$ then $a$ is the image of
$[f]_x+\mathfrak{m}_x$ for any germ $[f]_x$ for which $f(x)=a$\@.  To show
injectivity, if $[f]_x+\mathfrak{m}_x$ maps to $0$ then clearly $f(x)=0$ and
so $f\in\mathfrak{m}_x$\@.

\eqref{pl:stalk-fibre2} The map is clearly linear, so we verify that it is an
isomorphism.  Let $v_x\in\man{E}_x$\@.  Then $v_x$ is the image of
$[\xi]_x+\mathfrak{m}_x\gsections[r]{x}{\man{E}}$ for any germ $[\xi]_x$ for
which $\xi(x)=v_x$\@.  Also suppose that
$[\xi]_x+\mathfrak{m}_x\gsections[r]{x}{\man{E}}$ maps to zero.  Then
$\xi(x)=0$\@.  Since $\ssections[r]{\man{E}}$ is locally free (see the next
section in case the meaning here is not patently obvious), it follows that we
can write
\begin{equation*}
\xi(y)=f_1(y)\eta_1(y)+\dots+f_m(y)\eta_m(y)
\end{equation*}
for sections $\eta_1,\dots,\eta_m$ of class $\C^r$ in a neighbourhood of $x$
and for functions $f_1,\dots,f_m$ of class $\C^r$ in a neighbourhood of
$x$\@.  Moreover, the sections may be chosen such that
$\ifam{\eta_1(y),\dots,\eta_m(y)}$ is a basis for $\man{E}_y$ for every $y$
in some suitably small neighbourhood of $x$\@.  Thus
\begin{equation*}
\xi(x)=0\quad\implies\quad f_1(x)=\dots=f_m(x)=0,
\end{equation*}
giving $\xi\in\mathfrak{m}_x\gsections[r]{x}{\man{E}}$\@, as desired.

\eqref{pl:stalk-fibre3} The $\real$-linearity of the stated map is clear, and
the fact that the map is an isomorphism follows from the final assertion of
Proposition~\ref{prop:local-vector-space}\@.
\end{proof}
\end{proposition}

The preceding result relates stalks to fibres.  Subsequently, specifically in
Theorem~\ref{the:locally-free-sheaf}\@, we shall take a more global view
towards relating vector bundles and sheaves.

In the preceding result we were able to rebuild the fibre of a vector bundle
from the germs of sections.  There is nothing keeping one from making this
construction for a general sheaf.
\begin{definition}
Let $r\in\integernn\cup\{\infty,\omega\}$\@, let $\man{M}$ be a smooth or
real analytic manifold, as required, and let $\sF$ be a sheaf of
$\sfunc[r]{\man{M}}$-modules.  The \defn{fibre} of $\sF$ is the
$\real$-vector space $\man{E}(\sF)_x=\sF_x/\mathfrak{m}_x\sF_x$\@.\oprocend
\end{definition}

By Proposition~\ref{prop:stalk-fibre}\@, the fibres of
$\ssections[r]{\man{E}}$ are isomorphic to the usual fibres.  However, if
$\sF$ is merely a subsheaf of the sheaf of sections of a vector bundle, the
relationship to the usual notion of fibre is not generally what one expects,
as the following example illustrates.
\begin{example}\label{eg:loc-gen-sheaf1}
Let $r\in\integernn\cup\{\infty,\omega\}$\@.  Let us take $\man{M}=\real$ and
define a presheaf $\sI^r_0=\ifam{I_0(\nbhd{U})}_{\nbhd{U}\,\textrm{open}}$ by
\begin{equation*}
I^r_0(\nbhd{U})=\begin{cases}\func[r]{\nbhd{U}},&0\not\in\nbhd{U},\\
\setdef{f\in\func[r]{\nbhd{U}}}{f(0)=0},&x\in\nbhd{U}.\end{cases}
\end{equation*}
One directly verifies that $\sI^r_0$ is a sheaf.  Moreover, $\sI_0^r$ is a
sheaf of $\sfunc[r]{\real}$-modules; this too is easily verified.  Let us
compute the fibres associated with this sheaf.  The germs of this sheaf at
$x\in\real$ are readily seen to be given by
\begin{equation*}
\sI^r_{0,x}=\begin{cases}\gfunc[r]{x}{\real},&x\not=0,\\
\mathfrak{m}_0=\setdef{[f]_0\in\gfunc[r]{0}{\real}}{f(x)=0},&x=0.\end{cases}
\end{equation*}
Thus we have
\begin{equation*}
\man{E}(\sI_0^r)_x=\begin{cases}\gfunc[r]{x}{\real}/
\mathfrak{m}_x\gfunc[r]{x}{\real}\simeq\real,&x\not=0,\\
\mathfrak{m}_0/\mathfrak{m}_0^2\simeq\real,&x=0.\end{cases}
\end{equation*}
Note that the fibre at $0$ is ``bigger'' than we expect it to be.\oprocend
\end{example}

The next result shows that the sheaf fibre is always larger than the usual
fibre, a fact that will be useful for us in our proof of
Theorem~\ref{the:global-sections} below.
\begin{proposition}\label{prop:fibres}
Let\/ $r\in\integernn\cup\{\infty,\omega\}$\@, let\/
$\map{\pi}{\man{E}}{\man{M}}$ be a smooth or real analytic vector bundle, as
required, and let\/ $\sF$ be a subsheaf of\/ $\ssections[r]{\man{E}}$\@.
Then, for each\/ $x\in\man{M}$\@, there exists a natural epimorphism from the
fibre $\sF_x/\mathfrak{m}_x\sF_x$ onto the fibre\/ $\dist{F}(\sF)_x$\@.  If
the function\/ $x\mapsto\dim(\dist{F}(\sF)_x)$ is locally constant at\/
$x$\@, then the natural epimorphism is an isomorphism.
\begin{proof}
Let $x\in\man{M}$\@, let $\nbhd{N}$ be a neighbourhood of $x$ and let
$\ifam{\xi_j}_{j\in J}$ be local generators for $\dist{F}$ about $x$\@,
defined on $\nbhd{N}$\@.  Let us define a morphism $\Psi$ from $\oplus_{j\in
J}\sfunc[r]{\nbhd{N}}$ to $\ssections[r]{\man{E}|\nbhd{N}}$\@, a morphism
$\Phi$ from $\oplus_{j\in J}\sfunc[r]{\nbhd{N}}$ to
$\ssections[r]{\dist{F}}|\nbhd{N}$\@, and a morphism $\iota$ from
$\ssections[r]{\dist{F}}|\nbhd{N}$ to $\ssections[r]{\man{E}|\nbhd{N}}$ by
\begin{align*}
\Psi_{\nbhd{U}}\left(\oplus_{j\in J}f^j\right)=&\;
\sum_{j\in J}f^j\hat{\xi}_j,\\
\Phi_{\nbhd{U}}\left(\oplus_{j\in J}f^j\right)=&\;\sum_{j\in J}f^j\xi_j,\\
\iota_{\nbhd{U}}\left(\sum_{j\in J}f^j\xi_j\right)=&\;
\sum_{j\in J}f^j\hat{\xi}_j,
\end{align*}
for $\nbhd{U}\subset\nbhd{N}$ open.  Here $\xi_j$ denotes an element of
$\sections[r]{\dist{F}|\nbhd{U}}$ and $\hat{\xi}_j$ denotes the element of
$\sections[r]{\man{E}|\nbhd{U}}$ which is image of $\xi_j$ under the obvious
inclusion.  At the stalk level we thus have the commutative diagram
\begin{equation*}
\xymatrix{{\oplus_{j\in J}\gfunc[r]{x}{\man{M}}}
\ar[rr]^(0.55){\Psi_x}\ar[dr]_{\Phi_x}&&{\gsections[r]{x}{\man{E}}}\\
&{\sF_x}\ar[ru]_{\iota_x}}
\end{equation*}
We now take the tensor product of the diagram of
$\gfunc[r]{x}{\man{M}}$-modules with
$\gfunc[r]{x}{\man{M}}/\mathfrak{m}_x$\@, and use the commuting of tensor
product with direct sums~\cite[Proposition~II.3.7.7]{NB:89a} and
Proposition~\ref{prop:local-vector-space} to give the commutative diagram
\begin{equation*}
\xymatrix{{\oplus_{j\in J}\real}\ar[rr]^(0.55){\tilde{\Psi}_x}
\ar[dr]_{\tilde{\Phi}_x}&&{\man{E}_x}\\
&{\sF_x/\mathfrak{m}_x\sF_x}\ar[ru]_{\tilde{\iota}_x}}
\end{equation*}
Here the homomorphisms $\tilde{\Psi}_x$\@, $\tilde{\Phi}_x$\@, and
$\tilde{\iota}_{a,x}$ are defined by
\begin{align*}
\tilde{\Psi}_x\left(\oplus_{b\in B_a}\alpha_b\right)=&\;
\sum_{b\in B_a}\alpha_b\xi_b(x),\\
\tilde{\Phi}_x\left(\oplus_{b\in B_a}\alpha_b\right)=&\;
\sum_{b\in B_a}\alpha_b([\xi_b]_x+\mathfrak{m}_x\sF_{a,x}),\\
\tilde{\iota}_{a,x}\left(\sum_{b\in B_a}\alpha_b([\xi_b]_x+
\mathfrak{m}_x\sF_{a,x})\right)=&\;\sum_{b\in B_a}\alpha_b\xi_b(x).
\end{align*}
Note that $\image(\tilde{\Psi}_x)=\dist{F}(\sF)_x$ by definition of
$\dist{F}(\sF)_x$\@.  By the commuting of the preceding diagram,
$\image(\tilde{\iota}_x)=\dist{F}(\sF)_x$\@, and so $\tilde{\iota}_x$ is an
epimorphism.

For the final assertion, note that if $x\mapsto\dim(\dist{F}(\sF)_x)$ is
locally constant, then $\dist{F}(\sF)_x$ are the fibres of a subbundle of
$\man{E}$\@.  In this case, the final assertion then follows from the
correspondence between the notions of fibre for sheaves of sections of vector
bundles.
\end{proof}
\end{proposition}

\subsection{Subsheaves of sections of a real analytic vector
bundle}\label{subsec:analytic-alg}

Now we turn to the consideration of sheaves associated to real analytic
generalised subbundles.  It is here that we will rely on some of the
important known results about coherent real analytic sheaves.

We begin with some a few more or less standard observations about germs of
functions and sections.
\begin{proposition}\label{prop:analytic-germs}
Let\/ $\map{\pi}{\man{E}}{\man{M}}$ be an analytic vector bundle.  For\/
$x\in\man{M}$ the following statements hold:
\begin{compactenum}[(i)]
\item the ring\/ $\gfunc[\omega]{x}{\man{M}}$ is Noetherian,~\ie~if we have a
sequence\/ $\ifam{\alg{I}_j}_{j\in\integerp}$ of ideals of\/
$\gfunc[\omega]{x}{\man{M}}$ for which\/ $\alg{I}_j\subset\alg{I}_{j+1}$\@,\/
$j\in\integerp$\@, then there exists\/ $N\in\integerp$ such that\/
$\alg{I}_j=\alg{I}_N$ for\/ $j\ge N$\@;
\item the module\/ $\gsections[\omega]{x}{\man{E}}$ is Noetherian,~\ie~if we
have a sequence\/ $\ifam{\alg{A}_j}_{j\in\integerp}$ of submodules of\/
$\gsections[\omega]{x}{\man{E}}$ for which\/
$\alg{A}_j\subset\alg{A}_{j+1}$\@,\/ $j\in\integerp$\@, then there exists\/
$N\in\integerp$ such that\/ $\alg{A}_j=\alg{A}_N$ for\/ $j\ge N$\@;
\end{compactenum}
\begin{proof}
In the holomorphic case, the first part of the result is well-known and can
be found in most any text on several complex
variables~\cite[\eg][Theorem~6.3.3]{LH:66}\@.  The proof is an inductive one
based on an induction on $\dim(\man{M})$ using the Weierstrass Preparation
Theorem.  The real analytic Weierstrass Preparation Theorem is proved by
\citet[Theorem~6.1.3]{SGK/HRP:02}\@, and with this the holomorphic proof of
the Noetherian property applies in the real analytic case.  The second part
of the result follows since every finitely generated module over a Noetherian
ring is Noetherian~\cite[Theorem~VIII.1.8]{TWH:74}\@.
\end{proof}
\end{proposition}

General properties of Noetherian modules then give the following result.
\begin{corollary}
Let\/ $\map{\pi}{\man{E}}{\man{M}}$ be an analytic vector bundle and let\/
$\sF$ be a subsheaf of\/ $\ssections[\omega]{\man{E}}$\@.  Then, for\/
$x\in\man{M}$\@,\/ $\sF_x$ is finitely generated.
\begin{proof}
This follows from the fact that submodules of finitely generated Noetherian
modules are finitely generated~\cite[Theorem~VIII.1.9]{TWH:74}\@.
\end{proof}
\end{corollary}

This property of stalks of real analytic subsheaves is often not used
properly.  Specifically, this Noetherian property of stalks of real analytic
subsheaves is used as a standin for the much stronger property of being
locally finitely generated (see Definition~\ref{def:locfingen}).  It is
actually this latter property we will use in this paper, and that is most
useful in general.  So we now turn to this question of local finite
generation.

\subsection{Locally finitely generated subsheaves of modules}

First let us recall the correspondence between vector bundles, or more
generally regular generalised subbundles, with locally free, locally finitely
generated sheaves.  We need, therefore, to complement
Definition~\ref{def:locfingen} of locally finitely generated with the
following.
\begin{definition}
Let $\man{M}$ be a smooth or real analytic manifold, as required, let
$r\in\integernn\cup\{\infty,\omega\}$\@, and let $\sF$ be a sheaf of\/
$\sfunc[r]{\man{M}}$-modules over $\man{M}$\@.  The sheaf $\sF$ is
\defn{locally free} if, for each $x_0\in\man{M}$\@, there exists a
neighbourhood $\nbhd{U}$ of $x_0$ such that $\sF|\nbhd{U}$ is isomorphic to a
direct sum $\oplus_{a\in A}(\sR|\nbhd{U})$\@. \oprocend
\end{definition}

The following theorem then gives the correspondence we are after.  The result
is often found in the holomorphic case in texts on algebraic
geometry~\cite[\eg][Proposition~7.6.5]{JLT:02}\@, but is seldom given in the
setting here (but see \cite[\S2.2]{SR:05} for some discussion).
\begin{theorem}\label{the:locally-free-sheaf}
Let\/ $r\in\integernn\cup\{\infty,\omega\}$\@, let\/
$\map{\pi}{\man{E}}{\man{M}}$ be a smooth or real analytic vector bundle, as
required, and let\/ $\dist{F}\subset\man{E}$ be a regular generalised
subbundle of class $\C^r$\@.  Then\/ $\ssections[r]{\man{E}}$ is a locally
free, locally finitely generated sheaf of\/ $\sfunc[r]{\man{M}}$-modules.

Conversely, if\/ $\sF$ is a locally free, locally finitely generated sheaf
of\/ $\sfunc[r]{\man{M}}$-modules, then there exists a smooth or real
analytic vector bundle\/ $\map{\pi}{\man{E}}{\man{M}}$\@, as required, and a
regular generalised subbundle\/ $\dist{F}\subset\man{E}$ of class\/ $\C^r$
such that\/ $\sF$ is isomorphic to\/ $\ssections[r]{\dist{F}}$\@.
\begin{proof}
Note that regular generalised subbundles of class $\C^r$ are vector bundles
of class $\C^r$\@.  Thus we shall suppose, without loss of generality that
$\dist{F}=\dist{E}$ and that $\man{E}$ is not smooth or real analytic, but of
class $\C^r$\@.

First let $\map{\pi}{\man{E}}{\man{M}}$ be a vector bundle of class $\C^r$
and let $x_0\in\man{M}$\@.  Let $(\nbhd{V},\psi)$ be a vector bundle chart
such that the corresponding chart $(\nbhd{U},\phi)$ for $\man{M}$ contains
$x_0$\@.  Suppose that $\psi(\nbhd{V})=\phi(\nbhd{U})\times\real^m$ and let
$\eta_1,\dots,\eta_m\in\sections[r]{\man{E}|\nbhd{U}}$ satisfy
$\psi(\eta_j(x))=(\phi(x),\vect{e}_j)$ for $x\in\nbhd{U}$ and
$j\in\{1,\dots,m\}$\@.  Let us arrange the components $\eta^k_j$\@,
$j,k\in\{1,\dots,m\}$\@, of the sections $\eta_1,\dots,\eta_m$ in an $m\times
m$ matrix:
\begin{equation*}
\mat{\eta}(x)=\begin{bmatrix}\eta^1_1(x)&\cdots&\eta^1_m(x)\\
\vdots&\ddots&\vdots\\\eta^m_1(x)&\cdots&\eta^m_m(x)\end{bmatrix}.
\end{equation*}
Now let $\xi\in\sections[r]{\man{E}|\nbhd{U}}$\@, let the components of $\xi$
be $\xi^k$\@, $k\in\{1,\dots,k\}$\@, and arrange the components in a vector
\begin{equation*}
\vect{\xi}(x)=\begin{bmatrix}\xi^1(x)\\\vdots\\\xi^m(x)\end{bmatrix}.
\end{equation*}
Now fix $x\in\nbhd{U}$\@.  We wish to solve the equation
\begin{equation*}
\xi(x)=f^1(x)\eta_1(x)+\dots+f^m(x)\eta_m(x)
\end{equation*}
for $f^1(x),\dots,f^m(x)\in\real$\@.  Let us write
\begin{equation*}
\vect{f}(x)=\begin{bmatrix}f^1(x)\\\vdots\\f^m(x)\end{bmatrix}.
\end{equation*}
Writing the equation we wish to solve as a matrix equation we have
\begin{equation*}
\vect{\xi}(x)=\mat{\eta}(x)\vect{f}(x).
\end{equation*}
Therefore,
\begin{equation*}
\vect{f}(x)=\mat{\eta}^{-1}(x)\vect{\xi}(x),
\end{equation*}
noting that $\vect{\eta}(x)$ is invertible since the vectors
$\eta_1(x),\dots,\eta_m(x)$ are linearly independent.  By Cramer's Rule, or
some such, the components of $\mat{\eta}^{-1}$ are $\C^r$-functions of
$x\in\nbhd{U}$\@, and so $\xi$ is a $\func[r]{\nbhd{U}}$-linear combination
of $\eta_1,\dots,\eta_m$\@, showing that $\sections[r]{\man{E}|\nbhd{U}}$ is
finitely generated.  To show that this module is free, it suffices to show
that $\ifam{\eta_1,\dots,\eta_m}$ is linearly independent over
$\func[r]{\nbhd{U}}$\@.  Suppose that there exists
$f^1,\dots,f^m\in\func[r]{\nbhd{U}}$ such that
\begin{equation*}
f^1\eta_1+\dots+f^m\eta_m=0_{\sections[r]{\man{E}}}.
\end{equation*}
Then, for every $x\in\nbhd{U}$\@,
\begin{equation*}
f^1(x)\eta_1(x)+\dots+f^m(x)\eta_m(x)=0_x\quad\implies
\quad f^1(x)=\dots=f^m(x)=0,
\end{equation*}
giving the desired linear independence.

Next suppose that $\sF$ is a locally free, locally finitely generated sheaf
of $\sfunc[r]{\man{M}}$-modules.  Let us first define the total space of our
vector bundle.  For $x\in\man{M}$ define
\begin{equation*}
\man{E}_x=\sF_x/\mathfrak{m}_x\sF_x.
\end{equation*}
By Propositions~\ref{prop:local-vector-space} and~\ref{prop:stalk-fibre}\@,
$\man{E}_x$ is a $\real$-vector space.  We take
$\man{E}=\disjointunion_{x\in\man{M}}\man{E}_x$\@.  Let $x\in\man{M}$ and let
$\nbhd{U}_x$ be a neighbourhood of $x$ such that $F(\nbhd{U}_x)$ is a free
$\func[r]{\nbhd{U}_x}$-module.  By shrinking $\nbhd{U}_x$ if necessary, we
suppose that it is the domain of a coordinate chart $(\nbhd{U}_x,\phi_x)$\@.
Let $s_1,\dots,s_m\in F(\nbhd{U}_x)$ be such that $\ifam{s_1,\dots,s_m}$ is a
basis for $F(\nbhd{U}_x)$\@.  Note that $\ifam{[s_1]_y,\dots,[s_m]_x}$ is a
basis for $\sF_y$ for each $y\in\nbhd{U}_x$\@.  It is straightforward to show
that
\begin{equation*}
\ifam{[s_1]_y+\mathfrak{m}_y\sF_y,\dots,[s_m]_y+\mathfrak{m}_y\sF_y}
\end{equation*}
is then a basis for $\man{E}_y$\@.  For $y\in\nbhd{U}$ the map
\begin{equation*}
a^1([s_1]_y+\mathfrak{m}_y)+\dots+a^m([s_m]_y+\mathfrak{m}_y)\mapsto
(a^1,\dots,a^m)
\end{equation*}
is clearly an isomorphism.  Now define
$\nbhd{V}_x=\disjointunion_{y\in\nbhd{U}_x}\man{E}_y$ and define
$\map{\psi_x}{\nbhd{V}_x}{\phi_x(\nbhd{U})\times\real^m}$ by
\begin{equation*}
\psi_x(a^1([s_1]_y+\mathfrak{m}_y)+\dots+a^m([s_m]_y+\mathfrak{m}_y))=
(\phi_x(y),(a^1,\dots,a^m)).
\end{equation*}
This is clearly a vector bundle chart for $\man{E}$\@.  Moreover, this
construction furnishes a covering of $\man{E}$ by vector bundle charts.

Next we show that two overlapping vector bundle charts satisfy the
appropriate overlap condition.  Thus let $x,y\in\man{M}$ be such that
$\nbhd{U}_x\cap\nbhd{U}_y$ is nonempty.  Let $\ifam{s_1,\dots,s_m}$ and
$\ifam{t_1,\dots,t_m}$ be bases for $F(\nbhd{U}_x)$ and $F(\nbhd{U}_y)$\@,
respectively.  (Note that the cardinality of these bases agrees since, for
$z\in\nbhd{U}_x\cap\nbhd{U}_y$\@, $\ifam{[s_1]_z,\dots,[s_m]_z}$ and
$\ifam{[t_1]_z,\dots,[t_m]_z}$ are both bases for
$\sF_z$\@,~\cf~\cite[Corollary~IV.2.12]{TWH:74}\@.)  Note that
\begin{equation*}
r_{\nbhd{U}_x,\nbhd{U}_x\cap\nbhd{U}_x}(s_j)=
\sum_{k=1}^mf^k_jr_{\nbhd{U}_y,\nbhd{U}_x\cap\nbhd{U}_y}(t_k)
\end{equation*}
for $f^k_j\in\func[r]{\nbhd{U}_x\cap\nbhd{U}_y}$\@, $j,k\in\{1,\dots,m\}$\@.
At the stalk level we have
\begin{equation*}
[s_j]_z=\sum_{k=1}^m[f^k_j]_z[t_k]_z,
\end{equation*}
from which we conclude that
\begin{equation*}
([s_j]_z+\mathfrak{m}_z\sF_z)=
\sum_{k=1}^mf^k_j(z)([t_k]_z+\mathfrak{m}_z\sF_z),
\end{equation*}
From this we conclude that the matrix
\begin{equation*}
\mat{f}(z)=\begin{bmatrix}f^1_1(z)&\cdots&f^1_m(z)\\
\vdots&\ddots&\vdots\\f^m_1(z)&\cdots&f^m_m(z)\end{bmatrix}
\end{equation*}
is invertible, being the change of basis matrix for the two bases for
$\man{E}_z$\@.  Moreover, the change of basis formula gives
\begin{equation*}
\psi_y\scirc\psi_x^{-1}(\vect{z},(a^1,\dots,a^m))=
\left(\phi_y\scirc\phi_x^{-1}(\vect{z}),\left(\sum_{j=1}^ma^jf_j^1(z),\dots,
\sum_{j=1}^ma^jf_j^m(z)\right)\right)
\end{equation*}
for every $z\in\nbhd{U}_x\cap\nbhd{U}_y$\@, where $\vect{z}=\phi_x(z)$\@.
Thus we see that the covering by vector bundle charts has the proper overlap
condition to define a vector bundle structure for $\man{E}$\@.

It remains to show that $\ssections[r]{\man{E}}$ is isomorphic to $\sF$\@.
Let $\nbhd{U}\subset\man{M}$ be open and define
$\map{\Phi_{\nbhd{U}}}{F(\nbhd{U})}{\sections[r]{\man{E}|\nbhd{U}}}$ by
\begin{equation*}
\Phi_{\nbhd{U}}(s)(x)=[s]_x+\mathfrak{m}_x\sF_x.
\end{equation*}
For this definition to make sense, we must show that $\Phi_{\nbhd{U}}(s)$ is
of class $\C^r$\@.  Let $y\in\nbhd{U}$ and, using the above constructions,
let $\ifam{s_1,\dots,s_m}$ be a basis for $F(\nbhd{U}_y)$\@.  Let us
abbreviate $\nbhd{V}=\nbhd{U}\cap\nbhd{U}_y$\@.  Note that
$\ifam{r_{\nbhd{U},\nbhd{V}}(s_1),\dots,r_{\nbhd{U},\nbhd{V}}(s_m)}$ is a
basis for $F(\nbhd{V})$\@.  (To see that this is so, one can identify
$F(\nbhd{U})$ with a section of the \'etale space $\Et(\sF)$ over
$\nbhd{U}$\@,~\cf~the discussion at the end of Section~\ref{subsec:etale}\@,
and having done this the assertion is clear.)  We thus write
\begin{equation*}
r_{\nbhd{U},\nbhd{V}}(s)=f^1r_{\nbhd{U},\nbhd{V}}(s_1)+\dots+
f^mr_{\nbhd{U},\nbhd{V}}(s_m).
\end{equation*}
In terms of stalks we thus have
\begin{equation*}
[s]_z=[f^1]_z[s_1]_z+\dots+[f^m]_z[s_m]_z
\end{equation*}
for each $z\in\nbhd{V}$\@.  Therefore,
\begin{equation*}
\Phi_{\nbhd{U}}(s)(z)=f^1(z)([s_1]_z+\mathfrak{m}_z\sF_z)+\dots+
f^m(z)([s_m]_z+\mathfrak{m}_z\sF_z),
\end{equation*}
which (recalling that $\nbhd{U}_y$\@, and so also $\nbhd{V}$\@, is a chart
domain) gives the local representative of $\Phi_{\nbhd{U}}(s)$ on $\nbhd{V}$
as
\begin{equation*}
\vect{z}\mapsto(\vect{z},(f^1\scirc\phi_y^{-1}(\vect{z}),\dots,
f^m\scirc\phi_y^{-1}(\vect{z}))).
\end{equation*}
Since this local representative is of class $\C^r$ and since this
construction can be made for any $y\in\nbhd{U}$\@, we conclude that
$\Phi_{\nbhd{U}}(s)$ is of class $\C^r$\@.

Now, to show that the family of mappings $\ifam{\Phi_{\nbhd{U}}}_{\nbhd{U}\
\textrm{open}}$ is an isomorphism, by~\cite[Proposition~II.1.1]{RH:77a} it
suffices to show that the induced mapping on stalks is an isomorphism.  Let
us denote the mapping of stalks at $x$ by $\Phi_x$\@.  We again use our
constructions from the first part of this part of the proof and let
$\ifam{s_1,\dots,s_m}$ be a basis for $F(\nbhd{U}_x)$\@.  Let us show that
$\Phi_x$ is surjective.  Let $[\xi]_x\in\gsections[r]{x}{\man{M}}$\@,
supposing that $\xi\in\sections[r]{\man{E}|\nbhd{U}}$\@.  Let
$\nbhd{V}=\nbhd{U}\cap\nbhd{U}_x$\@.  Let the local representative of $\xi$
on $\nbhd{V}$ in the chart $(\nbhd{V}_x,\psi_x)$ be given by
\begin{equation*}
\vect{y}\mapsto(\vect{y},(f^1\scirc\phi_x^{-1}(\vect{y}),\dots,
f^m\scirc\phi_x^{-1}(\vect{y})))
\end{equation*}
for $f^1,\dots,f^m\in\func[r]{\nbhd{V}}$\@.  Then, if
\begin{equation*}
[s]_x=[f^1]_x[s_1]_x+\dots+[f^m]_x[s_m]_x,
\end{equation*}
we have $\Phi_x([s]_x)=[\xi]_x$\@.  To prove injectivity of $\Phi_x$\@,
suppose that $\Phi_x([s_x])=0$\@.  This means that $\Phi_x([s]_x)$ is the
germ of a section of $\man{E}$ over some neighbourhood $\nbhd{U}$ of $x$ that
is identically zero.  We may without loss of generality assume that
$\nbhd{U}\subset\nbhd{U}_x$\@.  We also assume without loss of generality (by
restriction of necessary) that $s\in F(\nbhd{U})$\@.  We thus have
\begin{equation*}
\Phi_{\nbhd{U}}(s)(y)=0,\qquad y\in\nbhd{U}.
\end{equation*}
Since
$\ifam{r_{\nbhd{U}_x,\nbhd{U}}(s_1),\dots,r_{\nbhd{U}_x,\nbhd{U}}(s_m)}$ is a
basis for $F(\nbhd{U})$ we write
\begin{equation*}
s=f^1r_{\nbhd{U}_x,\nbhd{U}}(s_1)+\dots+
f^mr_{\nbhd{U}_x,\nbhd{U}}(s_m).
\end{equation*}
for some uniquely defined $f^1,\dots,f^m\in\func[r]{\nbhd{U}}$\@.  We have
\begin{equation*}
\Phi_{\nbhd{U}}(s)(y)=f^1(y)([s_1]_y+\mathfrak{m}_y\sF_y)+\dots+
f^m(y)([s_m]_y+\mathfrak{m}_y\sF_y)
\end{equation*}
for each $y\in\nbhd{U}$\@.  Since
\begin{equation*}
\ifam{[s_1]_y+\mathfrak{m}_y\sF_y,\dots,[s_m]_y+\mathfrak{m}_y\sF_y}
\end{equation*}
is a basis for $\man{E}_y$\@, we must have $f^1(y)=\dots=f^m(y)=0$ for each
$y\in\nbhd{U}$\@, giving $[s]_x=0$\@.
\end{proof}
\end{theorem}

If we relax the assumption of regularity, it becomes more difficult, in fact
generally not possible, to establish any general results concerning the
finite generatedness of the sheaf.  However, in the patchy real analytic
case, we have the following important result.
\begin{theorem}\label{the:analytic-fingen}
Let\/ $\map{\pi}{\man{E}}{\man{M}}$ be a real analytic vector bundle and
let\/ $\sF$ be a patchy subsheaf of\/ $\ssections[\omega]{\man{E}}$\@.  Then,
for any\/ $x_0\in\man{M}$\@, there exists a neighbourhood\/ $\nbhd{U}$ of\/
$x_0$ and sections\/ $\xi_1,\dots,\xi_k$ of\/ $\sF$ over\/ $\nbhd{U}$ such
that, if\/ $\xi$ is a section of\/ $\sF$ over\/ $\nbhd{U}$\@, there exists\/
$f^1,\dots,f^k\in\func[\omega]{\nbhd{U}}$ such that
\begin{equation*}
\xi=f^1\xi_1+\dots+f^k\xi_k.
\end{equation*}
In particular,\/ $\sF$ is locally finitely generated.
\begin{proof}
A complete account of all that is needed to prove this result actually takes
a lot of effort.  We content ourselves by pointing the interested reader to
the textual literature for detailed proofs.

We shall need the notion of a germ of a function, not at a point, but at a
closed set $A\subset\man{M}$\@.  In this case, we consider pairs
$(f,\nbhd{U})$ where $\nbhd{U}\subset\man{M}$ is open with $A\subset\nbhd{U}$
and where $f\in\func[\omega]{\nbhd{U}}$\@.  Two such pairs $(f_1,\nbhd{U}_1)$
and $(f_2,\nbhd{U}_2)$ are equivalent if there exists
$\nbhd{V}\subset\nbhd{U}_1\cap\nbhd{U}_2$ with $A\subset\nbhd{V}$ and such
that $f_1|\nbhd{V}=f_2|\nbhd{V}$\@.  The set of equivalence classes we denote
by $\gfunc[\omega]{A}{\man{M}}$\@.  If $A,B\subset\man{M}$ are closed sets
such that $B\subset A$ we have the map
$\map{r_{A,B}}{\gfunc[\omega]{A}{\man{M}}}{\gfunc[\omega]{B}{\man{M}}}$ given
by restriction.  If $\alg{A}\subset\gfunc[\omega]{B}{\man{M}}$ we define
\begin{equation*}
r_{A,B}^{-1}(\alg{A})=\setdef{[f]_A\in\gfunc[\omega]{A}{\man{M}}}
{r_{A,B}([f]_A)\in\alg{A}}.
\end{equation*}
We abbreviate $r_{A,x_0}=r_{A,\{x_0\}}$\@.  If $\alg{A}$ is a subring
(\resp~ideal), one readily verifies that $r_{A,B}^{-1}(\alg{A})$ is also a
subring (\resp~ideal).  Similar constructions hold for germs of sections of a
vector bundle on closed sets.

The next lemma contains the hard technicalities needed to prove the theorem.
\begin{prooflemma}\label{plem:fingen}
Let \/ $\map{\pi}{\man{E}}{\man{M}}$ be a real analytic vector bundle, let\/
$x_0\in\man{M}$\@, and let\/ $\alg{A}$ be a submodule of\/
$\gsections[\omega]{x_0}{\man{E}}$\@.  Let\/ $\xi_1,\dots,\xi_k$ be sections
of\/ $\man{E}$ defined on some neighbourhood\/ $\nbhd{U}$ of\/ $x_0$ and such
that\/ $[\xi_1]_{x_0},\dots,[\xi_k]_{x_0}$ generate\/ $\alg{A}$ (since\/
$\gsections[\omega]{x_0}{\man{E}}$ is Noetherian).  Then there exists a
compact set\/ $C\subset\man{M}$ such that $x_0\in\interior(C)$ and such
that\/ $[\xi_1]_C,\dots,[\xi_k]_C$ generate the module\/
$r_{C,x_0}^{-1}(\alg{A})$\@.
\begin{subproof}
In the holomorphic case this result is proved as
\cite[Theorem~H.8]{RCG:90b}\@.  The result, as many of these finite
generation results, uses induction on $\dim(\man{M})$ and the Weierstrass
Preparation Theorem.  The result in the real analytic case can be proved in a
similar manner to the holomorphic case, using the real analytic Weierstrass
Preparation Theorem~\cite[Theorem~6.1.3]{SGK/HRP:02}\@.  The result stated in
\cite{RCG:90b} includes an estimate that is not required for our purposes,
and, moreover, does not hold in the real analytic case.  However, the
algebraic constructions we need can be extracted after chasing off the
details of the proof of the estimate.
\end{subproof}
\end{prooflemma}

Let $\sU=\ifam{\nbhd{U}_a}_{a\in A}$ be an open cover of $\man{M}$ and let
$\ifam{\sF_a}_{a\in A}$ be a family sheaves defining a patchy subsheaf
$\sF_{\sU}$ as in Definition~\ref{def:patchy-subsheaf}\@.  Let
$x_0\in\man{M}$ and let $a\in A$ be such that $x_0\in\nbhd{U}_a$\@.  Note
that, for any compact set $C\subset\nbhd{U}_a$ and for any
$\xi\in F_a(\nbhd{U}_a)$\@, we have $r_{C,x_0}([\xi]_C)=[\xi]_{x_0}$\@.
Therefore, letting
\begin{equation*}
\alg{A}=\setdef{[\xi]_{x_0}}{\xi\in F_a(\nbhd{U}_a)},
\end{equation*}
we see that
\begin{equation*}
r_{C,x_0}^{-1}(\alg{A})=\setdef{[\xi]_C}{\xi\in F_a(\nbhd{U}_a)}.
\end{equation*}
Suppose that $\alg{A}$ is generated by germs
$[\xi_1]_{x_0},\dots,[\xi_k]_{x_0}$\@, this by virtue of the fact that
$\gsections[\omega]{x_0}{\man{E}}$ is Noetherian.  By the lemma let
$C\subset\nbhd{U}_a$ be a compact set such that $x_0\in\interior(C)$ and such
that $[\xi_1]_C,\dots,[\xi_k]_C$ generate $r_{C,x_0}^{-1}(\alg{A})$\@.  Now
let $\nbhd{U}$ be a neighbourhood of $x_0$ such that $\nbhd{U}\subset C$\@.
Let $\xi\in F_a(\nbhd{U})$ so $\xi$ is the restriction to $\nbhd{U}$ of a
$\hat{\xi}$ section over $\nbhd{U}_a$\@, this because $\sF$ is patchy.  Note
that $[\hat{\xi}]_C\in r_{C,x_0}^{-1}(\alg{A})$ and so there exists
$[f^1]_C,\dots,[f^k]_C\in\gfunc[\omega]{C}{\man{M}}$ such that
\begin{equation*}
[\hat{\xi}]_C=[f^1]_C[\xi_1]_C+\dots+[f^k]_C[\xi_k]_C.
\end{equation*}
Therefore, we have
\begin{equation*}
\xi=(f^1|\nbhd{U})\xi_1|\nbhd{U}+\dots+
(f^k|\nbhd{U})\xi_k|\nbhd{U},
\end{equation*}
which is the first part of the result.  The second assertion follows
from Proposition~\pldblref{prop:locfingen}{pl:locfingen4}\@.
\end{proof}
\end{theorem}

From the theorem, we have the following corollaries, the third of which is
that with which we are presently concerned, but the first two of which are
also extremely important.
\begin{corollary}\label{cor:analytic-coherent}
If\/ $\map{\pi}{\man{E}}{\man{M}}$ is a real analytic vector bundle and if\/
$\sF$ is a patchy subsheaf of\/ $\ssections[\omega]{\man{E}}$\@, then\/ $\sF$
is coherent.
\begin{proof}
By Theorem~\ref{the:analytic-fingen} we know that $\sF$ is locally finitely
generated.  By Oka's Theorem, $\ssections[\omega]{\dist{E}}$ is coherent.
Since finitely generated subsheaves of coherent sheaves are
coherent~\cite[page~235]{HG/RR:84}\@, the result follows.
\end{proof}
\end{corollary}

The following corollary is one that is often used.
\begin{corollary}\label{cor:submod-fingen}
If\/ $\map{\pi}{\man{E}}{\man{M}}$ is a real analytic vector bundle and if\/
$\sM\subset\sections[\omega]{\man{E}}$ is a submodule of global sections,
then\/ $\sM$ is locally finitely generated.
\begin{proof}
This follows from Theorem~\ref{the:analytic-fingen} using the fact that the
sheaf $\sF_{\sM}$ (see Definition~\ref{def:module->sheaf} for the notation)
is obviously patchy; indeed, it is a patchy sheaf defined by the trivial open
covering.
\end{proof}
\end{corollary}

The following result is also interesting and useful.
\begin{corollary}\label{cor:subbunfingen}
If\/ $\map{\pi}{\man{E}}{\man{M}}$ is a real analytic vector bundle and if\/
$\dist{F}$ is an analytic generalised subbundle of\/ $\man{E}$\@, then\/
$\ssections[\omega]{\dist{F}}$ is locally finitely generated.
\begin{proof}
Let $x_0\in\man{M}$\@.  There then exists a neighbourhood $\nbhd{U}$ of $x_0$
and local generators $\ifam{\xi_a}_{a\in A}$ for $\dist{F}$ defined on
$\nbhd{U}$\@.  Let $\sM$ be the submodule of
$\sections[\omega]{\man{E}|\nbhd{U}}$ generated by the local generators.  We
then have the sheaf $\sF_{\sM}$ over $\nbhd{U}$\@, which is obviously patchy
(it has a single patch).  By Corollary~\ref{cor:submod-fingen} it follows
that $\sM$ is locally finitely generated.  Thus, by shrinking $\nbhd{U}$ we
arrive at a neighbourhood of $x_0$ with local sections
$\eta_1,\dots,\eta_k\in\sections[\omega]{\man{E}|\nbhd{U}}$ such that
\begin{equation*}
\dist{F}_x=\vecspan[\real]{\xi_a(x)|\enspace a\in A}=
\vecspan[\real]{\eta_j(x)|\enspace j\in\{1,\dots,k\}},
\end{equation*}
\cf~Proposition~\ref{prop:subbundlespan1}\@.
\end{proof}
\end{corollary}

Note that neither of the preceding two corollaries is a consequence of the
Noetherian property of germs of real analytic sheaves, but a consequence of
the deeper Theorem~\ref{the:analytic-fingen}\@.  This fact is routinely
misunderstood~\cite[\eg][Corollary~5.3]{AAA/YS:04}\@.  Let us make this
explicit by using one of our by now stock examples.
\begin{example}
We take $\man{M}=\real$ and consider the trivial vector bundle
$\man{E}=\real\times\real$ with projection $\pi(x,v)=x$\@.  We take
$S=\{0\}\cup\setdef{\frac{1}{j}}{j\in\integerp}$ and define a subsheaf
$\sF_S$ of $\ssections[\omega]{\man{E}}$ by
\begin{equation*}
\sF_S(\nbhd{U})=\setdef{\xi\in\sections[\omega]{\man{E}|\nbhd{U}}}
{\xi(x)=0\ \textrm{for all}\ x\in\nbhd{U}\cap S}.
\end{equation*}
The stalks $\sF_{S,x}$\@, $x\in\real$\@, are finitely generated since they
are Noetherian.  However, as we saw in Example~\ref{eg:!patchy}\@, $\sF_S$ is
not patchy, and, therefore, not locally finitely generated by
Theorem~\ref{the:analytic-fingen}\@.  Specifically, $\sF_S$ is not locally
finitely generated about $0$\@.\oprocend
\end{example}

In the smooth case, if we drop the assumption of regularity, then the sheaf
of sections of a singular generalised subbundle may be locally finitely
generated or not, depending on the precise example.  Let us illustrate 
what can happen.
\begin{example}\label{eg:singular-fingen}
We consider $\man{M}=\real^2$ and the distribution $\dist{D}$ given by
\begin{equation*}
\dist{D}_{(x_1,x_2)}=\begin{cases}\tb[(x_1,x_2)]{\real^2},&x_1\not=0,\\
\vecspan[\real]{\pderiv{}{x_1}},&x_1=0.\end{cases}
\end{equation*}
This distribution is generated by any pair of vector fields
\begin{equation*}
X_1(x_1,x_2)=\pderiv{}{x_1},\quad X_2(x_1,x_2)=f(x_1)\pderiv{}{x_2},
\end{equation*}
where $f\in\func[\infty]{\real}$ satisfies $f^{-1}(0)=\{0\}$\@.  Thus
$\dist{D}$ is a smooth distribution.  Moreover, we claim that if $X$ is any
section of $\dist{D}$ then $X=f^1X_1+f^2X_2$ for some
$f^1,f^2\in\func[\infty]{\real^2}$\@, provided we take $f$ to be defined by
$f(x)=x$\@.  Indeed, let us write
\begin{equation*}
X=g^1\pderiv{}{x_1}+g^2\pderiv{}{x_2},
\end{equation*}
noting that we must have $g^2(0,x_2)=0$ for every $x_2\in\real$\@.  We write
\begin{equation*}
g^2(x_1,x_2)=\int_0^{x_1}\pderiv{g^2}{x_1}(\xi,x_2)\,\d{\xi}=
x_1\int_0^1\pderiv{g^2}{x_1}(x_1\eta,x_2)\,\d{\eta}.
\end{equation*}
Thus our claim follows by taking
\begin{equation*}
f^1(x_1,x_2)=g^1(x_1,x_2),\quad
f^2(x_1,x_2)=\int_0^1\pderiv{g^2}{x_1}(x_1\eta,x_2)\,\d{\eta}.
\end{equation*}
This shows that, not only does $\dist{D}$ have a finite number of generators,
but also that $\sections[\infty]{\dist{D}}$ is free and finitely generated.
It follows that $\ssections[\infty]{\dist{D}}$ is free and finitely
generated.

By choosing more pathological smooth generators,~\eg~by taking
\begin{equation*}
f(x)=\begin{cases}\eul^{-1/x^2},&x\not=0,\\0,&x=0,\end{cases}
\end{equation*}
one imagines that the algebraic properties of the distribution should
deteriorate.  However, this is not seen in this case until one looks at the
Lie algebra generated by the generators.  The reader can see this clearly in
Example~\ref{eg:LinftyX-D(X)}\@.\oprocend
\end{example}

Be careful to understand that the example does not contradict
Theorem~\ref{the:locally-free-sheaf}\@.  For the previous example, the
theorem merely says that the sheaf $\ssections[\infty]{\dist{D}}$ of the
example is isomorphic to the sheaf of sections of some subbundle, in this
case the trivial bundle with two-dimensional fibre.

There is actually a simpler example of a subbundle that illustrates the main
point of the preceding example.  Indeed, the distribution $\dist{D}$ on
$\real$ defined by
\begin{equation*}
\dist{D}_x=\begin{cases}\tb[x]{\real},&x\not=0,\\\{0\},&x=0\end{cases}
\end{equation*}
has the desired property of being a singular subbundle whose subsheaf of
sections is free and finitely generated.  This can be proved along the lines
of the preceding example.  However, we prefer the example we work out because
we shall encounter it again in a different context, and it is illustrative to
be able to compare the various properties of the distribution.

The next example we consider shows that it is possible to have a distribution
whose submodule of sections is not finitely generated.
\begin{example}\label{eg:singular!fingen}
On $\man{M}=\real^2$ we consider the distribution $\dist{D}$ generated by the
two smooth vector fields
\begin{equation*}
X_1(x_1,x_2)=\pderiv{}{x_1},\quad X_2(x_1,x_2)=f(x_1)\pderiv{}{x_2},
\end{equation*}
where
\begin{equation*}
f(x)=\begin{cases}\eul^{-1/x^2},&x\in\realp,\\0,&x\in\realnp.\end{cases}
\end{equation*}
We claim that $\sections{\dist{D}}$ is not locally finitely generated.  To
see this, we refer ahead to Example~\ref{eg:converse-frobenius} where we show
that $\dist{D}$ is involutive but not integrable.  By Frobenius's Theorem,
Theorem~\ref{the:frobenius}\@, it follows that $\sections{\dist{D}}$ is not
locally finitely generated.\oprocend
\end{example}

Just as we saw after Example~\ref{eg:singular-fingen}\@, it is possible for
the preceding example to illustrate the point with a simpler distribution.
Indeed, the distribution $\dist{D}$ on $\real$ defined by
\begin{equation*}
\dist{D}_x=\begin{cases}\tb[x]{\real},&x\in\realp,\\
\{0\},&x\in\realnp,\end{cases}
\end{equation*}
has the property that it is involutive but not integrable.  Thus, again
according to Frobenius's Theorem, $\sections{\dist{D}}$ is not locally
finitely generated.

\subsection{Relationship between various notions of ``locally finitely
generated''}\label{subsec:fingen}

Note that we now have \emph{three} versions of ``locally finitely generated''
floating about, and they do not apply to the same thing, and the
relationships between them are not perfectly clear.  To be precise, we have
the notion of a locally finitely generated sheaf of
$\sfunc[r]{\man{M}}$-modules from Definition~\ref{def:locfingen}\@, the
notion of a locally finitely generated generalised subbundle from
Definition~\ref{def:subbundle}\@, and the notion of a locally finitely
generated submodule of sections from Definition~\ref{def:finite-module}\@.

The following property of locally finitely generated sheaves is of value.
\begin{lemma}\label{lem:local-generators}
Let\/ $\man{M}$ be a smooth or real analytic manifold, as required, let\/
$r\in\integernn\cup\{\infty,\omega\}$\@, and let\/
$\sF=\ifam{F(\nbhd{U})}_{\nbhd{U}\,\textrm{open}}$ be a locally finitely
generated sheaf of\/ $\sfunc[r]{\man{M}}$-modules.  If, for\/
$x_0\in\man{M}$\@,\/ $[s_1]_{x_0},\dots,[s_k]_{x_0}$ are generators for the\/
$\gfunc[\omega]{x_0}{\man{M}}$-module\/ $\sF_{x_0}$\@, then there exists a
neighbourhood\/ $\nbhd{U}$ of\/ $x_0$ such that\/ $[s_1]_x,\dots,[s_k]_x$ are
generators for\/ $\sF_x$ for each\/ $x\in\nbhd{U}$\@.
\begin{proof}
By hypothesis, there exists a neighbourhood $\nbhd{V}$ of $x_0$ and sections
$t_1,\dots,t_k\in F(\nbhd{V})$ such that $[t_1]_x,\dots,[t_m]_x$ generate
$\sF_x$ for all $x\in\nbhd{V}$\@.  Since $[s_1]_{x_0},\dots,[s_k]_{x_0}$
generate $\sF_{x_0}$\@,
\begin{equation*}
[t_l]_{x_0}=\sum_{j=1}^k[a^j_l]_{x_0}[s_j]_{x_0},\qquad l\in\{1,\dots,m\},
\end{equation*}
for germs $[a^j_l]_{x_0}\in\sR_{x_0}$\@.  We can assume, possibly by
shrinking $\nbhd{V}$\@, that $s_1,\dots,s_k\in F(\nbhd{V})$ and $a^j_l\in
R(\nbhd{V})$\@, $j\in\{1,\dots,k\}$\@, $l\in\{1,\dots,m\}$\@.  By definition
of germ, there exists a neighbourhood $\nbhd{U}\subset\nbhd{V}$ of $x_0$ such
that
\begin{equation*}
r_{\nbhd{V},\nbhd{U}}(t_l)=\sum_{j=1}^kr_{\nbhd{V},\nbhd{U}}(a^j_l)
r_{\nbhd{V},\nbhd{U}}(s_j),\qquad l\in\{1,\dots,m\}.
\end{equation*}
Taking germs at $x\in\nbhd{U}$ shows that the generators
$[t_1]_x,\dots,[t]_m$ of $\sF_x$ are linear combinations of
$[s_1]_x,\dots,[s_k]_x$\@, as desired.
\end{proof}
\end{lemma}

Note that the property of being locally finitely generated is one about
stalks, not one about local sections.  That is to say, if
$[s_1]_{x_0},\dots,[s_k]_{x_0}$ generate $\sF_{x_0}$ and if $\nbhd{U}$ is a
neighbourhood of $x_0$ for which $[s_1]_x,\dots,[s_k]_x$ generate $\sF_x$ for
each $x\in\nbhd{U}$\@, it is not clearly the case that $s_1,\dots,s_k$
generate $F(\nbhd{U})$\@.  What we can state is the following result, which
is a consequence of the vanishing of the cohomology groups of the sheaves in
these cases.
\begin{lemma}\label{lem:locfingen1}
Let\/ $r\in\integernn\cup\{\infty,\omega\}$\@, let\/ $\man{M}$ be a smooth or
real analytic manifold, as required, and let\/
$\sF=\ifam{F(\nbhd{U})}_{\nbhd{U}\,\textrm{open}}$ be a sheaf of\/
$\sfunc[r]{\man{M}}$-modules.  Assume one of the two cases:
\begin{compactenum}[(i)]
\item $r=\infty$\@;
\item $r=\omega$ and $\sF$ is coherent;
\end{compactenum}
Let\/ $x_0\in\man{M}$\@.  If\/
$[(s_1,\nbhd{U})]_{x_0},\dots,[(s_k,\nbhd{U})]_{x_0}$ generate\/
$\sF_{x_0}$\@, then there exists a neighbourhood\/ $\nbhd{W}\subset\nbhd{U}$
of\/ $x_0$ such that\/
$r_{\nbhd{U},\nbhd{V}}(s_1),\dots,r_{\nbhd{U},\nbhd{V}}(s_k)$ generate\/
$F(\nbhd{V})$ for every open set\/ $\nbhd{V}\subset\nbhd{W}$\@.
\begin{proof}
From the proof of Lemma~\ref{lem:local-generators} we see that there exists a
neighbourhood $\nbhd{W}$ of $x_0$ such that $\ifam{[s_1]_x,\dots,[s_k]_x}$
generate $\sF_x$ for every $x\in\nbhd{W}$\@.  If $\nbhd{V}\subset\nbhd{W}$\@,
we then have a presheaf morphism
$\Phi=\ifam{\Phi_{\nbhd{V}'}}_{\nbhd{V}'\subset\nbhd{V}\,\textrm{open}}$ from
$(\sfunc[r]{\nbhd{V}})^k$ to $\sF|\nbhd{V}$ given by
\begin{equation*}
\Phi_{\nbhd{V}'}(f^1,\dots,f^k)=f^1r_{\nbhd{U},\nbhd{V}'}(s_1)+\dots+
f^kr_{\nbhd{U},\nbhd{V}'}(s_k),
\end{equation*}
Note that the sequence
\begin{equation}\label{eq:locfingen1}
\xymatrix{{(\sfunc[r]{\nbhd{V}})^k}\ar[r]^{\Phi}&{\sF|\nbhd{V}}\ar[r]&{0}}
\end{equation}
is exact, by which we mean that it is exact on stalks.  If $s\in
F(\nbhd{V})$\@, exactness of~\eqref{eq:locfingen1} implies that, for
$x\in\nbhd{V}$\@,
\begin{equation*}
[s]_x=[g^1]_x[s_1]_x+\dots+[g^k]_x[s_k]_x
\end{equation*}
for $[g^1]_1,\dots,[g^k]_x\in\gfunc[r]{x}{\man{M}}$\@.  Since the preceding
expression involves only a finite number of germs, there exists a
neighbourhood $\nbhd{V}_x\subset\nbhd{V}$ of $x$ such that
\begin{equation*}
r_{\nbhd{V},\nbhd{V}_x}(s)=g^1_xr_{\nbhd{U},\nbhd{V}_x}(s_1)+\dots+
g^k_xr_{\nbhd{U},\nbhd{V}_x}(s_k)
\end{equation*}
for $g^1_x,\dots,g^k_x\in\func[r]{\nbhd{V}_x}$\@.  Let
$\sV=\ifam{\nbhd{V}_x}_{x\in\nbhd{V}}$\@.  If
$\nbhd{V}_x\cap\nbhd{V}_y\not=\emptyset$\@, define
$g^j_{xy}\in\func[r]{\nbhd{V}_x\cap\nbhd{V}_y}$ by
\begin{equation*}
g^j_{xy}=g^j_x|\nbhd{V}_x\cap\nbhd{V}_y-g^j_y|\nbhd{V}_x\cap\nbhd{V}_y,
\qquad j\in\{1,\dots,k\},
\end{equation*}
and note that $\ifam{(g^1_{xy},\dots,g^k_{xy})}_{x,y\in\nbhd{V}}\in
\cohomker{}{1}(\sV,\ker(\Phi))$\@.  We now note the following:
\begin{enumerate}
\item if $r=\infty$ then $\cohom{}{1}(\sV;\ker(\Phi))=0$ by
Theorem~\ref{the:smooth-cohom}\@;
\item if $r=\omega$ then $\ker(\Phi)$ is
coherent~\cite[page~237]{HG/RR:84}\@, and so $\cohom{}{1}(\sV;\ker(\Phi))=0$
by Cartan's Theorem~B.
\end{enumerate}
Since $\cohom{}{1}(\sV;\ker(\Phi))=0$\@, for each $x\in\nbhd{V}$ there exists
$\ifam{(h^1_x,\dots,h^k_x)}_{x\in\nbhd{V}_x}\in\cochain{}{1}(\sV;\ker(\Phi))$ such that
\begin{equation*}
h^j_y|\nbhd{V}_x\cap\nbhd{V}_y-h^j_x|\nbhd{V}_x\cap\nbhd{V}_y=
g^j_{xy}=g^j_x|\nbhd{V}_x\cap\nbhd{V}_y-g^j_y|\nbhd{V}_x\cap\nbhd{V}_y,
\qquad j\in\{1,\dots,k\}.
\end{equation*}
Define $f^j_x\in\func[r]{\nbhd{V}_x}$ by $f^j_x=g^j_x+h^j_x$\@, and note that
\begin{equation*}
f^j_x|\nbhd{V}_x\cap\nbhd{V}_y=f^j_y|\nbhd{V}_x\cap\nbhd{V}_y,
\qquad j\in\{1,\dots,k\}.
\end{equation*}
Thus there exists $f^j\in\func[r]{\nbhd{V}}$ such that $f^j|\nbhd{V}_x=f^j_x$
for each $j\in\{1,\dots,k\}$ and $x\in\nbhd{V}$\@.  Moreover, since
\begin{equation*}
h^1_xr_{\nbhd{U},\nbhd{V}_x}(s_1)+\dots+h^k_xr_{\nbhd{U},\nbhd{V}_x}(s_k)=0,
\end{equation*}
we have
\begin{equation*}
f^1r_{\nbhd{U},\nbhd{V}}(s_1)+\dots+f^kr_{\nbhd{U},\nbhd{V}}(s_k)=s,
\end{equation*}
as desired.
\end{proof}
\end{lemma}

We can now state the following result which connects the notions of ``locally
finitely generated.''
\begin{proposition}\label{prop:locfingen}
Let\/ $r\in\integernn\cup\{\infty,\omega\}$ and let\/
$\map{\pi}{\man{E}}{\man{M}}$ be a smooth or analytic vector bundle, as
required.  Let\/ $\sF$ be a subsheaf of\/ $\ssections[r]{\man{E}}$ and let\/
$\sM\subset\sections[r]{\man{E}}$ be a submodule.  The following statements
hold:
\begin{compactenum}[(i)]
\item \label{pl:locfingen1} if the subsheaf\/ $\sF$ is locally finitely
generated, then the generalised subbundle\/ $\dist{F}(\sF)$ is locally
finitely generated;
\item \label{pl:locfingen2} if the subsheaf\/ $\sF$ is locally finitely
generated, then, for each\/ $x\in\man{M}$\@, there exists a neighbourhood\/
$\nbhd{U}$ of\/ $x$ such that\/ $F(\nbhd{U})$ is a finitely generated\/ $\func[r]{\nbhd{U}}$-module;
\item \label{pl:locfingen3} if the submodule\/ $\sM$ is locally finitely
generated, then the generalised subbundle\/ $\dist{F}(\sM)$ is locally
finitely generated;
\item \label{pl:locfingen4} if the submodule\/ $\sM$ is locally finitely
generated, then the subsheaf\/ $\sF_{\sM}$ is locally finitely generated.
\end{compactenum}
\begin{proof}
\eqref{pl:locfingen1} Let $x_0\in\man{M}$\@.  Since $\sF$ is locally finitely
generated, there exists a neighbourhood $\nbhd{U}$ of $x_0$ and sections
$\xi_1,\dots,\xi_k\in F(\nbhd{U})$ such that
$\ifam{[\xi_1]_x,\dots,[\xi_k]_x}$ generate $\sF_x$ for each
$x\in\nbhd{U}$\@.  If $e_x\in\dist{F}(\sF)_x$ for $x\in\nbhd{U}$\@, by
definition of $\dist{F}(\sF)$ we have $e_x=\xi(x)$ for $[\xi]_x\in\sF_x$\@.
It follows that $\ifam{\xi_1(x),\dots,\xi_k(x)}$ span $\dist{F}(\sF)_x$\@,
showing that $\dist{F}(\sF)$ is locally finitely generated.

\eqref{pl:locfingen2} This follows from Lemma~\ref{lem:locfingen1} and the
Oka Coherence Theorem, $\sF$ being a finitely generated subsheaf of the
coherent sheaf $\ssections[r]{\man{E}}$\@, and so
coherent~\cite[page~235]{HG/RR:84}\@.

\eqref{pl:locfingen3} Let $x_0\in\man{M}$\@.  Since $\sM$ is locally finitely
generated, there exists a neighbourhood $\nbhd{U}$ of $x_0$ and
$\xi_1,\dots,\xi_k\in\sM$ which generate $F_{\sM}(\nbhd{U})$ as a
$\func[r]{\nbhd{U}}$-module.  If $e_x\in\dist{F}(\sM)_x$ for
$x\in\nbhd{U}$\@, by definition of $\dist{F}(\sM)$ we have $e_x=\xi(x)$ for
$\xi\in\sM$\@.  It follows that $\ifam{\xi_1(x),\dots,\xi_k(x)}$ span
$\dist{F}(\sM)_x$\@, showing that $\dist{F}(\sM)$ is locally finitely
generated.

\eqref{pl:locfingen4} Let $x_0\in\man{M}$\@.  Since $\sM$ is locally finitely
generated, there exists a neighbourhood $\nbhd{U}$ of $x_0$ and
$\xi_1,\dots,\xi_k\in\sM$ which generate $F_{\sM}(\nbhd{U})$ as a
$\func[r]{\nbhd{U}}$-module.  Let $x\in\nbhd{U}$ and let
$[(\xi,\nbhd{V})]_x\in\sF_{\sM,x}$ for some $\nbhd{V}\subset\nbhd{U}$\@.
Thus
\begin{equation*}
\xi=f^1(\eta_1|\nbhd{V})+\dots+f^m(\eta_m|\nbhd{V})
\end{equation*}
for $f^1,\dots,f^m\in\func[r]{\nbhd{V}}$ and $\eta_1,\dots,\eta_m\in\sM$\@.
Let us write
\begin{equation*}
\eta_a|\nbhd{U}=g_a^1\xi_1+\dots+g_a^k\xi_k
\end{equation*}
for $g_a^j\in\func[r]{\nbhd{U}}$\@, $a\in\{1,\dots,m\}$\@,
$j\in\{1,\dots,k\}$\@.  Thus
\begin{equation*}
\xi=\sum_{j=1}^k\sum_{a=1}^mf^j(g^a_j|\nbhd{V})(\xi_a|\nbhd{V}),
\end{equation*}
giving
\begin{equation*}
[\xi]_x=\sum_{j=1}^k\sum_{a=1}^m[f^j]_x[g^a_j]_x[\xi_a]_x,
\end{equation*}
so showing that $\ifam{[\xi_1]_x,\dots,[\xi_k]_x}$ generate $\sF_{\sM,x}$ for
every $x\in\nbhd{U}$\@.
\end{proof}
\end{proposition}

Note that the implications for a generalised subbundle being locally finitely
generated are not present.  This is because this is not interesting in the
following sense.  An analytic generalised subbundle $\dist{F}$ always has the
property that $\ssections[\omega]{\dist{F}}$ is locally finitely generated;
this is Corollary~\ref{cor:subbunfingen} above.  On the other hand, a smooth
generalised subbundle $\dist{F}$ is globally finitely generated by
Theorem~\ref{the:global-span} below, but $\ssections[\infty]{\dist{F}}$ may
not be locally finitely generated as can be seen from
Example~\ref{eg:singular!fingen}\@.

\section{Global generators for generalised
subbundles}\label{sec:global-sections}

One of the potential problems with our definition of a generalised subbundle
is that it relies on generators that are only locally defined.  In the smooth
or finitely differentiable case, these local generators can very often (for
example, when the neighbourhood on which the generators are defined is not
dense) be extended using the Tietsze Extension
Theorem~\cite[Theorem~5.5.9]{RA/JEM/TSR:88} to give global generators.
However, in the real analytic case, matters are more complicated, as
generally a local section simply \emph{cannot} be extended to be globally
defined.  Here one must thus attenuate one's objectives, and merely wonder
whether globally defined generators exist at all.  In this section we address
these questions, considering separately the smooth and finitely
differentiable cases, and the real analytic case.

\subsection{Global generators for $\C^r$-generalised subbundles}

Let us first look at the smooth and finitely differentiable case.  The result
here is due to \citet{HJS:08}\@, and we give the proof here since the result
is an important one.
\begin{theorem}\label{the:global-span}
Let\/ $r\in\integernn\cup\{\infty\}$ and let\/ $\map{\pi}{\man{E}}{\man{M}}$
be a smooth vector bundle whose fibres have bounded dimension and for which\/
$\man{M}$ is a smooth paracompact Hausdorff manifold of bounded dimension.
If\/ $\man{F}\subset\man{E}$ is a generalised subbundle then the following
statements are equivalent:
\begin{compactenum}[(i)]
\item \label{pl:global-span1} $\dist{F}$ is of class\/ $\C^r$\@;
\item \label{pl:global-span2} for each\/ $x_0\in\man{M}$ and each\/
$v_{x_0}\in\dist{F}_{x_0}$\@, there exists a neighbourhood\/ $\nbhd{N}$ of\/
$x_0$ and a\/ $\C^r$-section\/ $\xi\in\sections[r]{\man{E}}$ such that\/
$\xi(x_0)=v_{x_0}$ and\/ $\xi(x)\in\dist{F}_x$ for each\/ $x\in\nbhd{N}$\@;
\item \label{pl:global-span3} there exists a family\/
$\ifam{\xi_1,\dots,\xi_k}$ of\/ $\C^r$-sections on\/ $\man{M}$ such that
\begin{equation*}
\dist{F}_x=\vecspan[\real]{\xi_1(x),\dots,\xi_k(x)}
\end{equation*}
for each\/ $x\in\man{M}$\@.
\end{compactenum}
\begin{proof}
\eqref{pl:global-span1}$\implies$\eqref{pl:global-span2} Suppose that
$\dist{F}$ is of class $\C^r$\@, let $x_0\in\man{M}$\@, and let
$v_{x_0}\in\dist{F}_{x_0}$\@.  Let $\nbhd{N}$ be a neighbourhood of $x_0$ and
let $\ifam{\xi_j}_{j\in J}$ be a family of sections on $\nbhd{N}$ of class
$\C^r$ such that
\begin{equation*}
\dist{F}_x=\vecspan[\real]{\xi_j(x)|\enspace j\in J}
\end{equation*}
for $x\in\nbhd{N}$\@.  Let $j_1,\dots,j_k\in J$ be such that
$\ifam{\xi_{j_1}(x_0),\dots,\xi_{j_k}(x_0)}$ is a basis for $\dist{F}_{x_0}$\@.
Then
\begin{equation*}
v_{x_0}=c_1\xi_{j_1}(x_0)+\dots+c_k\xi_{j_k}(x_0)
\end{equation*}
for some uniquely defined $c_1,\dots,c_k\in\real$\@.  The section
$\xi=c_1\xi_{j_1}+\dots+c_k\xi_{j_k}$ defined on $\nbhd{N}$ is then of class
$\C^r$\@, is $\dist{F}$-valued on $\nbhd{N}$\@, and satisfies
$\xi(x_0)=v_{x_0}$\@.

\eqref{pl:global-span2}$\implies$\eqref{pl:global-span3} Since $\man{E}$ has
bounded fibre dimension there exists a least integer $n_0\in\integernn$ such
that $\rank(\dist{F}_x)\le n_0$ for every $x\in\man{M}$\@.  Use the notation
$\rank(\dist{F}_x)=\rank_{\dist{F}}(x)$\@.  For $k\in\{0,1,\dots,n_0+1\}$
denote
\begin{equation*}
\nbhd{U}_k=\setdef{x\in\man{M}}{\rank_{\dist{F}}(x)\ge k}.
\end{equation*}
By Proposition~\ref{prop:rank-semicont}\@, $\rank_{\dist{F}}$ is lower
semicontinuous, and so $\nbhd{U}_k$ is open for each
$k\in\{0,1,\dots,n_0+1\}$\@.  Moreover,
\begin{equation*}
\emptyset=\nbhd{U}_{n_0+1}\subset\nbhd{U}_{n_0}\subset\dots\subset
\nbhd{U}_1\subset\nbhd{U}_0=\man{M}.
\end{equation*}
We wish to define a certain open cover of each of the open sets
$\nbhd{U}_k$\@, $k\in\{0,1,\dots,n_0+1\}$\@, with a certain property.  We do
this inductively.

The following two general lemmata will be key in our inductive construction.
Note that the notation of the lemmata may or may nor correspond to the
notation of the theorem and its proof.  So beware.
\begin{prooflemma}\label{plem:global-span1}
Let\/ $\man{M}$ be a smooth, paracompact, Hausdorff manifold all of whose
connected components have dimension bounded by\/ $n\in\integernn$ and let\/
$\ifam{\nbhd{W}_j}_{j\in J}$ be an open cover of\/ $\man{M}$\@.  Then there
exists an open cover\/ $\ifam{\nbhd{V}_a}_{a\in A}$ of\/ $\man{M}$ with the
following properties:
\begin{compactenum}[(i)]
\item $\ifam{\nbhd{V}_a}_{a\in A}$ is a refinement of\/
$\ifam{\nbhd{W}_j}_{j\in J}$ (\ie~for each\/ $j\in J$ there exists\/ $a\in A$
such that\/ $\nbhd{W}_j\subset\nbhd{V}_a$);
\item there exist subsets\/ $A_1,\dots,A_{n+1}\subset A$ such that\/
$A=\disjointunion{}_{l=1}^{n+1}A_l$ and such that, whenever\/ $a_1,a_2\in
A_l$ for some\/ $l\in\{1,\dots,n+1\}$\@, it holds that\/
$\nbhd{V}_{a_1}\cap\nbhd{V}_{a_2}=\emptyset$\@.
\end{compactenum}
\begin{subproof}
Since $\man{M}$ is paracompact it possesses a Riemannian metric by
Corollary~5.5.13 of~\cite{RA/JEM/TSR:88}\@.  Therefore, we may assume that
$\man{M}$ is a metric space, and we denote the metric by $\d$\@.  A smooth
manifold can be triangulated by which we mean that there exists a
homeomorphism $\map{\Phi}{S}{\man{M}}$ where $S$ is a union of simplices
which intersect only at their boundaries~\cite[Theorem~8.4]{JRM:66}\@.  By
successive barycentric subdivisions of the simplices of $S$ we can assume
that all simplices $S_b$\@, $b\in B$\@, comprising $S$ are such that
$\Phi(S_b)\subset\nbhd{W}_j$ for some $j\in J$\@.  Let us adopt the usual
slight abuse of terminology and say that $\Phi(S_b)$\@, $b\in B$\@, is a
simplex.  Let $m\in\{0,1,\dots,n\}$ and define
\begin{equation*}
\sF_m=\setdef{F\subset\man{M}}
{F\ \textrm{is an open}\ m\textrm{-dimensional face of some simplex}}
\end{equation*}
and denote $\snorm{\sF_m}=\cup_{F\in\sF_m}F$\@.  For $F\in\sF_m$ denote
\begin{equation*}
C(F)=\closure(\cup\setdef{F'\in\sF_m}{F'\not=F}).
\end{equation*}
Note that $F$ is both open and closed in the relative topology of
$\snorm{\sF_m}$\@.  Also, $C(F)$ is closed in the relative topology of
$\snorm{\sF_m}$ since its intersection with any compact set is a finite union
of closed sets, and so closed.  This implies that $F\cap C(F)=\emptyset$ by
virtue of $F$ being relatively open.  Now let $x\in F$\@.  Then $\{x\}$ and
$C(F)$ are disjoint closed sets, and so
\begin{equation*}
\d(x,C(F))=\inf\setdef{\d(x,y)}{y\in C(F)}
\end{equation*}
is positive.  (Indeed, suppose otherwise.  Then there exists a sequence
$\ifam{y_k}_{k\in\integerp}$ in $C(F)$ converging to $x$\@.  Thus $x\in C(F)$
since $C(F)$ is relatively closed.  This contradicts the fact that $\{x\}$
and $C(F)$ are disjoint.)  Denote by
$\oball{x}{\frac{1}{2}\d(x,C(F))}$ the open ball in $\man{M}$ of
radius $\frac{1}{2}\d(x,C(F))$ and define $B(F)=\cup_{x\in
F}\oball{x}{\frac{1}{2}\d(x,C(F))}$\@.

We claim that if $F_1,F_2\in\sF_m$ are disjoint, then $B(F_1)$ and $B(F_2)$
are disjoint.  Suppose otherwise and let $x\in B(F_1)\cap B(F_2)$\@, Let
$x_1\in F_1$ and $x_2\in F_2$ be such that
$y\in\oball{x_j}{\frac{1}{2}\d(x_j,C(F_j))}$\@, $j\in\{1,2\}$\@.  Then
\begin{align*}
\d(x_1,x_2)\le&\;\d(x_1,x)+\d(x_2,x)<\tfrac{1}{2}\d(x_1,C(F_1))+
\tfrac{1}{2}\d(x_2,C(F_2))\\
\le&\;\max\{\d(x_1,C(F_1)),\tfrac{1}{2}\d(x_2,C(F_2))\}.
\end{align*}
Thus either
\begin{equation*}
\d(x_1,x_2)<\d(x_1,C(F_1))\quad\textrm{or}\quad\d(x_1,x_2)<\d(x_2,C(F_2)).
\end{equation*}
In the first case we have $x_2\in F_1$\@, contradicting the fact that $x_2\in
C(F_1)$\@, and in the second case we have $x_1\in F_2$\@, contradicting the
fact that $x_1\in C(F_2)$\@.  Thus $B(F_1)$ and $B(F_2)$ are disjoint if
$F_1,F_2\in\sF_m$ are disjoint.

For $F\in\sF_m$ let $j_F\in J$ be such that $F\subset\nbhd{W}_{j_F}$\@, this
being possible since our triangulation of $\man{M}$ was chosen in precisely
this manner.  Define an open set $\nbhd{V}_F=B(F)\cap\nbhd{W}_{j_F}$\@.
Note that
\begin{equation*}
F\subset\nbhd{V}_F\subset\nbhd{W}_{j_F}.
\end{equation*}
Moreover, since $B(F_1)$ and $B(F_2)$ are disjoint for $F_1,F_2\in\sF_m$
disjoint, it follows that $\nbhd{V}_{F_1}$ and $\nbhd{V}_{F_2}$ are disjoint
for $F_1,F_2\in\sF_m$ disjoint.  If $x\in\man{M}$ then $x$ belongs to some
open $m$-dimensional face from the triangulation of $\man{M}$ for some
$m\in\{0,1,\dots,n\}$\@.  Therefore,
\begin{equation*}
\man{M}=\bigcup_{m=0}^n\bigcup_{F\in\sF_m}\nbhd{V}_F.
\end{equation*}
Thus we have an open cover of $\man{M}$ that refines $\ifam{\nbhd{W}_j}_{j\in
J}$\@.  The index set for the open cover is the set
\begin{equation*}
A=\setdef{F}{F\in\sF_m,\ m\in\{0,1,\dots,n\}}.
\end{equation*}

It remains to show that the index set for the open cover satisfies the second
condition in the statement of the lemma.  For $m\in\{1,\dots,n+1\}$ let
$A_m=\sF_{m-1}$ so that $A$ is the disjoint union of $A_1,\dots,A_{n+1}$\@.
As we have shown above, for each $m\in\{1,\dots,n+1\}$\@, the family of open
sets $\ifam{\nbhd{V}_F}_{F\in A_m}$ is a pairwise disjoint union of open
sets, just as is asserted in the statement of the lemma.
\end{subproof}
\end{prooflemma}

The following technical lemma will also be useful.  The result is well-known,
but we were not able to find a proof for it, so we provide one here.
\begin{prooflemma}\label{plem:global-span2}
If\/ $\nbhd{U}$ is an open subset of a smooth, paracompact, Hausdorff
manifold\/ $\man{M}$ then there exists\/ $f\in\func{\man{M}}$ such that\/
$f(x)\in\realp$ for all\/ $x\in\nbhd{U}$ and\/ $f(x)=0$ for all\/
$x\in\man{M}\setminus\nbhd{U}$\@.
\begin{subproof}
We shall construct $f$ as the limit of a sequence of smooth functions
converging in the weak $\C^\infty$-topology.  We equip $\man{M}$ with a
Riemannian metric $\metric$\@, this by paracompactness of
$\man{M}$~\cite[Corollary~5.5.13]{RA/JEM/TSR:88}\@.  We denote by $\nabla$
the Levi-Civita connection of $\metric$\@.  Let $g\in\func{\man{M}}$\@.  If
$K\subset\man{M}$ is compact and if $r\in\integernn$\@, we define
\begin{equation*}
\dnorm{g}_{r,K}=\sup\setdef{\dnorm{\nabla^jg(x)}}
{x\in K,\ j\in\{0,1,\dots,r\}},
\end{equation*}
where $\dnorm{\cdot}$ indicates the norm induced on tensors by the norm
associated with the Riemannian metric.  One readily sees that the family of
seminorms $\dnorm{\cdot}_{r,K}$\@, $r\in\integernn$\@, $K\subset\man{M}$
compact, defines a locally convex topology agreeing with other definitions of
the weak topology.  Thus, if a sequence $\ifam{g_j}_{j\in\integerp}$
satisfies
\begin{equation*}
\lim_{j\to\infty}\dnorm{g-g_j}_{r,K}=0,\qquad
r\in\integernn,\ K\subset\man{M}\ \textrm{compact},
\end{equation*}
then $g$ is infinitely differentiable~\cite[\S4.3]{PWM:80}\@.

We suppose that $\man{M}$ is connected since, if it is not, we can construct
$f$ for each connected component, which suffices to give $f$ on $\man{M}$\@.
If $\man{M}$ is paracompact, connectedness allows us to conclude that
$\man{M}$ is second countable~\cite[Proposition~5.5.11]{RA/JEM/TSR:88}\@.
Using Lemma~2.76 of~\cite{CDA/KCB:06}\@, we let $\ifam{K_j}_{j\in\integerp}$
be a sequence of compact subsets of $\nbhd{U}$ such that
$K_j\subset\interior(K_{j+1})$ for $j\in\integerp$ and such that
$\cup_{j\in\integerp}K_j=\nbhd{U}$\@.  For $j\in\integerp$ let
$\map{g_j}{\man{M}}{\interval[0,1]}$ be a smooth function such that
$g_j(x)=1$ for $x\in K_j$ and $g_j(x)=0$ for $x\in\man{M}\setminus
K_{j+1}$\@; see~\cite[Proposition~5.5.8]{RA/JEM/TSR:88}\@.  Let us define
$\alpha_j=\dnorm{g_j}_{j,K_{j+1}}$ and take $\epsilon_j\in\realp$ to satisfy
$\epsilon_j<(\alpha_j2^j)^{-1}$\@.  We define $f$ by
\begin{equation*}
f(x)=\sum_{j=1}^\infty\epsilon_jg_j(x),
\end{equation*}
and claim that $f$ as defined satisfies the conclusions of the lemma.  First
of all, since each of the functions $g_j$ takes values in $\interval[0,1]$ we
have
\begin{equation*}
\snorm{f(x)}\le\sum_{j=1}^\infty\snorm{\epsilon_jg_j(x)}\le
\sum_{j=1}^\infty\epsilon_j\dnorm{g_j}_{0,K_{j+1}}\le
\sum_{j=1}^\infty\epsilon_j\dnorm{g_j}_{j,K_{j+1}}\le
\sum_{j=1}^\infty\frac{1}{2^j}\le1,
\end{equation*}
and so $f$ is well-defined.  If $x\in\nbhd{U}$ then there exists
$N\in\integerp$ such that $x\in K_N$\@.  Thus $g_N(x)=1$ and so
$f(x)\in\realp$\@.  If $x\in\man{M}\setminus\nbhd{U}$ then $g_j(x)=0$ for all
$j\in\integerp$ and so $f(x)=0$\@.  All that remains to show is that $f$ is
infinitely differentiable.

Let $x\in\man{M}$\@, let $m\in\integerp$\@, and let $j\in\integernn$ be such
that $j\le m$\@.  If $x\not\in K_{m+1}$ then $g_m$ is zero in a neighbourhood
of $x$\@, and so $\dnorm{\nabla^jg_m(x)}=0$\@.  If $x\in K_{m+1}$ then
\begin{align*}
\dnorm{\nabla^jg_m(x)}\le&\;\sup\setdef{\dnorm{\nabla^jg_m(x')}}
{x'\in K_{m+1}}\\
\le&\;\sup\setdef{\dnorm{\nabla^jg_m(x')}}
{x'\in K_{m+1},\ j\in\{0,1,\dots,m\}}=\alpha_m.
\end{align*}
Thus, whenever $j\le m$ we have $\dnorm{\nabla^jg_m(x)}\le\alpha_m$ for every
$x\in\man{N}$\@.

Let us define $f_m\in\func{\man{M}}$ by
\begin{equation*}
f_m(x)=\sum_{j=1}^m\epsilon_jg_j(x).
\end{equation*}
Let $K\subset\man{M}$ be compact, let $r\in\integernn$\@, and let
$\epsilon\in\realp$\@.  Take $N\in\integerp$ sufficiently large that
\begin{equation*}
\sum_{m=m_1+1}^{m_2}\frac{1}{2^m}<\epsilon,
\end{equation*}
for $m_1,m_2\ge N$ with $m_1<m_2$\@, this being possible by convergence of
$\sum_{j=1}^\infty\frac{1}{2^j}$\@.  Then, for $m_1,m_2\ge N$\@,
\begin{align*}
\dnorm{f_{m_1}-f_{m_2}}_{r,K}=&\;
\sup\setdef{\dnorm{\nabla^jf_{m_1}(x)-\nabla^jf_{m_2}(x)}}
{x\in K,\ j\in\{0,1,\dots,r\}}\\
=&\;\sup\asetdef{\adnorm{\sum_{m=m_1+1}^{m_2}\epsilon_m\nabla^jg_m(x)}}
{x\in K,\ j\in\{0,1,\dots,r\}}\\
\le&\;\sup\asetdef{\sum_{m=m_1+1}^{m_2}\epsilon_m\dnorm{\nabla^jg_m(x)}}
{x\in K,\ j\in\{0,1,\dots,r\}}
\le\sum_{m_1+1}^{m_2}\frac{1}{2^m}<\epsilon.
\end{align*}
Thus, for every $r\in\integernn$ and $K\subset\man{M}$ compact,
$\ifam{f_m}_{m\in\integerp}$ is a Cauchy sequence in the norm
$\dnorm{\cdot}_{r,K}$\@.  Completeness of the weak $\C^\infty$-topology
implies that the sequence $\ifam{f_m}_{m\in\integerp}$ converges to a
function that is infinitely differentiable.
\end{subproof}
\end{prooflemma}

With the use of the preceding lemmata, we can prove the following, reverting
now to the notation of the theorem and its proof.
\begin{prooflemma}\label{plem:global-span3}
There exist\/ $\C^r$-sections\/ $\xi^k_m$\@,\/ $k\in\{0,1,\dots,n_0\}$\@,\/
$m\in\{1,\dots,n+1\}$\@, on\/ $\man{M}$\@, taking values in $\dist{F}$\@,
such that, for every\/ $k\in\{0,1,\dots,n_0\}$ and for every
$x\in\nbhd{U}_k$\@,
\begin{equation*}
\dim(\vecspan[\real]{\setdef{\xi^j_m(x')}
{j\in\{0,1,\dots,k\},\ m\in\{1,\dots,n+1\}}})\ge k.
\end{equation*}
\begin{subproof}
We prove this by induction on $k$\@.  For $k=0$ the assertion is obvious.
Indeed, one merely takes the sections $\xi^0_1,\dots,\xi^0_{n+1}$ to all be
zero, and the conclusion holds in this case.  Now assume that the conclusions
hold for $k\in\{0,1,\dots,s\}$\@.  Let $x\in\nbhd{U}_{s+1}$ and define
\begin{equation*}
\alg{U}^s_x=\vecspan[\real]{\setdef{\xi^k_m(x)}
{k\in\{0,1,\dots,s\},\ m\in\{1,\dots,n+1\}}},
\end{equation*}
noting that $\dim(\alg{U}^s_x)\ge s$ by the induction hypothesis.  Since
$x\in\nbhd{U}_x$ there exists $v_x\in\dist{F}_x$ such that
\begin{equation*}
\dim(\vecspan[\real]{\alg{U}^s_x\cup\{v_x\}})\ge s+1.
\end{equation*}
By the hypotheses of part~\eqref{pl:global-span2} of the theorem, there
exists a neighbourhood $\nbhd{W}^{s+1}_x$ of $x$ and a $\C^r$-section $\xi_x$
on $\nbhd{W}^{s+1}_x$\@, taking values in $\dist{F}$\@, such that
$\xi_x(x)=v_x$\@.  By multiplying $\xi_x$ by a smooth function that is
positive on $\nbhd{W}^{s+1}_x$ and with compact support
\cite[see][Theorem~5.5.9]{RA/JEM/TSR:88} we can extend $\xi_x$ to a
$\dist{F}$-valued $\C^r$-section on all of $\man{M}$\@.  By shrinking
$\nbhd{W}^{s+1}_x$ if necessary, we can suppose that
\begin{equation*}
\dim(\vecspan[\real]{\alg{U}^s_{x'}\cup\{\xi_x(x')\}})\ge s+1,
\qquad x'\in\nbhd{W}^{s+1}_x,
\end{equation*}
by lower semicontinuity of the rank of a generalised subbundle
(Proposition~\ref{prop:rank-semicont}).  With $\nbhd{W}^{s+1}_x$ so specified
we have $\nbhd{W}^{s+1}_x\subset\nbhd{U}_{s+1}$\@.  Note that
$\ifam{\nbhd{W}^{s+1}_x}_{x\in\nbhd{U}_{s+1}}$ is then an open cover of
$\nbhd{U}_k$\@.

By Lemma~\ref{plem:global-span1} let $\ifam{\nbhd{V}^{s+1}_a}_{a\in A}$ be a
refinement of the open cover $\ifam{\nbhd{W}^{s+1}_x}_{x\in\nbhd{U}_{s+1}}$
of $\nbhd{U}_{s+1}$ such that the index set $A$ is a disjoint union of sets
$A_1,\dots,A_{n+1}$ such that
$\nbhd{V}^{s+1}_{a_1}\cap\nbhd{V}^{s+1}_{a_2}=\emptyset$ whenever
$a_1,a_2\in A_l$ are distinct for some $l\in\{1,\dots,n+1\}$\@.  Let us
denote $\nbhd{V}^{s+1}_l=\cup_{a\in A_l}\nbhd{V}^{s+1}_a$\@,
$l\in\{1,\dots,n+1\}$\@.  For $l\in\{1,\dots,m+1\}$ and $a\in A_l$\@, let
$x_{l,a}\in\nbhd{U}_k$ be such that
$\nbhd{V}^{s+1}_a\subset\nbhd{W}^{s+1}_{x_{l,a}}$\@.  Define a $\C^r$-section
$\eta^{s+1}_l$ on $\nbhd{V}^{s+1}_l$ by asking that, if
$x\in\nbhd{V}^{s+1}_a$ for some (necessarily unique) $a\in A_l$\@, then
$\eta^{s+1}_l(x)=\xi_{x_{l,a}}(x)$\@.  The section $\eta^{s+1}_l$ will then
have the property that
\begin{equation}\label{eq:global-span1}
\dim(\vecspan[\real]{\alg{U}^s_x\cup\{\eta^{s+1}_l(x)\}})\ge s+1,
\qquad x\in\nbhd{V}_l,\ l\in\{1,\dots,n+1\}.
\end{equation}
For each $l\in\{1,\dots,n+1\}$\@, by Lemma~\ref{plem:global-span2} let
$f_l\in\func{\man{M}}$ be such that $f_l(x)\in\realp$ for $x\in\nbhd{V}_l$
and such that $f_l(x)=0$ for $x\in\man{M}\setminus\nbhd{V}_l$\@.  Then define
$\xi^{s+1}_l=f_l\eta^{s+1}_l$ so that
\begin{equation}\label{eq:global-span2}
\dim(\vecspan[\real]{\alg{U}^s_x\cup\{\xi^{s+1}_l(x)\}})\ge s+1,
\qquad x\in\nbhd{V}_l,\ l\in\{1,\dots,n+1\},
\end{equation}
by~\eqref{eq:global-span1}\@, since $\xi^{s+1}_l(x)$ is a nonzero multiple of
$\eta^{s+1}_l(x)$ for all $x\in\nbhd{V}_l$\@.  To complete the proof of the
lemma, let $x\in\nbhd{U}_{s+1}$ and let $a\in A$ be such that
$x\in\nbhd{V}_a$\@.  Then $a\in A_l$ for some unique $l\in\{1,\dots,n+1\}$
and so $x\in\nbhd{V}_l$\@.  Therefore, by~\eqref{eq:global-span2} and
recalling the definition of $\alg{U}^s_x$\@,
\begin{align*}
&\vecspan[\real]{\alg{U}^s_x\cup\{\xi^{s+1}_l(x)\}}\subset
\vecspan[\real]{\setdef{\xi^j_m(x)}
{j\in\{0,1,\dots,s+1\},\ m\in\{1,\dots,n+1\}}}\\
\implies\quad&\dim(\vecspan[\real]{\setdef{\xi^j_m(x)}
{j\in\{0,1,\dots,s+1\},\ m\in\{1,\dots,n+1\}}})\ge s+1,
\end{align*}
as desired.
\end{subproof}
\end{prooflemma}

We now conclude the proof of this part of the theorem by showing that the
$\dist{F}$-valued $\C^r$-sections
\begin{equation*}
\ifam{\xi_m^j|\enspace j\in\{1,\dots,n_0\},\ m\in\{1,\dots,n+1\}}
\end{equation*}
from Lemma~\ref{plem:global-span3} generate $\dist{F}$\@.  Let $x\in\man{M}$
and let $k=\rank_{\dist{F}}(x)$ so that $x\in\nbhd{U}_k$\@.  By
Lemma~\ref{plem:global-span3} we have
\begin{multline*}
\dim(\vecspan[\real]{\setdef{\xi^j_m(x)}
{j\in\{1,\dots,n_0\},\ m\in\{1,\dots,n+1\}}})\\
\ge\dim(\vecspan[\real]{\setdef{\xi^j_m(x)}
{j\in\{1,\dots,k\},\ m\in\{1,\dots,n+1\}}})\ge k.
\end{multline*}
However, since $\dim(\dist{F}_x)=k$ and since all sections $\xi^j_m$\@,
$j\in\{1,\dots,n_0\}$\@, $m\in\{1,\dots,n+1\}$\@, are $\dist{F}$-valued, we
conclude that
\begin{equation*}
\vecspan[\real]{\setdef{\xi^j_m(x)}
{j\in\{1,\dots,n_0\},\ m\in\{1,\dots,n+1\}}}=\dist{F}_x,
\end{equation*}
as desired.

\eqref{pl:global-span3}$\implies$\eqref{pl:global-span1} This is obvious.
\end{proof}
\end{theorem}

Note that the theorem does not say that the sections $\xi_1,\dots,\xi_k$ from
part~\eqref{pl:global-span3} generate the stalks of the sheaf
$\ssections[r]{\dist{F}}$\@.  Indeed, as in
Example~\ref{eg:singular!fingen}\@, the stalks of a subsheaf of
$\C^r$-sections, $r\in\integernn\cup\{\infty\}$\@, are not generally finitely
generated.

\subsection{Global generators for real analytic generalised subbundles}

In the real analytic case, we cannot prove that there are generally finitely
many global generators.  However, we can use Cartan's Theorem~A, stated in
Theorem~\ref{the:theoremA}\@, to give local generation by global sections.
This is, of course, highly nontrivial since generally a local section of a
real analytic vector bundle cannot be extended to a global section.
\begin{theorem}\label{the:global-sections}
Let\/ $\map{\pi}{\man{E}}{\man{M}}$ be a real analytic vector bundle for
which\/ $\man{M}$ is paracompact and Hausdorff.  If\/
$\dist{F}\subset\man{E}$ is a generalised subbundle, then the following
statements are equivalent:
\begin{compactenum}[(i)]
\item \label{pl:analytic-gen1} $\dist{F}$ is real analytic;
\item \label{pl:analytic-gen2} for each\/ $x_0\in\man{M}$ and each\/
$v_{x_0}\in\dist{F}_{x_0}$\@, there exists a real analytic section\/
$\xi\in\sections[\omega]{\dist{F}}$ over\/ $\man{M}$ such that\/
$\xi(x_0)=v_{x_0}$\@;
\item \label{pl:analytic-gen3} for each\/ $x_0\in\man{M}$ there exists a
neighbourhood\/ $\nbhd{N}$ of\/ $x_0$ and a real analytic sections\/
$\xi_1,\dots,\xi_k\in\sections[\omega]{\dist{F}}$ over\/ $\man{M}$ such that
\begin{equation*}
\dist{F}_x=\vecspan[\real]{\xi_1(x),\dots,\xi_k(x)}
\end{equation*}
for each\/ $x\in\nbhd{N}$\@.
\end{compactenum}
\begin{proof}
\eqref{pl:analytic-gen1}$\implies$\eqref{pl:analytic-gen2} By
Proposition~\ref{prop:fibres} let
$V_{x_0}\in(\gfunc[\omega]{x_0}{\man{M}}/\mathfrak{m}_{x_0})
\otimes_{\gfunc[\omega]{x_0}{\man{M}}}\gsections[\omega]{x_0}{\dist{F}}$ be
such that $\tilde{\iota}_{x_0}(V_{x_0})=v_{x_0}$\@, where
$\tilde{\iota}_{x_0}$ is as in the proof of Proposition~\ref{prop:fibres}\@.
We write
\begin{equation*}
V_{x_0}=\sum_{j=1}^k([f^j]_{x_0}+\mathfrak{m}_{x_0})\otimes[\xi_j]_{x_0}
\end{equation*}
for some $[f^1]_{x_0},\dots,[f^k]_{x_0}\in\gfunc[\omega]{x_0}{\man{M}}$ and
$[\xi_1]_{x_0},\dots,[\xi_k]_{x_0}\in\gsections[\omega]{x_0}{\dist{F}}$\@.
By Corollary~\ref{cor:analytic-coherent} we know that
$\ssections[\omega]{\dist{F}}$ is coherent.  Therefore, by Cartan's Theorem~A
(Theorem~\ref{the:theoremA}), there exist global sections
$\ifam{\beta_a}_{a\in A}$ of $\dist{F}$ such that the germs $[\beta_a]$\@,
$a\in A$\@, generate $\gsections[\omega]{x_0}{\dist{F}}$ as a
$\gfunc[\omega]{x_0}{\man{M}}$-module.  It follows that we can write
\begin{equation*}
V_{x_0}=\sum_{a\in A}([g^a]_{x_0}+\mathfrak{m}_{x_0})\otimes[\beta_a]_{x_0},
\end{equation*}
where all but finitely many terms in the sum are zero.  Now define
\begin{equation*}
\xi=\sum_{a\in A}g^a(x_0)\xi_a,
\end{equation*}
again noting that all but finitely many terms in the sum are zero.  Since
\begin{equation*}
v_{x_0}=\tilde{\iota}_{x_0}(V_{x_0})=
\tilde{\iota}_{x_0}\left(\sum_{a\in A}([g^a]_{x_0}+\mathfrak{m}_{x_0})\otimes
[\beta_a]_{x_0}\right)=\sum_{a\in A}g^a(x_0)\xi_a(x_0)
\end{equation*}
(recalling from the proof of Proposition~\ref{prop:fibres} the definition of
$\tilde{\iota}_{x_0}$), it follows that $\xi(x_0)=v_{x_0}$\@, as desired.

\eqref{pl:analytic-gen2}$\implies$\eqref{pl:analytic-gen3} Let us take the
index set $A=\dist{F}$ and for $a\in A$ let $\xi_a$ be a real analytic global
section of $\dist{F}$ such that $\xi_a(x)=a$\@.  Let
$\alg{A}\subset\gsections[\omega]{x_0}{\dist{F}}$ be the submodule generated
by $[\xi_a]_{x_0}$\@, $a\in A$\@.  By Lemma~\ref{plem:fingen} from the proof
of Theorem~\ref{the:analytic-fingen}\@, there exists $a_1,\dots,a_k\in A$ and
a compact set $C\subset\man{M}$ such that $x_0\in\interior(C)$ and such that
$[\xi_{a_1}]_C,\dots,[\xi_{a_k}]_C$ generate $\rho^{-1}_{C,x_0}(\alg{A})$\@.
Then, as we saw in the proof of Theorem~\ref{the:analytic-fingen}\@,
$[\xi_{a_1}]_x,\dots,[\xi_{a_k}]_x$ generate
$\gsections[\omega]{x}{\dist{F}}$ for each $x\in\nbhd{N}$ for any
neighbourhood $\nbhd{N}\subset C$\@.  It follows that
\begin{equation*}
\dist{F}_x=\vecspan[\real]{\xi_{a_1}(x),\dots,\xi_{a_k}(x)}
\end{equation*}
for each $x\in\nbhd{N}$\@.

\eqref{pl:analytic-gen3}$\implies$\eqref{pl:analytic-gen1} This is obvious.
\end{proof}
\end{theorem}

\subsection{Swan's Theorem for regular generalised subbundles}

In this section we consider a different sort of result for global generators
when the subbundle is regular.  In this case, regularity provides additional
algebraic structure for the space of sections.  The result was proved for
$r=0$ and for a compact base by \citet{RGS:62}\@.  An analogous result for
vector bundles over an algebraic variety over an algebraically closed
field,~\eg~in the case of a holomorphic algebraic variety, was proved by
\citet{JPS:55}\@.  The result we give is in the $\C^r$-case for
$r\in\integernn\cup\{\infty,\omega\}$\@.

Let us first recall some notions from commutative algebra.  We let $\alg{R}$
be a commutative unit ring and let $\alg{A}$ be a unital $\alg{R}$-module.
The module $\alg{A}$ is \defn{projective} if there exists~(1)~a free module
$\alg{B}=\oplus_{i\in I}\alg{R}$\@,~(2)~a submodule
$\alg{C}\subset\alg{B}$\@, and~(3)~an injective homomorphism
$\map{\phi}{\alg{A}}{\alg{B}}$ for which
$\alg{B}=\image(\phi)\oplus\alg{C}$\@.  In brief, $\alg{A}$ is a direct
summand of a free module.

With this algebraic definition, we have the following result.
\begin{theorem}\label{the:swan}
Let\/ $r\in\integernn\cup\{\infty,\omega\}$ and let\/
$\map{\pi}{\man{E}}{\man{M}}$ be a smooth or real analytic vector bundle, as
required, whose fibres have bounded dimension and for which\/ $\man{M}$ is a
smooth paracompact Hausdorff manifold of bounded dimension.  The following
statements hold:
\begin{compactenum}[(i)]
\item \label{pl:swan1} if\/ $\dist{F}$ is a regular generalised subbundle
of\/ $\man{E}$ of class\/ $\C^r$\@, then\/ $\sections[r]{\dist{F}}$ is a
finitely generated projective module over\/ $\func[r]{\man{M}}$\@; that is to
say,\/ $\sections[r]{\dist{F}}$ is a direct summand of a finitely generated
free module over\/ $\func[r]{\man{M}}$\@;
\item \label{pl:swan2} if\/ $\sM$ is a finitely generated projective module
over\/ $\func[r]{\man{M}}$ then\/ $\sM$ is isomorphic to the module\/
$\sections[r]{\dist{F}}$ of\/ $\C^r$-sections of a\/ $\C^r$-generalised
subbundle of\/ $\man{E}$\@.
\end{compactenum}
\begin{proof}
\eqref{pl:swan1} Since $\dist{F}$ is a vector bundle of class $\C^r$\@, we
shall without loss of generality suppose that $\dist{F}=\man{E}$ and that
$\man{E}$ is of class $\C^r$\@.  The proof of this part of the theorem breaks
into two parts.  Were we to have access to an embedding theorem for
$\C^0$-manifolds, this could be averted, but we are not aware of such a
theorem.  In any case, we give separate proofs for $r\in\integernn$ and for
$r\in\{\infty,\omega\}$\@.  We comment that both proofs work for
$r\in\integerp\cup\{\infty\}$\@, but the first part applies to $r=0$ and the
second part applies to $r=\omega$\@.

We first consider $r\in\integernn$\@.  By Theorem~\ref{the:global-span}\@,
let $\xi_1,\dots,\xi_k$ be globally defined generators for $\man{E}$\@.  Let
$\real^k_{\man{M}}$ denote the trivial vector bundle $\man{M}\times\real^k$
and define a vector bundle map $\map{\Psi}{\real^k_{\man{M}}}{\man{E}}$ by
\begin{equation*}
\Psi(x,(v_1,\dots,v_k))=v_1\xi_1(x)+\dots+v_k\xi_k(x).
\end{equation*}
Clearly $\Psi$ is surjective and $\ker(\Psi)$ is a $\C^r$-subbundle of
$\real^k_{\man{M}}$~\cite[Proposition~3.4.18]{RA/JEM/TSR:88}\@.  Let
$\inprod{\cdot}{\cdot}$ be the standard inner product on $\real^k$ which we
think of as a vector bundle metric on $\real^k_{\man{M}}$\@.  Define
$\dist{G}_x$ to be the orthogonal complement to $\ker(\Psi_x)$\@.  Note that
$\dist{G}$ is then a $\C^r$-subbundle of $\real^k_{\man{M}}$\@.  We claim
that $\Psi|\dist{G}$ is a $\C^r$-vector bundle isomorphism onto $\man{E}$\@.
Certainly $\Psi$ is a $\C^r$-vector bundle map.  We claim that
$\Psi|\dist{G}$ is onto $\man{E}$\@.  Indeed, let $x\in\man{M}$ and let
$v_x\in\man{E}_x$\@.  Then, since $\image(\Psi)=\man{E}$\@, there exists
$(x,\vect{v})\in\real^k_{\man{M}}$ such that $\Psi(x,\vect{v})=v_x$\@.  Write
$(x,\vect{v})=(x,\vect{u}+\vect{w})$ for $(x,\vect{u})\in\ker(\Psi)$ and
where $\vect{w}$ is orthogonal to
$\vect{u}$\@,~\ie~$(x,\vect{w})\in\dist{G}_x$\@.  Then
\begin{equation*}
v_x=\Psi(x,\vect{v})=\Psi(x,\vect{u})+\Psi(x,\vect{w})=\Psi(x,\vect{w}),
\end{equation*}
and so $\Psi|\dist{G}$ is onto $\man{E}$\@.  To show that $\Psi$ is
injective, suppose that $\Psi(x,\vect{v}_1)=\Psi(x,\vect{v}_2)$ for
$(x,\vect{v}_1),(x,\vect{v}_2)\in\dist{G}_x$\@.  Then
\begin{equation*}
\Psi(x,\vect{v}_2-\vect{v}_1)=0_x\quad\implies
\quad(x,\vect{v}_2-\vect{v}_1)\in
\ker(\Psi),\quad\implies\quad\vect{v}_2-\vect{v}_1=\vect{0},
\end{equation*}
as desired.  Now we recall that $\C^r$-vector bundles over $\man{M}$ are
isomorphic if and only if their sets of sections are isomorphic as
$\func[r]{\man{M}}$-modules,~\cf~\cite[\S6]{EN:67}\@.  Thus
$\sections[r]{\man{E}}$ and $\sections[r]{\dist{G}}$ are isomorphic.

To complete this part of the proof, note that
$\sections[r]{\real^k_{\man{M}}}$ is a finitely generated free module over
$\func[r]{\man{M}}$\@.  To see this, one can easily show that the sections
$x\mapsto(x,\vect{e}_j)$\@, $j\in\{1,\dots,k\}$\@, form a basis for all
sections, where $\vect{e}_j\in\real^k$ is the $j$th standard basis vector.
Moreover, any section $\xi$ of $\real^k_{\man{M}}$ with
$\xi(x)=(x,\vect{v}(x))$ can be written uniquely as
\begin{equation*}
\xi(x)=(x,\vect{u}(x))+(x,\vect{w}(x))
\end{equation*}
with $(x,\vect{u}(x))\in\ker(\Psi_x)$ and where $\vect{w}(x)$ is orthogonal
to $\vect{u}(x)$\@.  That is, $\xi=\eta+\zeta$ for
$\eta\in\sections[r]{\ker(\Psi)}$ and $\zeta\in\sections[r]{\dist{E}}$\@, and
this decomposition is unique.  That is to say,
$\sections[r]{\real^k_{\man{M}}}=\sections[r]{\ker(\Psi)}\oplus
\sections[r]{\dist{E}}$\@, the direct sum being one of
$\func[r]{\man{M}}$-modules.  This gives this part of the theorem for
$r\in\integernn$\@.

For $r\in\{\infty,\omega\}$ we first note that there is a $\C^r$-proper
embedding $\map{\Phi}{\man{E}}{\real^N}$ into Euclidean space for some
$N\in\integerp$\@.  For $r=\infty$ this is the Whitney Embedding
Theorem~\cite[Theorem~6.9]{JML:03}\@.  In the real analytic case, this is the
Grauert Embedding Theorem~\cite{HG:58}\@.  Note that
$\tf{\Phi}$ is an injective vector bundle morphism over $\Phi$\@.  Now let
$\zsec{\man{E}}\subset\man{E}$ be the zero section, and note that
$\zsec{\man{E}}$ is canonically diffeomorphic to $\man{M}$\@.  Thus the
restricted vector bundle $\tb{\man{E}}|\zsec{\man{E}}$ is to be regarded as a
vector bundle over $\man{M}$\@.  Since the fibres of $\man{E}$ intersect
$\zsec{\man{E}}$ transversally at points $0_x$ of $\zsec{\man{E}}$\@, we have
\begin{equation*}
\tb[0_x]{\man{E}}\simeq\tb[x]{\man{M}}\oplus\tb[0_x]{\man{E}_x}\simeq
\tb[x]{\man{M}}\oplus\man{E}_x.
\end{equation*}
Thus $\tb{\man{E}}|\zsec{\man{E}}\simeq\tb{\man{M}}\oplus\man{E}$ and so
$\tf{\Phi}$\@, when restricted to $\tb{\man{E}}|\zsec{\man{E}}$\@, has image
as a vector bundle over the submanifold $\Phi(\man{M})$\@.  This is,
moreover, a subbundle of the trivial bundle $\Phi(\man{M})\times\real^N$\@.
Thus we have $\man{E}$ as a subbundle of the trivial bundle
$\real^N_{\man{M}}\eqdef\man{M}\times\real^N$\@.  Let $\inprod{\cdot}{\cdot}$
be the standard inner product on $\real^N$ which we think of as a vector
bundle metric on $\real^N_{\man{M}}$\@.  Define $\dist{G}_x$ to be the
orthogonal complement to $\man{E}_x$\@, noting that $\dist{G}$ is then a
$\C^r$-subbundle of $\real^N_{\man{M}}$ and that
$\real^N_{\man{M}}=\man{E}\oplus\man{G}$\@.  Let
$\map{\pi_1}{\real^N_{\man{M}}}{\man{E}}$ and
$\map{\pi_2}{\real^N_{\man{M}}}{\man{G}}$ be the projections, thought of as
vector bundle morphisms.  Note that $\sections[r]{\real^N_{\man{M}}}$ is
isomorphic, as a $\func[r]{\man{M}}$-module, to $\func[r]{\man{M}}^N$\@.
Moreover, the map from $\sections[r]{\real^N_{\man{M}}}$ to
$\sections[r]{\man{E}}\oplus\sections[r]{\man{G}}$ given by
\begin{equation*}
\vect{\xi}\mapsto(\pi_1\scirc\vect{\xi})\oplus(\pi_2\scirc\vect{\xi})
\end{equation*}
can be directly verified to be an isomorphism of
$\func[r]{\man{M}}$-modules.  In particular, $\sections[r]{\man{E}}$ is a
summand of the free, finitely generated $\func[r]{\man{M}}$-module $\sections[r]{\real^N_{\man{M}}}$\@.

\eqref{pl:swan2} By definition, there exists a module $\sN$ over
$\func[r]{\man{M}}$ such that
\begin{equation*}
\sM\oplus\sN\simeq\underbrace{\func[r]{\man{M}}\oplus\dots\oplus
\func[r]{\man{M}}}_{k\ \textrm{factors}}.
\end{equation*}
The direct sum on the right is naturally isomorphic to the set of sections of
the trivial vector bundle $\real^k_{\man{M}}=\man{M}\times\real^k$\@.  Thus
we can write $\sM\oplus\sN=\sections[r]{\real^k_{\man{M}}}$\@.  For
$a\in\{1,2\}$\@, let
$\map{\Pi_a}{\sections[r]{\real^k_{\man{M}}}}
{\sections[r]{\real^k_{\man{M}}}}$ be the projection onto the $a$th factor.
As per~\cite[\S6]{EN:67} (essentially), associated with $\Pi_a$ is a vector
bundle map $\map{\pi_a}{\real^k_{\man{M}}}{\real^k_{\man{M}}}$\@.  Since
$\Pi_a\scirc\Pi_a=\Pi_a$ (by virtue of $\Pi_a$ being a projection),
$\pi_a\scirc\pi_a=\pi_a$\@.  To show that $\sM$ is the set of sections of a
vector subbundle of $\real^k_{\man{M}}$ it suffices to show that $\pi_1$ has
locally constant rank.  Following along the lines of the proof of
Proposition~\ref{prop:rank-semicont} one can show that
$x\mapsto\rank(\pi_{a,x})$ is lower semicontinuous for $a\in\{1,2\}$\@.
However, since $\rank(\pi_{1,x})=\rank(\pi_{2,x})=k$ for all $x\in\man{M}$\@,
if $x\mapsto\rank(\pi_{1,x})$ is lower semicontinuous at $x_0$\@, then
$x\mapsto\rank(\pi_{2,x})$ is upper semicontinuous at $x_0$\@.  Thus we
conclude that both of these functions must be continuous at $x_0$\@.  Since
$x\mapsto\rank(\pi_{1,x})$ is integer-valued, it must therefore be locally
constant.
\end{proof}
\end{theorem}

\section{Differential constructions associated to
distributions}\label{sec:differential}

In this section we turn, for the first time, exclusively to
distributions,~\ie~generalised subbundles of the tangent bundle
$\tb{\man{M}}$ of a manifold $\man{M}$\@.  Here one has some particular
features of the tangent bundle that come into play, mainly associated with
the fact that sections of the tangent bundle are vector fields, which have
associated to them a great deal of additional and important structure.  It is
these features that we focus on in this section.  We begin by indicating how
standard vector field operations interact with our sheaf formalism for
distributions.  Then we look at two structures associated to
distributions:~invariant distributions and distributions arising from Lie
algebras of vector fields.

\subsection{Local diffeomorphisms and flows of vector
fields}\label{subsec:local-diffeos}

A standard assumption made in the differential geometry literature when
dealing with vector fields is that they are complete,~\ie~if $X$ is a vector
field on $\man{M}$\@, the flow $x\mapsto\flow{X}{t}(x)$ is defined for every
$(t,x)\in\real\times\man{M}$\@.  In order for the results in the paper to be
as general as possible, we shall not make this assumption.  Thus we have to
introduce some technicalities to deal with this, following \citet{HJS:73}\@.
The first is the notion of a local diffeomorphism, a notion which we will
encounter again in our treatment of the Orbit Theorem in
Section~\ref{sec:orbit-theorem}\@.
\begin{definition}
Let $r\in\{\infty,\omega\}$ and let $\man{M}$ be a $\C^r$-manifold.
\begin{compactenum}[(i)]
\item A \defn{$\C^r$-local diffeomorphism} on $\man{M}$ is a pair
$(\Phi,\nbhd{U})$ where $\nbhd{U}\subset\man{M}$ is (a possibly empty) open
subset called the \defn{domain} and where
$\map{\Phi}{\nbhd{U}}{\Phi(\nbhd{U})}$ is a $\C^r$-diffeomorphism.  The
\defn{image} of $(\Phi,\nbhd{U})$ is the open set $\Phi(\nbhd{U})$\@.
\item If $(\Phi,\nbhd{U})$ and $(\Psi,\nbhd{V})$ are $\C^r$-local
diffeomorphisms, their \defn{composition}
$(\Psi,\nbhd{V})\scirc(\Phi,\nbhd{U})$ is the $\C^r$-local diffeomorphism
$(\Psi\scirc\Phi|\Phi^{-1}(\nbhd{V}),\Phi^{-1}(\nbhd{V}))$\@.
\item If $(\Phi,\nbhd{U})$ is a $\C^r$-local diffeomorphism, its
\defn{inverse} $(\Phi,\nbhd{U})^{-1}$ is the $\C^r$-local diffeomorphism
$(\Phi^{-1},\Phi(\nbhd{U}))$\@.
\item A \defn{group of $\C^r$-local diffeomorphisms} is a family $\sG$ of
$\C^r$-local diffeomorphisms such that, if
$(\Phi,\nbhd{U}),(\Psi,\nbhd{V})\in\sG$ then
$(\Psi,\nbhd{V})\scirc(\Psi,\nbhd{U})\in\sG$ and
$(\Phi,\nbhd{U})^{-1}\in\sG$\@.\oprocend
\end{compactenum}
\end{definition}

The notion of a local diffeomorphism with an empty domain is possibly
confusing, but is a technical convenience.  We shall often make a slight
abuse of notation by denoting a local diffeomorphism by $\Phi$ rather than
$(\Phi,\nbhd{U})$\@.  In cases when we do this, we believe there will be no
loss in clarity.

Next we consider local diffeomorphisms generated by flows of vector fields.
We recall that if $X\in\sections[1]{\tb{\man{M}}}$ then the flow
$\flow{X}{t}(x)$ is defined for $(t,x)$ in an open subset of
$\real\times\man{M}$ that we denote by $D(X)$\@.  For a vector field $X$\@,
we denote by $I(X,x_0)\subset\real$ the domain of the maximal integral curve
of $X$ through $x_0$\@.  For $t\in\real$ we denote by $\nbhd{U}(X,t)$ the
largest (possibly empty) open subset of $\man{M}$ such that
$(\flow{X}{t},\nbhd{U}(X,t))$ is a local diffeomorphism.  We note that, given
$(t_0,x_0)\in\real\times\man{M}$ there exists a neighbourhood $\nbhd{U}$ of
$x_0$ such that $\flow{X}{t}(x)$ is defined for
$t\in\interval({-t_0},{t_0})$\@.

\subsection{Lie brackets, (local) diffeomorphisms, and sheaves}

In this section we indicate how two standard constructions normally defined
for vector fields can be applied to sheaves.  This has the advantage of
systematically handling situations where one wishes to deal with objects that
are not globally defined.

First let us consider the Lie bracket.  For $r\in\{\infty,\omega\}$ and for a
$\C^r$-manifold $\man{M}$\@, we note that $\sections[r]{\tb{\man{M}}}$ has
the structure of a Lie algebra.  This Lie algebra structure can be defined
most conveniently by recalling that the set of vector fields is in one-to-one
correspondence with the set of derivations of $\func[r]{\man{M}}$
(see~\cite[Theorem~4.2.16]{RA/JEM/TSR:88} for the classical $r=\infty$ case
and~\cite{JG:81} for the $r=\omega$ case).  The derivation associated to a
vector field $X$, denoted by $\lieder{X}{}$\@, is defined simply by
$\lieder{X}{f}(x)=\natpair{\der{f}(x)}{X(x)}$\@.  For
$X,Y\in\sections[r]{\tb{\man{M}}}$ we note that
\begin{equation*}
f\mapsto \lieder{X}{\lieder{Y}{f}}-\lieder{Y}{\lieder{X}{f}}
\end{equation*}
is a derivation (one can just check this).  Thus, associated to this
derivation is a unique vector field that we denote by $[X,Y]$ which is the
\defn{Lie bracket} of $X$ and $Y$\@.  This Lie bracket is easily verified to
have the properties that render $\sections[r]{\tb{\man{M}}}$ a $\real$-Lie
algebra.

Let us now extend the above construction to sheaves.  Let
$X\in\sections[r]{\tb{\man{M}}}$ be a smooth or real analytic vector field.
Let $\nbhd{U}\subset\man{M}$ be open.  We denote by $\lieder{X,\nbhd{U}}{}$
the derivation on $\func[r]{\nbhd{U}}$ defined by $X|\nbhd{U}$\@.  Note that
the diagram
\begin{equation*}
\xymatrix{{\func[r]{\nbhd{U}}}\ar[r]^{\lieder{X,\nbhd{U}}{}}
\ar[d]_{r_{\nbhd{U},\nbhd{V}}}&{\func[r]{\nbhd{U}}}\ar[d]^{r_{\nbhd{U},\nbhd{V}}}\\
{\func[r]{\nbhd{V}}}\ar[r]_{\lieder{X,\nbhd{V}}{}}&{\func[r]{\nbhd{V}}}}
\end{equation*}
commutes for open sets $\nbhd{U},\nbhd{V}\subset\man{M}$ with
$\nbhd{V}\subset\nbhd{U}$\@.  Thus we have a well-defined map
$\map{\lieder{X,x}{}}{\gfunc[r]{x}{\man{M}}}{\gfunc[r]{x}{\man{M}}}$ which is
a derivation of the $\real$-algebras.  Thus we can think of $\lieder{X}{}$ as
being a morphism of sheaves that is a derivation.  The set of such
derivations then forms a $\real$-Lie algebra with bracket
\begin{equation*}
[\lieder{X,x}{},\lieder{Y,x}{}]([f]_x)=\lieder{X,x}(\lieder{Y,x}([f]_x))-
\lieder{Y,x}(\lieder{X,x}([f]_x)).
\end{equation*}
Note that $\lieder{X,x}{}$ depends only on the germ of $X$ at $x$\@.

Now we consider diffeomorphisms acting on functions and vector fields.  Thus
we let $r\in\{\infty,\omega\}$\@, let $\Phi\in\mappings[r]{\man{M}}{\man{M}}$
be a diffeomorphism, let $f\in\func[r]{\man{M}}$ be a function, and let
$X\in\sections[r]{\tb{\man{M}}}$\@.  The \defn{pull-back}
(\resp~\defn{push-forward}) of $f$ by $\Phi$ is defined by
$\Phi^*f=f\scirc\Phi$ (\resp~$\Phi_*f=f\scirc\Phi^{-1}$) and the
\defn{pull-back} (\resp~\defn{push-forward} of $X$ by $\Phi$ is defined by
$\Phi_*X=\tf{\Phi}\scirc X\scirc\Phi^{-1}$
(\resp~$\Phi_*X=\tf{\Phi^{-1}}\scirc X\scirc\Phi$).

Let us formulate this in sheaf language, allowing for local diffeomorphisms.
Thus we let $r\in\{\infty,\omega\}$\@, let $\man{M}$ be a $\C^r$-manifold,
and let $(\Phi,\nbhd{U})$ be a local diffeomorphism with
$\nbhd{U}\not=\emptyset$\@.  Let us first consider the action of
$(\Phi,\nbhd{U})$ on functions.  We recall from Definition~\ref{def:di-ii}
that $\Phi$ induces the direct image sheaf $\Phi_*\sfunc[r]{\nbhd{U}}$ and
the inverse image sheaf $\Phi^{-1}\sfunc[r]{\nbhd{U}}$ defined by
\begin{equation*}
\Phi_*\sfunc[r]{\nbhd{U}}(\nbhd{W})=\func[r]{\Phi^{-1}(\nbhd{W})},\qquad
\nbhd{W}\in\Phi(\nbhd{U}),
\end{equation*}
and
\begin{equation*}
\Phi^{-1}\sfunc[r]{\nbhd{U}}(\nbhd{V})=\func[r]{\Phi(\nbhd{V})},\qquad
\nbhd{V}\subset\nbhd{U}.
\end{equation*}

The pull-back morphism is the sheaf morphism $\Phi^*$ from
$\Phi^{-1}\sfunc[r]{\Phi(\nbhd{U})}$ to $\sfunc[r]{\nbhd{U}}$ defined by
asking that
$\map{\Phi^*_{\nbhd{V}}}{\func[r]{\Phi(\nbhd{V})}}{\func[r]{\nbhd{V}}}$ be
given by $\Phi^*_{\nbhd{V}}(g)=g\scirc(\Phi|\nbhd{V})$\@.  Note that the
diagram
\begin{equation*}
\xymatrix{{\func[r]{\Phi(\nbhd{V})}}\ar[r]^(0.6){\Phi^*_{\nbhd{V}}}
\ar[d]_{r_{\Phi(\nbhd{V}),\Phi(\nbhd{W})}}&
{\func[r]{\nbhd{V}}}\ar[d]^{r_{\nbhd{V},\nbhd{W}}}\\
{\func[r]{\Phi(\nbhd{W})}}\ar[r]_(0.6){\Phi^*_{\nbhd{W}}}&
{\func[r]{\nbhd{W}}}}
\end{equation*}
commutes for each pair $\nbhd{V},\nbhd{W}\subset\nbhd{U}$ of open sets for
which $\nbhd{W}\subset\nbhd{V}$\@.  Thus we have an induced mapping
$\Phi^*_x$ from the stalk $\Phi^{-1}\gfunc[r]{x}{\nbhd{U}}$ to the stalk
$\gfunc[r]{x}{\nbhd{U}}$ given explicitly by
\begin{equation*}
\Phi^*_x([g]_{\Phi(x)})=[g\scirc\Phi]_x
\end{equation*}
for $x\in\nbhd{U}$\@.

The push-forward morphism is the sheaf morphism $\Phi_*$ from
$\sfunc[r]{\nbhd{U}}$ to $\Phi^{-1}\sfunc[r]{\nbhd{U}}$ defined by asking
that $\map{\Phi_{*,\nbhd{V}}}{\func[r]{\nbhd{V}}}{\func[r]{\Phi(\nbhd{V})}}$
be given by $\Phi_{*,\nbhd{V}}(f)=f\scirc(\Phi^{-1}|\Phi(\nbhd{V}))$\@.  Note
that the diagram
\begin{equation*}
\xymatrix{{\func[r]{\nbhd{V}}}\ar[r]^(0.4){\Phi_{*,{\nbhd{V}}}}
\ar[d]_{r_{\nbhd{W},\nbhd{V}}}&
{\func[r]{\Phi(\nbhd{V})}}\ar[d]^{r_{\Phi(\nbhd{V}),(\nbhd{W})}}\\
{\func[r]{\nbhd{W}}}\ar[r]_(0.4){\Phi_{*,\nbhd{W}}}&
{\func[r]{\Phi(\nbhd{W})}}}
\end{equation*}
commutes for every pair $\nbhd{W},\nbhd{V}\subset\nbhd{U}$ of open sets for
which $\nbhd{V}\subset\nbhd{W}$\@.  Thus we have an induced mapping
$\Phi_{*,x}$ from the stalk $\gfunc[r]{x}{\nbhd{U}}$ to the stalk
$\Phi^{-1}\gfunc[r]{x}{\nbhd{U}}$ given explicitly by
\begin{equation*}
\Phi_{*,x}([f]_x)=[f\scirc\Phi^{-1}]_{\Phi(x)}
\end{equation*}
for $x\in\nbhd{U}$\@.

Now let us see how this can be used to give the sheaf version of the
pull-back of a vector field by a local diffeomorphism.  Let $\man{M}$ and
$(\Phi,\nbhd{U})$ be as above, let $\nbhd{V}\subset\nbhd{U}$ be open, and let
$f\in\func[r]{\nbhd{V}}$ be a local section of $\sfunc[r]{\nbhd{U}}$ over
$\nbhd{V}$\@.  Let $Y$ be a local section of $\ssections[r]{\tb{\man{M}}}$
over $\Phi(\nbhd{V})$\@,~\ie~an element of the inverse image sheaf
$\Phi^{-1}\ssections{\tb{\man{M}}|\nbhd{U}}(\nbhd{V})$\@.  Then, for
$y\in\Phi(\nbhd{V})$\@, the tangent vector $\tf[y]{\Phi^{-1}}(Y(y))$ can be
used to derive the function $f$ at $\Phi^{-1}(y)$\@.  Thus we can define a
morphism $\map{\slieder{\Phi^*}{Y,\nbhd{V}}{}}{\func[r]{\nbhd{V}}}
{\func[r]{\Phi^{-1}(\nbhd{V})}}$ by
\begin{equation*}
\slieder{\Phi^*}{Y,\nbhd{V}}{f}(y)=
\natpair{\der{f}(\Phi^{-1}(y))}{\tf[y]{\Phi^{-1}}(Y(y))},\qquad y\in\Phi(\nbhd{V}).
\end{equation*}
Alternatively, we can regard $\slieder{\Phi^*}{Y,\nbhd{V}}{}$ as a morphism
from $\sfunc[r]{\nbhd{U}}(\nbhd{V})=\func[r]{\nbhd{V}}$ to itself by defining
\begin{equation*}
\slieder{\Phi^*}{Y,\nbhd{V}}{f}(x)=
\natpair{\der{f}(x)}{\tf[\Phi(x)]{\Phi^{-1}}(Y(\Phi(x)))},\qquad x\in\nbhd{V}.
\end{equation*}
This map is a derivation of $\real$-algebras,~\ie~it is $\real$-linear and
satisfies
\begin{equation*}
\slieder{\Phi^*}{Y,\nbhd{V}}{fg}=f(\slieder{\Phi^*}{Y,\nbhd{V}}{g})+
(\slieder{\Phi^*}{Y,\nbhd{V}}{f})g
\end{equation*}
for every $f,g\in\func[r]{\nbhd{V}}$\@.  Thus there exists a unique vector
field on $\nbhd{V}$\@, which we denote by $\Phi^*_{\nbhd{V}}Y$\@, such that
\begin{equation*}
\lieder{\Phi^*_{\nbhd{V}}Y}{f}=\slieder{\Phi^*}{Y,\nbhd{V}}{f}
\end{equation*}
for every $f\in\func[r]{\nbhd{V}}$\@.  Thus
$\Phi^*_{\nbhd{V}}Y\in\ssections[r]{\tb{\man{M}}|\nbhd{U}}(\nbhd{V})$\@.
Clearly the diagram
\begin{equation*}
\xymatrix{{\func[r]{\nbhd{V}}}\ar[r]^{\slieder{\Phi^*}{Y,\nbhd{V}}{}}
\ar[d]_{r_{\nbhd{V},\nbhd{W}}}&{\func[r]{\nbhd{V}}}
\ar[d]^{r_{\nbhd{V},\nbhd{W}}}\\
{\func[r]{\nbhd{W}}}\ar[r]_{\slieder{\Phi^*}{Y,\nbhd{W}}{}}&
{\func[r]{\nbhd{W}}}}
\end{equation*}
commutes for open sets $\nbhd{V},\nbhd{W}\subset\nbhd{U}$ satisfying
$\nbhd{W}\subset\nbhd{V}$\@.  Thus, for each $x\in\nbhd{U}$ we have a
well-defined morphism of stalks:
\begin{equation*}
\map{\slieder{\Phi^*}{Y,x}{}}{\gfunc[r]{x}{\nbhd{U}}}{\gfunc[r]{x}{\nbhd{U}}}.
\end{equation*}
Moreover, this map is a derivation.  That is, it is $\real$-linear and
satisfies
\begin{equation*}
\slieder{\Phi^*}{Y,x}{[f]_x[g]_x}=[f]_x(\slieder{\Phi^*}{Y,x}{[g]_x})+
(\slieder{\Phi^*}{Y,x}{[f]_x})[g]_x
\end{equation*}
for $[f]_x,[g]_x\in\gfunc[r]{x}{\nbhd{U}}$\@.  Thus there exists a unique
vector field germ, which we denote by $\Phi^*_x[Y]_{\Phi(x)}$\@, such that
\begin{equation*}
\lieder{\Phi^*_x[Y]_{\Phi(x)}}{[f]_x}=\slieder{\Phi^*}{Y,x}{[f]_x}
\end{equation*}
for every $[f]_x\in\gfunc[r]{x}{\nbhd{U}}$\@.  The upshot of the preceding
development is that, associated with the local diffeomorphism
$(\Phi,\nbhd{U})$\@, is a sheaf morphism $\Phi^*$ from
$\Phi^{-1}\ssections[r]{\tb{\man{M}}|\nbhd{U}}$ to
$\ssections[r]{\nbhd{U}}$\@.

Finally, we define the sheaf version of the push-forward of a vector field by
a local diffeomorphism.  Let $\man{M}$ and $(\Phi,\nbhd{U})$ be as above, let
$\nbhd{V}\subset\nbhd{U}$ be open, and let $g\in\func[r]{\Phi(\nbhd{V})}$ be
a local section of $\Phi^{-1}\sfunc[r]{\nbhd{U}}$ over $\nbhd{V}$\@.  Let
$X\in\sections[r]{\nbhd{V}}$ be a local section of
$\ssections[r]{\tb{\man{M}}}$ over $\nbhd{V}$\@,~\ie~an element of
$\ssections[r]{\tb{\man{M}}|\nbhd{U}}(\nbhd{V})$\@.  Then, for
$x\in\nbhd{V}$\@, the tangent vector
$\tf[x]{\Phi}(X(x))\in\tb[\Phi(x)]{\man{M}}$ can be used to derive the
function $g$ at $\Phi(x)$\@.  Thus, we can define a morphism
$\map{\slieder{\Phi_*}{X,\nbhd{V}}{}}{\func[r]{\Phi(\nbhd{V})}}
{\func[r]{\nbhd{V}}}$ by
\begin{equation*}
\slieder{\Phi_*}{X,\nbhd{V}}{g}(x)=
\natpair{\der{g}(\Phi(x))}{\tf[x]{\Phi}(X(x))},\qquad x\in\nbhd{V}.
\end{equation*}
Alternatively, we can regard $\slieder{\Phi_*}{X,\nbhd{V}}{}$ as a morphism
from $\Phi^{-1}\sfunc[r]{\nbhd{U}}(\nbhd{V})=\func[r]{\Phi(\nbhd{V})}$ to
itself by defining
\begin{equation*}
\slieder{\Phi_*}{X,\nbhd{V}}{g}(y)=
\natpair{\der{g}(y)}{\tf[\Phi^{-1}(y)]{\Phi}(X(\Phi^{-1}(y)))},\qquad y\in\Phi(\nbhd{V}).
\end{equation*}
This map is a derivation of $\real$-algebras,~\ie~it is $\real$-linear and
satisfies
\begin{equation*}
\slieder{\Phi_*}{X,\nbhd{V}}{gh}=g(\slieder{\Phi_*}{X,\nbhd{V}}{h})+
(\slieder{\Phi_*}{X,\nbhd{V}}{g})h
\end{equation*}
for every $g,h\in\func[r]{\Phi(\nbhd{V})}$\@.  Thus there exists a unique
vector field on $\Phi(\nbhd{V})$\@, which we denote by
$\Phi_{*,{\nbhd{V}}}X$\@, such that
\begin{equation*}
\lieder{\Phi_{*,{\nbhd{V}}}X}{g}=\slieder{\Phi_*}{X,\nbhd{V}}{g}
\end{equation*}
for every $g\in\func[r]{\Phi(\nbhd{V})}$\@.  Thus
$\Phi_{*,\nbhd{V}}X\in\Phi^{-1}\ssections[r]{\tb{\man{M}}|\nbhd{U}}(\nbhd{V})$\@.
Clearly the diagram
\begin{equation*}
\xymatrix{{\func[r]{\Phi(\nbhd{V})}}\ar[r]^{\slieder{\Phi_*}{X,\nbhd{V}}{}}
\ar[d]_{r_{\Phi(\nbhd{V}),\Phi(\nbhd{W})}}&{\func[r]{\Phi(\nbhd{V})}}
\ar[d]^{r_{\Phi(\nbhd{V}),\Phi(\nbhd{W})}}\\
{\func[r]{\Phi(\nbhd{W})}}\ar[r]_{\slieder{\Phi_*}{X,\nbhd{W}}{}}&
{\func[r]{\Phi(\nbhd{W})}}}
\end{equation*}
commutes for open sets $\nbhd{V},\nbhd{W}\subset\nbhd{U}$ satisfying
$\nbhd{W}\subset\nbhd{V}$\@.  Thus, for each $y\in\Phi(\nbhd{U})$ we have a
well-defined morphism of stalks:
\begin{equation*}
\map{\slieder{\Phi_*}{X,y}{}}{\Phi^{-1}\gfunc[r]{y}{\nbhd{U}}}
{\Phi^{-1}\gfunc[r]{y}{\nbhd{U}}}.
\end{equation*}
Moreover, this map is a derivation.  That is, it is $\real$-linear and
satisfies
\begin{equation*}
\slieder{\Phi_*}{X,y}{[g]_y[h]_y}=[g]_y(\slieder{\Phi_*}{X,y}{[h]_y})+
(\slieder{\Phi_*}{X,y}{[g]_y})[h]_y
\end{equation*}
for $[g]_y,[h]_y\in\Phi^{-1}\gfunc[r]{y}{\nbhd{U}}$\@.  Thus there exists a
unique vector field germ, which we denote by
$\Phi_{*,x}[X]_{\Phi^{-1}(y)}$\@, such that
\begin{equation*}
\lieder{\Phi_{*,x}[X]_{\Phi^{-1}(y)}}{[g]_y}=\slieder{\Phi_*}{X,y}{[g]_y}
\end{equation*}
for every $[g]_y\in\Phi^{-1}\gfunc[r]{y}{\nbhd{U}}$\@.  The upshot of the
preceding development is that, associated with the local diffeomorphism
$(\Phi,\nbhd{U})$\@, is a sheaf morphism $\Phi_*$ from
$\ssections[r]{\tb{\man{M}}|\nbhd{U}}$ to
$\Phi^{-1}\ssections[r]{\tb{\man{M}}|\nbhd{U}}$\@.

The preceding discussion is lengthy and notation-laden.  Therefore, it is
worth summarising the constructions we have made.  We let $(\Phi,\nbhd{U})$
be a local diffeomorphism, and note that our preceding developments define
the following sheaf morphisms:
\begin{minipage}[t]{0.4\linewidth}
\begin{enumerate}
\item $\map{\Phi^*}{\Phi^{-1}\sfunc[r]{\nbhd{U}}}{\sfunc[r]{\nbhd{U}}}$\@;
\item $\map{\Phi_*}{\sfunc[r]{\nbhd{U}}}{\Phi^{-1}\sfunc[r]{\nbhd{U}}}$\@;\savenum
\end{enumerate}
\end{minipage}
\begin{minipage}[t]{0.4\linewidth}
\begin{enumerate}\resumenum
\item $\map{\Phi^*}{\Phi^{-1}\ssections[r]{\tb{\man{M}}|\nbhd{U}}}
{\ssections[r]{\tb{\man{M}}|\nbhd{U}}}$\@;
\item $\map{\Phi_*}{\ssections[r]{\tb{\man{M}}|\nbhd{U}}}
{\Phi^{-1}\ssections[r]{\tb{\man{M}}|\nbhd{U}}}$\@.
\end{enumerate}
\end{minipage}

\subsection{Lie subalgebras of vector fields and
subsheaves}\label{subsec:lie-algebras}

In this section we consider Lie subalgebras of the Lie algebras
$\sections[r]{\tb{\man{M}}}$ and $\ssections[r]{\tb{\man{M}}}$ for
$r\in\{\infty,\omega\}$\@.  Let us first define the objects of interest.
\begin{definition}
Let $r\in\{\infty,\omega\}$\@, let $\man{M}$ be a manifold of class
$\C^r$\@, and let $x\in\man{M}$\@.
\begin{compactenum}[(i)]
\item A \defn{Lie subalgebra} of vector fields is a Lie subalgebra of the
$\real$-Lie algebra $\sections[r]{\tb{\man{M}}}$\@.
\item A \defn{Lie subalgebra} of $\ssections[r]{\tb{\man{M}}}$ is an
assignment to each open set\/ $\nbhd{U}\subset\man{M}$ a Lie subalgebra
$L(\nbhd{U})$ of vector fields on $\nbhd{U}$ with the property that
$L(\nbhd{V})=r_{\nbhd{U},\nbhd{V}}(L(\nbhd{U}))$ for every pair
$\nbhd{U},\nbhd{V}$ of open sets for which $\nbhd{V}\subset\nbhd{U}$\@.
\item A \defn{Lie subalgebra} of germs of vector fields at $x$ is a Lie
subalgebra of the $\real$-Lie algebra $\gsections[r]{x}{\tb{\man{M}}}$\@.
\item If $\sX\subset\sections[r]{\tb{\man{M}}}$ then the \defn{Lie subalgebra
generated by $\sX$} is the smallest Lie subalgebra of vector fields
containing $\sX$\@.  The Lie subalgebra is denoted by $\sL^{(\infty)}(\sX)$\@.
\item If $\sX$ is a subsheaf of sets of the sheaf
$\ssections[r]{\tb{\man{M}}}$\@\textemdash{}\ie~an assignment to each open
set $\nbhd{U}\subset\man{M}$ a subset
$X(\nbhd{U})\subset\sections[r]{\tb{\man{M}}|\nbhd{U}}$ with the assignment
satisfying $X(\nbhd{V})=r_{\nbhd{U},\nbhd{V}}(X(\nbhd{U}))$ for every pair of
open sets $\nbhd{U},\nbhd{V}$ for which
$\nbhd{V}\subset\nbhd{V}$\@\textemdash{}the \defn{Lie subalgebra of generated
by $\sX$} is the Lie subalgebra of $\ssections[r]{\tb{\man{M}}}$ defined by
assigning to the open set $\nbhd{U}$ the Lie subalgebra of vector fields on
$\nbhd{U}$ generated by $X(\nbhd{U})$\@.  This Lie subalgebra is denoted by
$\sL^{(\infty)}(\sX)=\ifam{\sL^{(\infty)}(X(\nbhd{U}))}_{\nbhd{U}\,\textrm{open}}$\@.
\item If $\sX_x\subset\gsections[r]{x}{\tb{\man{M}}}$\@, the \defn{Lie
subalgebra generated by $\sX_x$} is the smallest Lie subalgebra of germs of
vector fields at $x$ containing $\sX_x$\@.  This Lie subalgebra is denoted by
$\sL^{(\infty)}(\sX_x)$\@.\oprocend
\end{compactenum}
\end{definition}

One can easily give a precise characterisation of a Lie subalgebra generated
by a set in the three cases of the preceding definition.
\begin{proposition}\label{prop:liealggen}
Let\/ $r\in\{\infty,\omega\}$ and let\/ $\man{M}$ be a\/ $\C^r$-manifold.
Then the following statements hold:
\begin{compactenum}[(i)]
\item if\/ $\sX\subset\sections[r]{\tb{\man{M}}}$\@, the Lie subalgebra
generated by\/ $\sX$ is generated by finite\/ $\real$-linear combinations of
vector fields of the form
\begin{equation*}
[X_k,[X_{k-1},\dots,[X_2,X_1]\cdots]],\qquad k\in\integerp,\
X_1,\dots,X_k\in\sX;
\end{equation*}
\item if\/ $\sX=\ifam{X(\nbhd{U})}_{\nbhd{U},\textrm{open}}$ is a subsheaf of
sets of\/ $\ssections[r]{\tb{\man{M}}}$\@, then the Lie subalgebra of\/
$\ssections[r]{\tb{\man{M}}}$ generated by\/ $\sX$ is such that\/
$\sL^{(\infty)}(X(\nbhd{U}))$ is generated by finite\/ $\real$-linear
combinations of vector fields of the form
\begin{equation*}
[X_k,[X_{k-1},\dots,[X_2,X_1]\cdots]],\qquad k\in\integerp,\
X_1,\dots,X_k\in X(\nbhd{U});
\end{equation*}
\item if\/ $x\in\man{M}$ and if\/
$\sX_x\subset\gsections[r]{x}{\tb{\man{M}}}$\@, the Lie subalgebra generated
by\/ $\sX_x$ is generated by finite\/ $\real$-linear combinations of germs of
vector fields of the form
\begin{equation*}
[[X_k]_x,[[X_{k-1}]_x,\dots,[[X_2]_x,[X_1]_x]\cdots]],\qquad k\in\integerp,\
[X_1]_x,\dots,[X_k]_x\in\sX.
\end{equation*}
\end{compactenum}
\begin{proof}
We prove the first statement only.  The second follows from this and the
third follows via an entirely similar computation.  For vector fields
$X_1,\dots,X_k\in\sX$\@, since $\sL^{(\infty)}(\sX)$ is a Lie subalgebra of
$\sections[r]{\tb{\man{M}}}$\@, it follows by induction that
$[X_k,[X_{k-1},\dots,[X_2,X_1]\cdots]]\in\sL^{(\infty)}(\sX)$\@.  Since
$\sL^{(\infty)}(\sX)$ is a subspace of the $\real$-vector space
$\sections[r]{\tb{\man{M}}}$\@, it also follows that all $\real$-linear
combinations of such vector fields are in $\sL^{(\infty)}(\sX)$\@.

To prove the opposite inclusion, it suffices to show\textemdash{}since
$\sL^{(\infty)}(\sX)$ is the smallest Lie subalgebra containing
$\sX$\textemdash{}that the set of all $\real$-linear combinations in the
statement of the proposition forms a Lie algebra.  If we have two
$\real$-linear combinations of vector fields of the form stated in the
proposition, their Lie bracket will be in $\sL^{(\infty)}(\sX)$ if and only
if the Lie bracket of each of the summands is in $\sL^{(\infty)}(\sX)$ (by
linearity of the Lie bracket).  Consider two vector fields of the form stated
in the proposition:
\begin{gather*}
X=[X_k,[X_{k-1},\dots,[X_2,X_1]\cdots]]\\
Y=[Y_l,[Y_{l-1},\dots,[Y_2,Y_1]\cdots]].
\end{gather*}
We shall prove by induction that $[X,Y]\in\sL^{(\infty)}(\sX)$ for any $k$
and $l$\@.  Note that $[X,Y]\in\sL^{(\infty)}(\sX)$ for any $Y$ and $l$\@,
and for $k=1$\@.  Now suppose this is true for $k=1,\dots,m$\@.  Then, taking
$k=m+1$\@, we have
\begin{equation*}
[X,Y]=[[X_{m+1},X^1],Y]
\end{equation*}
where $X^1=[X_m,\dots,[X_2,X_1]\cdots]$\@.  By the Jacobi identity we have
\begin{equation*}
[[X_{m+1},X^1],Y]+[[Y,X_{m+1}],X^1]+[[X^1,Y],X_{m+1}]=
0_{\sections[r]{\tb{\man{M}}}}.
\end{equation*}
This gives
\begin{equation*}
[X,Y]=[X^1,[Y,X_{m+1}]]+[X_{m+1},[X^1,Y]].
\end{equation*}
By the induction hypothesis, $[X^1,[X_{m+1},Y]]\in\sL^{(\infty)}(\sX)$ since
$X^1$ is a bracket of length $m$\@.  Also $[X^1,Y]\in\sL^{(\infty)}(\sX)$ so
the second term on the right is in $\sL^{(\infty)}(\sX)$\@.  Thus the set of
linear combinations of the form stated in the proposition forms a Lie
subalgebra, giving the result.
\end{proof}
\end{proposition}

Let $\sX\subset\sections[r]{\tb{\man{M}}}$ be a family of vector fields.
Note that since $\sL^{(\infty)}(\sX)$ is a subspace, one immediately has
\begin{equation*}
\sL^{(\infty)}(\sX)=\sL^{(\infty)}(\vecspan[\real]{\sX}).
\end{equation*}
However, it is not generally the case that $\sL^{(\infty)}(\sX)$ is a
submodule.  Similarly flavoured statements hold for $\sX$ a subsheaf of sets
of $\ssections[r]{\tb{\man{M}}}$ or for $\sX_x$ a subset of germs from
$\gsections[r]{x}{\tb{\man{M}}}$\@.  In this respect, however, the following
result is useful.
\begin{proposition}\label{prop:liealgmodule}
If\/ $r\in\{\infty,\omega\}$\@, if\/ $\man{M}$ is a\/ $\C^r$-manifold, and
if\/ $x\in\man{M}$\@, then the following statements hold:
\begin{compactenum}[(i)]
\item if\/ $\sM\subset\sections[r]{\tb{\man{M}}}$ is a submodule of vector
fields, then\/ $\sL^{(\infty)}(\sM)$ is also a submodule of vector fields;
\item if\/ $\sF$ is a subsheaf\/ $\sfunc[r]{\man{M}}$-modules of\/
$\ssections[r]{\tb{\man{M}}}$\@, then\/ $\sL^{(\infty)}(\sF)$ is also a
subsheaf\/ $\sfunc[r]{\man{M}}$-modules of\/ $\ssections[r]{\tb{\man{M}}}$\@.
\item if\/ $\sF_x\subset\gsections[r]{x}{\tb{\man{M}}}$ is a submodule of
germs of vector fields, then\/ $\sL^{(\infty)}(\sF_x)$ is also a submodule of
germs of vector fields.
\end{compactenum}
\begin{proof}
We prove the first statement only; the second follows directly from this and
the third is proved in exactly the same way.  By
Proposition~\ref{prop:liealggen} it suffices to show that, for any
$f\in\func[r]{\man{M}}$ and for any $X_1,\dots,X_k\in\sM$\@,
\begin{equation*}
f[X_k,[X_{k-1},\dots,[X_2,X_1]\cdots]]\in\sL^{(\infty)}(\sM).
\end{equation*}
We prove this by induction on $k$\@, it clearly being true for $k=1$\@.
Assume now that the statement holds for $k\in\{1,\dots,m+1\}$\@.  Then
\begin{equation*}
f[X_{m+1},[X_m,\dots,[X_2,X_1]\cdots]]=
[fX_{m+1},[X_m,\dots,[X_2,X_1]\cdots]]+
(\lieder{[X_m,\dots,[X_2,X_1]\cdots]}{f})X_k.
\end{equation*}
By the induction hypothesis,
$(\lieder{[X_m,\dots,[X_2,X_1]\cdots]}{f})X_k\in\sL^{(\infty)}(\sM)$\@.
Since $fX_{m+1}\in\sM$\@, by Proposition~\ref{prop:liealggen} it follows that
$[fX_{m+1},[X_m,\dots,[X_2,X_1]\cdots]]\in\sL^{(\infty)}(\sM)$\@, giving the
result.
\end{proof}
\end{proposition}

Let $\sX$ be a family of smooth or real analytic vector fields on a smooth or
real analytic manifold $\man{M}$\@.  Following the notation of
Definition~\ref{def:generated}\@, associated with the family of vector fields
$\sL^{(\infty)}(\sX)$ is the distribution $\dist{D}(\sL^{(\infty)}(\sX))$
which we abbreviate by $\dist{L}^{(\infty)}(\sX)$\@.  The following result
simplifies some parts of the subsequent discussion.  For the statement, we
refer to Definition~\ref{def:generated} for the notation for the submodule of
vector fields generated by a family of vector fields.
\begin{proposition}\label{prop:Linfty-equiv1}
Let\/ $r\in\{\infty,\omega\}$\@, let\/ $\man{M}$ be a\/ $\C^r$-manifold, and
let\/ $\sX\subset\sections[r]{\tb{\man{M}}}$\@.  Then the distributions
\begin{compactenum}[(i)]
\item $\dist{L}^{(\infty)}(\sX)$\@,
\item $\dist{L}^{(\infty)}(\vecspan[\real]{\sX})$\@,
\item $\dist{L}^{(\infty)}(\modgen{\sX})$\@, and
\item $\dist{D}(\modgen{\sL^{(\infty)}(\sX)})$
\end{compactenum}
agree.
\begin{proof}
It is clear that
\begin{equation*}
\dist{L}^{(\infty)}(\sX)\subset\dist{L}^{(\infty)}(\vecspan[\real]{\sX})
\subset\dist{L}^{(\infty)}(\modgen{\sX}).
\end{equation*}
We will show that
$\dist{L}^{(\infty)}(\modgen{\sX})\subset\dist{L}^{(\infty)}(\sX)$\@.  By
Proposition~\ref{prop:liealggen} it suffices to show that
\begin{equation*}
[Y_k,[Y_{k-1},\dots,[Y_2,Y_1]\cdots]](x)\in\dist{L}^{(\infty)}(\sX)_x
\end{equation*}
for every $Y_1,\dots,Y_k\in\modgen{\sX}$ and for every $x\in\man{M}$\@.

We prove this by first showing that
\begin{equation*}
[Y_k,[Y_{k-1},\dots,[Y_2,Y_1]\cdots]]=f^1Z_1+\dots+f^sZ_s
\end{equation*}
for every $Y_1,\dots,Y_k\in\modgen{\sX}$\@, and where
$f^1,\dots,f^s\in\func[r]{\man{M}}$ and
\begin{equation*}
Z_j=[X_{j,l},[X_{j,l-1},\dots,[X_{j,2},X_{j,1}]\cdots]]
\end{equation*}
for $X_{j,1},\dots,X_{j,l}\in\sX$ with $l\in\{1,\dots,k\}$\@.  This we prove
by induction on $k$\@.  It is clearly true for $k=1$\@, so suppose it holds
for $k\in\{1,\dots,m\}$ and let $Y_1,\dots,Y_m,Y_{m+1}\in\modgen{\sX}$\@.
Write
\begin{equation*}
Y_{m+1}=f^1X_1+\dots+f^sX_s,\qquad f^1,\dots,f^s\in\func[r]{\man{M}},\
X_1,\dots,X_s\in\sX.
\end{equation*}
Then
\begin{equation*}
[Y_{m+1},[Y_m,\dots,[Y_2,Y_1]\cdots]]=
\sum_{j=1}^s(f^j[X_j,[Y_m,\dots,[Y_2,Y_1]]]-
(\lieder{[Y_m,\dots,[Y_2,Y_1]\cdots]}{f^j})X_j).
\end{equation*}
By the induction hypothesis,
\begin{equation*}
[Y_m,\dots,[Y_2,Y_1]\cdots]=g^1Z_1+\dots+g^dZ_d
\end{equation*}
for $g^1,\dots,g^d\in\func[r]{\man{M}}$ and where
\begin{equation*}
Z_a=[X_{a,l_a},[X_{a,l_a-1},\dots,[X_{a,2},X_{a,1}]\cdots]]
\end{equation*}
for $X_{a,1},\dots,X_{a,l_a}\in\sX$ and where $l_a\in\{1,\dots,m\}$ for each
$a\in\{1,\dots,d\}$\@.  Then
\begin{align*}
[X_j,[Y_m,\dots,[Y_2,Y_1]\cdots]]=&\;
\sum_{a=1}^d[X_j,g^a[X_{a,l_a},[X_{a,l_a-1},\dots,
[X_{a,2},X_{a,1}]\cdots]]]\\
=&\;\sum_{a=1}^d(g^a[X_j,[X_{a,l_a},[X_{a,l_a-1},\dots,
[X_{a,2},X_{a,1}]\cdots]]]\\
&\;+(\lieder{X_j}{g^a})[X_{a,l_a},[X_{a,l_a-1},\dots,
[X_{a,2},X_{a,1}]\cdots]]).
\end{align*}
This proves that $[Y_{m+1},[Y_m,\dots,[Y_2,Y_1]\cdots]]$ has the desired
form.

From this it immediately follows from Proposition~\ref{prop:liealggen} that
\begin{equation*}
[Y_k,[Y_{k-1},\dots,[Y_2,Y_1]\cdots]](x)\in\dist{L}^{(\infty)}(\sX)_x
\end{equation*}
for every $Y_1,\dots,Y_k\in\modgen{\sX}$ and for every $x\in\man{M}$\@, and
so the first three distributions in the statement of the proposition are
equal.  The equality of these distributions with the fourth distribution in
the statement of the proposition follows from
Proposition~\ref{prop:subbundlespan1}\@.
\end{proof}
\end{proposition}

There is a corresponding sheaf version of the preceding result which we state
for completeness; it follows immediately from the result above.  Following
the notation of Definition~\ref{def:setsheaves}\@, we denote
$\dist{L}^{(\infty)}(\sX)=\dist{D}(\sL^{(\infty)}(\sX))$ for a subsheaf of
sets $\sX$ of $\ssections[r]{\tb{\man{M}}}$\@.
\begin{proposition}\label{prop:Linfty-equiv2}
Let\/ $r\in\{\infty,\omega\}$\@, let\/ $\man{M}$ be a\/ $\C^r$-manifold, and
let\/ $\sX=\ifam{X(\nbhd{U})}_{\nbhd{U}\,\textrm{open}}$ be a subsheaf of
sets of $\ssections[r]{\tb{\man{M}}}$\@.  Then the distributions
\begin{compactenum}[(i)]
\item $\dist{L}^{(\infty)}(\sX)$\@,
\item $\dist{L}^{(\infty)}(\modgen{\sX})$\@, and
\item $\dist{D}(\modgen{\sL^{(\infty)}(\sX)})$
\end{compactenum}
agree.
\end{proposition}

The developments with Lie algebras thus far in this section have had to do
with families of vector fields and subsheaves of sets of
$\ssections[r]{\tb{\man{M}}}$\@.  We now turn to Lie algebraic constructions
in the case when $\sX=\sections[r]{\dist{D}}$ is the submodule of sections or
where $\sX=\ssections[r]{\dist{D}}$ is the sheaf of submodules of sections of
a distribution $\dist{D}$\@.

Let us first consider the submodule case.  Thus we let
$\sM\subset\sections[r]{\tb{\man{M}}}$ be a submodule of vector fields, and
we recall that $\sM\subset\sections[r]{\dist{D}(\sM)}$\@, but that the
inclusion is, in general, strict,~\cf~Example~\ref{eg:bad-generator}\@.
Moreover, unlike some of the other anomalies we have encountered and will
encounter in the paper, this one is not a result of a lack of analyticity.
However, there are interesting conclusions that hold in the analytic case,
and indeed more generally in the locally finitely generated case.
\begin{theorem}\label{the:LinftyX-D(X)}
Let\/ $r\in\{\infty,\omega\}$\@, let\/ $\man{M}$ be a\/ $\C^r$-manifold, and
let\/ $\sX\subset\sections[r]{\tb{\man{M}}}$ be a family of vector fields
such that
\begin{compactenum}[(i)]
\item $\modgen{\sX}$ is a locally finitely generated submodule of\/
$\sections[r]{\tb{\man{M}}}$ and
\item the module (by Proposition~\ref{prop:liealgmodule})\/
$\sL^{(\infty)}(\modgen{\sX})$ is locally finitely generated.
\end{compactenum}
Then\/
$\dist{L}^{(\infty)}(\sX)=\dist{L}^{(\infty)}(\sections[r]{\dist{D}(\sX)})$\@.
\begin{proof}
Our proof will rely on the Orbit Theorem and various constructions and
results from the theory of control systems and differential inclusions.

Let $x_0\in\man{M}$\@.  Denote $\sM=\modgen{\sX}$\@.  By hypothesis, there
exists a neighbourhood $\nbhd{U}$ of $x_0$\@, a finite subset, say
$\sX'=\ifam{X_1,\dots,X_k}$\@, of $\sX$\@, and $m\in\integerp$ such that
$F_{\sX'}(\nbhd{U})=F_{\sM}(\nbhd{U})$ (see
Definition~\ref{def:module->sheaf} for the notation) and such that the vector
fields
\begin{equation*}
[X_{a_1},[X_{a_2},\dots,[X_{a_{l-1}},X_{a_l}]\cdots]]|\nbhd{U},\qquad
l\in\{1,\dots,m\},\ a_1,\dots,a_l\in\{1,\dots,k\},
\end{equation*}
generate $\sL^{(\infty)}(F_{\sM}(\nbhd{U}))$\@, using
Proposition~\ref{prop:liealggen}\@.  Thus both $\sL^{(\infty)}(F_{\sX'})$ and
$\sL^{(\infty)}(F_{\sM}(\nbhd{U}))$ generate locally finitely generated
modules and, moreover, these modules agree.  Thus, by
Theorem~\ref{the:fingen-orbit-theorem}\@,
\begin{equation*}
\Orb(x_0,\sX'|\nbhd{U})=\Orb(x_0,\sM|\nbhd{U}),
\end{equation*}

Consider now two control systems defined on $\nbhd{U}$\@:
\begin{align*}
\Sigma\colon&\quad\xi'(t)=\sum_{j=1}^k\mu^j(t)X_j(\xi(t)),\qquad
\vect{\mu}(t)\in U\eqdef\{\vect{e}_1,\dots,\vect{e}_k,
-\vect{e}_1,\dots,-\vect{e}_k\},\\
\cohull(\Sigma)\colon&\quad\xi'(t)=\sum_{j=1}^k\mu^j(t)X_j(\xi(t)),\qquad
\vect{\mu}(t)\in\cohull(U).
\end{align*}
In each case, we consider controls to be locally integrable functions taking
values in the control set.  We have the associated differential inclusions
defined by the set-valued right-hand sides
\begin{align*}
F(x)=&\;\asetdef{\sum_{j=1}^ku^jX_j(x)}{\vect{u}\in U},\\
\cohull(F)(x)=&\;\asetdef{\sum_{j=1}^ku^jX_j(x)}{\vect{u}\in\cohull(U)},
\end{align*}
respectively.

We now state a lemma which is often used, but for which we were unable to
locate a proof.
\begin{prooflemma}
Let\/ $(\ms{M},\d_{\ms{M}})$ and\/ $(\ms{N},\d_{\ms{N}})$ be metric spaces,
let\/ $\ts{U}$ be a compact topological space, and let\/
$\map{f}{\ms{M}\times\ms{U}}{\ms{N}}$ be a continuous map for which the map\/
$x\mapsto f(x,u)$ is locally Lipschitz for each\/ $u\in\ts{U}$\@.  Then the
set-valued map
\begin{equation*}
x\mapsto F(x)\eqdef\setdef{f(x,u)}{u\in\ts{U}}
\end{equation*}
is locally Lipschitz,~\ie~for each\/ $x\in\ms{M}$ there exists\/ $L\in\realp$
and a neighbourhood\/ $\nbhd{X}$ of\/ $x$ such that
\begin{equation*}
F(x_1)\subset\bigcup_{y\in F(x_2)}\oball[\ms{N}]{y}{L\d_{\ms{M}}(x_1,x_2)}
\end{equation*}
for each\/ $x_1,x_2\in\nbhd{X}$\@, where\/ $\oball[\ms{N}]{y}{r}$ is the
ball of radius\/ $r$ centred at\/ $y\in\ms{N}$\@.
\begin{subproof}
Since $f$ is locally Lipschitz in its first argument, for each $x\in\ms{M}$
and $u\in\ts{U}$\@, there exists $L_u\in\realp$ and a neighbourhood
$\nbhd{X}_u$ of $x$ such that
\begin{equation*}
\d_{\ms{N}}(f(x_1,u),f(x_2,u))\le2L_u\d_{\ms{M}}(x_1,x_2)
\end{equation*}
for every $x_1,x_2\in\nbhd{X}_u$\@.  Continuity of $f$ and the metric ensures
that there exists a neighbourhood $\nbhd{Z}_u\subset\ts{U}$ of $u$ such that
\begin{equation*}
\d_{\ms{N}}(f(x_1,v_1),f(x_2,v_2))<L_u\d_{\ms{M}}(x_1,x_2)
\end{equation*}
for every $x_1,x_2\in\nbhd{X}_u$ and every $v_1,v_2\in\nbhd{Z}_u$\@, possibly
by also shrinking $\nbhd{X}_u$\@.  Let
\begin{equation*}
\nbhd{X}=\cap_{j=1}^k\nbhd{X}_{u_j},\quad L=\max\{L_{u_1},\dots,L_{u_k}\}.
\end{equation*}

Now let $x_1,x_2\in\nbhd{X}$ and let $f(x_1,u)\in F(x_1)$ for some
$u\in\ts{U}$\@.  Then $u\in\nbhd{Z}_j$ for some $j\in\{1,\dots,k\}$ and so
\begin{equation*}
\d_{\ms{N}}(f(x_1,u),f(x_2,u))<L_{u_j}\d_{\ms{M}}(x_1,x_2)\le L\d_{\ms{M}}(x_1,x_2),
\end{equation*}
giving $f(x_1,u)\in\oball_{\ms{N}}(f(x_2,u),L\d_{\ms{M}}(x_1,x_2))$\@.  We
thus conclude that
\begin{equation*}
F(x_1)\subset\bigcup_{y\in F(x_2)}\oball_{\ms{N}}(y,L\d_{\ms{M}}(x_1,x_2)),
\end{equation*}
as desired.
\end{subproof}
\end{prooflemma}

Let us define the family of vector fields
\begin{equation*}
\cohull(\sX')=\asetdef{\sum_{j=1}^ku^jX_j}{\vect{u}\in\cohull(U)}.
\end{equation*}
Since the vector fields $X_1,\dots,X_k$ are locally Lipschitz, by the lemma
the differential inclusion defined by $F$ is locally Lipschitz.  It is also
clearly compact-valued.  By the relaxation theorem of \citet{AFF:67} and
\citet{TW:62}\@, it follows that the reachable set for $\Sigma$ is dense in
the reachable set for $\cohull(\Sigma)$\@.  Since both $\Sigma$ and
$\cohull(\Sigma)$ are symmetric,~\ie~if $v_x\in F(x)$ then $-v_x\in F(x)$ and
similarly for $\cohull(F)$\@, it follows that the reachable set and the orbit
agree (see Proposition~4.3 in \cite{JB:01}\@, for example).  Thus
$\Orb(x_0,\sX')$ is dense in $\Orb(x_0,\cohull(\sX'))$\@.  By the Orbit
Theorem, $\Orb(x_0,\sX')$ is an immersed submanifold of the immersed
submanifold $\Orb(x_0,\cohull(\sX'))$\@.  Since
$x_0\in\interior(\Orb(x_0,\sX'))$ and
$x_0\in\interior(\Orb(x_0,\cohull(\sX')))$ (interior being taken in the orbit
topology in each case), this implies that the tangent spaces to the two
orbits agree at $x_0$\@.  By shrinking $\nbhd{U}$ we can ensure that
$\Orb(x_0,\sX')=\Orb(x_0,\cohull(\sX'))$\@.

Now note that, if $X\in\sections[r]{\dist{D}(\sX)|\nbhd{U}}$\@, then
$X(x)\in\vecspan[\real]{X_1(x),\dots,X_k(x)}$ for each $x\in\nbhd{U}$ since
$\ifam{X_1,\dots,X_k}$ generate $F_{\sM}(\nbhd{U})$\@.  Therefore, points in
$\Orb(x_0,\sections[r]{\dist{D}(\sX)|\nbhd{U}})$ are endpoints of
concatenations of curves whose tangent vectors are positive multiples of
tangent vectors in $\cohull(F)$\@.  Thus points in
$\Orb(x_0,\sections[r]{\dist{D}(\sX)|\nbhd{U}}))$ are endpoints of
concatenations of curves tangent to $\Orb(x_0,\cohull(\sX'))$\@.  Thus
$\Orb(x_0,\sections[r]{\dist{D}(\sX)|\nbhd{U}})\subset
\Orb(x_0,\cohull(\sX'))$\@.  Since the opposite inclusion is obvious, we have
$\Orb(x_0,\sections[r]{\dist{D}(\sX)|\nbhd{U}})=\Orb(x_0,\cohull(\sX'))$\@.

Putting the above arguments together, the tangent spaces at $x_0$ of
$\Orb(\sections[r]{\dist{D}(\sX)|\nbhd{U}})$ and $\Orb(x_0,\sX|\nbhd{U})$
agree.  By Theorems~\ref{the:Linfty-orbit}
and~\ref{the:fingen-orbit-theorem} we have
\begin{equation*}
\dist{L}^{(\infty)}(\sections[r]{\dist{D}(\sX)})_{x_0}\subset
\tb[x_0]{\Orb(x_0,\sections[r]{\dist{D}(\sX)|\nbhd{U}})}=
\tb[x_0]{\Orb(x_0,\sX|\nbhd{U})}=\dist{L}^{(\infty)}(\sX)_{x_0}.
\end{equation*}
Since the inclusion $\dist{L}^{(\infty)}(\sX)_{x_0}\subset
\dist{L}^{(\infty)}(\sections[r]{\dist{D}(\sX)})_{x_0}$ is clear and since
$x_0$ is arbitrary, the theorem follows.
\end{proof}
\end{theorem}

The result has two interesting and often applicable corollaries.
\begin{corollary}
Let\/ $\man{M}$ be a\/ $\C^\infty$-manifold, let\/
$\sX\subset\sections{\tb{\man{M}}}$\@, and let\/ $x$ be a regular point of\/
$\dist{D}(\sX)$ and of\/ $\dist{L}^{(\infty)}(\sX)$\@.  Then\/
$\dist{L}^{(\infty)}(\sX)_x=
\dist{L}^{(\infty)}(\sections[\infty]{\dist{D}(\sX)})_x$ for every\/
$x\in\man{M}$\@.
\begin{proof}
This follows from Theorem~\ref{the:LinftyX-D(X)}\@, along with
Theorem~\ref{the:locally-free-sheaf}\@.
\end{proof}
\end{corollary}

\begin{corollary}\label{cor:LinftyX-D(X)-analytic}
Let\/ $\man{M}$ be a\/ $\C^\omega$-manifold and let\/
$\sX\subset\sections[\omega]{\tb{\man{M}}}$\@.  Then\/
$\dist{L}^{(\infty)}(\sX)_x=
\dist{L}^{(\infty)}(\sections[\omega]{\dist{D}(\sX)})_x$\@.
\begin{proof}
This follows from Theorem~\ref{the:LinftyX-D(X)}\@, along with
Theorem~\ref{the:analytic-fingen} (noting that $\sF_{\sX}$ is obviously
patchy).
\end{proof}
\end{corollary}

Let us consider a few examples that illustrate the subtlety of the preceding
results.
\begin{examples}\label{eg:LinftyX-D(X)}
We consider $\man{M}=\real^2$ and the distribution $\dist{D}$ given by
\begin{equation*}
\dist{D}_{(x_1,x_2)}=\begin{cases}\tb[(x_1,x_2)]{\real^2},&x_1\not=0,\\
\vecspan[\real]{\pderiv{}{x_1}},&x_1=0.\end{cases}
\end{equation*}
This distribution is generated by any pair of vector fields
\begin{equation*}
X_1(x_1,x_2)=\pderiv{}{x_1},\quad X_2(x_1,x_2)=f(x_1)\pderiv{}{x_2},
\end{equation*}
where $f\in\func{\real}$ satisfies $f^{-1}(0)=\{0\}$\@.  Thus $\dist{D}$ is a
smooth distribution.  Moreover, we claim that if $X$ is any section of
$\dist{D}$ then $X=f^1X_1+f^2X_2$ for some $f^1,f^2\in\func{\real^2}$\@,
provided we take $f$ defined by $f(x)=x$\@.  Indeed, let us write
\begin{equation*}
X=g^1\pderiv{}{x_1}+g^2\pderiv{}{x_2},
\end{equation*}
noting that we must have $g^2(0,x_2)=0$ for every $x_2\in\real$\@.  We write
\begin{equation*}
g^2(x_1,x_2)=\int_0^{x_1}\pderiv{g^2}{x_1}(\xi,x_2)\,\d{\xi}=
x_1\int_0^1\pderiv{g^2}{x_1}(x_1\eta,x_2)\,\d{\eta}.
\end{equation*}
Thus our claim follows by taking
\begin{equation*}
f^1(x_1,x_2)=g^1(x_1,x_2),\quad
f^2(x_1,x_2)=\int_0^1\pderiv{g^2}{x_1}(x_1\eta,x_2)\,\d{\eta}.
\end{equation*}
This shows that, not only does $\dist{D}$ have a finite number of generators
(as per Theorem~\ref{the:global-span}), but also that $\sections{\dist{D}}$
is finitely generated.  Note that it is possible to choose generators for
$\dist{D}$ that do not generate $\sections{\dist{D}}$\@,~\eg~by taking
$f(x)=x^2$ we see that $X_1$ and $X_2$ as above have this
property,~\cf~Example~\ref{eg:!stalk-generate}\@.  By choosing more
pathological generators,~\eg~by taking
\begin{equation*}
f(x)=\begin{cases}\eul^{-1/x^2},&x\not=0,\\0,&x=0,\end{cases}
\end{equation*}
one imagines that the algebraic properties of the distribution should
deteriorate.  We shall see now that this is true as concerns the Lie algebra
generated by the generators.

We take $\sX=\ifam{X_1,X_2}$ and consider a few $f$'s.
\begin{compactenum}
\item First let us consider $f(x)=x$\@.  We compute
\begin{equation*}
[X_1,X_2](x_1,x_2)=\pderiv{}{x_2}.
\end{equation*}
Therefore, $\dist{L}^{(\infty)}(\sX)=\tb{\real^2}$\@.  Thus we must have
$\dist{L}^{(\infty)}(\sX)=\dist{L}^{(\infty)}(\sections[r]{\dist{D}(\sX)})$\@.
Since $f$ is analytic, this is in agreement with
Corollary~\ref{cor:LinftyX-D(X)-analytic}\@.

\item Next consider $f(x)=x^2$\@.  In this case we have
\begin{equation*}
[X_1,X_2](x_1,x_2)=2x_1\pderiv{}{x_2},\quad
[X_1,[X_1,X_2]](x_1,x_2)=2\pderiv{}{x_2}.
\end{equation*}
Thus we can again conclude that $\dist{L}^{(\infty)}(\sX)=\tb{\real^2}$\@,
implying that
$\dist{L}^{(\infty)}(\sX)=\dist{L}^{(\infty)}(\sections[r]{\dist{D}(\sX)})$\@.
This again is consistent with Corollary~\ref{cor:LinftyX-D(X)-analytic}\@.
Note, however, the distribution $\dist{L}^{(\infty)}(\sX)$ is generated by
different brackets than was the case when we took $f(x)=x$\@.  Thus the fact
that
$\dist{L}^{(\infty)}(\sX)=\dist{L}^{(\infty)}(\sections[r]{\dist{D}(\sX)})$
is less obvious in this case.

\item The final case we consider is
\begin{equation*}
f(x)=\begin{cases}\eul^{-1/x^2},&x\not=0,\\0,&x=0.\end{cases}
\end{equation*}
We note that
\begin{equation*}
\dist{D}(\sX)_{x_1,x_2}=\begin{cases}\tb[(x_1,x_2)]{\real^2},&x_1\not=0,\\
\vecspan[\real]{\pderiv{}{x_1}},&x_1=0.\end{cases}
\end{equation*}
Thus, for example, the vector fields
\begin{equation*}
X'_1(x_1,x_2)=\pderiv{}{x_1},\quad X'_2(x_1,x_2)=x_1\pderiv{}{x_2}
\end{equation*}
generate $\dist{D}(\sX')$\@.  As above, one can compute
$[X_1,X_2]=\pderiv{}{x_2}$\@, and so we have
$\dist{L}^{(\infty)}(\sections[r]{\dist{D}(\sX)})=\tb{\real^2}$\@.  However,
one can easily show that $\dist{L}^{(\infty)}(\sX)=\dist{D}(\sX)$ and so
$\dist{L}^{(\infty)}(\sX)\subsetneq
\dist{L}^{(\infty)}(\sections[r]{\dist{D}(\sX)}q)$\@.  By
Theorem~\ref{the:LinftyX-D(X)} we conclude that
$\sL^{(\infty)}(\modgen{\sX})$ is not locally finitely generated.  This is a
problem with the generators $\sX$ for $\dist{D}(\sX)$ being smooth but not
analytic.\oprocend
\end{compactenum}
\end{examples}

Theorem~\ref{the:LinftyX-D(X)} and its corollaries can also be adapted to
subsheaves instead of submodules.
\begin{theorem}
Let\/ $r\in\{\infty,\omega\}$\@, let\/ $\man{M}$ be a\/ $\C^r$-manifold, and
let\/ $\sX=\ifam{X(\nbhd{U})}_{\nbhd{U}\,\textrm{open}}$ be a subsheaf of
sets of the sheaf\/ $\ssections[r]{\tb{\man{M}}}$ such that
\begin{compactenum}[(i)]
\item $\modgen{\sX}$ is a locally finitely generated subsheaf of\/
$\ssections[r]{\tb{\man{M}}}$ and
\item the subsheaf (by Proposition~\ref{prop:liealgmodule})\/
$\sL^{(\infty)}(\modgen{\sX})$ is locally finitely generated.
\end{compactenum}
Then\/
$\dist{L}^{(\infty)}(\sX)=\dist{L}^{(\infty)}(\ssections[r]{\dist{D}(\sX)})$\@.
\begin{proof}
This follows from Theorem~\ref{the:LinftyX-D(X)}\@.
\end{proof}
\end{theorem}

\begin{corollary}
let\/ $\man{M}$ be a\/ $\C^\infty$-manifold, let\/
$\sX=\ifam{X(\nbhd{U})}_{\nbhd{U}\,\textrm{open}}$ be a subsheaf of sets of
the sheaf\/ $\ssections[r]{\tb{\man{M}}}$\@, and let\/ $x$ be a regular point
of\/ $\dist{D}(\sX)$ and of\/ $\dist{L}^{(\infty)}(\sX)$\@.  Then\/
$\dist{L}^{(\infty)}(\sX)=\dist{L}^{(\infty)}(\ssections[r]{\dist{D}(\sX)})$\@.
\end{corollary}

\begin{corollary}
Let\/ $\man{M}$ be a\/ $\C^\omega$-manifold and let\/
$\sX=\ifam{X(\nbhd{U})}_{\nbhd{U}\,\textrm{open}}$ be a subsheaf of sets of
the sheaf\/ $\ssections[r]{\tb{\man{M}}}$ for which\/ $\modgen{\sX}$ is
patchy.  Then\/
$\dist{L}^{(\infty)}(\sX)=\dist{L}^{(\infty)}(\ssections[r]{\dist{D}(\sX)})$\@.
\end{corollary}

\subsection{Distributions and subsheaves invariant under vector fields and
diffeomorphisms}\label{subsec:invariance}

Using the constructions of the preceding section, in this section we
introduce the notion of distributions and subsheaves that are invariant under
vector fields and diffeomorphisms.  In this section, when we say ``subsheaf''
of $\ssections[r]{\tb{\man{M}}}$ we mean a subsheaf of
$\sfunc[r]{\man{M}}$-modules.
\begin{definition}
Let $r\in\{\infty,\omega\}$\@, let $\man{M}$ be a manifold of class $\C^r$\@,
let $\dist{D}$ be a distribution of class $\C^r$\@, let $\sF$ be a subsheaf
of $\ssections[r]{\tb{\man{M}}}$\@, let $X$ be a $\C^r$-vector field, and let
$(\Phi,\nbhd{U})$ be a $\C^r$-local diffeomorphism.  The distribution
$\dist{D}$
\begin{compactenum}[(i)]
\item is \defn{invariant} under $X$ if $[X,Y]\in\sections[r]{\dist{D}}$ for
every $Y\in\sections[r]{\dist{D}}$ and
\item is \defn{invariant} under $(\Phi,\nbhd{U})$ if
$\Phi^*Y\in\sections[r]{\dist{D}|\Phi(\nbhd{U})}$ for every
$Y\in\sections[r]{\dist{D}|\nbhd{U}}$\@.\savenum
\end{compactenum}
The subsheaf $\sF$
\begin{compactenum}[(i)]\resumenum
\item is \defn{invariant} under $X$ at $x$ if $[[X,Y]]_x\in\sF_x$ for every
$[Y]_x\in\sF_x$\@,
\item is \defn{invariant} under $(\Phi,\nbhd{U})$ at $x\in\nbhd{U}$ if
$\Phi^*_x[X]_x\in\sF_{\Phi(x)}$ for every $[X]_x\in\sF_x$\@,
\item is \defn{invariant} under $X$ if it is invariant under $X$ at $x$ for
each $x\in\man{M}$\@, and
\item is \defn{invariant} under $(\Phi,\nbhd{U})$ if it is invariant under
$\Phi$ at $x$ for each $x\in\nbhd{U}$\@.\oprocend
\end{compactenum}
\end{definition}

One would like to think that invariance under a vector field and invariance
under its flow are equivalent.  This is true under suitable hypotheses.  Let
us first look at the case of invariant sheaves.  Here the statement requires
us to exploit the topologies on stalks of subsheaves of vector bundles from
Section~\ref{subsec:stalk-topology}\@.
\begin{theorem}\label{the:invariant-sheaves}
Let\/ $r\in\{\infty,\omega\}$\@, let\/ $\man{M}$ be a\/ $\C^r$-manifold,
let\/ $\sF=\ifam{F(\nbhd{U})}_{\nbhd{U}\,\textrm{open}}$ be a subsheaf of\/
$\ssections[r]{\tb{\man{M}}}$\@, let\/ $X\in\sections[r]{\tb{\man{M}}}$\@,
and let\/ $x\in\man{M}$\@.  Consider the following two statements:
\begin{compactenum}[(i)]
\item \label{pl:Finvariant1} $\sF$ is invariant under\/ $X$ at\/ $x$\@;
\item \label{pl:Finvariant2} for each\/ $T\in\realp$ there exists a
neighbourhood\/ $\nbhd{U}$ of\/ $x$ such that $\sF$ is invariant under\/
$(\flow{X}{t},\nbhd{U})$ at\/ $x$ for each\/ $t\in\interval[-T,T]$\@.
\end{compactenum}
Then~\eqref{pl:Finvariant2}$\implies$\eqref{pl:Finvariant1} if\/ $\sF_x$ is a
closed submodule of\/ $\gsections[r]{x}{\tb{\man{M}}}$\@, and
\eqref{pl:Finvariant1}$\implies$\eqref{pl:Finvariant2} if there exists a
neighbourhood\/ $\nbhd{V}$ of\/ $x$ such that\/ $F(\nbhd{V})$ is a finitely
generated\/ $\func[r]{\nbhd{V}}$-module.
\begin{proof}
\eqref{pl:Finvariant1}$\implies$\eqref{pl:Finvariant2} Let $x\in\man{M}$ and
let $X\in\sections[r]{\tb{\man{M}}}$ satisfy
\begin{equation*}
\setdef{[X,Y]_x}{[Y]_x\in\sF_x}\subset\sF_x.
\end{equation*}
Suppose that $Y_1,\dots,Y_k$ generate $F(\nbhd{V})$ for some neighbourhood
$\nbhd{V}$ of $x$\@.  By hypothesis,
\begin{equation*}
[X,Y_j]_x=\sum_{i=1}^k[f^i_j]_x[Y_i]_x,\qquad j\in\{1,\dots,k\},
\end{equation*}
for some $[f^i_j]_x\in\gfunc[r]{x}{\man{M}}$\@, $i,j\in\{1,\dots,k\}$\@.  As
this expression involves only finitely many germs, we may assume $\nbhd{V}$
sufficiently small that
\begin{equation*}
[X,Y_j](y)=\sum_{i=1}^kf^i_j(y)Y_i(y),\qquad y\in\nbhd{V}.
\end{equation*}
Let $\nbhd{U}\subset\nbhd{V}$ be sufficiently small that
$\flow{X}{t}(y)\in\nbhd{V}$ for every $t\in\interval[{-T},T]$ and every
$y\in\nbhd{U}$\@.  For $t\in\interval[{-T},T]$ and $y\in\nbhd{U}$ define
$v_j(t,y)=(\flow{X}{t})^*Y_j(y)$\@, $j\in\{1,\dots,k\}$\@, so that $t\mapsto
v_j(t,y)$ is a curve in $\tb[y]{\man{M}}$\@.  By
\cite[Theorem~4.2.19]{RA/JEM/TSR:88}
\begin{align*}
\deriv{}{t}v_j(t,y)=&\;\deriv{}{t}(\flow{X}{t})^*Y_j(y)=
(\flow{X}{t})^*[X,Y_j](y)=(\flow{X}{t})^*\left(\sum_{i=1}^kf^i_jY_i\right)(y)\\
=&\;\sum_{i=1}^k(\flow{X}{t})^*f^i_j(y)(\flow{X}{t})^*Y_i(y)=
\sum_{i=1}^k(\flow{X}{t})^*f^i_j(y)v_i(t).
\end{align*}
Define $\mat{A}_y(t)\in\real^{k\times k}$ by
\begin{equation*}
A^i_{y,j}(t)=(\flow{X}{t})^*f^i_j(y)
\end{equation*}
and let $\map{\mat{\Psi}_y}{\real}{\real^{k\times k}}$ be the solution to the
matrix initial value problem
\begin{equation*}
\deriv{}{t}\mat{\Psi}_y(t)=\mat{A}_y(t)\mat{\Psi}_y(t),
\qquad\mat{\Psi}_y(0)=\mat{I}_k.
\end{equation*}
We claim that
\begin{equation*}
v_j(t,y)=\sum_{i=1}^k\Psi^i_{y,j}(t)Y_i(y).
\end{equation*}
Indeed, 
\begin{align*}
\deriv{}{t}\left(\sum_{i=1}^k\Psi^i_{y,j}(t)Y_i(y)\right)=&\;
\sum_{i=1}^k\deriv{}{t}\Psi^i_{y,j}(t)Y_i(y)
=\sum_{i,l=1}^kA^l_{y,j}(t)\Psi^i_{y,l}(t)Y_i(y)\\
=&\;\sum_{l=1}^k(\flow{X}{t})^*f^l_j(y)
\left(\sum_{i=1}^k\Psi^i_{y,l}(t)Y_i(y)\right).
\end{align*}
Moreover,
\begin{equation*}
\sum_{i=1}^k\Psi^i_{y,j}(0)Y_i(y)=Y_j(y),\qquad v_j(0,y)=Y_j(y).
\end{equation*}
Thus
\begin{equation*}
t\mapsto v_j(t,y)\enspace\textrm{and}\enspace
t\mapsto\sum_{i=1}^k\Psi^i_{y,j}(t)Y_i(y)
\end{equation*}
satisfy the same differential equation with the same initial condition.  Thus
they are equal.  This gives
\begin{equation*}
(\flow{X}{t})^*Y_j(y)=\sum_{i=1}^k\Psi^i_{y,j}(t)Y_i(y)
\end{equation*}
for every $t\in\interval[{-T},T]$ and $y\in\nbhd{U}$\@.  Now let
$[Y]_x\in\sF_x$ and suppose that $Y$ is a local section over
$\nbhd{W}\subset\nbhd{U}$\@.  Following the proof of
Proposition~\pldblref{prop:locfingen}{pl:locfingen4}\@, we can write
\begin{equation*}
[Y]_x=\sum_{j=1}^k[\eta^j]_x[Y_j]_x
\end{equation*}
for $[\eta^j]_x\in\gfunc[r]{x}{\man{M}}$\@, $j\in\{1,\dots,k\}$\@.
Therefore, possibly after shrinking $\nbhd{W}$\@, we can write
\begin{equation*}
Y=\eta^1(Y_j|\nbhd{W})+\dots+\eta^k(Y_k|\nbhd{W})
\end{equation*}
for some $\eta^1,\dots,\eta^k\in\func[r]{\nbhd{W}}$\@.  Therefore, for
$y\in\nbhd{W}$ and $t\in\interval[-T,T]$\@,
\begin{equation}\label{eq:invariant-sheaf2}
(\flow{X}{t})^*Y(y)=\sum_{j=1}^k\eta^j(y)(\flow{X}{t})^*Y_j(y)=
\sum_{j=1}^k\eta^j(y)\sum_{i=1}^k\Psi^i_{y,j}(t)Y_i(y),
\end{equation}
and so $[(\flow{X}{t})^*Y]_x\in\sF_x$\@.  This gives this part of the
theorem.

\eqref{pl:Finvariant2}$\implies$\eqref{pl:Finvariant1} This part of the proof
we carry out separately for $r=\infty$ and $r=\omega$\@.

Let us first consider the case $r=\infty$\@.  Let $[Y]_x\in\sF_x$\@, let
$\nbhd{V}$ be a sufficiently small neighbourhood that $Y$ can be taken as
being defined on $\nbhd{V}$\@, and let $\nbhd{U}\subset\nbhd{V}$ and
$T\in\realp$ be sufficiently small that $\flow{X}{t}(y)\in\nbhd{V}$ for every
$y\in\nbhd{U}$ and $t\in\interval[-T,T]$\@.  We can without loss of
generality assume that $\nbhd{V}$ is the domain of an admissible chart, and
so deal with local representatives.  Let us state a lemma that will be
helpful here and in the second part of the proof where we treat the real
analytic case.
\begin{prooflemma}
Let\/ $\nbhd{U}_1\subset\real^{n_1}$ and\/ $\nbhd{U}_2\subset\real^{n_2}$ be
open sets, let\/ $\alg{V}$ be a finite-dimensional normed vector space, and
let\/ $\map{f}{\nbhd{U}_1\times\nbhd{U}_2}{\alg{V}}$ be continuously
differentiable.  Then, for any compact set\/ $K\subset\nbhd{U}_2$\@, any\/
$\vect{x}_1\in\nbhd{U}_1$\@, and any\/ $\vect{v}\in\real^{n_1}$\@,
\begin{equation*}
\lim_{t\to0}\left(\sup\asetdef{
\adnorm{\frac{1}{t}(f(\vect{x}_1+t\vect{v},\vect{x}_2)-
f(\vect{x}_1,\vect{x}_2))-\plinder{1}{f}(\vect{x}_1,\vect{x}_2)\cdot\vect{v}}}
{\vect{x}_2\in K}\right)=0.
\end{equation*}
\begin{subproof}
Define
\begin{equation*}
g_{\vect{x}_1,\vect{v}}(t,\vect{x}_2)=f(\vect{x}_1+t\vect{v},\vect{x}_2)-
f(\vect{x}_1,\vect{x}_2)-t\plinder{1}{f}(\vect{x}_1,\vect{x}_2)\cdot\vect{v}.
\end{equation*}
Note that
$(t,\vect{x}_2)\mapsto\frac{1}{t}g_{\vect{x}_1,\vect{v}}(t,\vect{x}_2)$ is
continuous for $t$ sufficiently near zero (taking its value when $t=0$ to be
zero).  This implies that, for every $\epsilon\in\realp$\@, there exists
$\delta\in\realp$ such that
$\dnorm{\frac{1}{t}g_{\vect{x}_1,\vect{v}}(t,\vect{x}_2)}<\epsilon$ for every
$t\in\interval[-\delta,\delta]$ and $\vect{x}_2\in K$\@.  From this the
result follows.
\end{subproof}
\end{prooflemma}

Applying the lemma successively to all derivatives, and recalling the
notation from Section~\ref{subsec:stalk-topology}\@, we have
\begin{equation*}
\lim_{t\to0}\adnorm{\frac{1}{t}((\flow{X}{t})^*Y-Y)-[X,Y]}_{r,K}=0
\end{equation*}
for every compact $K\subset\nbhd{U}$ and every $r\in\integernn$\@, using
\cite[Theorem~4.2.19]{RA/JEM/TSR:88}\@.  This shows that
\begin{equation*}
\lim_{t\to0}\frac{1}{t}((\flow{X}{t})^*Y-Y)=[X,Y]
\end{equation*}
in $\sections[\infty]{\tb{\man{M}}|\nbhd{U}}$ in the weak
$\C^\infty$-topology described in Section~\ref{subsec:stalk-topology}\@.
Since we are assuming that $\sF_x$ is closed, $r_{\nbhd{U},x}^{-1}(\sF_x)$ is
closed (since $r_{\nbhd{U},x}$ is continuous as in \cite[\S19.5]{GK:69}).  By
hypothesis, $\frac{1}{t}((\flow{X}{t})^*Y-Y)\in r_{\nbhd{U},x}^{-1}(\sF_x)$
for every sufficiently small $t$\@, it follows that
$[X.Y]\in r_{\nbhd{U},x}^{-1}(\sF_x)$\@.  Thus
$r_{\nbhd{U},x}([X,Y])\in\sF_x$ for every sufficiently small neighbourhood
$\nbhd{U}$ of $x$ and so $[[X,Y]]_x\in\sF_x$\@, which is the result in the
smooth case.

Now we consider the real analytic case.  We let $[Y]_x\in\sF_x$ and suppose
that $Y$ is defined in a neighbourhood $\nbhd{V}$ of $x$\@.  As above, let
$\nbhd{U}\subset\nbhd{V}$ be a neighbourhood of $x$ and let $T\in\realp$ be
such that $\flow{X}{t}(y)\in\nbhd{V}$ for every $y\in\nbhd{V}$ and
$t\in\interval[-T,T]$\@.  We suppose that $\nbhd{V}$ is the domain of a
coordinate chart and work locally in Euclidean space.  We let $\ol{\nbhd{Z}}$
be a neighbourhood in the complexification of
$\interval(-T,T)\times\nbhd{V}$\@.  Let $\ol{x}$ be the image of $x$ in this
complexification.  By shrinking $\ol{\nbhd{Z}}$ we extend the vector fields
$X$ and $Y$ to holomorphic vector fields $\ol{X}$ and $\ol{Y}$\@,
respectively, on $\ol{\nbhd{Z}}$\@.  We can extend the real analytic map
\begin{equation}\label{eq:invariant-sheaf1}
\interval(-T,T)\times\nbhd{U}\ni(t,y)\mapsto\flow{X}{t}(y)\in\nbhd{V}
\end{equation}
to a holomorphic map defined in a neighbourhood of
$\interval(-T,T)\times\nbhd{V}$ in $\ol{\nbhd{Z}}$\@.  Let us denote the
variables in the complexification by $(\tau,z)$\@.  Let us write the
complexification of the map~\eqref{eq:invariant-sheaf1} as
$(\tau,z)\mapsto\ol{\Phi}^{\ol{X}}_\tau(z)$\@.  Note that
$\tderivatzero{}{\tau}(\ol{\Phi}^{\ol{X}}_\tau)^*\ol{Y}(z)$ is the
holomorphic extension of $[X,Y]$\@, since it agrees with $[X,Y]$ at points in
$\nbhd{V}$ using \cite[Theorem~4.2.19]{RA/JEM/TSR:88}\@.  Let us denote this
extension by $[\ol{X},\ol{Y}]$\@.  Using the lemma above (noting that complex
differentiation with respect to $\tau=t+\imag s$ is represented by the
Jacobian of the corresponding real mapping) and recalling the notation from
Section~\ref{subsec:stalk-topology}\@,
\begin{equation*}
\lim_{\tau\to0}\adnorm{\frac{1}{\tau}((\ol{\Phi}^{\ol{X}}_\tau)^*\ol{Y}-\ol{Y})-
[\ol{X},\ol{Y}]}_K=0
\end{equation*}
for every compact set $K\subset\nbhd{V}$\@.  We now recall the notation
introduced preceding the proof of Theorem~\ref{the:analytic-closed-module}\@.
Since every submodule of $\gsections[\omega]{x}{\man{E}}$ is closed by
Theorem~\ref{the:analytic-closed-module}\@, the definition of the topology on
$\gsections[\omega]{x}{\man{E}}$ ensures that
$r_{\ol{\nbhd{Z}},\ol{x}}^{-1}(\rho_{\ol{x},x}^{-1}(\sF_x))$ is
closed,~\cf~\cite[\S19.5]{GK:69}\@.  Thus
$[\ol{X},\ol{Y}]\in
r_{\ol{\nbhd{Z}},\ol{x}}^{-1}(\rho_{\ol{x},x}^{-1}(\sF_x))$ since, by
hypothesis,
$\frac{1}{\tau}((\ol{\Phi}^{\ol{X}}_\tau)^*\ol{Y}-\ol{Y})\in
r_{\ol{\nbhd{Z}},\ol{x}}^{-1}(\rho_{\ol{x},x}^{-1}(\sF_x))$ for $\tau$
sufficiently near zero.  It follows, therefore, that $[[X,Y]]_x\in\sF_x$\@,
giving the theorem.
\end{proof}
\end{theorem}

The theorem then has the following immediate corollary, when combined with
Theorem~\ref{the:analytic-closed-module} and
Proposition~\ref{prop:analytic-germs}\@.
\begin{corollary}
Let\/ $\man{M}$ be a real analytic manifold, let\/ $\sF$ be a subsheaf of\/
$\ssections[\omega]{\tb{\man{M}}}$\@, let\/
$X\in\sections[r]{\tb{\man{M}}}$\@, and let\/ $x\in\man{M}$\@.  Then the
following two statements are equivalent:
\begin{compactenum}[(i)]
\item $\sF$ is invariant under\/ $X$ at\/ $x$\@;
\item for each\/ $T\in\realp$ there exists a neighbourhood\/ $\nbhd{U}$ of\/
$x$ such that $\sF$ is invariant under\/ $(\flow{X}{t},\nbhd{U})$ at\/ $x$
for each\/ $t\in\interval[-T,T]$\@.
\end{compactenum}
\end{corollary}

Next we consider the relationship between invariance of distributions under
vector fields and their flows.
\begin{theorem}\label{the:Dinvariant}
Let\/ $r\in\{\infty,\omega\}$\@, let\/ $\man{M}$ be a\/ $\C^r$-manifold,
let\/ $\dist{D}$ be a\/ $\C^r$-distribution, and let\/
$X\in\sections[r]{\tb{\man{M}}}$\@.  Consider the following two statements:
\begin{compactenum}[(i)]
\item \label{pl:Dinvariant1} $\dist{D}$ is invariant under\/ $X$\@;
\item \label{pl:Dinvariant2} for each\/ $T\in\realp$ and each\/ $x\in\man{M}$
there exists a neighbourhood\/ $\nbhd{U}$ of\/ $x$ such that\/
$\dist{D}|\nbhd{U}$ is invariant under\/ $(\flow{X}{t},\nbhd{U})$ for each\/
$t\in\interval[-T,T]$\@.
\end{compactenum}
Then~\eqref{pl:Dinvariant2}$\implies$\eqref{pl:Dinvariant1} always, and
\eqref{pl:Dinvariant1}$\implies$\eqref{pl:Dinvariant2} if\/
$\ssections[r]{\dist{D}}$ is locally finitely generated.
\begin{proof}
\eqref{pl:Dinvariant1}$\implies$\eqref{pl:Dinvariant2} This follows from the
corresponding part of Theorem~\ref{the:invariant-sheaves}\@, taking
particular note that in the proof of the theorem the
equation~\eqref{eq:invariant-sheaf2} holds in a neighbourhood of $x$\@.

\eqref{pl:Dinvariant2}$\implies$\eqref{pl:Dinvariant1} Let $x\in\man{M}$ and
let $\epsilon\in\realp$ be such that $\flow{X}{t}(x)$ exists for
$t\in\interval(-\epsilon,\epsilon)$\@.  Then, for each
$t\in\interval(-\epsilon,\epsilon)$\@, we have
\begin{equation*}
(\flow{X}{t})_*Y(x)\in\dist{D}_x.
\end{equation*}
Therefore, since $\dist{D}_x$ is a subspace, we
use~\cite[Theorem~4.2.19]{RA/JEM/TSR:88} to compute
\begin{equation*}
[Y,X](x)=\derivatzero{}{t}(\flow{X}{t})_*Y(x)\in\dist{D}_x,
\end{equation*}
as desired.
\end{proof}
\end{theorem}

Of course, the preceding theorem then has the following immediate corollary,
when combined with Proposition~\ref{prop:analytic-germs}\@.
\begin{corollary}
Let\/ $\man{M}$ be a real analytic manifold, let\/ $\dist{D}$ be an analytic
distribution on\/ $\man{M}$\@, and let\/ $X\in\sections[r]{\tb{\man{M}}}$\@.
Then the following two statements are equivalent:
\begin{compactenum}[(i)]
\item $\dist{D}$ is invariant under\/ $X$ at\/ $x$\@;
\item for each\/ $T\in\realp$ and each\/ $x\in\man{M}$ there exists a
neighbourhood\/ $\nbhd{U}$ of\/ $x$ such that\/ $\dist{D}|\nbhd{U}$ is
invariant under\/ $(\flow{X}{t},\nbhd{U})$ for each\/
$t\in\interval[-T,T]$\@.
\end{compactenum}
\end{corollary}

\section{The Orbit Theorem}\label{sec:orbit-theorem}

In this section we consider an important theorem that is not so immediately
connected to the theory of distributions, but, as we shall see, leads to
important theorems regarding special classes of distributions.  Contributions
to the Orbit Theorem have been made by
\citet{PS:74a,PS:74b,HJS:73,AJK:74,CL:70b,RH:62}\@.

\subsection{Partially defined vector fields}

In our development of the Orbit Theorem, we follow \citet{HJS:73} and
consider the quite general setting where vector fields are not only not
necessarily complete, but only defined on an open subset of the manifold.  In
this section we develop the notation and rules for dealing with these sorts
of vector fields.
\begin{definition}
Let $r\in\{\infty,\omega\}$ and let $\man{M}$ be a manifold of class
$\C^r$\@.
\begin{compactenum}[(i)]
\item A \defn{partially defined vector field of class $\C^r$} is a pair
$(X,\nbhd{U})$ where $\nbhd{U}\subset\man{M}$ is open and where $X$ is a
$\C^r$-vector field on $\nbhd{U}$\@.
\item A family $\sX=\ifam{(X_j,\nbhd{U}_j)}_{j\in J}$ of partially defined
vector fields of class $\C^r$ is \defn{everywhere defined} if, for every
$x\in\man{M}$\@, there exists $j\in J$ such that $x\in\nbhd{U}_j$\@.\oprocend
\end{compactenum}
\end{definition}

Let us outline, for clarity, how some of the standard vector field
constructions apply to partially defined vector fields.  The \defn{sum} of
two partially defined vector fields $(X,\nbhd{U})$ and $(Y,\nbhd{V})$ is the
partially defined vector field
\begin{equation*}
(X,\nbhd{U})+(Y,\nbhd{V})=(X|\nbhd{U}\cap\nbhd{V}+Y|\nbhd{U}\cap\nbhd{V},
\nbhd{U}\cap\nbhd{V}).
\end{equation*}
Given partially defined vector fields $(X,\nbhd{U})$ and $(Y,\nbhd{V})$\@,
their \defn{Lie bracket} is the partially defined vector field
\begin{equation*}
[(X,\nbhd{U}),(Y,\nbhd{V})]=
([X|\nbhd{U}\cap\nbhd{V},Y|\nbhd{U}\cap\nbhd{V}],\nbhd{U}\cap\nbhd{V}).
\end{equation*}

Motivated by the constructions of Section~\ref{subsec:lie-algebras}\@, let us
say that a family $\sX=\ifam{(X_j,\nbhd{U}_j)}_{j\in J}$ of partially defined
vector fields is a \defn{Lie algebra} of partially defined vector fields if
\begin{equation*}
(X_{j_1},\nbhd{U}_{j_1})+(X_{j_2},\nbhd{U}_{j_2}),
[(X_{j_1},\nbhd{U}_{j_1}),(X_{j_2},\nbhd{U}_{j_2})]\in\sX
\end{equation*}
for every $j_1,j_2\in J$\@.  For a family $\sX=\ifam{(X_j,\nbhd{U}_j)}_{j\in
J}$ of partially defined vector fields on $\man{M}$\@, we denote by
$\sL^{(\infty)}(\sX)$ the smallest Lie algebra of partially defined vector
fields containing $\sX$\@.

Associated in a natural way to a partially defined family of vector fields is
a distribution.
\begin{definition}
Let $r\in\{\infty,\omega\}$ and let $\man{M}$ be a manifold of class
$\C^r$\@.  Given a family $\sX=\ifam{(X_j,\nbhd{U}_j)}_{j\in J}$ of partially
defined vector fields of class $\C^r$\@, we can define a distribution
$\dist{D}(\sX)$ by
\begin{equation*}\eqoprocend
\dist{D}(\sX)_x=\vecspan[\real]{X_j(x)|\enspace x\in\nbhd{U}_j}.
\end{equation*}
\end{definition}

In the definition we did not assert anything about whether the distribution
$\dist{D}(\sX)$ is of class $\C^r$ if $\sX$ is a family of partially defined
$\C^r$-vector fields.  In the smooth case, the resulting distribution
\emph{is} smooth, provided that all vector fields are bounded.
\begin{proposition}\label{prop:smooth-partial}
Let\/ $\man{M}$ be a smooth paracompact Hausdorff manifold and let\/
$\metric$ be a smooth Riemannian metric on\/ $\man{M}$\@.  If\/
$\sX=\ifam{(X_j,\nbhd{U}_j)}_{j\in J}$ is a family of smooth partially
defined vector fields\/ $\man{M}$ such that
\begin{equation*}
\sup\setdef{\dnorm{X_j(x)}_{\metric}}{x\in\nbhd{U}_j}<\infty
\end{equation*}
for each $j\in J$\@, then the distribution $\dist{D}(\sX)$ is smooth.
(Here\/ $\dnorm{\cdot}_{\metric}$ denotes the norm on the tangent spaces
induced by the Riemannian metric\/ $\metric$\@.)
\begin{proof}
For each $j\in J$\@, we use Lemma~\ref{plem:global-span2} from the proof of
Theorem~\ref{the:global-span} to assert the existence of
$f_j\in\func[\infty]{\man{M}}$ such that $f_j(x)\in\realp$ for
$x\in\nbhd{U}_j$ and $f_j(x)=0$ for $x\in\man{M}\setminus\nbhd{U}_j$\@.  We
can then define a smooth vector field $\hat{X}_j$ on $\man{M}$ by
\begin{equation*}
\hat{X}_j(x)=\begin{cases}f_j(x)X_j(x),&x\in\nbhd{U}_j,\\
0,&x\in\man{M}\setminus\nbhd{U}_j.\end{cases}
\end{equation*}
Define the family $\hat{\sX}=\ifam{\hat{X}_j}_{j\in J}$ of vector fields on
$\man{M}$\@.  It is clear that $\dist{D}(\sX)=\dist{D}(\hat{\sX})$\@, giving
the result.
\end{proof}
\end{proposition}

For families of partially defined analytic vector fields, the corresponding
assertion does not generally hold.
\begin{example}\label{eg:analytic-partial}
Note that Example~\enumdblref{eg:orbits}{enum:partially-orbit} below shows
that families of partially defined $\C^\omega$-vector fields do not define
$\C^\omega$-distributions.  Indeed, in this example, the family of analytic
partially vector fields $\sX=\ifam{(X_j,\nbhd{U}_j)}_{j\in\{1,2\}}$ defines
the distribution
\begin{equation*}
\dist{D}(\sX)_{(x_1,x_2)}=
\begin{cases}\vecspan[\real]{\pderiv{}{x_1}},&x_1\le-1,\\
\vecspan[\real]{\pderiv{}{x_2}},&x_1\ge-1,\\
\tb[(x_1,x_2)]{\real^2},&x_1\in\interval(-1,1).\end{cases}
\end{equation*}
By Proposition~\ref{prop:analytic-rank}\@, this is not an analytic
distribution.  It \emph{is}\@, however, a smooth distribution.\oprocend
\end{example}

Let $\sX$ be a family of partially defined vector fields and, following our
notation of Section~\ref{subsec:lie-algebras}\@, let us abbreviate
$\dist{D}(\sL^{(\infty)}(\sX))$ by $\dist{L}^{(\infty)}(\sX)$\@.  By
Proposition~\ref{prop:liealggen}\@, $\dist{L}^{(\infty)}(\sX)_x$ is generated
by tangent vectors of the form
\begin{equation*}
[(X_{j_k},\nbhd{U}_{j_k}),[(X_{j_{k-1}},\nbhd{U}_{j_{k-1}}),\dots,
[(X_{j_2},\nbhd{U}_{j_2}),(X_{j_1},\nbhd{U}_{j_1})]\cdots]](x),
\end{equation*}
where $k\in\integerp$ and
$(X_{j_1},\nbhd{U}_{j_1}),\dots,(X_{j_k},\nbhd{U}_{j_k})\in\sX$ are such that
$x\in\nbhd{U}_{j_1}\cap\dots\cap\nbhd{U}_{j_k}$\@.

It is also convenient to introduce some notation regarding germs of partially
defined vector fields.  Thus we let $\sX=\ifam{(X_j,\nbhd{U}_j)}_{j\in J}$ be
a family of partially defined vector fields and we denote by
\begin{equation*}
\sX_x=\setdef{r_{\nbhd{U}_j,x}(X_j)}{x\in\nbhd{U}_j}
\end{equation*}
the set of germs at $x$ of vector fields from $\sX$\@.  Note that
\begin{equation}\label{eq:partially-dist}
\dist{D}(\sX)_x=\vecspan[\real]{X(x)|\enspace[X]_x\in\sX_x}.
\end{equation}

\subsection{Orbits}\label{subsec:orbits}

Now we turn to defining orbits for partially defined families of vector
fields.  The first step is to study a ``group'' of diffeomorphisms associated
with a partially defined family of vector fields.  For a partially defined
vector field $(X,\nbhd{U})$ and for $t\in\real$\@, we recall from
Section~\ref{subsec:local-diffeos} that $\nbhd{U}(X,t)\subset\nbhd{U}$ is the
open set such that $\flow{X}{t}(x)$ is defined for each
$x\in\nbhd{U}(X,t)$\@.
\begin{definition}
Let $r\in\{\infty,\omega\}$\@, let $\man{M}$ be a $\C^r$-manifold, let $\sD$
be a family of $\C^r$-local diffeomorphisms, and let $\sX$ be a family of
partially defined vector fields of class $\C^r$\@.
\begin{compactenum}[(i)]
\item A \defn{group of $\C^r$-local diffeomorphisms} is a family
$\sG=\ifam{(\Phi_j,\nbhd{U}_j)}_{j\in J}$ of $\C^r$-local diffeomorphisms
such that
$(\Phi_{j_1},\nbhd{U}_{j_1})\scirc(\Phi_{j_2},\nbhd{U}_{j_2})\in\sG$ and
$(\Phi_j,\nbhd{U}_j)^{-1}\in\sG$ for every $j,j_1,j_2\in J$\@.
\item A group $\sG=\ifam{(\Phi_j,\nbhd{U}_j)}_{j\in J}$ of $\C^r$-local
diffeomorphisms is \defn{everywhere defined} if, for each $x\in\man{M}$\@,
there exists $j\in J$ such that $x\in\nbhd{U}_j$\@.
\item The \defn{group of local diffeomorphisms} generated by $\sD$ is the
smallest group of $\C^r$-local diffeomorphisms containing $\sD$ and which is
closed under the operations of composition and inverse of local
diffeomorphisms.
\item The \defn{group of local diffeomorphisms} generated by $\sX$ is the
group of $\C^r$-local diffeomorphisms generated by those local
diffeomorphisms of the form $(\flow{X}{t},\nbhd{U}(X,t))$ for
$(X,\nbhd{U})\in\sX$ and for $t\in\real$\@.\oprocend
\end{compactenum}
\end{definition}

Note that a group of local diffeomorphisms is not a actually a group in the
usual sense of the word since, for example, local diffeomorphisms with empty
domain do not have unique inverses.

Let us obtain a concrete description of $\Diff(\sX)$\@.  First of all, to
simplify notation, since the open set $\nbhd{U}(X,t)$ associated with the
mapping $\flow{X}{t}$ is implicit, we shall write $\flow{X}{t}$ for the local
diffeomorphism $(\flow{X}{t},\nbhd{U}(X,t))$\@.  Then, if
$\vect{X}=\ifam{X_1,\dots,X_k}$ is a family of vector fields from $\sX$ and
if $\vect{t}=(t_1,\dots,t_k)\in\real^k$\@, then we denote
\begin{equation*}
\flow{\vect{X}}{\vect{t}}=\flow{X_k}{t_k}\scirc\cdots\scirc\flow{X_1}{t_1},
\end{equation*}
which we think of as a composition of local diffeomorphisms and so a local
diffeomorphism.  For $x\in\man{M}$ we shall also write
\begin{equation*}
\flow{\vect{X}}{\vect{t}}(x)=
\flow{X_k}{t_k}\scirc\cdots\scirc\flow{X_1}{t_1}(x),
\end{equation*}
with the understanding that this is defined if $x$ is in the domain of the
local diffeomorphism $\flow{\vect{X}}{\vect{t}}$\@.  The set of such $x$'s we
denote by $\nbhd{U}(\vect{X},\vect{t})$\@, noting that this is an open subset
of $\man{M}$\@.  One can then directly verify that
\begin{equation}\label{eq:Diff(X)}
\Diff(\sX)=\asetdef{\flow{\vect{X}}{\vect{t}}}
{\vect{X}\in\sX^k,\ \vect{t}\in\real^k,\ k\in\integerp}.
\end{equation}

The preceding discussion is greatly complicated by the fact that we allow
vector fields from $\sX$ to possibly not be complete and/or not globally
defined.  If all vector fields from $\sX$ are complete and globally defined,
then one easily sees that
\begin{equation*}
\Diff(\sX)=\asetdef{\flow{X_k}{t_k}\scirc\cdots\scirc\flow{X_1}{t_1}}
{X_1,\dots,X_k\in\sX,\ t_1,\dots,t_k\in\real}.
\end{equation*}
The reader would benefit by keeping this special case in mind.

We can now define what we mean by an orbit for a family of partially defined
vector fields.
\begin{definition}
Let $r\in\{\infty,\omega\}$\@, let $\man{M}$ be a $\C^r$-manifold, and let
$\sX=\ifam{(X_j,\nbhd{U}_j)}_{j\in J}$ be an everywhere defined family of
partially defined $\C^r$-vector fields.  The \defn{orbit} of $\sX$ through
$x_0\in\man{M}$ is the set
\begin{equation*}\eqoprocend
\Orb(x_0,\sX)=\asetdef{\flow{\vect{X}}{\vect{t}}(x_0)}
{\vect{X}\in\sX^k,\ \vect{t}\in\real^k,\ k\in\integerp}.
\end{equation*}
\end{definition}

Note that two distinct orbits are disjoint.  Thus the set of orbits defines a
partition of $\man{M}$\@.

Let us understand the concept of an orbit by using some examples.  Most of
our examples involve complete vector fields, so obviating some of the
complications of the constructions above.
\begin{examples}\label{eg:orbits}
\begin{compactenum}
\item \label{enum:singular-foliation} We take $\man{M}=\real^2$ and define
\begin{equation*}
X_1=x_1\pderiv{}{x_1},\quad X_2=x_2\pderiv{}{x_2}.
\end{equation*}
The flows of $X_1$ and $X_2$ are easily determined explicitly.  Using these
flows one can readily determine the orbits for $\sX=\ifam{X_1,X_2}$\@.  Let
us illustrate how to do this in two cases; the other cases follow in the same
manner.
\begin{compactenum}[(a)]
\item \emph{$\vect{x}_0=(x_{01},x_{02})$ with $x_{01}\in\realp$ and
$x_{20}=0$\@:} Since $X_2=0$ on the $x_1$-axis and since $X_1$ is tangent to
the $x_1$-axis, $\Orb(\vect{x}_0,\sX)$ will be contained in the $x_1$-axis.
Moreover, if $x_1\in\realp$ and if we define $t_1=\log\frac{x_1}{x_{01}}$\@,
then $\flow{X_1}{t_1}(\vect{x}_0)=(x_1,0)$\@.  Moreover, for any
$t\in\realp$\@,
\begin{equation*}
\flow{X_1}{t}(\vect{x}_0)\in\setdef{(x_1,0)}{x_1\in\realp}.
\end{equation*}
Thus we must have
\begin{equation*}
\Orb(\vect{x}_0,\sX)=\setdef{(x_1,0)}{x_1\in\realp}.
\end{equation*}

\item \emph{$\vect{x}_0=(x_{01},-x_{02})$ with $x_{01},x_{02}\in\realp$\@:}
Here we let $(x_1,-x_2)\in\real^2$ with $x_1,x_2\in\realp$\@.  We then define
$t_1=\log\frac{x_1}{x_{01}}$ and $t_2=\log\frac{x_2}{x_{02}}$ and note that
$\flow{X_2}{t_2}\scirc\flow{X_1}{t_1}(\vect{x}_0)=(x_1,x_2)$\@.  Moreover,
for $t_1,\dots,t_k\in\real$ and for $j_1,\dots,j_k\in\{1,2\}$\@,
\begin{equation*}
\flow{X_{j_k}}{t_k}\scirc\cdots\scirc\flow{X_{j_1}}{t_1}(\vect{x}_0)\in
\setdef{(x_1,-x_2)}{x_1,x_2\in\realp}.
\end{equation*}
This shows that
\begin{equation*}
\Orb(\vect{x}_0,\sX)=\setdef{(x_1,-x_2)}{x_1,x_2\in\realp}.
\end{equation*}
\end{compactenum}
In any case, it is easy to see that there are nine distinct orbits for the
family of vector fields $\sX=\ifam{X_1,X_2}$\@, and these are determined to
be
\begin{align*}
\Orb_1((0,0),\sX)=&\;\{(0,0)\},\\
\Orb_2((1,0),\sX)=&\;\setdef{(x_1,0)}{x_1\in\realp},\\
\Orb_3((-1,0),\sX)=&\;\setdef{(-x_1,0)}{x_1\in\realp},\\
\Orb_4((0,1),\sX)=&\;\setdef{(0,x_2)}{x_2\in\realp},\\
\Orb_5((0,-1),\sX)=&\;\setdef{(0,-x_2)}{x_2\in\realp},\\
\Orb_6((1,1),\sX)=&\;\setdef{(x_1,x_2)}{x_1,x_2\in\realp},\\
\Orb_7((-1,1),\sX)=&\;\setdef{(-x_1,x_2)}{x_1,x_2\in\realp},\\
\Orb_8((1,-1),\sX)=&\;\setdef{(x_1,-x_2)}{x_1,x_2\in\realp},\\
\Orb_9((-1,-1),\sX)=&\;\setdef{(-x_1,-x_2)}{x_1,x_2\in\realp}.
\end{align*}
We depict these orbits in Figure~\ref{fig:R2orbits}\@.
\begin{figure}[htbp]
\centering
\includegraphics{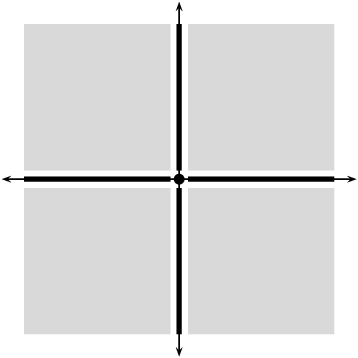}
\caption{Orbits}\label{fig:R2orbits}
\end{figure}%

\item \label{enum:!analytic-orbit} We let $\man{M}=\real^2$ and define
\begin{equation*}
X_1=\pderiv{}{x_1},\quad X_2=f(x_1)\pderiv{}{x_2},
\end{equation*}
where
\begin{equation*}
f(x)=\begin{cases}\eul^{-1/x^2},&x\in\realp,\\0,&x\in\realnp.\end{cases}
\end{equation*}
We take $\sX=\ifam{X_1,X_2}$ and claim that $\Orb(\vect{x},\sX)=\real^2$ for
every $\vect{x}\in\real^2$\@.  It suffices to show that, for example,
$\Orb(\vect{0},\sX)=\real^2$\@.  For $(x_1,x_2)\in\real^2$ with $x_1>0$ we
define $t_1=x_1$ and $t_2=\frac{x_2}{f(x_1)}$\@, and directly compute
\begin{equation*}
\flow{X_2}{t_2}\scirc\flow{X_1}{t_1}(\vect{0})=(x_1,x_2).
\end{equation*}
If $(x_1,x_2)\in\real^2$ with $x_1\le0$ we define $t_1=1$\@,
$t_2=\frac{x_2}{f(x_1)}$\@, and $t_3=-1+x_1$\@, and directly compute
\begin{equation*}
\flow{X_3}{t_3}\scirc\flow{X_2}{t_2}\scirc
\flow{X_1}{t_1}(\vect{0})=(x_1,x_2).
\end{equation*}
This shows that $\Orb(\vect{0},\sX)=\real^2$ as desired.

\item \label{enum:linear-orbit} For $u\in\real$ define a vector field $X_u$
on $\real^2$ by
\begin{equation*}
X_u(\vect{x})=(\vect{x},\mat{A}\vect{x}+\vect{b}u),
\end{equation*}
where
\begin{equation*}
\mat{A}=\begin{bmatrix}0&1\\0&0\end{bmatrix},\quad
\vect{b}=\begin{bmatrix}0\\1\end{bmatrix}.
\end{equation*}
We then consider the family of vector fields on $\real^2$ given by
$\sX=\ifam{X_u}_{u\in\real}$\@.  Let
$(x_{11},x_{12}),(x_{21},x_{22})\in\real^2$ and let
$T\in\real\setminus\{0\}$\@.  Define
\begin{equation*}
\begin{gathered}
u_1=\frac{-3Tx_{12}-Tx_{22}-4 x_{11}+4x_{21}}{T^2},\\
u_2=\frac{Tx_{12}+3 Tx_{22}+4 x_{11}-4x_{21}}{T^2}.
\end{gathered}
\end{equation*}
One can then directly verify (by solving the linear differential equations
defining the flow) that
\begin{equation}\label{eq:linear-orbitu}
(x_{21},x_{22})=\flow{X_{u_2}}{T/2}\scirc\flow{X_{u_1}}{T/2}(x_{11},x_{12}),
\end{equation}
showing that $\Orb(\vect{x},\sX)=\real^2$ for every
$\vect{x}\in\real^2$\@.
\item \label{enum:partially-orbit} The final example we consider is one where
the vector fields are partially defined, but not globally defined.  We take
$\man{M}=\real^2$ and define a family
$\sX=\ifam{(X_j,\nbhd{U}_j)}_{j\in\{1,2\}}$ of partially defined vector
fields by
\begin{equation*}
\nbhd{U}_1=\setdef{(x_1,x_2)}{x_1<1},\quad\nbhd{U}_2=\setdef{(x_1,x_2)}{x_1>-1}
\end{equation*}
and
\begin{equation*}
X_1=\pderiv{}{x_1},\quad X_2=\pderiv{}{x_2}.
\end{equation*}
It is easy to see that
\begin{equation*}
\Orb((x_{01},x_{02}),\sX)=
\begin{cases}\setdef{(x_1,x_2)}{x_1<1},&x_{01}<1,\\
\setdef{(x_{01},x_2)}{x_2\in\real},&x_{01}\ge1.\end{cases}
\end{equation*}
Note that the orbits are analytic submanifolds.\oprocend
\end{compactenum}
\end{examples}

\subsection{Fixed-time orbits}\label{subsec:Torbits}

In this section we consider a modification of the notion of an orbit as
defined in the previous section.  Let $r\in\{\infty,\omega\}$\@, let
$\man{M}$ be a $\C^r$-manifold, and let $\sX$ be a family of partially vector
fields of class $\C^r$\@.  Above we defined $\Diff(\sX)$ as the subgroup of
local diffeomorphisms defined by flows of vector fields from $\sX$\@.  In
this section we modify this construction slightly to give the orbit
corresponding to flows whose ``total time'' is fixed.  To make this
construction, we first consider flows whose ``total time'' is zero.

For convenience and to reestablish notation, we recall the explicit
characterisation of $\Diff(\sX)$ from above.  As above, for a vector field
$X$ we shall often denote the local diffeomorphism
$(\flow{X}{t},\nbhd{U}(X,t))$ simply by $\flow{X}{t}$\@.  Let
$\vect{X}=\ifam{X_1,\dots,X_k}$ be a finite family of vector fields from
$\sX$ and let $\vect{t}=\ifam{t_1,\dots,t_k}\in\real^k$\@.  Then we define a
local diffeomorphism
\begin{equation*}
\flow{\vect{X}}{\vect{t}}=\flow{X_k}{t_k}\scirc\cdots\scirc\flow{X_1}{t_1},
\end{equation*}
understanding implicitly that this is only interesting when the resulting
composition has nonempty domain.  With this notation,
\begin{equation*}
\Diff(\sX)=\asetdef{\flow{\vect{X}}{\vect{t}}}
{\vect{X}\in\sX^k,\ \vect{t}\in\real^k,\ k\in\integerp}.
\end{equation*}
Now let $T\in\real$\@.  Define
\begin{equation*}
\Diff_T(\sX)=\asetdef{\flow{\vect{X}}{\vect{t}}}
{\vect{X}\in\sX^k,\ \vect{t}\in\real^k,\ \sum_{j=1}^kt_k=T,\ k\in\integerp}.
\end{equation*}
The case where $T=0$ is particularly interesting, as we shall see below.  The
following properties of $\Diff_0(\sX)$ are useful in understanding some of
the subsequent constructions.
\begin{proposition}
Let\/ $r\in\{\infty,\omega\}$\@, let\/ $\man{M}$ be a\/ $\C^r$-manifold, and
let\/ $\sX\subset\sections[r]{\tb{\man{M}}}$\@.  Then the following
statements hold:
\begin{compactenum}[(i)]
\item \label{pl:Diff01} $\Diff_0(\sX)$ is a subgroup of the group\/
$\Diff(\sX)$ of local diffeomorphisms; that is, if\/
$(\Phi,\nbhd{U}),(\Psi,\nbhd{V})\in\Diff_0(\sX)$\@, then\/
$(\Phi,\nbhd{U})^{-1}\in\Diff_0(\sX)$ and\/
$(\Phi,\nbhd{U})\scirc(\Psi,\nbhd{V})\in\Diff_0(\sX)$\@;
\item \label{pl:Diff02} $\Diff_0(\sX)$ is a normal subgroup; that is, if
$(\Phi,\nbhd{U})\in\Diff_0(\sX)$ and if $(\Psi,\nbhd{V})\in\Diff(\sX)$\@,
then $(\Psi,\nbhd{V})\scirc(\Phi,\nbhd{U})\scirc(\Psi,\nbhd{V})^{-1}\in
\Diff_0(\sX)$\@.
\end{compactenum}
\begin{proof}
\eqref{pl:Diff01} Let $\vect{X}=\ifam{X_1,\dots,X_k}$ and
$\vect{Y}=\ifam{Y_1,\dots,Y_m}$ be families of vector fields and
$\vect{t}=\ifam{t_1,\dots,t_k}$ and $\vect{s}=\ifam{s_1,\dots,s_m}$ be
families of real numbers such that $\Phi=\flow{\vect{X}}{\vect{t}}$ and
$\Psi=\flow{\vect{Y}}{\vect{s}}$\@.  Thus
\begin{equation*}
\sum_{j=1}^kt_j=\sum_{l=1}^ms_l=0.
\end{equation*}
Then $(\Phi,\nbhd{U})^{-1}$ is defined by
\begin{equation*}
\Phi=\flow{X_1}{-t_1}\scirc\cdots\scirc\flow{X_k}{-t_k},
\end{equation*}
and so $(\Phi,\nbhd{U})^{-1}\in\Diff_0(\sX)$\@.  Similarly,
$(\Phi,\nbhd{U})\scirc(\Psi,\nbhd{V})$ is defined by
\begin{equation*}
\flow{X_k}{t_k}\scirc\cdots\scirc\flow{X_1}{t_1}\scirc
\flow{Y_m}{t_m}\scirc\cdots\flow{Y_1}{t_1},
\end{equation*}
and so $(\Phi,\nbhd{U})\scirc(\Psi,\nbhd{V})\in\Diff_0(\sX)$\@.

\eqref{pl:Diff02} Let $\vect{X}=\ifam{X_1,\dots,X_k}$ and
$\vect{Y}=\ifam{Y_1,\dots,Y_m}$ be families of vector fields and
$\vect{t}=\ifam{t_1,\dots,t_k}$ and $\vect{s}=\ifam{s_1,\dots,s_m}$ be
families of real numbers such that $\Phi=\flow{\vect{X}}{\vect{t}}$ and
$\Psi=\flow{\vect{Y}}{\vect{s}}$\@.  Thus
\begin{equation*}
\sum_{j=1}^kt_j=0.
\end{equation*}
Note that $(\Psi,\nbhd{V})\scirc(\Phi,\nbhd{U})\scirc(\Psi,\nbhd{V})^{-1}$ is
defined by
\begin{equation*}
\flow{Y_m}{s_m}\scirc\cdots\scirc\flow{Y_1}{s_1}\scirc
\flow{X_k}{t_k}\scirc\cdots\scirc\flow{X_1}{t_1}\scirc
\flow{Y_1}{-s_1}\scirc\cdots\scirc\flow{Y_m}{-s_m},
\end{equation*}
and so $(\Psi,\nbhd{V})\scirc(\Phi,\nbhd{U})\scirc(\Psi,\nbhd{V})^{-1}\in
\Diff_0(\sX)$\@, as desired.
\end{proof}
\end{proposition}

We can now make the following definition.
\begin{definition}
Let $r\in\{\infty,\omega\}$\@, let $\man{M}$ be a $\C^r$-manifold, let
$\sX\subset\sections[r]{\tb{\man{M}}}$ be a family of $\C^r$-vector fields,
and let $T\in\real$\@.  The \defn{$T$-orbit} of $\sX$ through $x_0\in\man{M}$
is the set
\begin{equation*}
\Orb_T(x_0,\sX)=\asetdef{\flow{\vect{X}}{\vect{t}}(x_0)}
{\vect{X}\in\sX^k,\ \vect{t}\in\real^k,\ \sum_{j=1}^kt_j=T,\ k\in\integerp}.
\end{equation*}
A \defn{fixed-time orbit} of $\sX$ through $x_0\in\man{M}$ is a set of the
form $\Orb_T(x_0,\sX)$ for some $T\in\real$\@.\oprocend
\end{definition}

Let us give some examples of fixed-time orbits.
\begin{examples}
We resume some of the examples whose orbits we studied in Example~\ref{eg:orbits}\@.
\begin{compactenum}
\item \label{enum:singular-fixed} Let $\man{M}=\real^2$ and take
$X_1=x_1\pderiv{}{x_1}$ and $X_2=x_2\pderiv{}{x_2}$\@.  We determined the
orbits of $\sX=\ifam{X_1,X_2}$ in
Example~\enumdblref{eg:orbits}{enum:singular-foliation}\@.  One readily
computes
\begin{equation*}
\Orb_T((x_1,0),\sX)=\{(x_1\eul^T,0)\}
\end{equation*}
and
\begin{equation*}
\Orb_T((0,x_2),\sX)=\{(0,x_2\eul^T)\}.
\end{equation*}
Thus the $T$-orbits for points on the coordinate axes are singletons.  With a
little more effort, or better, an application of
Theorem~\ref{the:fingen-orbit-theorem} below, we can determine the $T$-orbits
for the other points.  The $T$-orbits are the submanifolds
\begin{equation*}
\man{S}_c^+=\setdef{(x,cx^{-1})}{x\in\realp},\quad
\man{S}_c^-=\setdef{(x,cx^{-1})}{x\in\realp},
\end{equation*}
parameterised by the nonzero real number $c$\@.  In
Figure~\ref{fig:R2Torbits}
\begin{figure}[htbp]
\centering
\includegraphics[width=0.4\hsize]{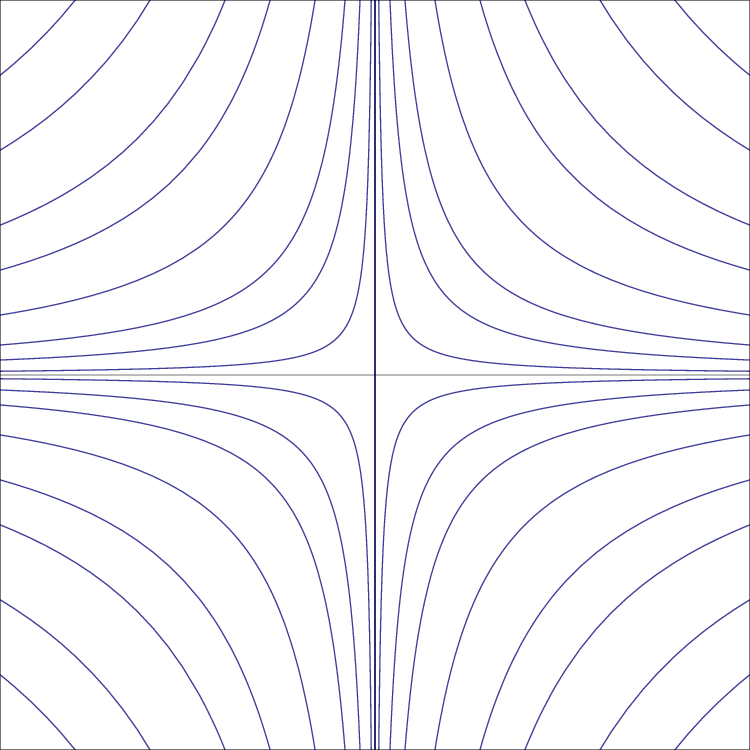}
\caption{Fixed-time orbits}\label{fig:R2Torbits}
\end{figure}%
we show the collection of these fixed-time orbits.  To determine which of the
submanifolds is a $T$-orbit through a point $(x_1,x_2)$ that is not on a
coordinate axis, one follows the integral curve of either $X_1$ or $X_2$ (it
matters not which) for time $T$ and the submanifold upon which one finds
oneself is the $T$-orbit for $(x_1,x_2)$\@.

\item We consider here the example first considered in
Example~\enumdblref{eg:orbits}{enum:linear-orbit}\@.  Thus we consider the
family $\ifam{X_u}_{u\in\real}$ of vector fields with
\begin{equation*}
X_u(\vect{x})=(\vect{x},\mat{A}\vect{x}+\vect{b}u),
\end{equation*}
where
\begin{equation*}
\mat{A}=\begin{bmatrix}0&1\\0&0\end{bmatrix},\quad
\vect{b}=\begin{bmatrix}0\\1\end{bmatrix}.
\end{equation*}
We saw in~\eqref{eq:linear-orbitu} that $\Orb_T(\vect{x},\sX)=\real^2$ for
every $\vect{x}\in\real^2$ and $T\in\real\setminus\{0\}$\@.
\item Let us now consider the partially defined vector fields
$\sX=\ifam{(X_1,\nbhd{U}_1),(X_2,\nbhd{U}_2)}$ from
Example~\enumdblref{eg:orbits}{enum:partially-orbit}\@.  As with
Example~\ref{enum:singular-fixed} above and by
Theorem~\ref{the:fingen-orbit-theorem}\@, it suffices to determine the
$0$-orbit through each point.  The $T$-orbit for a point is then determined
by following a curve that is a concatenation of integral curves for $X_1$ and
$X_2$ for a total time $T$\@.  The $T$-orbit with be the $0$-orbit through
the resulting point.  A few moments of thought will lead one to conclude that
\begin{equation*}
\Orb_0((x_1,x_2),\sX)=\begin{cases}
\setdef{(x_1+s,x_2-s)}{s\in\interval({-\infty},{1-x_1})},&x_1<1,\\
\{(x_1,x_2)\},&x_1\ge1.\end{cases}
\end{equation*}
The $0$-orbits are depicted in Figure~\ref{fig:partially-orbits}\@.
\begin{figure}[htbp]
\centering
\includegraphics{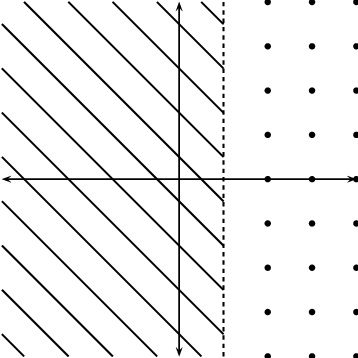}
\caption{The $0$-orbits for a family of partially defined vector
fields}\label{fig:partially-orbits}
\end{figure}%
Note that the dimensions of the $0$-orbits differ on disjoint open
sets.\oprocend
\end{compactenum}
\end{examples}

The following result gives a useful characterisation of $\Diff_T(\sX)$\@.
\begin{proposition}\label{prop:OrbTchar}
Let\/ $r\in\{\infty,\omega\}$\@, let\/ $\man{M}$ be a\/ $\C^r$-manifold, and
let\/ $\sX$ be a family of partially defined vector fields of class\/
$\C^r$\@.  For\/ $x_0\in\man{M}$ and\/ $T\in\real$ suppose that\/
$\Orb_T(x_0,\sX)\not=\emptyset$\@.  Then, for every\/
$x\in\Orb_T(x_0,\sX)$\@, we have
\begin{equation*}
\Orb_T(x_0,\sX)=\setdef{\Phi(x)}{\Phi\in\Diff_0(\sX)}.
\end{equation*}
\begin{proof}
Write $x=\flow{\vect{X}}{\vect{t}}(x_0)$ for $\vect{X}\in\sX^k$ and
$\vect{t}\in\real^k$ with $\sum_{j=1}^kt_j=T$\@.  If $\Phi\in\Diff_0(\sX)$
then we have $\Phi=\flow{\vect{Y}}{\vect{s}}$ for $\vect{Y}\in\sX^m$ and
$\vect{s}\in\real^m$ with $\sum_{l=1}^ms_l=0$\@.  Then,
\begin{equation*}
\Phi(x)=\flow{\vect{Y}}{\vect{s}}\scirc\flow{\vect{X}}{\vect{t}}(x_0),
\end{equation*}
and so it is evident that $\Phi(x)\in\Orb_T(x_0,\sX)$\@.  Thus
\begin{equation*}
\setdef{\Phi(x)}{\Phi\in\Diff_0(\sX)}\subset\Orb_T(x_0,\sX).
\end{equation*}

Conversely, let $y\in\Orb_T(x_0,\sX)$ and write
$y=\flow{\vect{Y}}{\vect{s}}(x_0)$ for $\vect{Y}\in\sX^m$ and
$\vect{s}\in\real^m$ with $\sum_{l=1}^ms_l=T$\@.  If we write
$\hat{\vect{t}}=(t_k,\dots,t_1)$\@, then
\begin{equation*}
y=\flow{\vect{Y}}{\vect{s}}(x_0)=\flow{\vect{Y}}{\vect{s}}
\scirc\flow{\vect{X}}{-\hat{\vect{t}}}(x),
\end{equation*}
from which we deduce that
\begin{equation*}
\Orb_T(x_0,\sX)\subset\setdef{\Phi(x)}{\Phi\in\Diff_0(\sX)},
\end{equation*}
as desired.
\end{proof}
\end{proposition}

The following corollary gives a particular simple form for the fixed-time
orbit in some cases.
\begin{corollary}
Let\/ $r\in\{\infty,\omega\}$\@, let\/ $\man{M}$ be a\/ $\C^r$-manifold, and
let\/ $\sX$ be a family of partially defined vector fields of class\/
$\C^r$\@.  Let\/ $X\in\sX$\@,\/ $x_0\in\man{M}$\@, and\/ $T\in\realp$
satisfy\/ $x\in\nbhd{U}(X,T)$\@.  Then
\begin{equation*}
\Orb_T(x_0,\sX)=\setdef{\Phi\scirc\flow{X}{T}(x_0)}{\Phi\in\Diff_0(\sX)}.
\end{equation*}
\end{corollary}

\subsection{The Orbit Theorem}

Before we state the Orbit Theorem we need some terminology and notation.

Let $r\in\{\infty,\omega\}$ and let $\man{M}$ be a $\C^r$-manifold.  For a
$\C^r$-local diffeomorphism $(\Phi,\nbhd{U})$ and for a partially defined
vector field $(X,\nbhd{V})$\@, denote by $\Phi_*X$ the partially defined
vector field with domain $\Phi(\nbhd{U}\cap\nbhd{V})$ and defined by
$\Phi_*X(x)=\tf[\Phi^{-1}(x)]{\Phi}\scirc X\scirc\Phi^{-1}(x)$\@.

With this terminology, we are ready to state the Orbit Theorem.
\begin{theorem}\label{the:orbit-theorem}
Let\/ $r\in\{\infty,\omega\}$\@, let\/ $\man{M}$ be a\/ $\C^r$-manifold, and
let\/ $\sX$ be an everywhere defined family of partially defined vector
fields of class\/ $\C^r$\@.  Then, for each\/ $x_0\in\man{M}$\@,
\begin{compactenum}[(i)]
\item \label{pl:orbit-theorem1} $\Orb(x_0,\sX)$ is a connected immersed\/
$\C^r$-submanifold of\/ $\man{M}$ and
\item \label{pl:orbit-theorem2} for each\/ $x\in\Orb(x_0,\sX)$\@,
\begin{equation*}
\tb[x]{\Orb(x_0,\sX)}=\vecspan[\real]{\Phi_*X(x)\mid\enspace
\Phi\in\Diff(\sX),\ X\in\sX}.
\end{equation*}
\end{compactenum}
Moreover,\/ $\man{M}$ is the disjoint union of the set of orbits.
\begin{proof}
Let us denote a family of partially defined vector fields by
\begin{equation*}
\sO(\sX)=\setdef{\Phi_*X}{\Phi\in\Diff(\sX),\ X\in\sX}.
\end{equation*}
We also define a distribution $\dist{O}$ by
\begin{equation*}
\dist{O}_x=\vecspan[\real]{\Phi_*X(x)|\enspace\Phi\in\Diff(\sX),\ X\in\sX}.
\end{equation*}
The following lemma gives a useful property of these subspaces.
\begin{prooflemma}\label{plem:Odistdim}
For\/ $x_0\in\man{M}$ and for\/ $x\in\Orb(x_0,\sX)$\@,\/
$\dim(\dist{O}_x)=\dim(\dist{O}_{x_0})$\@.
\begin{subproof}
If $x\in\Orb(x_0,\sX)$ then there exists $\Psi\in\Diff(\sX)$ such that
$x=\Psi(x_0)$\@.  Let $\Phi\in\Diff(\sX)$ be such that $x_0$ is in the image
of $\Phi$ and let $X\in\sX$ so that $\Phi_*X(x_0)\in\dist{O}_{x_0}$\@.  Then
\begin{align*}
\tf[x_0]{\Psi}(\Phi_*X(x_0))=&\;
\tf[x_0]{\Psi}\scirc\tf[\Phi^{-1}(x_0)]{\Phi}\scirc X\scirc\Phi^{-1}(x_0)\\
=&\;\tf[x_0]{\Psi}\scirc\tf[\Phi^{-1}(x_0)]{\Phi}\scirc X\scirc\Phi^{-1}
\scirc\Psi^{-1}\scirc\Psi(x_0)\\
=&\;\tf[(\Psi\scirc\Phi)^{-1}(x)]{(\Psi\scirc\Phi)}\scirc X\scirc
(\Psi\scirc\Phi)^{-1}(x)\\
=&\;(\Psi\scirc\Phi)_*X(x)\in\dist{O}_x
\end{align*}
since $\Psi\scirc\Phi\in\Diff(\sX)$\@.  Since $\tf[x_0]{\Psi}$ is an
isomorphism, this shows that $\dim(\dist{O}_x)\ge\dim(\dist{O}_{x_0})$\@.
This argument can be reversed to give the opposite inequality.
\end{subproof}
\end{prooflemma}

For $x\in\man{M}$ suppose that $\dim(\dist{O}_x)=m_x$\@.  Let
$\sY_x=\ifam{Y_1,\dots,Y_{m_x}}\subset\sO(\sX)$ be such that
$\ifam{Y_1(x),\dots,Y_{m_x}(x)}$ is a basis for $\dist{O}_x$\@.  Define
$\map{\phi_x}{\real^{m_x}}{\man{M}}$ by
\begin{equation*}
\phi_x(t_1,\dots,t_{m_x})=\flow{Y_{m_x}}{t_{m_x}}\scirc\cdots\scirc
\flow{Y_1}{t_1}(x).
\end{equation*}
Since $\pderiv{\phi_x}{t_j}$ is equal to $Y_j(x)$ when evaluated to
$t_1=\dots=t_{m_x}=0$\@, it follows that $\phi_x$ is an embedding in a
neighbourhood of $\vect{0}\in\real^{m_x}$\@; let us denote the image of this
neighbourhood under $\phi_x$ by $\nbhd{U}(\sY_x)$\@.  Thus $\nbhd{U}(\sY_x)$
is a $\C^r$-submanifold.
\begin{prooflemma}\label{plem:orbit-basis}
$\nbhd{U}(\sY_x)\subset\Orb(x,\sX)$\@.
\begin{subproof}
By definition we can write $Y_j=\Phi_{j*}X_j$ for $\Phi_j\in\Diff(\sX)$ and
$X_j\in\sX$\@, $j\in\{1,\dots,m_x\}$\@.  We recall from
\cite[Proposition~4.2.4]{RA/JEM/TSR:88} the formula
$\flow{Y}{t}=\Phi\scirc\flow{X}{t}\scirc\Phi^{-1}$ which holds if $Y=\Phi_*X$
for vector fields $X$ and $Y$ and a diffeomorphism $\Phi$\@.  Using this
formula we have, for $\vect{t}\in\phi_x^{-1}(\nbhd{U}(\sY_x))$\@,
\begin{equation*}
\phi_x(\vect{t})=\flow{Y_{m_x}}{t_{m_x}}\scirc\cdots\scirc\flow{Y_1}{t_1}(x)=
\Phi_{m_x}\scirc\flow{X_{m_x}}{t_{m_x}}\scirc
\Phi_{m_x}^{-1}\scirc\cdots\scirc
\Phi_1\scirc\flow{X_1}{t_1}\scirc\Phi_1^{-1}(x),
\end{equation*}
showing that $\nbhd{U}(\sY_x)\subset\Orb(x,\sX)$\@, as desired, since all
mappings in the above composition are in $\Diff(\sX)$\@.
\end{subproof}
\end{prooflemma}

\begin{prooflemma}\label{plem:orbit-tangent}
$\tb[y]{\nbhd{U}(\sY_x)}=\dist{O}_y$ for every\/ $y\in\nbhd{U}(\sY_x)$\@.
\begin{subproof}
Let $j\in\{1,\dots,m_x\}$ and let
\begin{equation*}
\Phi_{\vect{t}}=\flow{Y_{m_x}}{t_{m_x}}\scirc\cdots\scirc
\flow{Y_{j+1}}{t_{j+1}}.
\end{equation*}
Then
\begin{align*}
\pderiv{\phi_x}{t_j}=&\;\pderiv{}{t_j}\flow{Y_{m_x}}{t_{m_x}}\scirc\cdots
\scirc\flow{Y_1}{t_1}(x)\\
=&\;\pderiv{}{t_j}\Phi_{\vect{t}}\scirc\flow{Y_j}{t_j}\scirc
\Phi_{\vect{t}}^{-1}\scirc\flow{Y_{m_x}}{t_{m_x}}\scirc\cdots\scirc
\flow{Y_1}{t_1}(x)\\
=&\;\tf{\Phi_{\vect{t}}}\scirc Y_j\scirc\Phi_{\vect{t}}^{-1}\scirc
\phi_x(\vect{t})=\Phi_{\vect{t}*}Y_j(\phi_x(\vect{t}))\in
\dist{O}_{\phi_x(\vect{t})}.
\end{align*}
Since this holds for every $j\in\{1,\dots,m_x\}$\@, it follows that
$\image(\tf[\vect{t}]{\phi_x})\subset\dist{O}_{\phi_x(\vect{t})}$\@.  By
Lemma~\ref{plem:Odistdim} we have
\begin{equation*}
\dim(\dist{O}_{\phi_x(\vect{t})})=\dim(\dist{O}_x)=
\dim(\image(\tf[\vect{t}]{\phi_x})).
\end{equation*}
Therefore,
\begin{equation*}
\dist{O}_{\phi_x(\vect{t})}=\image(\tf[\vect{t}]{\phi_x})=
\tb[\phi_x(\vect{t})]{\nbhd{U}(\sY_x)},
\end{equation*}
as desired.
\end{subproof}
\end{prooflemma}

\begin{prooflemma}
The subsets\/ $\ifam{\nbhd{U}(\sY_x)}_{x\in\man{M}}$ for\/ $\sY_x$ as defined
above, form a basis for a topology on\/ $\man{M}$\@.
\begin{subproof}
By,~\eg~Theorem~5.3 of \cite{SW:70}\@, it suffices to show that for any pair
$\nbhd{U}(\sY_{x_1})$ and $\nbhd{U}(\sY_{x_2})$ of such subsets there exists
$\nbhd{U}(\sY_x)$ such that
\begin{equation*}
\nbhd{U}(\sY_x)\subset\nbhd{U}(\sY_{x_1})\cap\nbhd{U}(\sY_{x_2}).
\end{equation*}
Let $x\in\nbhd{U}(\sY_{x_1})\cap\nbhd{U}(\sY_{x_2})$ and let
$Y_1,\dots,Y_{m_x}\in\sO(\sX)$ be such that
\begin{equation*}
\dist{O}_x=\vecspan[\real]{Y_1(x),\dots,Y_{m_x}(x)}.
\end{equation*}
Let $\map{\phi_x}{\real^{m_x}}{\man{M}}$ be the map defined above.  By
Lemma~\ref{plem:orbit-tangent} it follows that
\begin{equation*}
Y_1(y),\dots,Y_{m_x}(y)\in\tb[y]{\nbhd{U}(\sY_{x_1})},\quad
Y_1(y),\dots,Y_{m_x}(y)\in\tb[y]{\nbhd{U}(\sY_{x_2})}
\end{equation*}
for every $y$ in a neighbourhood of $x$\@.  Therefore, the integral curves of
the vector fields $Y_1,\dots,Y_{m_x}$ with initial conditions in
$\nbhd{U}(\sY_{x_1})$ (\resp~$\nbhd{U}(\sY_{x_2})$) nearby $x$ remain in
$\nbhd{U}(\sY_{x_1})$ (\resp~$\nbhd{U}(\sY_{x_2})$).  Therefore,
concatenations of these integral curves nearby $x$ will also remain in
$\nbhd{U}(\sY_{x_1})$ (\resp~$\nbhd{U}(\sY_{x_2})$).  In short, for
$t_1,\dots,t_{m_x}$ sufficiently near zero,
\begin{equation*}
\flow{Y_{m_k}}{t_{m_k}}\scirc\cdots\scirc\flow{Y_1}{t_1}(x)
\in\nbhd{U}(\sY_{x_1})\quad (\textrm{\resp}\
\flow{Y_{m_k}}{t_{m_k}}\scirc\cdots\scirc\flow{Y_1}{t_1}(x)
\in\nbhd{U}(\sY_{x_2})).
\end{equation*}
Thus, by restricting $\phi_x$ to a small enough neighbourhood $\nbhd{N}$ of
$\vect{0}$\@, if we define $\nbhd{U}(\sY_x)=\phi_x(\nbhd{N})$\@, we have
\begin{equation*}\eqsubqed
\nbhd{U}(\sY_x)\subset\nbhd{U}(\sY_{x_1})\cap\nbhd{U}(\sY_{x_2}).
\end{equation*}
\end{subproof}
\end{prooflemma}

Let us call the topology on $\man{M}$ generated by the sets
$\ifam{\nbhd{U}(\sY_x)}_{x\in\man{M}}$ the \defn{orbit topology}\@.
\begin{prooflemma}
In the orbit topology, the orbits are connected, open, and closed.
\begin{subproof}
Let $X\in\sX$\@.  Since the integral curve $t\mapsto\flow{X}{t}(x)$ is a
continuous curve of $\man{M}$ and since it is tangent to $\nbhd{U}(\sY_x)$\@,
it follows that the curve is continuous in the relative topology on
$\nbhd{U}(\sY_x)$\@.  Since $\nbhd{U}(\sY_x)$ is a submanifold, the relative
topology is the same as the topology induced by the immersion $\phi_x$\@.
Thus the integral curve $t\mapsto\flow{X}{t}(x)$ is continuous in the orbit
topology.  Therefore, by definition of orbits, orbits are path connected and
so connected.

If $x\in\Orb(x_0,\sX)$ then every set $\nbhd{U}(\sY_x)$ is a subset of
$\Orb(x_0,\sX)$\@.  Since $\nbhd{U}(\sY_x)$ is open, it follows that
$\Orb(x_0,\sX)$ is open.

Note that $\man{M}$ is a disjoint union of orbits.  Therefore, the complement
of an orbit is a union of orbits.  Thus the complement of an orbit is a union
of open sets and so open.  Thus an orbit is closed.
\end{subproof}
\end{prooflemma}

This last lemma shows that the orbits are the connected components of
$\man{M}$ in the orbit topology.  This, in particular, gives the final
assertion of the theorem.
\begin{prooflemma}
For\/ $x_0\in\man{M}$\@,
\begin{equation*}
\Orb(x_0,\sX)=\bigcup_{x\in\Orb(x_0,\sX)}\nbhd{U}(\sY_x),
\end{equation*}
the union being over the neighbourhoods\/ $\nbhd{U}(\sY_x)$ constructed
above.
\begin{subproof}
It is trivial that
$\Orb(x_0,\sX)\subset\bigcup_{x\in\Orb(x_0,\sX)}\nbhd{U}(\sY_x)$\@.
The converse inclusion follows from Lemma~\ref{plem:orbit-basis}\@.
\end{subproof}
\end{prooflemma}

Since each of the sets $\nbhd{U}(\sY_x)$ is diffeomorphic to an open subset
of $\real^{m_x}$ by $\phi_x^{-1}$\@, it follows that
$(\nbhd{U}(\sY_x),\phi_x^{-1})$ is a chart for $\Orb(x_0,\sX)$ for every
$x\in\Orb(x_0,\sX)$\@.  The overlap map between intersection charts
$(\nbhd{U}(\sY_{x_1}),\phi_{x_1}^{-1})$ and
$(\nbhd{U}(\sY_{x_2}),\phi_{x_2}^{-1})$ are obtained by concatenations of
flows of vector fields from $\sX$\@, and so are diffeomorphisms.  This shows
that $\Orb(x_0,\sX)$ is an immersed submanifold as in
\eqref{pl:orbit-theorem1}\@.  Assertion~\eqref{pl:orbit-theorem2} follows
from the Lemma~\ref{plem:orbit-tangent} and the definition of the
differentiable structure on the orbits.
\end{proof}
\end{theorem}

\begin{remark}
In the proof of the Orbit Theorem we prescribed a topology on $\man{M}$ that
we will, on occasion, make reference to.  This topology is called the
\defn{orbit topology}\@, and is defined as follows, using the notation from
the proof of the Orbit Theorem.  For a family $\sX$ of partially defined
$\C^r$-vector fields, $r\in\{\infty,\omega\}$\@, and for $x\in\man{M}$\@, let
$\sY_x=\ifam{Y_1,\dots,Y_{m_x}}$ be vector fields from the family
\begin{equation*}
\sO(\sX)=\setdef{\Phi_*X}{\Phi\in\Diff(\sX),\ X\in\sX}
\end{equation*}
for which $\ifam{Y_1(x),\dots,Y_{m_x}(x)}$ are a basis for
$\tb[x]{\Orb(x,\sX)}$\@.  Define $\map{\phi_x}{\real^{m_x}}{\man{M}}$ by
\begin{equation*}
\phi_x(t_1,\dots,t_{m_x})=\flow{Y_{m_x}}{t_{m_x}}\scirc\cdots\scirc
\flow{Y_1}{t_1}(x),
\end{equation*}
and note that this map, restricted to a neighbourhood of
$\vect{0}\in\real^{m_x}$\@, is an embedding.  The image of this neighbourhood
under $\phi_x$ is denoted by $\nbhd{U}(\sY_x)$\@.  The sets $\nbhd{U}(\sY_x)$
form a basis for a topology, and this topology is the orbit topology.  In the
proof of the Orbit Theorem it was shown that the orbits are the connected
components of $\man{M}$ in the orbit topology.

Sometimes the following characterisation of the orbit topology is useful.
Let $x\in\man{M}$\@, let $k\in\integerp$\@, let
$\vect{X}=\ifam{X_1,\dots,X_k}\subset\sX$\@, and let $\nbhd{U}\subset\real^k$
be a neighbourhood of $\vect{0}$ be such that the map
\begin{equation*}
\nbhd{U}\ni\vect{t}\mapsto\flow{\vect{X}}{\vect{t}}(x)\in\man{M}
\end{equation*}
is defined.  One may then define the orbit topology as the final topology
induced by the above family of mappings,~\ie~the finest topology for which
all of these maps are continuous.\oprocend
\end{remark}

It is interesting to consider whether there are stronger conclusions that can
be drawn from the Orbit Theorem when the vector fields are globally defined.
In the smooth case when the partially defined vector fields are bounded, by
Proposition~\ref{prop:smooth-partial}\@, there is no extra structure present
when the vector fields from $\sX$ are globally defined.  In the analytic
case, however, if the vector fields are globally defined, there are
additional conclusions that can be drawn.  To state this clearly, the
following definition is convenient.
\begin{definition}
Let $r\in\{\infty,\omega\}$\@, let $\man{M}$ be a $\C^r$-manifold, and let
$\sX$ be a family of partially defined vector fields of class $\C^r$\@.  An
orbit $\man{O}\subset\man{M}$ for $\sX$ is
\begin{compactenum}[(i)]
\item \defn{regular} if, for each $x_0\in\man{O}$ there exists a neighbourhood
$\nbhd{U}$ of $x_0$ such that $\dim(\Orb(x,\sX))=\dim(\man{O})$ for each
$x\in\nbhd{U}$ and is
\item \defn{singular} if it is not regular.\oprocend
\end{compactenum}
\end{definition}

With this terminology, we have the following result, first for the smooth
case.
\begin{proposition}
Let\/ $\man{M}$ be a smooth manifold and let\/
$\sX$ be a family of smooth partially defined vector fields.  Then the union
of the regular orbits for\/ $\sX$ is an open dense subset of\/ $\man{M}$\@.
\begin{proof}
By the Orbit Theorem, the distribution whose subspace at $x$ is
$\tb[x]{\Orb(x,\sX)}$ is the distribution associated to a family of smooth
partially defined vector fields.  Thus, by
Proposition~\ref{prop:smooth-partial}\@, this is a smooth distribution.  The
result follows from Proposition~\ref{prop:rank-semicont}\@.
\end{proof}
\end{proposition}

\begin{proposition}
Let\/ $\man{M}$ be a real analytic manifold and let\/
$\sX\subset\sections[\omega]{\tb{\man{M}}}$ be a family of real analytic
vector fields.  Then the following statements hold:
\begin{compactenum}[(i)]
\item the set of singular orbits for\/ $\sX$ is a locally analytic subset
of\/ $\man{M}$\@;
\item if\/ $\man{M}$ is connected, then all regular orbits of\/ $\sX$ have
the same dimension.
\end{compactenum}
\begin{proof}
By Theorem~\ref{the:fingen-orbit-theorem} below it follows that the
distribution whose subspace at $x$ is $\tb[x]{\Orb(x,\sX)}$ is analytic.  The
result then follows from Proposition~\ref{prop:analytic-rank}\@.
\end{proof}
\end{proposition}

Note that in the analytic case, we really do require that the vector fields
be globally defined; Example~\enumdblref{eg:orbits}{enum:partially-orbit}
suffices to illustrate this.

The Orbit Theorem gives us an insightful description of the tangent spaces to
$\sX$-orbits.  However, computationally the description is not the most
useful since it requires that we know something about the group
$\Diff(\sX)$\@.  One can wonder whether there is a simpler ``infinitesimal''
description.  If one has some intuition about things analytic, one might
imagine that such an infinitesimal description is possible for the analytic
version of the Orbit Theorem.  We shall show that this is true.  We begin by
describing a subspace of the tangent spaces to the $\sX$-orbits.

The proof of the next theorem is adapted from that of
\citet[Proposition~4.15]{JB:01}\@.
\begin{theorem}\label{the:Linfty-orbit}
Let\/ $r\in\{\infty,\omega\}$\@, let\/ $\man{M}$ be a\/ $\C^r$-manifold, and
let\/ $\sX$ be a family of\/ $\C^r$-partially defined vector fields on\/
$\man{M}$\@.  Then
\begin{equation*}
\dist{L}^{(\infty)}(\sX)_{x_0}\subset\tb[x_0]{\Orb(x_0,\sX)}
\end{equation*}
for every\/ $x_0\in\man{M}$\@.
\begin{proof}
We abbreviate a local diffeomorphism $(\Phi,\nbhd{U})$ by $\Phi$ and a
partially defined vector field $(X,\nbhd{U})$ by $X$\@.

Let $x_0\in\man{M}$\@.  By Proposition~\ref{prop:liealggen} it suffices to
show that
\begin{equation}\label{eq:Linfty-orbit1}
[X_k,[X_{k-1},\dots,[X_2,X_1]\cdots]](x_0)\in\tb[x_0]{\Orb(x_0,\sX)}
\end{equation}
for any vector fields $X_1,\dots,X_k\in\sX$ whose domains contain $x_0$\@.
We prove this using some notation and three lemmata.  First the notation.
For local diffeomorphisms $\Phi$ and $\Psi$ of $\man{M}$ denote
\begin{equation*}
[\Phi,\Psi]=\Phi^{-1}\scirc\Psi^{-1}\scirc\Phi\scirc\Psi.
\end{equation*}
With this notation we have the following lemma.
\begin{prooflemma}\label{plem:Linfty-orbit1}
Let\/ $X_1,\dots,X_k\in\sections[r]{\tb{\man{M}}}$ and recursively define
\begin{align*}
\Psi_1(t_1)=&\;\flow{X_1}{t_1},\\
\Psi_2(t_1,t_2)=&\;[\flow{X_1}{t_1},\flow{X_2}{t_2}],\\
\Psi_3(t_1,t_2,t_3)=&\;[[\flow{X_1}{t_1},\flow{X_2}{t_2}],\flow{X_3}{t_3}],\\
\vdots\,&\\
\Psi_k(t_1,t_2,\dots,t_k)=&\;
[\cdots[\flow{X_1}{t_1},\flow{X_2}{t_2}],\dots,\flow{X_k}{t_k}],
\end{align*}
for\/ $t_1,\dots,t_k\in\real$ such that all flows are defined.  Then
\begin{equation*}
\left.\pderiv{}{t_1}\right|_{t_1=0}\Psi_k(t_1,\dots,t_k)(x_0)=-X_1(x_0)+
\left.\pderiv{}{t_1}\right|_{t_1=0}\flow{X_k}{t_k}\scirc
\Psi_{k-1}(t_1,\dots,t_{k-1})\scirc\flow{X_2}{t_2}(x_0).
\end{equation*}
\begin{subproof}
It is easy to see by induction that if any of the numbers $t_1,\dots,t_k$ are
zero, then $\Psi_k=\id_{\man{M}}$\@.  We shall use this fact frequently.

First note that differentiation of the relation
\begin{equation*}
\Psi_{k-1}(t_1,\dots,t_{k-1})\scirc\Psi_{k-1}(t_1,\dots,t_{k-1})^{-1}(x)=x
\end{equation*}
gives
\begin{multline*}
\pderiv{}{t_1}(\Psi_{k-1}(t_1,\dots,t_{k-1})^{-1}(x))\\
=-\tf{\Psi_{k-1}(t_1,\dots,t_{k-1})^{-1}}\cdot\left(\left(
\pderiv{}{t_1}\Psi_{k-1}(t_1,\dots,t_{k-1})\right)\scirc
\Psi_{k-1}(t_1,\dots,t_{k-1})^{-1}(x)\right).
\end{multline*}
Evaluating at $t_1=0$ and using the fact stated at the beginning of the proof
then gives
\begin{equation*}
\left.\pderiv{}{t_1}\right|_{t_1=0}
(\Psi_{k-1}(t_1,\dots,t_{k-1})^{-1}(x))=-X_1(x).
\end{equation*}
Using this fact along with the statement made at the beginning of the lemma,
we calculate
\begin{align*}
\left.\pderiv{}{t_1}\right|_{t_1=0}\Psi_k(t_1,\dots,t_k)(x_0)=&\;
\left.\pderiv{}{t_1}\right|_{t_1=0}
[\Psi_{k-1}(t_1,\dots,t_{k-1}),\flow{X_k}{t_k}](x_0)\\
=&\;\left.\pderiv{}{t_1}\right|_{t_1=0}\Psi_{k-1}(t_1,\dots,t_{k-1})^{-1}\scirc
\flow{X_2}{-t_2}\scirc\Psi_{k-1}(t_1,\dots,t_{k-1})\scirc\flow{X_k}{t_k}(x_0)\\
=&\;-X_1(x_0)+\left.\pderiv{}{t_1}\right|_{t_1=0}\flow{X_k}{t_k}\scirc
\Psi_{k-1}(t_1,\dots,t_{k-1})\scirc\flow{X_2}{t_2}(x_0),
\end{align*}
giving the result.
\end{subproof}
\end{prooflemma}

We also recall the definition of the pull-back of a vector field $X$ by a
diffeomorphism $\Phi$\@: $\Phi^*X=\tf{\Phi^{-1}}\scirc X\scirc\Phi$\@.  With
this notation we have the following lemma.
\begin{prooflemma}\label{plem:Linfty-orbit2}
With the notation from Lemma~\ref{plem:Linfty-orbit1}\@,
\begin{equation*}
\left.\pderiv{}{t_k}\cdots\pderiv{}{t_1}\right|_{t_1=\dots=t_k=0}
\Psi_k(t_1,\dots,t_k)(x_0)=
\left.\pderiv{}{t_k}\cdots\pderiv{}{t_2}\right|_{t_2=\dots=t_k=0}
((\flow{X_k}{t_k})^*\cdots(\flow{X_2}{t_2})^*X_1)(x_0).
\end{equation*}
\begin{subproof}
We prove this by induction on $k$\@.  For $k=2$ we use
Lemma~\ref{plem:Linfty-orbit1} to determine that
\begin{align*}
\left.\pderiv{}{t_2}\pderiv{}{t_1}\right|_{t_1=t_2=0}\Psi_2(t_1,t_2)(x_0)=&\;
\left.\pderiv{}{t_2}\pderiv{}{t_1}\right|_{t_1=t_2=0}
\flow{X_2}{-t_2}\scirc\flow{X_1}{t_1}\scirc\flow{X_2}{t_2}(x_0)\\
=&\;\left.\pderiv{}{t_2}\right|_{t_2=0}\tf{\flow{X_2}{-t_2}}\scirc
X_1(\flow{X_2}{t_2}(x_0))
=\left.\pderiv{}{t_2}\right|_{t_2=0}((\flow{X_s}{t_2})^*X_1)(x_0),
\end{align*}
giving the lemma for $k=2$\@.

Now suppose the lemma holds for $k\in\{1,\dots,m-1\}$\@.  An application of
Lemma~\ref{plem:Linfty-orbit1} and the induction hypothesis gives
\begin{align*}
\left.\pderiv{}{t_m}\cdots\pderiv{}{t_1}\right|_{t_1=\dots=t_m=0}&
\Psi_m(t_1,\dots,t_m)(x_0)\\
=&\;\left.\pderiv{}{t_m}\cdots\pderiv{}{t_1}\right|_{t_1=\dots=t_m=0}
\flow{X_m}{-t_m}\scirc\Psi_{m-1}(t_1,\dots,t_{m-1})\scirc\flow{X_m}{t_m}(x_0)\\
=&\;\left.\pderiv{}{t_m}\cdots\pderiv{}{t_2}\right|_{t_2=\dots=t_m=0}
\tf{\flow{X_m}{-t_m}}\cdot\left(\left.\pderiv{}{t_1}\right|_{t_1=0}
\Psi_{m-1}(t_1,\dots,t_{m-1})\scirc\flow{X_m}{t_m}(x_0)\right)\\
=&\;\left.\pderiv{}{t_m}\cdots\pderiv{}{t_2}\right|_{t_2=\dots=t_m=0}
(\flow{X_m}{t_m})^*
\left.\pderiv{}{t_1}\right|_{t_1=0}\Psi_{m-1}(t_1,\dots,t_{m-1})\\
=&\;\left.\pderiv{}{t_m}\right|_{t_m=0}(\flow{X_m}{t_m})^*
\left.\pderiv{}{{t_{m-1}}}\cdots\pderiv{}{t_2}\pderiv{}{t_1}
\right|_{t_1=\dots=t_{m-1}=0}\Psi_{m-1}(t_1,\dots,t_{m-1})\\
=&\;\left.\pderiv{}{t_m}\right|_{t_m=0}(\flow{X_m}{t_m})^*
\left.\pderiv{}{{t_{m-1}}}\cdots\pderiv{}{t_2}\right|_{t_2=\dots=t_{m-1}=0}
((\flow{X_{m-1}}{t_{m-1}})^*\cdots(\flow{X_2}{t_2})^*X_1)(x_0)\\
=&\;\left.\pderiv{}{t_m}\cdots\pderiv{}{t_2}\right|_{t_2=\dots=t_m=0}
(\flow{X_m}{t_m})^*\cdots(\flow{X_2}{t_2})^*X_1)(x_0),
\end{align*}
which is the result.  (Note that our freely swapping partial derivatives with
pull-backs is justified since we are differentiating the pull-back with
respect to its argument, and the pull-back is linear in its argument.)
\end{subproof}
\end{prooflemma}

Now we prove the key fact.
\begin{prooflemma}\label{plem:Linfty-orbit3}
We use the notation from Lemma~\ref{plem:Linfty-orbit2}\@.  For\/
$x_0\in\man{M}$\@, if we define\/ $\Psi_{x_0}(s)=\Psi_k(s,\dots,s)(x_0)$ for
all\/ $s\in\real$ such that the expression makes sense, then
\begin{equation*}
\left.\deriv{^j}{s^j}\right|_{s=0}\Psi_{x_0}(s)=0,\qquad
j\in\{0,1,\dots,k-1\},
\end{equation*}
and
\begin{equation*}
\left.\deriv{^k}{s^k}\right|_{s=0}\Psi_{x_0}(s)=
k![X_k,\dots,[X_2,X_1]\cdots](x_0).
\end{equation*}
\begin{subproof}
Now let $j\in\{0,1,\dots,k-1\}$\@.  By the Chain
Rule for high-order derivatives,~\cite[Supplement~2.4A]{RA/JEM/TSR:88}\@,
\begin{equation}\label{eq:Linfty-orbit2}
\derivatzero{^j}{s^j}\Psi_{x_0}(s)=
\sum_{\substack{j_1,\dots,j_k\in\{0,1,\dots,j\}\\j_1+\dots+j_k=j}}
\frac{j!}{j_1!\cdots j_k!}
\left.\pderiv{^{j_1}}{{t_1^{j_1}}}\cdots\pderiv{^{j_k}}{{t_k^{j_k}}}
\right|_{t_1=\dots=t_k=0}\Phi_k(t_1,\dots,t_k)(x_0).
\end{equation}
Note that each term in the above will have $j_a=0$ for some
$a\in\{1,\dots,k\}$\@.  The partial derivatives
in~\eqref{eq:Linfty-orbit2}\@, when evaluated at $t_1=\dots=t_k=0$\@, will
then necessarily be taken with $t_a=0$ for some $a\in\{1,\dots,k\}$\@.  By
our comment at the beginning of the proof of
Lemma~\ref{plem:Linfty-orbit1}\@, it follows that all such partial
derivatives will be zero.

By the same reasoning, in the expression
\begin{equation*}
\derivatzero{^k}{s^k}\Psi_{x_0}(s)=
\sum_{\substack{j_1,\dots,j_k\in\{0,1,\dots,k\}\\j_1+\dots+j_k=k}}
\frac{k!}{j_1!\cdots j_k!}
\left.\pderiv{^{j_1}}{{t_1^{j_1}}}\cdots\pderiv{^{j_k}}{{t_k^{j_k}}}
\right|_{t_1=\dots=t_k=0}\Phi_k(t_1,\dots,t_k)(x_0)
\end{equation*}
for the $k$th derivative, the only nonzero term in the sum occurs when
$j_1=\dots=j_k=1$\@, since otherwise at least one of the numbers
$j_1,\dots,j_k$ will be zero.  That is to say,
\begin{equation*}
\derivatzero{^k}{s^k}\Psi_{x_0}(s)=k!
\left.\pderiv{}{t_1}\cdots\pderiv{}{t_k}\right|_{t_1=\dots=t_k=0}
\Phi_k(t_1,\dots,t_k)(x_0).
\end{equation*}
Let us now turn to the proof of the fact that this expression is the iterated
Lie bracket in the statement of the lemma.

We prove this by showing that, for any $j\in\{2,\dots,k\}$\@,
\begin{multline}\label{eq:Linfty-orbit3}
\left.\pderiv{}{t_j}\right|_{t_j=0}((\flow{X_k}{t_k})^*\cdots
(\flow{X_j}{t_j})^*[X_{j-1},\dots,[X_2,X_1]\cdots])(x_0)\\
=((\flow{X_k}{t_k})^*\cdots(\flow{X_{j+1}}{t_{j+1}})^*
[X_j,\dots,[X_2,X_1]\cdots])(x_0).
\end{multline}
This we prove by induction on $j$\@.  For $j=2$ we have
\begin{align*}
\left.\pderiv{}{t_2}\right|_{t_2=0}((\flow{X_k}{t_k})^*\cdots(\flow{X_2}{t_2})^*
X_1)(x_0)=&\;((\flow{X_k}{t_k})^*\cdots(\flow{X_3}{t_3})^*
\left.\pderiv{}{t_2}\right|_{t_2=0}\flow{X_2}{t_2}X_1)(x_0)\\
=&\;((\flow{X_k}{t_k})^*\cdots(\flow{X_3}{t_3})^*[X_2,X_1])(x_0),
\end{align*}
using the well-known characterisation of the Lie bracket by
\begin{equation*}
[X,Y](x)=\derivatzero{}{t}(\flow{X}{t})^*Y(x)
\end{equation*}
\cite[see][Theorem~4.2.19]{RA/JEM/TSR:88}\@.  Now suppose
that~\eqref{eq:Linfty-orbit3} holds for $j\in\{2,\dots,m-1\}$\@.  Then we use
the induction hypotheses to get
\begin{align*}
\left.\pderiv{}{t_m}\right|_{t_m=0}&((\flow{X_k}{t_k})^*\cdots
(\flow{X_m}{t_m})^*[X_{m-1},\dots,[X_2,X_1]\cdots])(x_0)\\
=&\;((\flow{X_k}{t_k})^*\cdots(\flow{X_{m+1}}{t_{m+1}})^*
\left.\pderiv{}{t_m}\right|_{t_m=0}\flow{X_m}{t_m}
[X_{m-1},\dots,[X_2,X_1]\cdots])(x_0)\\
=&\;((\flow{X_k}{t_k})^*\cdots(\flow{X_{m+1}}{t_{m+1}})^*
[X_m,\dots,[X_2,X_1]\cdots])(x_0),
\end{align*}
giving~\eqref{eq:Linfty-orbit3}\@.

Now we use Lemma~\ref{plem:Linfty-orbit2} and recursively
apply~\eqref{eq:Linfty-orbit3} to give
\begin{align*}
\left.\pderiv{}{t_k}\cdots\pderiv{}{t_1}\right|_{t_1=\dots=t_k=0}
\Psi_k(t_1,\dots,t_k)(x_0)=&\;
\left.\pderiv{}{t_k}\cdots\pderiv{}{t_2}\right|_{t_2=\dots=t_k=0}
((\flow{X_k}{t_k})^*\cdots(\flow{X_2}{t_2})^*X_1)(x_0)\\
=&\;[X_k,[X_{k-1},\dots,[X_2,X_1]\cdots]](x_0)
\end{align*}
as desired.
\end{subproof}
\end{prooflemma}

We can now complete the proof of the theorem by
verifying~\eqref{eq:Linfty-orbit1}\@.  Recalling the distribution $\dist{O}$
from the proof of the Orbit Theorem, the vector fields $X_1,\dots,X_k$ are
$\dist{O}$-valued and so tangent to $\Orb(x_0,\sX)$\@.  Therefore, with the
notation of Lemma~\ref{plem:Linfty-orbit3}\@, $\Psi_{x_0}(s)\in\Orb(x_0,\sX)$
for every $s\in\real$\@.  Therefore, since the first $k-1$ derivatives of
$\Psi_{x_0}$ at $s=0$ vanish by Lemma~\ref{plem:Linfty-orbit3}\@, the $k$th
derivative of $\Psi_{x_0}$ is tangent to $\Orb(x_0,\sX)$\@.  That is to say,
by Lemma~\ref{plem:Linfty-orbit3}\@,
\begin{equation*}
\left.\deriv{^k}{s^k}\right|_{s=0}\Psi_{x_0}(s)=
k![X_k,\dots,[X_2,X_1]\cdots](x_0)\in\tb[x_0]{\Orb(x_0,\sX)},
\end{equation*}
showing that~\eqref{eq:Linfty-orbit1} holds.
\end{proof}
\end{theorem}

An immediate consequence of Theorems~\ref{the:orbit-theorem}
and~\ref{the:Linfty-orbit} is the following well-known result of
\citet{PKR:38} and \citet{WLC:39}\@.  The proof of \citeauthor{WLC:39} is
given in the Cartan-like framework of differential forms, whereas the proof
of \citeauthor{PKR:38} is, like our proof, a vector field proof.
\begin{corollary}
Let\/ $\man{M}$ be a connected\/ $\C^\infty$-manifold and let\/ $\sX$ be a
family of smooth partially defined vector fields.  If\/
$\dist{L}^{(\infty)}(\sX)=\tb{\man{M}}$ then, for\/ $x_1,x_2\in\man{M}$\@,
there exists\/ $X_1,\dots,X_k\in\sX$ and\/ $t_1,\dots,t_k\in\real$ such that
\begin{equation*}
x_2=\flow{X_k}{t_k}\scirc\cdots\scirc\flow{X_1}{t_1}(x_1).
\end{equation*}
\begin{proof}
Let $x_1\in\man{M}$\@.  By Theorem~\ref{the:Linfty-orbit} we have
$\dist{L}^{(\infty)}(\sX)_x\subset\tb[x]{\Orb(x_1,\sX)}$ for every
$x\in\Orb(x_1,\sX)$ and so
\begin{equation*}
\dist{L}^{(\infty)}(\sX)_x\subset\tb[x]{\Orb(x_1,\sX)}\subset\tb[x]{\man{M}}=
\dist{L}^{(\infty)}(\sX)_x,
\end{equation*}
giving $\tb[x]{\Orb(x_1,\sX)}=\tb[x]{\man{M}}$ for every
$x\in\Orb(x_1,\sX)$\@.  Thus $\Orb(x_1,\sX)$ is an open submanifold of
$\man{M}$\@.  Thus, recalling the basis for the orbit topology from the proof
of the Theorem~\ref{the:orbit-theorem}\@, open subsets of $\Orb(x_1,\sX)$ in
the orbit topology are open subsets in the relative topology on $\man{M}$\@.
Since $\man{M}$ is a disjoint union of its orbits, it is a disjoint union of
open sets.  Each component in this disjoint union is necessarily closed since
its complement is open, being a union of open sets.  Thus each orbit is a
connected component of $\man{M}$\@.  Since $\man{M}$ is assumed connected it
follows that $\Orb(x_1,\sX)=\man{M}$\@, which is the result.
\end{proof}
\end{corollary}

The converse of the Rashevsky\textendash{}Chow Theorem is generally false.
\begin{example}\label{eg:chow-converse}
Recall Example~\enumdblref{eg:orbits}{enum:!analytic-orbit} where
$\man{M}=\real^2$ and where $\sX=\ifam{X_1,X_2}$ is defined by
\begin{equation*}
X_1=\pderiv{}{x_1},\quad X_2=f(x_1)\pderiv{}{x_2},
\end{equation*}
where
\begin{equation*}
f(x)=\begin{cases}\eul^{-1/x^2},&x\in\realp,\\0,&x\in\realnp.\end{cases}
\end{equation*}
In Example~\enumdblref{eg:orbits}{enum:!analytic-orbit} we explicitly showed
that $\man{M}=\Orb(\vect{0},\sX)$\@.  However, one can also directly show
that
\begin{equation*}
\dist{L}^{(\infty)}(\sX)_{(x_1,x_2)}=\begin{cases}
\tb[(x_1,x_2)]{\real^2},&x_1>0,\\
\vecspan[\real]{\pderiv{}{x_1}},&x_1\le0.\end{cases}
\end{equation*}
Thus $\dist{L}^{(\infty)}(\sX)\subsetneq\tb{\man{M}}$\@.\oprocend
\end{example}

\subsection{The finitely generated Orbit Theorem}

Now we turn to characterising situations where the tangent spaces to the
orbits are \emph{exactly} the subspaces $\dist{L}^{(\infty)}(\sX)$\@.
\begin{theorem}\label{the:fingen-orbit-theorem}
Let\/ $r\in\{\infty,\omega\}$\@, let\/ $\man{M}$ be a\/ $\C^r$-manifold, and
let\/ $\sX$ be a family of partially defined vector fields of class\/ $\C^r$
such that\/ $\sL^{(\infty)}(\sX)_x$ is a finitely generated submodule of\/
$\gsections[r]{x}{\tb{\man{M}}}$ for each\/ $x\in\man{M}$\@.  Then, for
each\/ $x_0\in\man{M}$\@,
\begin{compactenum}[(i)]
\item $\Orb(x_0,\sX)$ is a connected immersed\/ $\C^r$-submanifold of\/
$\man{M}$ and
\item for each\/ $x\in\Orb(x_0,\sX)$\@,
$\tb[x]{\Orb(x_0,\sX)}=\dist{L}^{(\infty)}(\sX)_x$\@.
\end{compactenum}
Moreover,\/ $\man{M}$ is the disjoint union of the set of orbits.
\begin{proof}
From the Orbit Theorem and Theorem~\ref{the:Linfty-orbit} it only remains to
show that $\tb[x]{\Orb(x_0,\sX)}\subset\dist{L}^{(\infty)}(\sX)_x$\@.  For
$X\in\sX$ we have $[X,Y]\in\sL^{(\infty)}(\sX)$ for every
$Y\in\sL^{(\infty)}(\sX)$\@, this since $\sL^{(\infty)}(\sX)$ is a Lie
subalgebra.  Since $\sL^{(\infty)}(\sX)_x$ is assumed to be finitely
generated, Theorem~\ref{the:invariant-sheaves} and~\eqref{eq:partially-dist}
gives $(\flow{X}{t})_*Y(x)\in\dist{L}^{(\infty)}(\sX)_x$ for every $X\in\sX$
and $t\in\real$ such that $x\in\flow{X}{t}(\nbhd{U}(X,t))$\@.  A trivial
induction then gives
\begin{equation}\label{eq:fingenorbit2}
(\flow{X_k}{t_k})_*\cdots(\flow{X_1}{t_1})_*Y(x)=
(\flow{X_k}{t_k}\scirc\cdots\scirc\flow{X_1}{t_1})_*Y(x)
\in\dist{L}^{(\infty)}(\sX)_x
\end{equation}
for every suitable $X_1,\dots,X_k\in\sX$ and $t_1,\dots,t_k\in\real$\@, where
we use the fact that push-forward commutes with
composition~\cite[Proposition~4.2.3]{RA/JEM/TSR:88}\@.  However, by
Theorem~\ref{the:orbit-theorem}\@,
\begin{equation*}
\tb[x]{\Orb(x_0,\sX)}=\setdef{\Phi_*X(x)}{\Phi\in\Diff(\sX),\ X\in\sX},
\end{equation*}
and so~\eqref{eq:fingenorbit2} implies that
$\tb[x]{\Orb(x_0,\sX)}\subset\dist{L}^{(\infty)}(\sX)$\@.
\end{proof}
\end{theorem}

\begin{remark}
In the statement of the preceding theorem we asked that
$\sL^{(\infty)}(\sX)_x$ be a submodule of $\gsections[r]{x}{\tb{\man{M}}}$
for each $x\in\man{M}$\@.  For a general family $\sX$ of vector fields, it
will not be the case that $\sL^{(\infty)}(\sX)_x$ is a submodule.  However,
if all one is interested in is the tangent spaces to orbits, then, by
Propositions~\ref{prop:liealgmodule} and~\ref{prop:Linfty-equiv1}\@, one can
replace with $\sX_x$ with the module $\modgen{\sX_x}$ generated by $\sX_x$\@.
Equivalently, also by Proposition~\ref{prop:Linfty-equiv1}\@, one can replace
$\sL^{(\infty)}(\sX)_x$ with the module $\modgen{\sL^{(\infty)}(\sX)_x}$
generated by $\sL^{(\infty)}(\sX)_x$\@.  As long as the module
$\sL^{(\infty)}(\modgen{\sX_x})$ or $\modgen{\sL^{(\infty)}(\sX)_x}$ is
locally finitely generated, it will hold that
$\tb[x]{\Orb(x_0,\sX)}=\dist{L}^{(\infty)}(\sX)_x$ for all\/ $x_0\in\man{M}$
and $x\in\Orb(x_0,\sX)$\@.\oprocend
\end{remark}

Again, it is important to distinguish between a distribution generated by a
family of vector fields being finitely generated and the module generated by
a family of vector fields being finitely generated.

This gives the following important results for families of analytic vector
fields and certain families of smooth vector fields.
\begin{corollary}
Let\/ $\man{M}$ be a\/ $\C^\infty$-manifold and let\/ $\sX$ be a family of
partially defined smooth vector fields such that the distribution\/
$\dist{L}^{(\infty)}(\sX)$ is regular.  Then, for each\/ $x_0\in\man{M}$\@,
\begin{compactenum}[(i)]
\item $\Orb(x_0,\sX)$ is a connected immersed smooth submanifold of\/
$\man{M}$ and
\item for each\/ $x\in\Orb(x_0,\sX)$\@,
$\tb[x]{\Orb(x_0,\sX)}=\dist{L}^{(\infty)}(\sX)$\@.
\end{compactenum}
Moreover,\/ $\man{M}$ is the disjoint union of the set of orbits.
\begin{proof}
This follows from Theorem~\ref{the:fingen-orbit-theorem}\@, along with
Theorem~\ref{the:locally-free-sheaf}\@.
\end{proof}
\end{corollary}

\begin{corollary}\label{cor:analytic-orbit-theorem}
Let\/ $\man{M}$ be an analytic manifold and let\/ $\sX$ be a family of
partially defined analytic vector fields.  Then, for each\/
$x_0\in\man{M}$\@,
\begin{compactenum}[(i)]
\item $\Orb(x_0,\sX)$ is a connected immersed analytic submanifold of\/
$\man{M}$ and
\item for each\/ $x\in\Orb(x_0,\sX)$\@,
$\tb[x]{\Orb(x_0,\sX)}=\dist{L}^{(\infty)}(\sX)$\@.
\end{compactenum}
Moreover,\/ $\man{M}$ is the disjoint union of the set of orbits.
\begin{proof}
This follows from Theorem~\ref{the:fingen-orbit-theorem}\@, along with
Theorem~\ref{the:analytic-fingen}\@.
\end{proof}
\end{corollary}

The hypothesis of finite generation is necessary, and
Example~\ref{eg:chow-converse} serves to demonstrate this necessity.

\subsection{The fixed-time Orbit Theorem}

In this section we give the version of the Orbit Theorem corresponding to the
fixed-time orbits considered in Section~\ref{subsec:Torbits}\@.  In order to
understand the tangent spaces to the fixed-time orbits for a family
$\sX=\ifam{(X_j,\nbhd{U}_j)}_{j\in J}$ of partially defined vector fields, we
introduce a family of partially defined vector fields by
\begin{equation*}
\sX_0=\asetdef{\sum_{l=1}^k\lambda_l(X_{j_l},\nbhd{U}_{j_l})}
{j_1,\dots,j_k\in J,\ \lambda_1,\dots,\lambda_k\in\real,\
\sum_{l=1}^k\lambda_j=0,\ k\in\integerp}.
\end{equation*}
For brevity, in the preceding expression we have suppressed the domain of the
partially defined vector fields.

With the family of vector fields $\sX_0$ at our disposal, we can state the
fixed-time Orbit Theorem.
\begin{theorem}\label{the:orbitT-theorem}
Let\/ $r\in\{\infty,\omega\}$\@, let\/ $\man{M}$ be a\/ $\C^r$-manifold, and
let\/ $\sX$ be a family of partially defined vector fields of class\/
$\C^r$\@.  For\/ $x_0\in\man{M}$ and\/ $T\in\realp$ satisfying\/
$\Orb_T(x_0,\sX)\not=\emptyset$\@, it holds that
\begin{compactenum}[(i)]
\item \label{pl:ftorbit-theorem1} $\Orb_T(x_0,\sX)$ is a connected immersed\/
$\C^r$-submanifold of\/ $\man{M}$ and
\item \label{pl:ftorbit-theorem2} for each\/ $x\in\Orb_T(x_0,\sX)$\@,
\begin{equation*}
\tb[x]{\Orb_T(x_0,\sX)}=L(\affhull(\setdef{\Phi_*X(x)}
{\Phi\in\Diff(\sX),\ X\in\sX})),
\end{equation*}
where the right-hand side of this expression denotes the linear subspace
associated with the affine hull.
\end{compactenum}
As a result,\/ $\dim(\Orb(x_0,\sX))-\dim(\Orb_T(x_0,\sX))\in\{0,1\}$\@.
\begin{proof}
We shall consider the manifold $\real\times\man{M}$ and so let us introduce
some convenient notation for using this manifold.  For
$\nbhd{U}\subset\man{M}$ open, note that
$\tb{(\real\times\nbhd{U})}=\tb{\real}\times\tb{\nbhd{U}}\simeq
\real\times\real\times\tb{\nbhd{U}}$\@.  For a partially defined vector field
$(X,\nbhd{U})$ define a partially vector field
$(\hat{X},\real\times\nbhd{U})$ by
\begin{equation*}
\hat{X}(s,x)=(s,1,X(x))\in\real\times\real\times\tb{\nbhd{U}}.
\end{equation*}
The following simple lemma gives some of the useful
properties\textemdash{}the last of which we shall not use until the proof of
Theorem~\ref{the:fingen-orbitT-theorem} below\textemdash{}of this natural
extension of vector fields from $\man{M}$ to $\real\times\nbhd{U}$\@.
\begin{prooflemma}\label{plem:hattedvf}
For partially defined vector fields\/ $(X,\nbhd{U})$ and\/ $(Y,\nbhd{V})$\@,
the following statements hold when they make sense:
\begin{compactenum}[(i)]
\item \label{pl:hattedvf1} $\flow{\hat{X}}{t}(s,x)=(s+t,\flow{X}{t}(x))$\@;
\item \label{pl:hattedvf2}
$(\flow{\hat{X}}{t})_*\hat{Y}=\widehat{(\flow{X}{t})_*Y}$\@;
\item \label{pl:hattedvf3} $[\hat{X},\hat{Y}](s,x)=(s,0,[X,Y](x))$\@.
\end{compactenum}
\begin{subproof}
\eqref{pl:hattedvf1} We have
\begin{equation*}
\deriv{}{t}(s+t,\flow{X}{t}(x))=(s+t,1,X(\flow{X}{t}(x)))=
\hat{X}(s+t,\flow{X}{t}(x)),
\end{equation*}
from which the desired conclusion follows by definition of integral curves.

\eqref{pl:hattedvf2} Using~\eqref{pl:hattedvf1} we compute
\begin{align*}
(\flow{\hat{X}}{t})_*\hat{Y}(s,x)=&\;
\tf{\flow{\hat{X}}{t}}\scirc\hat{Y}\scirc\flow{\hat{X}}{-t}(s,x)=
\tf{\flow{\hat{X}}{t}}(s-t,1,Y(\flow{X}{-t}(x)))\\
=&\;(s,1,\tf{\flow{X}{t}}\scirc Y\scirc\flow{X}{-t}(x))=
\widehat{(\flow{X}{t})_*Y}(s,x),
\end{align*}
as desired.

\eqref{pl:hattedvf3} Let $f\in\func[r]{\real\times\nbhd{U}}$\@, denote by
$\map{f_s}{\nbhd{U}}{\real}$ and $\map{f^x}{\real}{\real}$ the functions
$f_s(x)=f^x(s)=f(s,x)$\@.  Then compute
\begin{equation*}
\hat{X}f(s,x)=\deriv{f^x}{s}+Xf_s(x),\quad
\hat{Y}f(s,x)=\deriv{f^x}{s}+Yf_s(x),
\end{equation*}
which gives
\begin{multline*}
(\hat{X}\hat{Y}-\hat{Y}\hat{X})f(s,x)=\deriv[1]{f^x}{x^2}+
Y\pderiv{f}{s}(s,x)+X\pderiv{f}{s}(s,x)+XYf_s(x)\\
-\deriv[1]{f^x}{x^2}-X\pderiv{f}{s}(s,x)-Y\pderiv{f}{s}(s,x)-YXf_s(x)=
[X,Y]f_s(x),
\end{multline*}
and our result follows from this.
\end{subproof}
\end{prooflemma}

Now let us define a family $\hat{\sX}$ of partially defined $\C^r$-vector
fields on $\real\times\man{M}$ by
\begin{equation*}
\hat{\sX}=\setdef{(\hat{X},\real\times\nbhd{U})}{(X,\nbhd{U})\in\sX}.
\end{equation*}
As usual, when it is convenient we will suppress the domain
$\real\times\nbhd{U}$ when we write the partially defined vector field
$(\hat{X},\real\times\nbhd{U})$\@.  We can describe the fixed-time orbit in
terms of the orbits of $\hat{\sX}$\@.
\begin{prooflemma}
For each\/ $s_0\in\real$\@,
\begin{equation*}
\Orb_T(x_0,\sX)=\setdef{x\in\man{M}}{(s_0+T,x)\in\Orb((s_0,x_0),\hat{\sX})}.
\end{equation*}
\begin{subproof}
Let $x\in\Orb_T(x_0,\sX)$ so that $x=\flow{\vect{X}}{\vect{t}}(x_0)$ where
$\vect{X}\in\sX^k$ and $\vect{t}\in\real^k$ satisfies $\sum_{j=1}^kt_j=T$\@.
If $\vect{X}=\ifam{X_1,\dots,X_k}$\@, let us denote
$\hat{\vect{X}}=\ifam{\hat{X}_1,\dots,\hat{X}_k}$\@.  Then, by
Lemma~\ref{plem:hattedvf}\@, note that
\begin{equation*}
\flow{\hat{\vect{X}}}{\vect{t}}(s_0,x_0)=
(s_0+T,\flow{\vect{X}}{\vect{t}}(x_0))=(s_0+T,x).
\end{equation*}
This shows that
\begin{equation*}
\Orb_T(x_0,\sX)\subset
\setdef{x\in\man{M}}{(s_0+T,x)\in\Orb((s_0,x_0),\hat{\sX})}.
\end{equation*}

Conversely, suppose that $(s_0+T,x)\in\Orb((s_0,x_0),\hat{\sX})$\@.  Thus
there exists $\hat{\vect{X}}\in\hat{\sX}{}^k$ and $\vect{t}\in\real^k$ such
that $(s_0+T,x)=\flow{\hat{\vect{X}}}{\vect{t}}(s_0,x_0)$\@.  By the form of
the flow of vector fields from $\hat{\sX}$ as given in
Lemma~\ref{plem:hattedvf} we must have $\sum_{j=1}^kt_j=T$ and
$x=\flow{\vect{X}}{\vect{t}}(x_0)$\@, giving
\begin{equation*}
\setdef{x\in\man{M}}{(s_0+T,x)\in\Orb((s_0,x_0),\hat{\sX})}\subset
\Orb_T(x_0,\sX),
\end{equation*}
as desired.
\end{subproof}
\end{prooflemma}

\eqref{pl:ftorbit-theorem1} In what follows, we make the identification
\begin{equation}\label{eq:orbitT-theorem1}
\Orb_T(x_0,\sX)=\setdef{(T,x)\in\real\times\man{M}}
{(T,x)\in\Orb((0,x_0),\hat{\sX})},
\end{equation}
according to the lemma.  We can now prove the first part of the theorem.  Let
$\man{M}_T=\{T\}\times\man{M}$\@, noting that $\man{M}_T$ is obviously a
submanifold of $\real\times\man{M}$\@.  We claim that
\begin{equation*}
\tb[(T,x)]{\Orb((0,x_0),\hat{\sX})}+\tb[(T,x)]{\man{M}_T}=
\tb[(T,x)]{(\real\times\man{M})}.
\end{equation*}
This follows since $\codim(\man{M}_T)=1$ and since, for each $X\in\sX$\@, the
vector field $\hat{X}$ is tangent to $\Orb((0,x_0),\hat{\sX})$ and is not
tangent to $\man{M}_T$\@.  Thus $\Orb_T(x_0,\sX)$ and $\man{M}_T$ intersect
transversely, and so it follows from~\citet[Corollary~3.5.13]{RA/JEM/TSR:88}
that
\begin{equation*}
\Orb_T(x_0,\sX)=\Orb((0,x_0),\hat{\sX})\cap\man{M}_T
\end{equation*}
is an immersed submanifold of class $C^r$ by the Orbit Theorem.  To prove the
first part of the theorem, it remains to show that $\Orb_T(x_0,\sX)$ is
connected.  Let $x\in\Orb_T(x_0,\sX)$\@, let $X_1,\dots,X_k\in\sX$\@, and let
$t_1,\dots,t_k\in\real$ be such that
\begin{equation*}
\flow{X_1}{t_1}\scirc\cdots\scirc\flow{X_k}{t_k}(x_0)=x,\quad
\sum_{j=1}^kt_j=0,
\end{equation*}
this being possible by Proposition~\ref{prop:OrbTchar}\@.  Then the curve
\begin{equation*}
\interval[0,1]\ni t\mapsto
\flow{X_1}{t_1t}\scirc\cdots\scirc\flow{X_k}{t_kt}(x_0)
\end{equation*}
is a curve in $\Orb_T(x_0,\sX)$ (again by Proposition~\ref{prop:OrbTchar})
connecting $x_0$ and $x$\@, giving path connectedness and so connectedness of
$\Orb_T(x_0,\sX)$\@.

\eqref{pl:ftorbit-theorem2} Given our identification of $\Orb_T(x_0,\sX)$ as
an immersed submanifold of $\real\times\man{M}$\@, we can describe its
tangent space accordingly:
\begin{equation}\label{eq:Torbit-identify}
\tb[(T,x)]{\Orb_T(x_0,\sX)}=
\setdef{(a,v_x)\in\tb[(T,x)]{(\real\times\man{M})}}
{a=0,\ (0,v_x)\in\tb[(T,x)]{\Orb((0,x_0),\hat{\sX})}}.
\end{equation}
By the Orbit Theorem,
\begin{equation*}
\tb[(T,x)]{\Orb((0,x_0),\hat{\sX})}=
\vecspan[\real]{\hat{\Phi}_*\hat{X}(T,x)\mid\enspace
\hat{\Phi}\in\Diff(\hat{\sX}),\ \hat{X}\in\hat{\sX}}.
\end{equation*}
By Lemma~\ref{plem:hattedvf}\@, the definition of $\Diff(\hat{\sX})$
from~\eqref{eq:Diff(X)}\@, and an elementary induction, for
$\hat{\Phi}\in\Diff(\hat{\sX})$ and $\hat{X}\in\hat{\sX}$ we have
\begin{equation}\label{eq:hatpf}
\hat{\Phi}_*\hat{X}(T,x)=\widehat{\Phi_*X}(T,x)=(T,1,\Phi_*X(x))
\end{equation}
for $\Phi\in\Diff(\sX)$\@.  Let $\Phi_1,\dots,\Phi_k\in\Diff(\hat{\sX})$\@,
let $\hat{X}_1,\dots,\hat{X}_k\in\hat{\sX}$\@, and let
$\lambda_1,\dots,\lambda_k\in\real$\@.  Then
\begin{align*}
&\sum_{j=1}^k\lambda_j\hat{\Phi}_{j*}\hat{X}_j(T,x)\in
\tb[(T,x)]{\Orb_T(x_0,\sX)}\\
\iff\quad&
\left(T,\sum_{j=1}^k\lambda_j,
\sum_{j=1}^k\lambda_j\Phi_{j*}X_j(x)\right)\in\tb[(T,x)]{\Orb_T(x_0,\sX)}
\end{align*}
using~\eqref{eq:hatpf}\@.  We immediately conclude
from~\eqref{eq:Torbit-identify} that $\sum_{j=1}^k\lambda_j=0$\@, and so
\begin{equation*}
\tb[x]{\Orb_T(x_0,\sX)}=L(\affhull(\setdef{\Phi_*X(x)}
{\Phi\in\Diff(\sX),\ X\in\sX})),
\end{equation*}
recalling that the linear part of the affine hull of a subset $S$ of a vector
space $\alg{V}$ is given by
\begin{equation*}
\sum_{j=1}^k\lambda_jv_j,\qquad k\in\integerp,\ v_1,\dots,v_k\in S,\
\lambda_1,\dots,\lambda_k\in\real,\ \sum_{j=1}^k\lambda_j=0
\end{equation*}
(this follows directly from the definition of the affine
hull,~\cf~Theorem~1.2.5 of~\cite{RJW:94} and Theorem~1.2 of~\cite{RTR:70}).

For the final assertion of the theorem, let us abbreviate
\begin{equation*}
\alg{V}_x=\vecspan[\real]{\setdef{\Phi_*X(x)}{\Phi\in\Diff(\sX),\ X\in\sX}}
\subset\tb[x]{\man{M}}
\end{equation*}
and
\begin{equation*}
\alg{A}_x=\affhull(\setdef{\Phi_*X(x)}{\Phi\in\Diff(\sX),\ X\in\sX})
\subset\tb[x]{\man{M}}.
\end{equation*}
We have two cases,~\cf~the discussion on the bottom of page~8
of~\cite{RJW:94}\@.
\begin{enumerate}
\item $\alg{A}_x=\alg{V}_x$\@: In this case we have
$\alg{V}_x=\alg{A}_x=L(\alg{A}_x)$ and so
\begin{equation*}
\dim(\Orb(x_0,\sX))=\dim(\Orb_T(x_0,\sX))
\end{equation*}
since the tangent spaces to the orbit and the fixed-time orbits agree.
\item $\alg{A}_x\subsetneq\alg{V}_x$\@: In this case there exists
$v_0\in\alg{V}_x\setminus\alg{A}_x$ and
\begin{equation*}
\alg{V}_x=\vecspan[\real]{\{v_0\}\cup\alg{A}_x},
\end{equation*}
which gives
\begin{equation*}
\dim(\Orb(x_0,\sX))=\dim(\Orb_T(x_0,\sX))+1
\end{equation*}
as desired.\qed
\end{enumerate}
\end{proof}
\end{theorem}

\begin{remark}
As with the proof of the Orbit Theorem, the proof of the fixed-time Orbit
Theorem prescribes a topology on $\man{M}$\@.  Let us extract this
\defn{fixed-time orbit topology} here.  We will use the notation from the
proof of the fixed-time Orbit Theorem.  In the proof we saw that
$\Orb_T(x_0,\sX)$ was naturally identified with an immersed submanifold of
$\real\times\man{M}$ and was also a subset of the (non-fixed-time) orbit
$\Orb((0,x_0),\hat{\sX})$\@.  Since $\Orb((0,x_0),\hat{\sX})$ has the orbit
topology, the natural topology on $\Orb_T(x_0,\sX)$ is that induced by this
orbit topology.

As with the orbit topology, one can describe the fixed-time orbit topology as
a final topology.  We do this as follows.  Denote
\begin{equation*}
\real^k_0=\setdef{\vect{t}\in\real^k}{t_1+\dots+t_k=0},
\end{equation*}
noting that $\real^k_0$ is a $(k-1)$-dimensional subspace.  Now let
$x\in\man{M}$\@, let $k\in\integerp$\@, let
$\vect{X}=\ifam{X_1,\dots,X_k}\subset\sX$\@, and let $\nbhd{U}\subset\real^k$
be a neighbourhood of $\vect{0}$ be such that the map
\begin{equation*}
\real^k_0\cap\nbhd{U}\ni\vect{t}\mapsto\flow{\vect{X}}{\vect{t}}(x)\in\man{M}
\end{equation*}
is defined.  One may then define the orbit topology as the final topology
induced by the above family of mappings.\oprocend
\end{remark}

As with Theorem~\ref{the:Linfty-orbit} for orbits, there is an easily
described subspace of the tangent spaces to the fixed-time orbits.  Let us
describe this here.  We let $\sD(\sX)$ be the derived algebra of
$\sL^{(\infty)}(\sX)$\@.  That is to say, $\sD(\sX)$ is the subspace of
$\sL^{(\infty)}(\sX)$ generated by vector fields of the form $[Y_1,Y_2]$
where $Y_1,Y_2\in\sL^{(\infty)}(\sX)$\@.  (Here again we are suppressing the
domain of partially defined vector fields.)  As with Lie algebras of globally
defined vector fields, we have the notion of an ideal in
$\sL^{(\infty)}(\sX)$\@.  Indeed, a subset $\sI\subset\sL^{(\infty)}(\sX)$ is
an \defn{ideal} if $[X,Y]\in\sI$ for every $X\in\sI$ and
$Y\in\sL^{(\infty)}(\sX)$\@.  Note that $\sI_x$ is then an ideal of
$\sL^{\infty)}(\sX)_x$ in the usual sense.

With the preceding comments, we have the following explicit characterisation
of the derived algebra.
\begin{proposition}\label{prop:derivedalg}
Let\/ $r\in\{\infty,\omega\}$\@, let\/ $\man{M}$ be a\/ $\C^r$-manifold, and
let\/ $\sX$ be a family of partially defined\/ $\C^r$-vector fields.  Then
the derived algebra\/ $\sL^{(\infty)}(\sX)$ is comprised of finite\/
$\real$-linear combinations of vector fields of the form
\begin{equation*}
[X_k,[X_{k-1},\dots,[X_2,X_1]\cdots]],\qquad X_1,\dots,X_k\in\sX,\ k\ge2.
\end{equation*}
\begin{proof}
We first claim that the derived algebra is an ideal of
$\sL^{(\infty)}(\sX)$\@, indeed the ideal generated by elements of the form
$[X_1,X_2]$ for $X_1,X_2\in\sX$\@.  First, it is clear from the definition
that $\sD(\sX)$ is an ideal.  It is also clear, since
$\sX\subset\sL^{(\infty)}(\sX)$\@, that $[X_1,X_2]\in\sD(\sX)$ for every
$X_1,X_2\in\sX$\@.  Thus $\sD(\sX)$ contains the ideal generated by brackets
from $\sX$\@.  Now consider an element from $\sD(\sX)$ of the form
$[Y_1,Y_2]$ for $Y_1,Y_2\in\sL^{(\infty)}(\sX)$\@, noting that all elements
of $\sD(\sX)$ are finite linear combinations of such elements.  By
Proposition~\ref{prop:liealggen} (more properly, its adaptation to partially
defined vector fields), $[Y_1,Y_2]$ is a finite linear combination of
brackets of the type in the statement of the result.  Moreover, if we
consider the proof of Proposition~\ref{prop:liealggen}\@, we can see that the
brackets involved will be of the form in the statement of the proposition
with $k\ge2$\@.  From this we conclude our claim that $\sD(\sX)$ is the ideal
of $\sL^{(\infty)}(\sX)$ generated by brackets $[X_1,X_2]$ for
$X_1,X_2\in\sX$\@.

From the fact that the derived algebra is an ideal and that it contains all
vector fields of the form $[X_1,X_2]$ for $X_1,X_2\in\sX$\@, it follows that
\begin{equation*}
[X_k,[X_{k-1},\dots,[X_2,X_1]\cdots]]\in\sD(\sX)
\end{equation*}
for every $X_1,\dots,X_k\in\sX$\@, $k\ge2$\@.  Conversely, the set of finite
$\real$-linear combinations of the form in the statement of the result is
easily shown to be an ideal of $\sL^{(\infty)}(\sX)$ and it clearly contains
the brackets $[X_1,X_2]$ for $X_1,X_2\in\sX$\@.  Thus $\sD(\sX)$ is contained
in the this set of linear combinations, which completes the proof.
\end{proof}
\end{proposition}

We then define
\begin{equation*}
\sI(\sX)=\vecspan[\real]{X+Y|\enspace X\in\sX_0,\ Y\in\sD(\sX)}.
\end{equation*}
The following characterisation of $\sI(\sX)$ is useful.
\begin{proposition}
Let\/ $r\in\{\infty,\omega\}$\@, let\/ $\man{M}$ be a\/ $\C^r$-manifold, and
let\/ $\sX$ be a family of partially defined\/ $\C^r$-vector fields.  Then
the following statements hold:
\begin{compactenum}[(i)]
\item \label{pl:Linfty01} $\sI(\sX)$ is an ideal of\/
$\sL^{(\infty)}(\sX)$\@;
\item \label{pl:Linfty02} the codimension of\/ $\sI(\sX)_x$ in
$\sL^{(\infty)}(\sX)_x$ is zero if\/
$\sI(\sX)_x\cap\sX_x\not=\emptyset$ and is one otherwise.
\end{compactenum}
\begin{proof}
\eqref{pl:Linfty01} If $Y\in\sL^{(\infty)}(\sX)$ and if $X\in\sI(\sX)$\@,
then $[Y,X]$ is obviously in the derived algebra of $\sL^{(\infty)}(\sX)$\@,
by definition of the derived algebra.  Since $\sD(\sX)\subset\sI(\sX)$\@,
this part of the result follows.

\eqref{pl:Linfty02} From Propositions~\ref{prop:liealggen}
and~\ref{prop:derivedalg} we have that any element of $\sL^{(\infty)}(\sX)$
can be written as
\begin{equation*}
\sum_{j=1}^k\lambda_jX_j+Y,\qquad X_1,\dots,X_k\in\sX,\
\lambda_1,\dots,\lambda_k\in\real,\ Y\in\sD(\sX)
\end{equation*}
(suppressing the domains of partially defined vector fields, as usual).  Thus
$\sL^{(\infty)}(\sX)_x$ is the sum of the subspaces $\sX_x$ and
$\sD(\sX)_x$\@.  Referring to Proposition~\ref{prop:derivedalg}\@,
$\sI(\sX)_x$ is the sum of the subspaces $L(\affhull(\sX_x))$ and
$\sD(\sX)_x$\@.  Note that $L(\affhull(\sX_x))$ is a subspace of $\sX_x$\@.
Moreover, as in the last step in the proof of
Theorem~\ref{the:orbitT-theorem}\@, $L(\affhull(\sX_x))=\sX_x$ if and only if
$\sX_x\cap L(\affhull(\sX_x))\not=\emptyset$\@.  Also as in the last step in
the proof of Theorem~\ref{the:orbitT-theorem}\@, if
$\sX_x\cap L(\affhull(\sX_x))=\emptyset$\@, then the codimension of
$L(\affhull(\sX_x))$ in $\sX_x$ is one.  This gives this part of the result.
\end{proof}
\end{proposition}

We define
\begin{equation*}
\dist{I}(\sX)_x=\setdef{X(x)}{X\in\sI(\sX)}
\end{equation*}
so that $\dist{I}(\sX)$ is a distribution on $\man{M}$\@.  Since $\sI(\sX)$
is an ideal of $\sL^{(\infty)}(\sX)$\@, it is also a Lie subalgebra, and so
is a Lie subalgebra of $\sections[r]{\tb{\man{M}}}$\@.  The picture one
should have in mind is that $\sI(\sX)$ is to $\Orb_T(x,\sX)$ what
$\sL^{(\infty)}(\sX)$ is to $\Orb(x,\sX)$\@.  For example, one should think
of $\sL^{(\infty)}(\sX)$ as being the ``Lie algebra'' of the ``Lie group''
$\Diff(\sX)$\@.  Upon doing so, one should think of $\sI(\sX)$ as the Lie
subalgebra (actually ideal) of $\sL^{(\infty)}(\sX)$ corresponding to the
``subgroup'' (actually, ``normal subgroup'') $\Diff_0(\sX)$ of
$\Diff(\sX)$\@.  Moreover, we have the following theorem.
\begin{theorem}\label{the:Linfty-Torbit}
Let\/ $r\in\{\infty,\omega\}$\@, let\/ $\man{M}$ be a\/ $\C^r$-manifold, and
let\/ $\sX$ be a family of partially defined\/ $\C^r$-vector fields on\/
$\man{M}$\@.  For\/ $x_0\in\man{M}$ and\/ $T\in\real$ satisfying\/
$\Orb_T(x_0,\sX)\not=\emptyset$\@, it holds that
\begin{equation*}
\dist{I}(\sX)_x\subset\tb[x]{\Orb_T(x_0,\sX)}
\end{equation*}
for every\/ $x\in\Orb_T(x_0,\sX)$\@.
\begin{proof}
Let $x\in\Orb_T(x_0,\sX)$\@.  Let $X_1,\dots,X_k\in\sX$ and let
$\lambda_1,\dots,\lambda_k\in\real$ satisfy $\sum_{j=1}^k\lambda_j=0$\@.
Then take $X=\sum_{j=1}^k\lambda_jX_j$ (we remind the reader once again that
we suppress the domains of partially defined vector fields).  We have
\begin{align*}
X(x)=&\;\lambda_1X_1(x)+\dots+\lambda_kX_k(x)=
\derivatzero{}{t}\flow{X}{t}(x)=
\derivatzero{}{t}\flow{\lambda_1X_1+\dots+\lambda_kX_k}{t}(x)\\
=&\;\derivatzero{}{t}\flow{\lambda_1X_1}{t}\scirc\cdots\scirc
\flow{\lambda_kX_k}{t}(x)=
\derivatzero{}{t}\flow{X_1}{\lambda_1t}\scirc\cdots\scirc
\flow{X_k}{\lambda_kt}(x).
\end{align*}
Since
\begin{equation*}
\flow{X_1}{\lambda_1t}\scirc\cdots\scirc\flow{X_k}{\lambda_kt}\in
\Diff_0(\sX),
\end{equation*}
it follows from Proposition~\ref{prop:OrbTchar} that
$X(x)\in\tb[x]{\Orb_T(x_0,\sX)}$\@.

Recalling the notation $[\Phi,\Psi]$ for local diffeomorphisms $\Phi$ and
$\Psi$ from the proof of Theorem~\ref{the:Linfty-orbit}\@, note that for
$t_1,t_2\in\real$ we have
$[\flow{X_1}{t_1},\flow{X_2}{t_2}]\in\Diff_0(\sX)$\@.  An induction then
gives
\begin{equation*}
[\cdots[\flow{X_1}{t_1},\flow{X_2}{t_2}],\dots,\flow{X_k}{t_k}]\in
\Diff_0(\sX)
\end{equation*}
for vector fields $X_1,\dots,X_k\in\sX$ and $t_1,\dots,t_k\in\real$\@.  Thus
the curve
\begin{equation*}
s\mapsto[\cdots[\flow{X_1}{s},\flow{X_2}{s}],\dots,\flow{X_k}{s}](x)
\end{equation*}
is a curve in $\Orb_T(x_0,\sX)$ by Proposition~\ref{prop:OrbTchar}\@.  From
Lemma~\ref{plem:Linfty-orbit3} from the proof of
Theorem~\ref{the:Linfty-orbit} we then have
\begin{equation*}
[X_k,\dots,[X_2,X_1]\cdots](x)\in\tb[x]{\Orb_T(x_0,\sX)}
\end{equation*}
for any $X_1,\dots,X_k\in\sX$ and for $x\in\Orb_T(x_0,\sX)$\@.  Combining
this conclusion with that from the preceding paragraph, and using
Proposition~\ref{prop:liealggen}\@, we have that
$\dist{I}(\sX)_x\subset\tb[x]{\Orb_T(x_0,\sX)}$\@, as desired.
\end{proof}
\end{theorem}

\subsection{The finitely generated fixed-time Orbit Theorem}

Our final version of the Orbit Theorem considers the fixed-time orbits, but
under the assumption that $\sI(\sX)$ is a locally finitely generated module.
\begin{theorem}\label{the:fingen-orbitT-theorem}
Let\/ $r\in\{\infty,\omega\}$\@, let\/ $\man{M}$ be a\/ $\C^r$-manifold, and
let\/ $\sX$ be a family of partially defined vector fields of class\/ $\C^r$
such that\/ $\sL^{(\infty)}(\sX)_x$ is a finitely generated submodule of\/
$\gsections[r]{x}{\tb{\man{M}}}$ for each\/ $x\in\man{M}$\@.  Then, for
each\/ $x_0\in\man{M}$ and\/ $T\in\real$ such that\/
$\Orb_T(x_0,\sX)\not=\emptyset$\@,
\begin{compactenum}[(i)]
\item $\Orb_T(x_0,\sX)$ is a connected immersed\/ $\C^r$-submanifold of\/
$\man{M}$ and
\item for each\/ $x\in\Orb_T(x_0,\sX)$\@,\/
$\tb[x]{\Orb_T(x_0,\sX)}=\dist{I}(\sX)_x$\@.
\end{compactenum}
\begin{proof}
From Proposition~\ref{prop:OrbTchar}\@, the fixed-time Orbit Theorem, and
Theorem~\ref{the:Linfty-Torbit}\@, it only remains to show that
$\tb[x]{\Orb_T(x_0,\sX)}\subset\dist{I}^{(\infty)}(\sX)_x$\@.

We shall adopt the notation and assume the setting of
Theorem~\ref{the:orbitT-theorem}\@.  In particular, we recall the
identifications~\eqref{eq:orbitT-theorem1} and~\eqref{eq:Torbit-identify}\@.
Note that since $\sL^{(\infty)}(\sX)_x$ is finitely generated, it follows
from Lemma~\ref{plem:hattedvf} from the proof of
Theorem~\ref{the:orbitT-theorem} that $\sL^{(\infty)}(\hat{\sX})_x$ is also
locally finitely generated.  Let $x\in\Orb_T(x_0,\sX)$\@.  From the finitely
generated Orbit Theorem we have
\begin{equation*}
\tb[(T,x)]{\Orb((0,x_0),\hat{\sX})}=\dist{L}^{(\infty)}(\hat{\sX})_{(T,x)}.
\end{equation*}
From lemma~\ref{plem:hattedvf} from the proof of
Theorem~\ref{the:orbitT-theorem} we have
\begin{equation*}
\dist{L}^{(\infty)}(\hat{\sX})_{(T,x)}=
\setdef{(T,0,v_x)\in\real\times\real\times\tb{\man{M}}}
{v_x\in\dist{L}^{(\infty)}(\sX)}.
\end{equation*}
By Proposition~\ref{prop:liealggen}, if $\hat{Y}\in\sL^{(\infty)}(\hat{\sX})$
then we can write
\begin{equation*}
\hat{Y}=\sum_{j=1}^k\lambda_j\hat{X}_j+X
\end{equation*}
for $X_1,\dots,X_k\in\sX$\@, for $\lambda_1,\dots,\lambda_k\in\real$ and for
$X\in\sD(\hat{\sX})$\@.  If $\hat{Y}$ further has the property that
$\hat{Y}(T,x)\in\tb[(T,x)]{\Orb_T(x_0,\sX)}$ (under the
identification~\eqref{eq:Torbit-identify}), then it follows that
$\sum_{j=1}^k\lambda_j=0$\@, using the fact that $X(T,x)=(T,0,X'(x))$ for
some $X'\in\sD(\sX)$\@.  Thus
$\tb[(T,x)]{\Orb_T(x_0,\sX)}\subset\dist{I}(\sX)_x$\@, as desired.
\end{proof}
\end{theorem}

As with the finitely generated Orbit Theorem, the assumption that
$\sL^{(\infty)}(\sX)_x$ is finitely generated is necessary, and
Example~\ref{eg:chow-converse} serves to demonstrate this necessity.

As is usual with these notions of finite generatedness, one has the following
two cases where finite generation is guaranteed.
\begin{corollary}
Let\/ $\man{M}$ be a\/ $\C^\infty$-manifold and let\/ $\sX$ be a family of
partially defined smooth vector fields such that the distribution\/
$\dist{L}^{(\infty)}(\sX)$ is regular.  Then, for each\/ $x_0\in\man{M}$
and\/ $T\in\real$ such that\/ $\Orb_T(x_0,\sX)\not=\emptyset$\@,
\begin{compactenum}[(i)]
\item $\Orb_T(x_0,\sX)$ is a connected immersed smooth submanifold of\/
$\man{M}$ and
\item for each\/ $x\in\Orb_T(x_0,\sX)$\@,\/
$\tb[x]{\Orb_T(x_0,\sX)}=\dist{I}(\sX)_x$\@.
\end{compactenum}
\begin{proof}
This follows from Theorem~\ref{the:fingen-orbitT-theorem}\@, along with
Theorem~\ref{the:locally-free-sheaf}\@.
\end{proof}
\end{corollary}

\begin{corollary}
Let\/ $\man{M}$ be an analytic manifold and let\/ $\sX$ be a family of
partially defined analytic vector fields.  Then, for each\/ $x_0\in\man{M}$
and\/ $T\in\real$ such that\/ $\Orb_T(x_0,\sX)\not=\emptyset$\@,
\begin{compactenum}[(i)]
\item $\Orb_T(x_0,\sX)$ is a connected immersed analytic submanifold of\/
$\man{M}$ and
\item for each\/ $x\in\Orb_T(x_0,\sX)$\@,\/
$\tb[x]{\Orb_T(x_0,\sX)}=\dist{I}(\sX)_x$\@.
\end{compactenum}
\begin{proof}
This follows from Theorem~\ref{the:fingen-orbitT-theorem}\@, along with
Theorem~\ref{the:analytic-fingen}\@.
\end{proof}
\end{corollary}

\section{Frobenius's Theorem for subsheaves and
distributions}\label{sec:frobenius}

We next use the Orbit Theorem to prove Frobenius's Theorem.  We give the
statement in terms of both distributions and subsheaves.

\subsection{Involutive distributions and subsheaves}

Frobenius's Theorem connects two concepts: integrability and involutivity.
Let us first consider involutivity.
\begin{definition}
Let\/ $r\in\{\infty,\omega\}$\@, let $\man{M}$ be a manifold of class
$\C^r$\@, let $\dist{D}$ be a distribution of class $\C^r$\@, and let
$\sF=\ifam{F(\nbhd{U})}_{\nbhd{U}\,\textrm{open}}$ be a subsheaf of
$\ssections[r]{\tb{\man{M}}}$\@.
\begin{compactenum}[(i)]
\item The subsheaf $\sF$ is \defn{involutive} if $[X,Y]\in F(\nbhd{U})$ for
every $X,Y\in F(\nbhd{U})$ and every open $\nbhd{U}\subset\man{M}$\@,~\ie~if
$\sF$ is a Lie subalgebra of $\ssections[r]{\tb{\man{M}}}$\@.
\item The distribution $\dist{D}$ is \defn{involutive} if
$\dist{L}^{(\infty)}(\ssections[r]{\dist{D}})=\dist{D}$\@.\oprocend
\end{compactenum}
\end{definition}

Let us explore the relationship between the two notions of involutivity.  We
start with two examples that show that, in general, there will be no exact
correspondence between involutivity of sheaves and involutivity of the
distributions they generate.
\begin{example}\label{eg:involutive1}
We consider one of the cases of Example~\ref{eg:LinftyX-D(X)}\@.
Specifically, we take $\man{M}=\real^2$ and take the two vector fields
\begin{equation*}
X_1=\pderiv{}{x_1},\quad X_2=f(x_1)\pderiv{}{x_2},
\end{equation*}
where
\begin{equation*}
f(x)=\begin{cases}\eul^{-1/x^2},&x\not=0,\\0,&x=0,\end{cases}
\end{equation*}
We let $\sX=\{X_1,X_2\}$ and let
$\sF=\ifam{F(\nbhd{U})}_{\nbhd{U}\,\textrm{open}}$ be the subsheaf defined by
$F(\nbhd{U})=\sL^{(\infty)}(\modgen{\sX|\nbhd{U}})$\@.  By
Proposition~\ref{prop:Linfty-equiv2} we have
\begin{equation*}
\dist{D}(\sF)=\dist{L}^{(\infty)}(\sX).
\end{equation*}
A simple inductive argument shows that the only nonzero Lie brackets of the
form in Proposition~\ref{prop:liealggen} are given by
\begin{equation*}
[\underbrace{X_1,[X_1,\cdots,[X_1}_{k\ \textrm{times}},X_2]\cdots]]=
f^{(k)}\pderiv{}{x_2}.
\end{equation*}
Thus we conclude from Proposition~\ref{prop:liealggen} that
\begin{equation*}
\dist{D}(\sF)_{(x_1,x_2)}=\begin{cases}\tb[(x_1,x_2)]{\real^2},&x_1\not=0,\\
\vecspan[\real]{\pderiv{}{x_1}},&x_1=0.\end{cases}
\end{equation*}
From Example~\ref{eg:LinftyX-D(X)} we see that $\dist{D}(\sF)$ is not
involutive, although $\sF$ clearly is.\oprocend
\end{example}

\begin{example}\label{eg:involutive2}
Let us consider the vector fields
\begin{equation*}
X_1(x_1,x_2)=(x_1^2+x_2^2)\pderiv{}{x_1},\quad
X_2(x_1,x_2)=(x_1^2+x_2^2)\pderiv{}{x_2}
\end{equation*}
on $\real^2$\@.  We let $\sX=\{X_1,X_2\}$ and take $\sF=\modgen{\sX}$\@.  By
Proposition~\ref{prop:subbundlespan1} the distribution generated by these
vector fields is
\begin{equation*}
\dist{D}(\sF)_{(x_1,x_2)}=
\begin{cases}\tb[(x_1,x_2)]{\real^2},&(x_1,x_2)\not=(0,0),\\
0,&(x_1,x_2)=(0,0).\end{cases}
\end{equation*}
By Frobenius's Theorem below, $\dist{D}(\sF)$ is involutive since it is
integrable (integrability is discussed in the next section).

We claim that $\sF$ is also involutive.  To see this, let
$X,Y\in\modgen{X_1,X_2}$ and write
\begin{equation*}
X=\alpha_XX_1+\beta_XX_2,\quad Y=\alpha_YX_1+\beta_XX_2
\end{equation*}
for smooth or analytic functions $\alpha_X$\@, $\alpha_Y$\@, $\beta_X$\@, and
$\beta_Y$\@.  A direct calculation then shows that
$[X,Y]\in\modgen{X_1,X_2}$\@.  Thus, not only is $\dist{D}(\sF)$\@, the
generators $\ifam{X_1,X_2}$ generate a subsheaf of vector fields that is
involutive.

Next consider the generators
\begin{equation*}
X'_1(x_1,x_2)=(x_1^2+x_2^2)\pderiv{}{x_1},\quad
X'_2(x_1,x_2)=(x_1^4+x_2^4)\pderiv{}{x_2}
\end{equation*}
for $\dist{D}(\sF)$\@.  We let $\sF'$ be the subsheaf generated by these
vector fields.  In this case we calculate
\begin{equation*}
[X'_1,X'_2](x_1,x_2)=
\underbrace{-\frac{2x_2(x_1^4+x_2^4)}{x_1^2+x_2^2}}_{\alpha}X'_1+
\underbrace{\frac{4x_1^3(x_1^2+x_2^2)}{x_1^4+x_2^4}}_{\beta}X'_2.
\end{equation*}
We calculate
\begin{equation*}
\pderiv{\beta}{x_1}=\frac{4(x_1^8-x_1^6x_2^2+5x_1^4+x_2^4+3x_1^2x_2^6)}
{(x_1^4+x_2^4)^2},
\end{equation*}
which is not continuous since
\begin{equation*}
\pderiv{\beta}{x_1}(x_1,0)=4,\quad\pderiv{\beta}{x_1}(0,x_2)=0.
\end{equation*}
Thus, while $\dist{D}(\sF')=\dist{D}(\sF)$ is involutive, the subsheaf $\sF'$
is not.\oprocend
\end{example}

Note that the second of these examples applies to the analytic case.  The
first example, however, is smooth but not real analytic.  Indeed, by
Theorem~\ref{the:LinftyX-D(X)}\@, if $\sF$ is an involutive analytic subsheaf
of vector fields, then $\dist{D}(\sF)$ is necessarily involutive.

\subsection{Integral manifolds}

A related notion to an orbit is the following.
\begin{definition}
Let $r\in\{\infty,\omega\}$\@, let $\man{M}$ be a $\C^r$-manifold, let\/
$\sF$ be a subsheaf of\/ $\ssections[r]{\tb{\man{M}}}$\@, and let $\dist{D}$
be a $\C^r$-distribution.
\begin{compactenum}[(i)]
\item An \defn{integral manifold} of $\dist{D}$ is a $\C^r$-immersed
submanifold $\man{S}$ of $\man{M}$ such that $\tb[x]{\man{S}}=\dist{D}_x$ for
every $x\in\man{S}$\@.
\item An integral manifold $\man{S}$ for $\dist{D}$ is \defn{maximal} if it
is connected, and if every connected integral manifold $\man{S}'$ for
$\dist{D}$ such that $\man{S}'\cap\man{S}\not=\emptyset$ is an open
submanifold of $\man{S}$\@.
\item The distribution $\dist{D}$ is \defn{integrable} if, for each
$x\in\man{M}$\@, there exists an integral manifold of $\dist{D}$ containing
$x$\@.
\item The subsheaf $\sF$ is \defn{integrable} if the distribution
$\dist{D}(\sF)$ is integrable.
\item A \defn{$\C^r$-foliation} of $\man{M}$ is a family
$\ifam{\man{S}_a}_{a\in A}$ of pairwise disjoint immersed $\C^r$-submanifolds
such that
\begin{compactenum}[(a)]
\item $\man{M}=\cup_{a\in A}\man{S}_a$ and
\item \label{pl:foliation2} for each $x_0\in\man{M}$\@, there exists a
neighbourhood $\nbhd{N}$ of $x$ and a family $\ifam{X_b}_{b\in B}$ of\/
$\C^r$-vector fields for which
$\tb[x]{\man{S}_a}=\vecspan[\real]{X_b(x)|\enspace b\in B}$\@.\oprocend
\end{compactenum}
\end{compactenum}
\end{definition}

Let us illustrate these definitions with examples.
\begin{examples}\label{eg:integral-manifolds}
\begin{compactenum}
\item In Example~\enumdblref{eg:orbits}{enum:singular-foliation} we
considered an example with $\man{M}=\real^2$ and define
\begin{equation*}
X_1=x_1\pderiv{}{x_1},\quad X_2=x_2\pderiv{}{x_2}.
\end{equation*}
By $\dist{D}$ we denote the distribution generated by the vector fields
$\sX=\ifam{X_1,X_2}$\@.  The orbits for $\sX$ are shown in
Figure~\ref{fig:R2orbits}\@, and we note that these are also the maximal
integral manifolds for $\dist{D}$\@.  Note that the dimension of the integral
manifolds passing through distinct points may have different dimensions.
Moreover, the family of maximal integral manifolds comprises a foliation.

\item Let $\man{M}=\real^3$ and define
\begin{equation*}
X_1=\pderiv{}{x_2},\quad X_2=\pderiv{}{x_1}+x_2\pderiv{}{x_3}.
\end{equation*}
We let $\dist{D}$ be the distribution generated by these vector fields.  We
shall see that Frobenius's Theorem provides an easy means of verifying that
this distribution does not possess integral manifolds.  However, let us
verify this ``by hand'' to possibly get some insight.  Let us fix $t\in\real$
and compute
\begin{align*}
\flow{X_1}{t}(0,0,0)=&\;(0,t,0),\\
\flow{X_2}{t}\scirc\flow{X_1}{t}(0,0,0)=&\;(t,t,t^2),\\
\flow{X_1}{-t}\scirc\flow{X_2}{t}\scirc\flow{X_1}{t}(0,0,0)=&\;(t,0,t^2),\\
\flow{X_2}{-t}\scirc\flow{X_1}{-t}\scirc\flow{X_2}{t}\scirc
\flow{X_1}{t}(0,0,0)=&\;(0,0,t^2).
\end{align*}
Now suppose that $\man{S}$ is an integral manifold for $\dist{D}$ containing
$\vect{0}=(0,0,0)$\@.  Thus $\man{S}$ must be two-dimensional.  Since
$X_1,X_2\in\sections{\dist{D}}$ and since $\tb{\man{S}}\subset\dist{D}$\@, it
must be the case the integral curves, and therefore concatenations of
integral curves, of $X_1$ and $X_2$ with initial conditions in $\man{S}$ must
remain in $\man{S}$\@.  Therefore, for each $t\in\real$\@, we must have
$(0,0,t^2)\in\man{S}$\@.  Therefore,
$\pderiv{}{x_3}\in\tb[\vect{0}]{\man{S}}$\@.  However, one readily checks
that $\ifam{X_1(\vect{0}),X_2(\vect{0}),\pderiv{}{x_3}}$ is linearly
independent, prohibiting $\man{S}$ from being two-dimensional.  Thus
$\dist{D}$ has no integral manifold passing through $\vect{0}$\@.  One can
show, in fact, that $\dist{D}$ possesses no integral manifolds passing
through \emph{any} point.

\item \label{enum:!analytic-intman} We next consider the example from
Example~\enumdblref{eg:orbits}{enum:!analytic-orbit}\@.  We take
$\man{M}=\real^2$ and define
\begin{equation*}
X_1(x_1,x_2)=\pderiv{}{x_1},\quad X_2(x_1,x_2)=f(x_1)\pderiv{}{x_2},
\end{equation*}
where
\begin{equation*}
f(x)=\begin{cases}\eul^{-1/x^2},&x\in\realp,\\0,&x\in\realnp.\end{cases}
\end{equation*}
In Example~\enumdblref{eg:orbits}{enum:!analytic-orbit} we showed that there
was one orbit, and this was all of $\real^2$\@.

If $(x_{01},x_{02})\in\real^2$ with $x_{01}>0$ we can see that
\begin{equation*}
\setdef{(x_1,x_2)\in\real^2}{x_1\in\realp}
\end{equation*}
is the unique maximal integral manifold of $\dist{D}$ through
$(x_{01},x_{02})$\@.  For $(x_{01},x_{02})\in\real^2$ with $x_{01}<0$\@, the
maximal integral manifold of $\dist{D}$ through $(x_{01},x_{02})$ is
\begin{equation*}
\setdef{(-x_1,x_{02})}{x_1\in\realp}.
\end{equation*}
Note that there are no integral manifolds through points on the
$x_2$-axis.

\item Our next example shows that integral manifolds can be
isolated.  We consider the vector fields
\begin{equation*}
X_1(x_1,x_2,x_3)=x_1x_3\pderiv{}{x_1}+\pderiv{}{x_2},\quad
X_3(x_1,x_2,x_2)=\pderiv{}{x_3}
\end{equation*}
on $\real^3$\@, and let $\dist{D}$ be the distribution generated by these
vector fields.  Note that
\begin{equation*}
\man{S}=\setdef{(x_1,x_2,x_3)}{x_1=0}
\end{equation*}
is an integral manifold for $\dist{D}$\@.  However, $\dist{D}$ is not
integrable.  Let us verify this, using Frobenius's Theorem below.  We
calculate
\begin{equation*}
[X_1,X_2](x_1,x_2,x_3)=-x_1\pderiv{}{x_1},
\end{equation*}
and note that $[X_1,X_2](x_1,x_2,x_3)\in\dist{D}_{(x_1,x_2,x_3)}$ if and only
if $x_1=0$\@.  Thus $\dist{D}|(\real^3\setminus\man{S})$ is not involutive
and so $\dist{D}$ possesses no integral manifolds other than
$\man{S}$\@.
\item Let us now consider the family
$\sX=\ifam{(X_j,\nbhd{U}_j)}_{j\in\{1,2\}}$ of partially defined vector
fields from Example~\enumdblref{eg:orbits}{enum:partially-orbit}\@.  In this
example, the orbits are also integral manifolds for $\dist{D}(\sX)$\@.  The
resulting family of immersed submanifolds defines a smooth foliation
(\cf~Proposition~\ref{prop:smooth-partial}) but not an analytic foliation
(\cf~Example~\ref{eg:analytic-partial}).\oprocend
\end{compactenum}
\end{examples}

\subsection{Frobenius's Theorem, examples, and counterexamples}

The smooth part of the following theorem was proved by \citet{GF:77} and the
analytic part was proved by \citet{TN:66}\@.  Contributions also come
from~\cite{RH:60}\@.
\begin{theorem}\label{the:frobenius}
Let\/ $r\in\{\infty,\omega\}$\@, let\/ $\man{M}$ be a\/ $\C^r$-manifold,
let\/ $\sF$ be a subsheaf of\/ $\ssections[r]{\tb{\man{M}}}$\@, and let\/
$\dist{D}$ be a\/ $\C^r$-distribution on\/ $\man{M}$\@.  Then the following
statements hold:
\begin{compactenum}[(i)]
\item \label{pl:frobenius1} if\/ $\dist{D}$ is integrable then it is
involutive;
\item \label{pl:frobenius2} if\/ $\ssections[r]{\dist{D}}$ is locally
finitely generated and if\/ $\dist{D}$ is involutive, then it is
integrable;
\item \label{pl:frobenius3} if\/ $\sF$ is locally finitely generated and
involutive, then it is integrable.\savenum
\end{compactenum}
In particular,
\begin{compactenum}[(i)]\resumenum
\item \label{pl:frobenius4} if\/ $r=\infty$ and if\/ $\rank_{\dist{D}}$ is
locally constant, then\/ $\dist{D}$ is integrable if and only if it is
involutive and
\item \label{pl:frobenius5} if\/ $r=\omega$ then\/ $\dist{D}$ is integrable
if and only if it is involutive.
\end{compactenum}
Moreover, in case the hypotheses are satisfied in either of the above cases,
the set of maximal integral manifolds forms a foliation of $\man{M}$\@.
\begin{proof}
In the proof of the theorem we shall use the Orbit Theorem for a certain
class of partially defined vector fields.  Specifically, given a distribution
$\dist{D}$\@, we shall consider the partially defined vector fields
$(X,\nbhd{U})$ where $\nbhd{U}\subset\man{M}$ is open and
$X\in\sections[r]{\dist{D}|\nbhd{U}}$\@.  Thus the class of partially defined
vector fields is exactly the collection of local sections of the subsheaf
$\ssections[r]{\dist{D}}$\@.  For this reason, we this family of partially
defined vector fields simply by $\ssections[r]{\dist{D}}$\@.

\eqref{pl:frobenius1} First suppose that $\dist{D}$ is integrable.  Let
$x_0\in\man{M}$ and let $\man{S}$ be the maximal integral manifold through
$x_0$\@.  Since every $\dist{D}$-valued vector field is tangent to $\man{S}$
since $\tb{\man{S}}\subset\dist{D}$\@, it follows that
$\nbhd{U}\subset\man{S}$ for some neighbourhood $\nbhd{U}$ of $x_0$ in the
orbit topology.  Thus $\tb{\nbhd{U}}\subset\tb{\man{S}}$\@.  Moreover, if
$x\in\nbhd{U}$ and if $v_x\in\tb[x]{\man{S}}=\dist{D}_x$\@, then let
$X\in\ssections[r]{\dist{D}}(\nbhd{U})$ be such that $v_x=X(x)$ (recalling
Theorem~\ref{the:global-span} if $r=\infty$ or
Theorem~\ref{the:global-sections} if $r=\omega$).  Then, since
$\flow{X}{t}(x)\in\nbhd{U}\subset\Orb(x_0,\ssections[r]{\dist{D}})$ for $t$
sufficiently small, it follows that $X(x)=v_x\in\tb[x]{\nbhd{U}}$\@.  That is
to say, for $x\in\nbhd{U}$\@,
$\tb[x]{\man{S}}\subset\tb[x]{\Orb(x_0,\ssections[r]{\dist{D}})}$\@.  Thus
$\tb[x]{\man{S}}=\tb[x]{\Orb(x_0,\ssections[r]{\dist{D}})}$ since we
obviously have
$\tb[x]{\Orb(x_0,\ssections[r]{\dist{D}})}\subset\tb[x]{\man{S}}$\@.  Since
$\dist{L}^{(\infty)}(\ssections[r]{\dist{D}})_x\subset
\tb[x]{\Orb(x_0,\ssections[r]{\dist{D}})}$ by
Theorem~\ref{the:Linfty-orbit}\@, it follows that
$\dist{L}^{(\infty)}(\ssections[r]{\dist{D}})_x=\dist{D}_x$\@, so $\dist{D}$
is involutive.

\eqref{pl:frobenius2} Conversely, suppose that $\dist{D}$ is involutive and
that $\ssections[r]{\dist{D}}$ is locally finitely generated.  Involutivity
of $\dist{D}$ implies that if $X,Y\in\sections[r]{\dist{D}|\nbhd{U}}$ for
some open set $\nbhd{U}\subset\man{M}$\@, then $[X,Y](x)\in\dist{D}_x$ for
every $x\in\nbhd{U}$\@.  Thus $[X,Y]\in\sections[r]{\dist{D}|\nbhd{U}}$ and
so $\ssections[r]{\dist{D}}$ is an involutive subsheaf.  Since
$\sL^{(\infty)}(\ssections[r]{\dist{D}})=\ssections[r]{\dist{D}}$\@, it
follows that $\sL^{(\infty)}(\ssections[r]{\dist{D}})$ is locally finitely
generated.  Therefore, by Theorem~\ref{the:fingen-orbit-theorem}\@,
$\tb[x]{\Orb(x,\ssections[r]{\dist{D}})}=\dist{D}_x$ for every
$x\in\man{M}$\@.  Thus $\Orb(x,\ssections[r]{\dist{D}})$ is an integral
manifold for $\dist{D}$ through $x$\@.

\eqref{pl:frobenius3} By Theorem~\ref{the:LinftyX-D(X)} and involutivity of
$\sF$\@,
$\dist{L}^{(\infty)}(\sF)=\dist{L}^{(\infty)}(\ssections[r]{\dist{D}})$\@.
This part of the result then follows from part~\eqref{pl:frobenius2}\@.

Parts~\eqref{pl:frobenius4} and~\eqref{pl:frobenius5} follow from
Theorems~\ref{the:locally-free-sheaf} and~\ref{the:analytic-fingen}\@,
respectively.

Now we verify the final assertion of the theorem.  Disjointness of the orbits
ensures that the orbits are maximal integral manifolds, and moreover shows
that the maximal integral manifolds form a partition of $\man{M}$\@.  That
this partition is a foliation follows since local generators for $\dist{D}$
will satisfy part~\eqref{pl:foliation2} in the definition of a foliation.
\end{proof}
\end{theorem}

Example~\ref{eg:involutive2} shows that part~\eqref{pl:frobenius1} of the
proceeding theorem is not generally true for subsheaves; that is, it may be
the case that an integrable locally finitely generated subsheaf may not be
involutive.  Note that
Example~\enumdblref{eg:integral-manifolds}{enum:!analytic-intman} shows that
any attempt to relax the constant rank condition in the $\C^\infty$-case will
be met with failure in general.  Let us clarify this.
\begin{example}\label{eg:converse-frobenius}
We consider the vector fields
\begin{equation*}
X_1(x_1,x_2)=\pderiv{}{x_1},\quad X_2(x_1,x_2)=f(x_1)\pderiv{}{x_2}
\end{equation*}
on $\real^2$\@, where
\begin{equation*}
f(x)=\begin{cases}\eul^{-1/x^2},&x\in\realp,\\0,&x\in\realnp.\end{cases}
\end{equation*}
As we have seen in
Example~\enumdblref{eg:integral-manifolds}{enum:!analytic-intman}\@, the
distribution $\dist{D}$ generated by $\ifam{X_1,X_2}$ is not integrable since
there is no integral manifold for $\dist{D}$ passing through the points of
the form $(0,x_2)$\@, $x_2\in\real$\@.  We claim that $\dist{D}$ is
involutive.  To see this, suppose that $X,Y\in\sections[\infty]{\dist{D}}$
and write
\begin{equation*}
X=\alpha_X\pderiv{}{x_1}+\beta_X\pderiv{}{x_2},\quad
Y=\alpha_Y\pderiv{}{x_1}+\beta_Y\pderiv{}{x_2}
\end{equation*}
for smooth functions $\alpha_X$\@, $\alpha_Y$\@, $\beta_X$\@, and
$\beta_Y$\@.  Then compute
\begin{multline*}
[X,Y](x_1,x_2)=\left(\pderiv{\alpha_X}{x_1}\alpha_Y+
\pderiv{\alpha_X}{x_2}\beta_Y-\alpha_X\pderiv{\alpha_Y}{x_1}-
\pderiv{\alpha_Y}{x_2}\beta_X\right)\pderiv{}{x_1}\\
-\left(\alpha_X\pderiv{\beta_Y}{x_1}-\alpha_Y\pderiv{\beta_X}{x_1}-
\pderiv{\beta_X}{x_2}\beta_Y+\beta_X\pderiv{\beta_Y}{x_2}\right)
\pderiv{}{x_2}.
\end{multline*}
We consider three cases.
\begin{enumerate}
\item $x_1\in\realn$\@: Here, in some neighbourhood of $(x_1,x_2)$\@,
$\beta_X$ and $\beta_Y$ are zero.  In this case, in this neighbourhood
$[X,Y]$ is collinear with $\pderiv{}{x_1}$ and so
$[X,Y](x_1,x_2)\in\dist{D}_{(x_1,x_2)}$\@.

\item $x_1=0$\@: In this case, the requirement that $X$ and $Y$ are
$\dist{D}$-valued implies that
\begin{equation*}
\beta_X(x_1,x_2)=\beta_Y(x_1,x_2)=0,\quad
\pderiv{\beta_X}{x_1}(x_1,x_2)=\pderiv{\beta_Y}{x_1}(x_1,x_2)=0.
\end{equation*}
Thus, in this case we again have $[X,Y](x_1,x_2)$ collinear with
$\pderiv{}{x_1}$\@, and so $[X,Y](x_1,x_2)\in\dist{D}_{(x_1,x_2)}$\@.

\item $x_1\in\realp$\@: Here we obviously have
$[X,Y](x_1,x_2)\in\dist{D}_{(x_1,x_2)}$ since
$\dist{D}_{(x_1,x_2)}=\tb[(x_1,x_2)]{\real^2}$ when $x_1\in\realp$\@.
\end{enumerate}
This gives the desired involutivity of $\dist{D}$\@.\oprocend
\end{example}

\printbibliography
\end{document}